\documentclass[a4paper,10pt]{article}
\usepackage{hyperref} 
\usepackage{amsmath}
\usepackage{amsfonts}
\usepackage{amssymb}
\usepackage{graphicx}
\usepackage{booktabs}
\usepackage[hang]{subfigure}
\usepackage{tabularx}
\usepackage{fancyhdr}
\usepackage{ifthen}
\usepackage{bm}
\usepackage{verbatim} 
\usepackage{amssymb}
\usepackage{caption}
\usepackage{multirow}
\usepackage[numbers,square,comma,sort&compress]{natbib}
\usepackage{color}
\usepackage{xcolor}
\usepackage{commath}
\usepackage{pdflscape}
\usepackage{rotating}
\usepackage{epstopdf}
\usepackage{geometry}
\usepackage{psfrag}
\usepackage{rotating}
\usepackage{relsize}
\usepackage{threeparttable}
\usepackage{todonotes}

\title{A Comprehensive $p$-VEM Framework for Advanced \\ Variable Stiffness Plates with Arbitrary Shapes}

\author{Paola Pia Foligno$^{1}$, Daniele Boffi$^{2,}$$^{3,}$$^{4}$, Fabio Credali$^{2,}$\footnote{Corresponding author.\\ \textsl{Email addresses:} \texttt{paolapia.foligno@polimi.it} (Paola Pia Foligno), \texttt{daniele.boffi@kaust.edu.sa} (daniele Boffi), \texttt{fabio.credali@kaust.edu.sa} (fabio Credali), \texttt{riccardo.vescovini@polimi.it} (Riccardo Vescovini)} ,  Riccardo Vescovini$^1$\\
\small{\textit{$^{1}$Dipartimento di Scienze e Tecnologie Aerospaziali, Politecnico di Milano,}}
\small{\textit{20156 Milano, Italy}} \\
\small{\textit{$^{2}$
CEMSE Division, King Abdullah University of Science and Technology,}} 
\small{\textit{23955 Thuwal, Saudi Arabia}} \\
\small{\textit{$^{3}$
Dipartimento di Matematica ``F. Casorati", Universit\`a degli Studi di Pavia,}}   
\small{\textit{27100 Pavia, Italy}} 
\\
\small{\textit{$^{4}$
IMATI ``E. Magenes", CNR,}}   
\small{\textit{27100 Pavia, Italy}} 
}

\begin{document}
\definecolor{bluepoli}{cmyk}{0.4,0.1,0,0.4}


\captionsetup[figure]{labelfont={color=black}} 
\captionsetup[table]{labelfont={color=black}} 
\captionsetup[algorithm]{labelfont={color=black}} 


\newcommand\T{\rule{0pt}{2.6ex}}
\newcommand\B{\rule[-1.2ex]{0pt}{0pt}}

\newcounter{algsubstate}
\renewcommand{\thealgsubstate}{\alph{algsubstate}}
\newenvironment{algsubstates}
  {\setcounter{algsubstate}{0}%
   \renewcommand{\STATE}{%
     \stepcounter{algsubstate}%
     \Statex {\small\thealgsubstate:}\space}}
  {}

\newcommand\numfontsize{\@setfontsize\Huge{200}{60}}





\newcommand{\hsp}{\hspace{0pt}}

\makeatletter
\renewcommand*\cleardoublepage{%
  \clearpage\if@twoside\ifodd\c@page\else
  \null
  \AddToShipoutPicture*{\BackgroundPic}
  \thispagestyle{empty}%
  \newpage
  \if@twocolumn\hbox{}\newpage\fi\fi\fi}
\makeatother


\definecolor{orangePoliMi}{rgb}{1,0.59,0}			

\newcommand{ \oo }[1]{\overline{#1}}
\newcommand{\plbr}[1]{ \left( #1 \right) }
\newcommand{\sqbr}[1]{ \left[ #1 \right] }
\newcommand{\tr}[0]{ ^{\mathrm{T}} }
\newcommand{ \red }[1]{\textcolor{red}{#1}}
\newcommand{ \blue }[1]{\textcolor{blue}{#1}}
\newcommand{ \green }[1]{\textcolor{green}{#1}}
\newcommand{\fracp}[2]{ \frac{\partial #1}{\partial #2} } 

\newcommand{ \rv }[1]{\textcolor{magenta}{#1}}

\newcommand{ \eq }[1]{Eq.~(\ref{eq:#1})}
\newcommand{ \eqs }[2]{Eqs.~(\ref{eq:#1}) and (\ref{eq:#2})}
\newcommand{ \eqss }[3]{Eqs.~(\ref{eq:#1}), (\ref{eq:#2}) and (\ref{eq:#3})}
\newcommand{ \eqto }[2]{Eqs.~(\ref{eq:#1})-(\ref{eq:#2})}
\newcommand{ \fig }[1]{Figure~\ref{fig:#1}}
\newcommand{ \figs }[2]{Figures~\ref{fig:#1} and \ref{fig:#2}}
\newcommand{ \figto }[2]{Figures~\ref{fig:#1} to \ref{fig:#2}}
\newcommand{ \tab }[1]{Table~\ref{tab:#1}}
\newcommand{ \tabs }[2]{Tables~\ref{tab:#1} and ~\ref{tab:#2}}

\newcommand{ \alg }[1]{Algorithm~\ref{alg:#1}}

\newcommand{ \chap }[1]{Chapter~\ref{chapter:#1}}
\newcommand{ \sect }[1]{Section~\ref{section:#1}}
\newcommand{ \subsect }[1]{Subsection~\ref{subsection:#1}}

\newcommand{\cubr}[1]{ \left\{ #1 \right\} }
\newcommand{\shbr}[1]{ \left| #1 \right| }
\newcommand{\nrbr}[1]{ \left\| #1 \right\| }
\newcommand{\anbr}[1]{ \langle #1 \rangle }

\newcommand{ \mr }[1]{\mathrm{#1}}
\newcommand{ \de }[0]{\, \mathrm{d}}

\newcommand{ \nl }[0]{{\mathrm{nl}}} 
\newcommand{ \lin }[0]{{\mathrm{l}}} 
\newcommand{ \imp }[0]{{\mathrm{imp}}} 
\newcommand{ \buck }[0]{{\mathrm{b}}} 
\newcommand{ \inter }[0]{{\mathrm{int}}} 
\newcommand{ \tang }[0]{{\mathrm{T}}} 
\newcommand{\diff}[0]{{\mathrm{\partial}}} 
\newcommand{\const}[0]{{\mathrm{\theta}}} 
\newcommand{\curl}[0]{{\mathrm{\times}}} 
\newcommand{\ther}[0]{{\mathrm{th}}} 
\newcommand{\curltext}[0]{\text{curl}} 
\newcommand{\diverg}[0]{\text{div}} 
\newcommand{\svd}[0]{\text{svd}} 

\newcommand{\Vspace}[0]{\mathcal{V}} 
\newcommand{\Vspaced}[0]{\mathcal{V}_h} 
\newcommand{\Vspacedenh}[0]{\mathcal{\tilde{V}}_h} 
\newcommand{\VspacedenhV}[0]{\mathcal{\tilde{V}}_h} 
\newcommand{\VspacedenhA}[0]{\mathcal{\tilde{V}}^{\mathcal{A}}_h} 
\newcommand{\Pspace}[0]{\mathcal{P}} 
\newcommand{\Qpolyspace}[0]{\mathcal{Q}} 
\newcommand{\Pspacecurv}[0]{\mathcal{\tilde{P}}} 
\newcommand{\Pspacevec}[0]{\bm{\mathcal{P}}} 
\newcommand{\Pspacevecther}[0]{\bm{\mathcal{P}}^*} 
\newcommand{\Pspacevecunc}[0]{\bm{\mathcal{P}}^\star} 
\newcommand{\Pspaceunc}[0]{\mathcal{P}} 
\newcommand{\kord}[0]{p} 
\newcommand{\laug}[0]{\ell_E} 
\newcommand{\kordaug}[0]{p+\ell_E} 
\newcommand{\Vspacedim}[0]{N_d \big|_E} 
\newcommand{\pspacedim}[0]{N_p \big|_E} 
\newcommand{\pspacedimone}[0]{N_{p-1} \big|_E} 
\newcommand{\pspacedimtwo}[0]{N_{p-2} \big|_E} 
\newcommand{\spacedimaug}[0]{N_{d+\ell_E} \big|_E} 
\newcommand{\spacedimaugV}[0]{N_{d+\ell_E} \big|_E} 
\newcommand{\pspacedimaug}[0]{N_{p+\ell_E} \big|_E} 
\newcommand{\pspacedimaugV}[0]{N_{p+\ell_E} \big|_E} 
\newcommand{\pspacedimaugVone}[0]{N_{p+\ell_E-1} \big|_E} 
\newcommand{\pspacedimaugone}[0]{N_{p+\ell_E-1} \big|_E} 
\newcommand{\globspacedim}[0]{N_d} 
\newcommand{\qpoly}[0]{q} 
\newcommand{\qpolyvec}[0]{\bm{q}} 
\newcommand{\qpolyvecther}[0]{\bm{q}^*} 
\newcommand{\qpolymtrx}[0]{\bm{Q}} 
\newcommand{\qpolyvecunc}[0]{\bm{q}^\star} 
\newcommand{\globdom}[0]{\mathcal{T}_h} 
\newcommand{\VEMspacedimglob}[0]{N_d} 
\newcommand{\POLYspacedimglob}[0]{N_p} 

\newcommand{\bila}[0]{a} 
\newcommand{\bilael}[0]{a^E} 
\newcommand{\bilad}[0]{a_h} 
\newcommand{\biladel}[0]{a_h^E} 
\newcommand{\force}[0]{f} 
\newcommand{\forcevec}[0]{\bm{f}} 
\newcommand{\forced}[0]{f_h} 
\newcommand{\forcedvec}[0]{\bm{f}_h} 
\newcommand{\bodyforcevec}[0]{\bm{\bar{b}}} 
\newcommand{\lineforcevec}[0]{\bm{\bar{t}}} 

\newcommand{\unku}[0]{u} 
\newcommand{\testu}[0]{\delta u} 
\newcommand{\unkud}[0]{u_h} 
\newcommand{\testud}[0]{\delta u_h} 

\newcommand{\surom}[0]{\Omega} 
\newcommand{\boundc}[0]{\Gamma} 
\newcommand{\el}[0]{E} 
\newcommand{\diammax}[0]{h} 
\newcommand{\boundel}[0]{\Gamma^\el} 
\newcommand{\bound}[0]{\Gamma^E_e} 
\newcommand{\boundcurv}[0]{\Tilde{\Gamma}^E_e} 
\newcommand{\normedgeel}[0]{\bm{\hat{n}}_{\Gamma^\el}} 
\newcommand{\normedge}[0]{\bm{\hat{n}}_{\Gamma_e}} 
\newcommand{\normedgecurv}[0]{\bm{\hat{n}}_{\Tilde{\Gamma}_e}} 
\newcommand{\normedgecurvtilde}[0]{\tilde{\bm{\hat{n}}}_{\Tilde{\Gamma}_e}} 

\newcommand{\nv}[0]{n_v} 
\newcommand{\nedges}[0]{n_e} 
\newcommand{\nedgescurv}[0]{\tilde{n}_e} 
\newcommand{\Nv}[0]{N_v} 
\newcommand{\Ne}[0]{N_e} 
\newcommand{\NE}[0]{N_E} 
\newcommand{\diam}[0]{h_E} 

\newcommand{\trialfcn}[0]{\psi} 

\newcommand{\sing}[0]{S} 

\newcommand{\Pinabla}[0]{\Pi^{\nabla}_p} 
\newcommand{\Piok}[0]{\Pi^{0}_p} 
\newcommand{\Piokone}[0]{\Pi^{0}_{\kord-1}} 
\newcommand{\Pioktwo}[0]{\Pi^{0}_{p-2}} 
\newcommand{\Bm}[0]{\bm{B}^{\nabla}} 
\newcommand{\Gm}[0]{\bm{G}^{\nabla}} 
\newcommand{\Dm}[0]{\bm{D}} 
\newcommand{\DmG}[0]{\bm{D}^\mathcal{G}} 
\newcommand{\DmT}[0]{\bm{D}^\mathcal{T}} 
\newcommand{\Qm}[0]{\bm{Q}} 
\newcommand{\Hm}[0]{\bm{H}} 
\newcommand{\Pinablam}[0]{\bm{\tilde{\Pi}}^{\nabla}_p} 
\newcommand{\Pim}[0]{\bm{\Pi}} 
\newcommand{\Piokm}[0]{\bm{\tilde{\Pi}}^{0}_p} 
\newcommand{\Pioktwom}[0]{\bm{\tilde{\Pi}}^{0}_{p-2}} 
\newcommand{\PinablaS}[1]{\Pi^{\mathrm{V#1}}_p} 
\newcommand{\PinablaStwo}[1]{\Pi^{\mathrm{V#1}}_{\kordaug-1}} 
\newcommand{\BmS}[1]{\bm{B}^{\mathrm{V#1}}} 
\newcommand{\GmS}[1]{\bm{G}^{\mathrm{V#1}}} 
\newcommand{\PinablamS}[1]{\bm{\tilde{\Pi}}^{\mathrm{V#1}}_p} 
\newcommand{\PiA}[0]{\Pi^{\mathcal{A}}_p} 
\newcommand{\BmA}[0]{\bm{B}^{\mathcal{A}}} 
\newcommand{\GmA}[0]{\bm{G}^{\mathcal{A}}} 

\newcommand{\dofs}[0]{{\bm{c}_h}} 
\newcommand{\bardofs}[0]{{\bar{\bm{c}}_h}} 
\newcommand{\globdofs}[0]{{\bm{C}_h}} 
\newcommand{\barglobdofs}[0]{{\bar{\bm{C}}_h}} 
\newcommand{\trialfcnm}[0]{\bm{\Psi}} 

\newcommand{\KE}[0]{\bm{K}^E} 
\newcommand{\KEc}[0]{\bm{K}_c^E} 
\newcommand{\KEs}[0]{\bm{K}_s^E} 
\newcommand{\KES}[1]{\bm{S}_{#1}^E} 
\newcommand{\KEScalar}[0]{S^E} 
\newcommand{\FE}[0]{\bm{f}^E} 
\newcommand{\KEg}[0]{\bm{K}_g^E} 
\newcommand{\KEtg}[0]{\bm{K}_\tang^E} 
\newcommand{\ME}[0]{\bm{M}^E} 
\newcommand{\bodyfE}[0]{\bm{f}_b^E} 
\newcommand{\therfE}[0]{\bm{f}_{\mathrm{th}}^E} 
\newcommand{\linefE}[0]{\bm{f}_{\mathrm{l}}^E} 
\newcommand{\KEnl}[0]{\bm{K}_\nl^E } 
\newcommand{\intfE}[0]{\bm{f}_\inter^E} 

\newcommand{\Acoeff}[0]{\mathcal{A}} 

\newcommand{\disp}[0]{\bm{u}} 
\newcommand{\dispther}[0]{\bm{u}^*} 
\newcommand{\dispd}[0]{\bm{u}_h} 
\newcommand{\bardispd}[0]{\bar{\bm{u}}_h} 
\newcommand{\dispdther}[0]{\bm{u}^*_h} 
\newcommand{\testdisp}[0]{\delta \bm{u}} 
\newcommand{\testdispd}[0]{\delta \bm{u}_h} 
\newcommand{\Deltadispd}[0]{\Delta \bm{u}_h} 
\newcommand{\testdispdther}[0]{\delta \bm{u}^*_h} 
\newcommand{\strain}[0]{\bm{\varepsilon}} 
\newcommand{\stress}[0]{\bm{S}} 
\newcommand{\straind}[0]{\bm{\varepsilon}_d} 
\newcommand{\strainl}[0]{\bm{\varepsilon}_\ell} 
\newcommand{\cost}[0]{\mathbb{C}}  
\newcommand{\costs}[0]{\mathbb{D}_s}  
\newcommand{\mass}[0]{\mathbb{M}}  
\newcommand{\masss}[0]{\mathbb{M}_s}  

\newcommand{\Vspacevec}[0]{\bm{\mathcal{V}}} 
\newcommand{\Vspacedvec}[0]{\bm{\mathcal{V}}_h} 
\newcommand{\Vspacedvecenh}[0]{\bm{\mathcal{\tilde{V}}}_h} 
\newcommand{\VspacedvecenhC}[0]{\bm{\mathcal{\tilde{V}}}^{\mathbb{C}}_h} 
\newcommand{\VspacedvecenhN}[0]{\bm{\mathcal{\tilde{V}}}^{\mathbb{N}}_h} 
\newcommand{\VspacedvecenhR}[0]{\bm{\mathcal{\tilde{V}}}^{\mathbb{R}}_h} 
\newcommand{\VspacedvecenhV}[0]{\bm{\mathcal{\tilde{V}}}_h} 
\newcommand{\VspacedvecenhCV}[0]{\bm{\mathcal{\tilde{V}}}^{\mathbb{C},SS}_h} 

\newcommand{\PiK}[0]{\bm{\Pi}^{\bm{\varepsilon}}_p} 
\newcommand{\BmK}[0]{\bm{B}^{\bm{\varepsilon}}} 
\newcommand{\GmK}[0]{\bm{G}^{\bm{\varepsilon}}} 
\newcommand{\PiKm}[0]{\bm{\tilde{\Pi}}^{\bm{\varepsilon}}_p} 
\newcommand{\PiokK}[0]{\bm{\Pi}_p^{0}} 
\newcommand{\PiokKtwo}[0]{\bm{\Pi}_{p-2}^{0}} 
\newcommand{\PiokKone}[0]{\bm{\Pi}_{p-1}^{0}} 
\newcommand{\PiokKther}[0]{\bm{\Pi}_p^{0,*}} 
\newcommand{\PiokKtilde}[0]{\tilde{\bm{\Pi}}_p^{0}} 
\newcommand{\PiokKtwotilde}[0]{\tilde{\bm{\Pi}}_{p-2}^{0}} 
\newcommand{\PiokKonetilde}[0]{\tilde{\bm{\Pi}}_{p-1}^{0}} 
\newcommand{\PiokKthertilde}[0]{\tilde{\bm{\Pi}}_p^{0,*}} 
\newcommand{\PiKS}[1]{\bm{\Pi}^{\mathrm{V#1}}_p} 
\newcommand{\PiKenhascalar}[0]{\Pi^{\mathcal{K}}_\kord} 
\newcommand{\PiKenha}[0]{\bm{\Pi}^{\mathcal{K}}_\kord} 
\newcommand{\PiKstnd}[0]{\bm{\Pi}^{\mathcal{\hat{K}}}_\kord} 
\newcommand{\PiKVC}[0]{\bm{\Pi}^{{\mathcal{K}_\mathrm{VC}}}_\kord} 
\newcommand{\PiKenhatilde}[0]{\tilde{\bm{\Pi}}^{\mathcal{K}}_\kord} 
\newcommand{\PiKstndtilde}[0]{\tilde{\bm{\Pi}}^{\mathcal{\hat{K}}}_\kord} 
\newcommand{\PiKVCtilde}[0]{\tilde{\bm{\Pi}}^{{\mathcal{K}_\mathrm{VC}}}_\kord} 
\newcommand{\PiGenha}[0]{\bm{\Pi}^{\mathcal{G}}_\kord} 
\newcommand{\PiGstnd}[0]{\bm{\Pi}^{\mathcal{\hat{G}}}_\kord} 
\newcommand{\PiGVC}[0]{\bm{\Pi}^{{\mathcal{G}_\mathrm{VC}}}_\kord} 
\newcommand{\PiGenhatilde}[0]{\tilde{\bm{\Pi}}^{\mathcal{G}}_\kord} 
\newcommand{\PiGstndtilde}[0]{\tilde{\bm{\Pi}}^{\mathcal{\hat{G}}}_\kord} 
\newcommand{\PiGVCtilde}[0]{\tilde{\bm{\Pi}}^{{\mathcal{G}_\mathrm{VC}}}_\kord} 
\newcommand{\PiMenha}[0]{\bm{\Pi}^{\mathcal{M}}_\kord} 
\newcommand{\PiMstnd}[0]{\bm{\Pi}^{\mathcal{\hat{M}}}_\kord} 
\newcommand{\PiMenhatilde}[0]{\tilde{\bm{\Pi}}^{\mathcal{M}}_\kord} 
\newcommand{\PiMstndtilde}[0]{\tilde{\bm{\Pi}}^{\mathcal{\hat{M}}}_\kord} 
\newcommand{\PiTenha}[0]{\bm{\Pi}^{\mathcal{T}}_\kord} 
\newcommand{\PiTstnd}[0]{\bm{\Pi}^{\mathcal{\hat{T}}}_\kord} 
\newcommand{\PiTVC}[0]{\bm{\Pi}^{{\mathcal{T}_\mathrm{VC}}}_\kord} 
\newcommand{\PiTenhatilde}[0]{\tilde{\bm{\Pi}}^{\mathcal{T}}_\kord} 
\newcommand{\PiTstndtilde}[0]{\tilde{\bm{\Pi}}^{\mathcal{\hat{T}}}_\kord} 
\newcommand{\PiTVCtilde}[0]{\tilde{\bm{\Pi}}^{{\mathcal{T}_\mathrm{VC}}}_\kord} 
\newcommand{\PiKenhaS}[0]{\bm{\Pi}^{\mathcal{K}}_{\kordaug}} 
\newcommand{\PiKstndS}[0]{\bm{\Pi}^{\mathcal{\hat{K}}}_{\kordaug}} 
\newcommand{\PiKVCS}[0]{\bm{\Pi}^{{\mathcal{K}_\mathrm{VC}}}_{\kordaug}} 
\newcommand{\PiKenhatildeS}[0]{\tilde{\bm{\Pi}}^{\mathcal{K}}_{\kordaug}} 
\newcommand{\PiKstndtildeS}[0]{\tilde{\bm{\Pi}}^{\mathcal{\hat{K}}}_{\kordaug}} 
\newcommand{\PiKVCtildeS}[0]{\tilde{\bm{\Pi}}^{{\mathcal{K}_\mathrm{VC}}}_{\kordaug}} 
\newcommand{\PiKenhaSone}[0]{\bm{\Pi}^{\mathcal{K}}_{\kordaug-1}} 
\newcommand{\PiKstndSone}[0]{\bm{\Pi}^{\mathcal{\hat{K}}}_{\kordaug-1}} 
\newcommand{\PiKVCSone}[0]{\bm{\Pi}^{{\mathcal{K}_\mathrm{VC}}}_{\kordaug-1}} 
\newcommand{\PiKenhatildeSone}[0]{\tilde{\bm{\Pi}}^{\mathcal{K}}_{\kordaug-1}} 
\newcommand{\PiKstndtildeSone}[0]{\tilde{\bm{\Pi}}^{\mathcal{\hat{K}}}_{\kordaug-1}} 
\newcommand{\PiKVCtildeSone}[0]{\tilde{\bm{\Pi}}^{{\mathcal{K}_\mathrm{VC}}}_{\kordaug-1}} 
\newcommand{\BmKenha}[0]{\bm{B}^{\mathcal{K}}} 
\newcommand{\GmKenha}[0]{\bm{G}^{\mathcal{K}}} 
\newcommand{\GmKenhatilde}[0]{\tilde{\bm{G}}^{\mathcal{K}}} 
\newcommand{\BmKstnd}[0]{\bm{B}^{\mathcal{\hat{K}}}} 
\newcommand{\GmKstnd}[0]{\bm{G}^{\mathcal{\hat{K}}}} 
\newcommand{\GmKstndtilde}[0]{\tilde{\bm{G}}^{\mathcal{\hat{K}}}} 
\newcommand{\BmKVC}[0]{\bm{B}^{{\mathcal{K}_\mathrm{VC}}}} 
\newcommand{\GmKVC}[0]{\bm{G}^{{\mathcal{K}_\mathrm{VC}}}} 
\newcommand{\GmKVCtilde}[0]{\tilde{\bm{G}}^{{\mathcal{K}_\mathrm{VC}}}} 
\newcommand{\BmGenha}[0]{\bm{B}^{\mathcal{G}}} 
\newcommand{\GmGenha}[0]{\bm{G}^{\mathcal{G}}} 
\newcommand{\GmGenhatilde}[0]{\tilde{\bm{G}}^{\mathcal{G}}} 
\newcommand{\BmGstnd}[0]{\bm{B}^{\mathcal{\hat{G}}}} 
\newcommand{\GmGstnd}[0]{\bm{G}^{\mathcal{\hat{G}}}} 
\newcommand{\GmGstndtilde}[0]{\tilde{\bm{G}}^{\mathcal{\hat{G}}}} 
\newcommand{\BmGVC}[0]{\bm{B}^{{\mathcal{G}_\mathrm{VC}}}} 
\newcommand{\GmGVC}[0]{\bm{G}^{{\mathcal{G}_\mathrm{VC}}}} 
\newcommand{\GmGVCtilde}[0]{\tilde{\bm{G}}^{{\mathcal{G}_\mathrm{VC}}}} 
\newcommand{\QmMenha}[0]{\bm{Q}^{\mathcal{M}}} 
\newcommand{\HmMenha}[0]{\bm{H}^{\mathcal{M}}} 
\newcommand{\QmMstnd}[0]{\bm{Q}^{\mathcal{\hat{M}}}} 
\newcommand{\HmMstnd}[0]{\bm{H}^{\mathcal{\hat{M}}}} 
\newcommand{\Qmok}[0]{\bm{Q}_p^{0}} 
\newcommand{\Hmok}[0]{\bm{H}_p^{0}} 
\newcommand{\Qmoktwo}[0]{\bm{Q}_{p-2}^{0}} 
\newcommand{\Hmoktwo}[0]{\bm{H}_{p-2}^{0}} 
\newcommand{\Qmokone}[0]{\bm{Q}_{p-1}^{0}} 
\newcommand{\Hmokone}[0]{\bm{H}_{p-1}^{0}} 
\newcommand{\BmTenha}[0]{\bm{B}^{\mathcal{T}}} 
\newcommand{\GmTenha}[0]{\bm{G}^{\mathcal{T}}} 
\newcommand{\GmTenhatilde}[0]{\tilde{\bm{G}}^{\mathcal{T}}} 
\newcommand{\BmTstnd}[0]{\bm{B}^{\mathcal{\hat{T}}}} 
\newcommand{\GmTstnd}[0]{\bm{G}^{\mathcal{\hat{T}}}} 
\newcommand{\GmTstndtilde}[0]{\tilde{\bm{G}}^{\mathcal{\hat{T}}}} 
\newcommand{\BmTVC}[0]{\bm{B}^{{\mathcal{T}_\mathrm{VC}}}} 
\newcommand{\GmTVC}[0]{\bm{G}^{{\mathcal{T}_\mathrm{VC}}}} 
\newcommand{\GmTVCtilde}[0]{\tilde{\bm{G}}^{{\mathcal{T}_\mathrm{VC}}}} 
\newcommand{\BmKenhaS}[0]{\bm{B}^{\mathcal{K}\star}} 
\newcommand{\GmKenhaS}[0]{\bm{G}^{\mathcal{K}\star}} 
\newcommand{\GmKenhaStilde}[0]{\tilde{\bm{G}}^{\mathcal{K}\star}} 
\newcommand{\BmKstndS}[0]{\bm{B}^{\mathcal{\hat{K}\star}}} 
\newcommand{\GmKstndS}[0]{\bm{G}^{\mathcal{\hat{K}\star}}} 
\newcommand{\GmKstndStilde}[0]{\tilde{\bm{G}}^{\mathcal{\hat{K}\star}}} 
\newcommand{\BmKVCS}[0]{\bm{B}^{{\mathcal{K}_\mathrm{VC}\star}}} 
\newcommand{\GmKVCS}[0]{\bm{G}^{{\mathcal{K}_\mathrm{VC}\star}}} 
\newcommand{\GmKVCStilde}[0]{\tilde{\bm{G}}^{{\mathcal{K}_\mathrm{VC}\star}}} 

\newcommand{\bilKc}[0]{\mathcal{K}} 
\newcommand{\bilKenha}[0]{\mathcal{K}^\el} 
\newcommand{\bilKstnd}[0]{\mathcal{\hat{K}}^\el} 
\newcommand{\bilKVC}[0]{{\mathcal{K}_\mathrm{VC}}^\el} 
\newcommand{\bilKVCstar}[0]{{\mathcal{K}_\mathrm{VC}}_*^\el} 
\newcommand{\bilKdVC}[0]{{\mathcal{K}_\mathrm{VC}}_h^\el} 
\newcommand{\bilGc}[0]{\mathcal{G}} 
\newcommand{\bilGenha}[0]{\mathcal{G}^\el} 
\newcommand{\bilGstnd}[0]{\mathcal{\hat{G}}^\el} 
\newcommand{\bilGVC}[0]{{\mathcal{G}_\mathrm{VC}}^\el} 
\newcommand{\bilGVCstar}[0]{{\mathcal{G}_\mathrm{VC}}_*^\el} 
\newcommand{\bilGdVC}[0]{{\mathcal{G}_\mathrm{VC}}_h^\el} 
\newcommand{\bilMc}[0]{\mathcal{M}} 
\newcommand{\bilMenha}[0]{\mathcal{M}^\el} 
\newcommand{\bilMstnd}[0]{\mathcal{\hat{M}}^\el} 
\newcommand{\bilMdenha}[0]{\mathcal{M}_h^\el} 
\newcommand{\bilF}[0]{\mathcal{F}^\el} 
\newcommand{\bilFd}[0]{\mathcal{F}_h^\el} 
\newcommand{\bilTc}[0]{\mathcal{T}} 
\newcommand{\bilTenha}[0]{\mathcal{T}^\el} 
\newcommand{\bilTstnd}[0]{\mathcal{\hat{T}}^\el} 
\newcommand{\bilTVC}[0]{{\mathcal{T}_\mathrm{VC}}^\el} 
\newcommand{\bilTVCstar}[0]{{\mathcal{T}_\mathrm{VC}}_*^\el} 
\newcommand{\bilTdVC}[0]{{\mathcal{T}_\mathrm{VC}}_h^\el} 
\newcommand{\bilKVCnl}[0]{{\mathcal{K}_\mathrm{VC}}_{\nl}^\el} 
\newcommand{\bilKdVCnl}[0]{{\mathcal{K}_\mathrm{VC}}_{\nl,h}^\el} 
\newcommand{\bilFVCint}[0]{\mathcal{\bar{F}}_{\inter}^\el} 
\newcommand{\bilFdVCint}[0]{\mathcal{\bar{F}}_{\inter,h}^\el} 

\newcommand{\trialfcnvec}[0]{\bm{\psi}} 
\newcommand{\trialfcnvecther}[0]{\bm{\psi}^*} 

\newcommand{\vel}[0]{\bm{u}} 
\newcommand{\testvel}[0]{\delta\bm{u}} 
\newcommand{\pres}[0]{\rho} 
\newcommand{\testpres}[0]{\delta \rho} 
\newcommand{\veld}[0]{\bm{u}_h} 
\newcommand{\presd}[0]{\rho_h} 

\newcommand{\Qspace}[0]{Q} 
\newcommand{\Qspaced}[0]{Q_h} 
\newcommand{\Gspace}[0]{\bm{\mathcal{G}}} 
\newcommand{\Gspaceperp}[0]{\bm{\mathcal{G}^{\perp}}} 
\newcommand{\gpolyperp}[0]{\bm{g}^{\perp}} 

\newcommand{\bilb}[0]{b} 
\newcommand{\bilbd}[0]{b} 

\newcommand{\Pinablavec}[0]{\bm{\Pi}^{\bm{\nabla}}_p} 
\newcommand{\BmN}[0]{\bm{B}^{\bm{\nabla}}} 
\newcommand{\GmN}[0]{\bm{G}^{\bm{\nabla}}} 
\newcommand{\Pinablavecm}[0]{\bm{\tilde{\Pi}}^{\bm{\nabla}}_p} 
\newcommand{\PinablavecS}[1]{\bm{\Pi}^{\mathrm{V#1}}_p} 

\newcommand{ \fixme }[1]{\textcolor{red}{FIXME: #1}}

\date{}
\maketitle

\begin{abstract}

This paper presents a comprehensive, high-order ($\kord$-version) Virtual Element Method (VEM) framework for the structural analysis of innovative variable stiffness plates. VEM is particularly suited for complex configurations due to its ability to handle arbitrary polygonal meshes, including curved edges. However, its mathematical formulation may hinder its spread in the engineering community. This work illustrates a formulation with an accessible implementation using well-known FEM notation and integrating at the same time a set of new advanced capabilities. Specifically, both standard stabilized and advanced self-stabilized strategies are adopted. To further improve the robustness of VEM in the presence of variable coefficients, polynomial projections taking into account the coefficients are employed. This approach is referred to as Variable Coefficients-VEM approach ($\mathrm{VC}$-VEM). This unified framework is applied to linear static, free-vibration, and buckling analyses, and validated against analytical solutions and numerical benchmarks. In particular, plates with cutouts and problems featuring high-gradient solutions are investigated, demonstrating that the proposed comprehensive approach provides a flexible and ready-to-implement tool for advanced structural design.

\end{abstract}

\section{Introduction}
 
Innovative structural configurations feature plates characterized by a variable stiffness skin, where the fibers follow curvilinear paths. These variable stiffness (VS) plates offer significant potential for weight minimization and improved efficiency compared to classical designs. Thus, increasing attention has been devoted to their study.

The literature on these innovative configurations is largely limited to the use of classical numerical methods, namely the Finite Element Method and the Ritz method. Recent applications in finite element frameworks include~\cite{raju2015buckling}, where shear buckling and postbuckling responses were compared to classical straight fiber configurations, and~\cite{manickam2018thermal}, which addressed thermal buckling. 
On the other hand, applications of the Ritz method can be found for the mechanical~\cite{wu2012buckling,coburn2016buckling} and thermal~\cite{vescovini2018thermal} buckling, where the use of simple geometries, such as rectangular domains, is due to the global nature of the trial functions. 
Interesting exceptions where the Ritz method is constructed on more complex geometries include~\cite{janssens2021semi,jing2023discrete}, in which plates with arbitrary cutouts have been investigated, and~\cite{jing2024free}, which accounted for arbitrary geometries. Innovative Ritz-based formulations for the linear and geometrically nonlinear responses of arbitrary domains have been proposed in~\cite{vescovini2023ritz,vescovini2025geometrically}. 

It is therefore clear that the study of VS plates requires advanced and effective numerical methods capable of tackling the challenges arising in this field. 
In this regard, the Virtual Element Method (VEM) represents a promising alternative, as it employs general polygonal elements (even with curved edges~\cite{daveiga2019virtual}), thereby significantly simplifying the mesh generation for complex geometries.

Pioneering works~\cite{daveiga2013basic,daveiga2014hitchhiker} introduced the innovative formulation of the method and its implementation for the Laplace problem, while its accuracy for general second-order elliptic problems has been assessed in~\cite{daveiga2016virtual}. 

The mathematical foundations of VEM have been extensively investigated over the last decade. However, most of the existing literature is formulated within a rigorous mathematical framework that may hinder its spread in the engineering community, this consideration being exacerbated by a relatively steep learning curve. Investigations of elasticity problems within the VEM framework can be found in~\cite{daveiga2013virtual,daveiga2015virtual,artioli2017arbitrary,mengolini2019engineering,daltri2020error}, while Kirchhoff-Love and Reissner-Mindlin plate models have been addressed in~\cite{brezzi2013virtual,mora2018virtual} and~\cite{daveiga2019virtual2,daltri2022first}, respectively. Applications to buckling problems can be found in~\cite{mora2020virtual}.
A key limitation of the VEM is the need for a stabilization term to handle the non-polynomial residual component of the trial functions and ensure the well-posedness of the discrete problem. The stability of the method has been investigated in~\cite{daveiga2017stability,mascotto2018ill}, where several stabilization formulas have also been introduced. Generally, in the standard VEM, the stabilization term does not conform with the physics of the problem under consideration and it is scaled by a user-defined tuning parameter. This may reduce the generality of the method and lead to over-stabilization, resulting in degraded solution quality. This is particularly relevant for anisotropic problems~\cite{berrone2022comparison} and eigenvalue problems~\cite{boffi2020approximation}. Alternative approaches have therefore been proposed to avoid arbitrary stabilization terms, such as strategies based on higher-order polynomial projections~\cite{berrone2025lowest}, divergence-free polynomial spaces~\cite{berrone2025stabilization} and enriched VEM spaces with additional internal degrees of freedom (DOFs)~\cite{lamperti2023hu}. The effect of the stabilization parameter has been studied in~\cite{fujimoto2024study}, showing that choosing it based on the bending strain energy density yields more accurate results compared to classical approaches. Self-stabilized formulations for mixed linear elasticity have been proposed in~\cite{liguori2024stabilization}.

The treatment of variable coefficients in the bilinear form is another relevant aspect of the VEM that requires further investigation, particularly in advanced applications. 
In~\cite{daveiga2016virtual}, the standard $L^2$ projection was adopted and the variable coefficients were introduced only in the construction of the final discrete bilinear form. This approach was shown to perform better in the presence of variable coefficients compared to other alternatives.
A VEM application to composites with spatially varying fiber directions was investigated in~\cite{reddy2019virtual}. The authors demonstrated that approximating the fiber direction as a constant equal to the average at the nodal values yields more accurate results than using the element's centroid value. Spatially varying material properties have been studied within a Hellinger–Reissner VEM framework in 2D~\cite{artioli2018family} and in 3D~\cite{visinoni2024family}.
A recent work~\cite{foligno2026benchmarking} introduces a novel formulation, denoted as $\mathrm{VC}$-VEM, in which the coefficient is fully incorporated into the definition of the projection operators, both in stabilized and self-stabilized settings. 

This paper aims to bridge the gap between theoretical VEM formulations and engineering practice by presenting a robust, high-order $\kord$-VEM computational framework with focus on the analysis of variable stiffness plates with complex geometries. For this purpose, the main objective of this work is to provide a self-contained reference, illustrating VEM formulation with an implementation-oriented matrix notation, familiar to the engineering computational mechanics community. 
In addition to the above mentioned objectives, the proposed advanced formulation exploits the geometric flexibility of VEM to consider polygonal elements with curved edges and hanging nodes, hence simplifying the mesh generation for plates with complex boundaries and cutouts. To address the critical choice of the stabilization term, both stabilized and self-stabilized strategies are included and discussed. Furthermore, to enhance numerical robustness in the presence of non-uniform elasticity properties, a novel Variable Coefficient-VEM ($\mathrm{VC}$-VEM) approach is adopted to integrate spatial variability within the projection operators.

The manuscript is organized as follows. \sect{structural_modeling} presents the structural modeling of the variable stiffness plate. In \sect{VEM}, the VEM formulation is introduced, including the stabilized and self-stabilized strategies, as well as the $\mathrm{VC}$-VEM approach. \sect{results} provides the numerical results, in which the developed formulation is validated against analytical solutions and numerical results from the literature. Lastly, the conclusions are drawn in \sect{conclusions}.
\section{Structural modeling} \label{section:structural_modeling}

The structures under investigation are variable stiffness (VS) plates, where the fiber orientation varies across the domain, rendering the stiffness a function of the planar position. A two-dimensional model is employed, with a Cartesian reference system $xyz$, whose origin is located on the plate midsurface. The $x$ and $y$ axes are directed in the longitudinal and transverse directions, respectively, and the $z$ axis is in the thickness direction. The geometry is arbitrary, with maximum dimensions equal to $a$ and $b$ and thickness $t$. 
The plate domain is denoted by $\surom$ and its Lipschitz boundary $\boundc$ consists of a finite number of smooth curves $\cubr{\boundc_i}_{i=1,\dots,\Ne}$, with $\Ne$ being the number of edges. Each curve $\boundc_i$ is of class $\mathcal{C}^{m+1}$ for $m \geq 0$ and it is parametrized by an invertible $\mathcal{C}^{m+1}$ map $\gamma_i: I_i  = \sqbr{a_i,b_i} \rightarrow \boundc_i$~\cite{daveiga2019virtual}.

\subsection{Variational statement} \label{section:PVW}

The formulation is developed within a displacement-based variational framework.

Linear static, buckling, and free-vibration analyses are of concern. Hence, the variational statement can be expressed in a unified form as~\cite{liew1993pb}:
\begin{equation}
    \plbr{\beta_1 + \beta_2 + \beta_3} \delta W_i + \beta_2 \delta W_b + \beta_3 \delta W_k = \beta_1 \delta W_e,
    \label{eq:variational_statement}
\end{equation}
where $W_\alpha$, with $\alpha=i,b,k,e$, are used to denote the internal virtual work, the pre-buckling energy contribution, the first variation of the kinetic energy contribution and the external virtual work, respectively. The Boolean flags $\beta_i$, with $i=1,2,3$, are chosen dependently on the analysis of interest, following the summary reported in \tab{analysis_coefficients}.
\begin{table}[htbp]
\centering
\begin{tabular}{ccccc}
\toprule
{Analysis type} & $\beta_1$ & $\beta_2$ & $\beta_3$\\
\midrule
{Linear static}      & 1 & 0 & 0\\
{Buckling}           & 0 & 1 & 0\\
{Free-vibration}     & 0 & 0 & 1\\
\bottomrule
\end{tabular}
\caption{Variational statement coefficients for the different analysis types.}
\label{tab:analysis_coefficients}
\end{table}

\subsection{Plate model} \label{section:plate}

Hereafter, the two-dimensional plate model is presented. First, the kinematics, the strain measure, and the constitutive law, which accounts for the spatially varying fiber orientations, are detailed. Lastly, the energy terms are derived within the variational framework.

\subsubsection{Kinematics}

The plate kinematics is modeled according to the First-order Shear Deformation Theory (FSDT), which allows thin and relatively thick panels to be considered. The displacement of a generic point on the plate is expressed as~\cite{reddy2004mechanics}:
\begin{equation}
\begin{aligned}
\bm{d}\plbr{x,y,z} = 
\begin{Bmatrix}
    d_x\plbr{x,y,z}\\
    d_y\plbr{x,y,z}\\
    d_z\plbr{x,y,z}
\end{Bmatrix}
&=
\begin{Bmatrix}
    u\plbr{x,y}\\
    v\plbr{x,y}\\
    w\plbr{x,y}
\end{Bmatrix}
+ z
\begin{bmatrix}
    1 & 0\\
    0 & 1\\
    0 & 0\\
\end{bmatrix}
\begin{Bmatrix}
    \theta_x\plbr{x,y}\\
    \theta_y\plbr{x,y}\\
    0
\end{Bmatrix}= \bm{u}^0\plbr{x,y} + z \bm{L}^{\bm{\theta}} \, \bm{\theta}\plbr{x,y}\\
&= 
\begin{bmatrix}
\bm{I} & z\,\bm{L}^{\bm{\theta}}
\end{bmatrix}
\begin{Bmatrix}
\bm{u}^0\plbr{x,y}\\
\bm{\theta}\plbr{x,y}
\end{Bmatrix}
= 
\begin{bmatrix}
\bm{I} & z\,\bm{L}^{\bm{\theta}}
\end{bmatrix}
\bm{u},
\end{aligned}
    \label{eq:FSDT}
\end{equation}
where $\bm{u}^0$ and $\bm{\theta}$ are the generalized displacements and rotations of the midsurface, and $\bm{u}$ is the vector collecting them. The conventions of the plate kinematics are shown in \fig{plate_model}. 

The strains $\hat{\bm{\varepsilon}} = \begin{Bmatrix} \varepsilon_{xx},\; \varepsilon_{yy},\; \gamma_{xy} \end{Bmatrix}^T$ and $\bm{\gamma} = \begin{Bmatrix} \gamma_{yz},\; \gamma_{xy} \end{Bmatrix}^T$ are expressed in terms of the displacements as:
\begin{equation}
\begin{aligned}
\hat{\bm{\varepsilon}} = 
\bm{\varepsilon}^0 \plbr{\bm{u}} + z \, \bm{k} \plbr{\bm{u}}, \qquad \bm{\gamma} = 
\bm{\gamma} \plbr{\bm{u}},
\end{aligned}
\label{eq:GL_strains}
\end{equation}
where $\bm{\varepsilon}^0$ denotes the membrane strains, $\bm{k}$ the curvatures, and $\bm{\gamma}$ the transverse shear strains, and their expression reads:
\begin{equation}
\begin{split}
\bm{\varepsilon}^0 \plbr{\bm{u}} = 
\begin{Bmatrix}
u_{,x} \\
v_{,y} \\
u_{,y} + v_{,x}
\end{Bmatrix},
\qquad
\bm{k} \plbr{\bm{u}} = 
\begin{Bmatrix}
\theta_{x,x} \\
\theta_{y,y} \\
\theta_{x,y} + \theta_{y,x}
\end{Bmatrix},
\qquad
\bm{\gamma} \plbr{\bm{u}} = 
\begin{Bmatrix}
\theta_{y} + w_{,y} \\
\theta_{x} + w_{,x} \\
\end{Bmatrix},
\label{eq:fsdt_strains_lin}
\end{split}
\end{equation}
where $(\cdot)_{,x}$ and $(\cdot)_{,y}$ are the derivatives with respect to the in-plane coordinates. 

In view of future developments, the generalized strains can be conveniently collected into a single vector:
\begin{equation}
\strain \plbr{\bm{u}} = 
\begin{Bmatrix}
\bm{\varepsilon}^0 \plbr{\bm{u}} \\
\bm{k}\plbr{\bm{u}} \\
\bm{\gamma}\plbr{\bm{u}}
\end{Bmatrix}.
\label{eq:elin_enl}
\end{equation}

\subsubsection{Constitutive law}

\begin{figure}[!htbp]
    \centering
    \subfigure[Conventions for generalized displacements.\label{fig:plate_model}]{
    \includegraphics[width=0.45\textwidth]{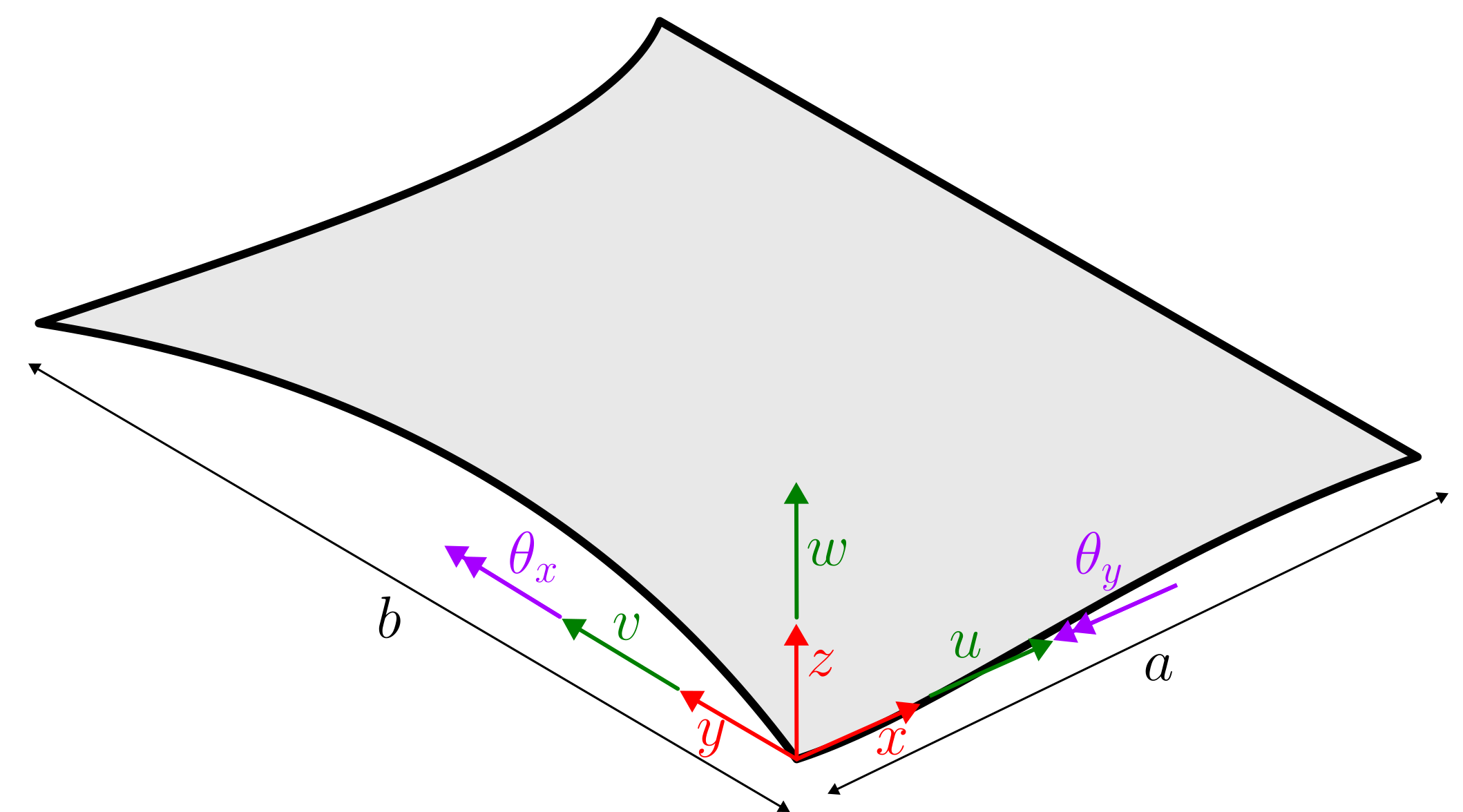}
    }
    \subfigure[Fiber orientations.\label{fig:plate_VS}]{
    \includegraphics[width=0.45\textwidth]{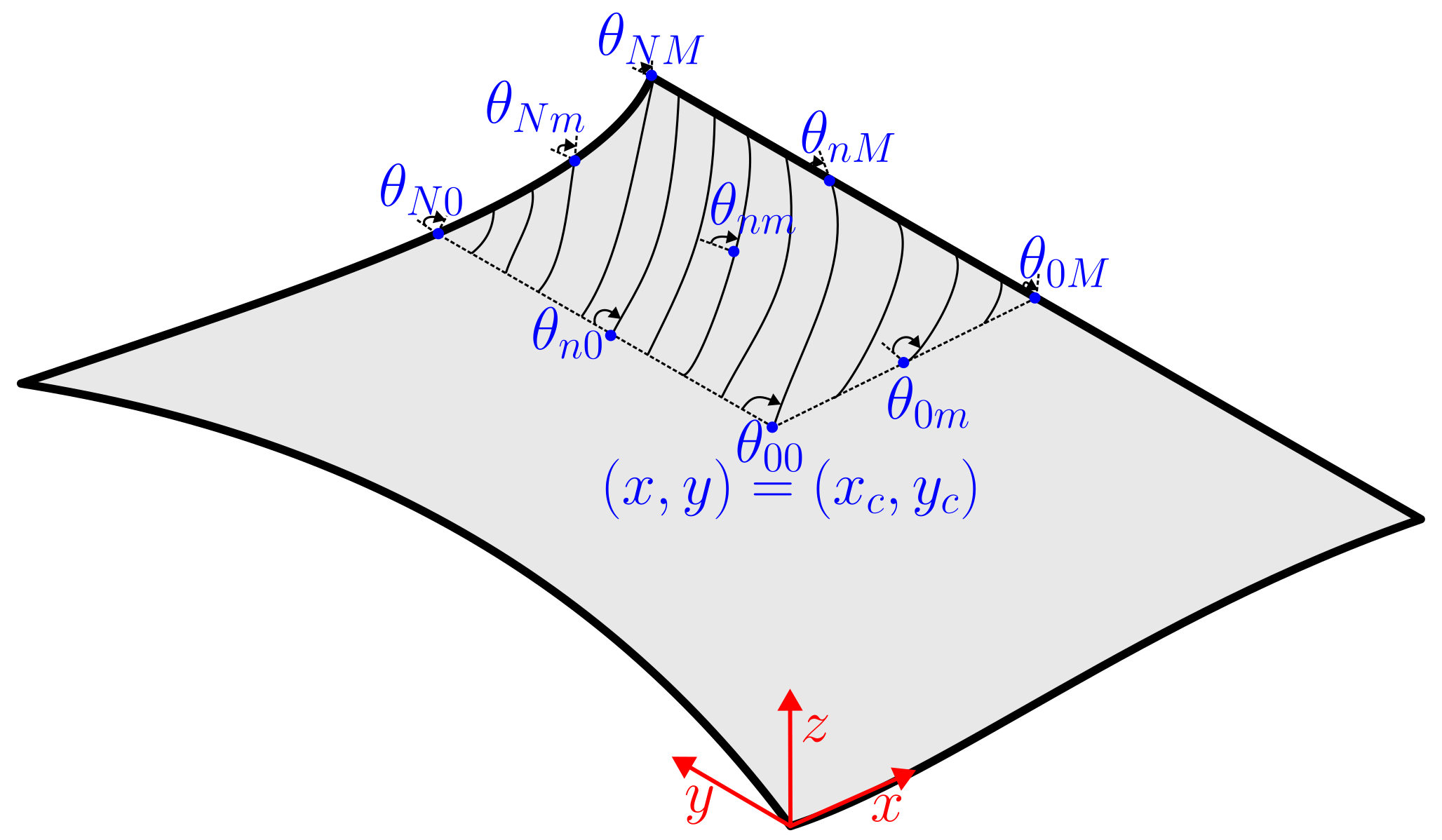}
    }
    \caption{Plate model.}
\end{figure}

In variable stiffness plates, the fiber orientation angles are a function of the position. Different representations have been proposed in the literature. Linear variations~\cite{olmedo1993buckling}, Lobatto distributions~\cite{alhajahmad2008design} and NURBS~\cite{nagendra1995optimization} are examples. In this work, the angles $\theta_{mn}$ are specified on a grid of M $\times$ N points over one quarter of the plate domain, and the angles at a generic point of the domain are retrieved via Lagrange polynomials interpolation~\cite{wu2012buckling,wu2018optimization,zhao2019prestressed,vescovini2020semi}:
\begin{equation}
\theta \plbr{x,y} = \sum_{m=0}^{M-1} \sum_{n=0}^{N-1} \theta_{mn} \prod_{m \neq i} \plbr{ \frac{|x| - x_i}{x_m - x_i } }
\prod_{n \neq j} \plbr{ \frac{|y| - y_j}{y_n - y_j } }.
    \label{eq:lagrange_th}
\end{equation}
The proposed distribution is illustrative of a specific choice, but other choices could be easily accommodated within the present framework. 
From \eq{lagrange_th}, the angle $\theta \plbr{x,y}$ is nonlinear with respect to the coordinates $x$ and $y$. The linear variation is retrieved as a special case by considering only two points, e.g. the center and the edge of the plate. The fiber orientation distribution is illustrated in \fig{plate_VS}.

For each ply $k$, the constitutive law is expressed in global coordinates as:
\begin{equation}
\begin{aligned}
&\bm{\sigma}_k = \bar{\bm{Q}}_k\plbr{x,y} \,\left(\hat{\bm{\varepsilon}}_k - \bm{\varepsilon}_k^\ther \right),
&\bm{\tau}_k = \bar{\bm{Q}}_{nk}\plbr{x,y} \, \bm{\gamma},
\end{aligned}
    \label{eq:plate_therm_constlaw}
\end{equation}
where $\bar{\bm{Q}}_k\plbr{x,y}$ and $\bar{\bm{Q}}_{nk}\plbr{x,y}$ are the constitutive matrices expressed in laminate axes. The vector $\hat{\bm{\varepsilon}}_k$ is the total deformation, and $\bm{\varepsilon}_k^\ther$ is the thermal contribution, accounted only for the in-plane behavior, as out-of-plane thermal stresses are typically negligible for the plates under consideration~\cite{vescovini2020semi}.
The thermal contribution $\bm{\varepsilon}_k^\ther$ is expressed as in~\cite{vescovini2020semi}: 
\begin{equation}
\bm{\varepsilon}_k^\ther = \bar{\bm{\alpha}}_k \plbr{x,y} \Delta T,
    \label{eq:plate_therm_def}
\end{equation}
where $\bar{\bm{\alpha}}_k$ is the ply thermal expansion coefficient in global coordinates and $\Delta T$ is the temperature gradient. The coefficients of $\bar{\bm{\alpha}}_k$ are assumed to be temperature-independent.

The thermo-elastic constitutive law is then obtained as:
\begin{equation}
\bm{S} = 
\begin{Bmatrix}
\bm{N} \\
\bm{M} \\
\bm{Q} 
\end{Bmatrix}
=
\begin{bmatrix}
\bm{A}\plbr{x,y} & \bm{B}\plbr{x,y} & \bm{0} \\
\bm{B}\plbr{x,y} & \bm{D}\plbr{x,y} & \bm{0} \\
\bm{0} & \bm{0} & \bm{A}^{\mathrm{s}}\plbr{x,y} \\
\end{bmatrix}
\strain - \begin{Bmatrix}
\hat{\bm{N}} \\
\hat{\bm{M}} \\
\bm{0} 
\end{Bmatrix} \Delta T
=
\cost\plbr{x,y} \strain - \hat{\bm{R}}\Delta T,
\label{eq:fsdt_cost_platemodel}
\end{equation}
where $\bm{N}$, $\bm{M}$ and $\bm{Q}$ are the forces and moments per unit length and $\bm{S}$ is the vector collecting them. The resulting matrices $\bm{A}$, $\bm{B}$, $\bm{D}$ and $\bm{A}^{\mathrm{s}}$ are the laminate stiffness matrices, and $\hat{\bm{N}}$ and $\hat{\bm{M}}$ are the thermal forces and moments per unit length, defined as~\cite{vescovini2020semi}:
\begin{equation}
\begin{aligned}
&\hat{\bm{N}} \plbr{x,y} = \sum_{k=1}^N \int_{t_k}^{t_{k+1}} \bar{\bm{Q}}_k \, \bar{\bm{\alpha}}_k \plbr{x,y}  \de z,
&\hat{\bm{M}} \plbr{x,y} = \sum_{k=1}^N \int_{t_k}^{t_{k+1}} z\, \bar{\bm{Q}}_k \, \bar{\bm{\alpha}}_k \plbr{x,y}  \de z.
    \label{eq:plate_thermal_components}
\end{aligned}
\end{equation}

\subsubsection{Energy terms}

Having defined the kinematics, the strain measure, and the constitutive relation, it is now convenient to derive the energy terms in the Principle of Virtual Work framework.
In particular, the internal virtual work expression reads:
\begin{equation}
\delta W_i = \int_\Omega \strain \plbr{\delta \bm{u}}^T \cost \plbr{x,y} \strain \plbr{\bm{u}} \de \Omega.
    \label{eq:int_virt_work_plate_linear}
\end{equation}

The buckling contribution is expressed as:
\begin{equation}
\delta W_b = \int_\Omega \strain_\buck \plbr{\delta \bm{u}}^T 
    \begin{bmatrix}
        N_{xx} & N_{xy}\\
        N_{xy} & N_{yy} 
    \end{bmatrix}
    \strain_\buck \plbr{\bm{u}} \de \Omega = \int_\Omega \strain_\buck \plbr{\delta \bm{u}}^T \mathbb{N} \plbr{x,y} \strain_\buck \plbr{\bm{u}} \de \Omega,
    \label{eq:buck_virt_work_plate}
\end{equation}
with:
\begin{equation}
\strain_\buck \plbr{\cdot} = 
\begin{bmatrix}
    0 & 0 & ,x & 0 & 0\\
    0 & 0 & ,y & 0 & 0
\end{bmatrix}.
    \label{eq:plate_epsilonbuck}
\end{equation}
For the free-vibration analysis, the contribution due to inertial forces is:
\begin{equation}
\begin{aligned}
\delta W_k &= \int_\Omega \delta \bm{u}^T \int_t \rho_0
\begin{bmatrix}
   \bm{I} & z\bm{L}\\
    z\bm{L} &  z^2\bm{L}
\end{bmatrix}
\de z \bm{\ddot{u}} \de \Omega
= \int_\Omega \delta \bm{u}^T 
\begin{bmatrix}
   \bm{I}_0 & \bm{I}_1\\
    \bm{I}_1 &  \bm{I}_2
\end{bmatrix}
\bm{\ddot{u}} \de \Omega= \int_\Omega \delta \bm{u}^T  \mass \bm{\ddot{u}} \de \Omega,
\end{aligned}
    \label{eq:inertial_work_plate}
\end{equation}
where $\bm{\ddot{u}}$ denotes the second time derivative of $\bm{u}$, $\bm{I}_0,\bm{I}_1,\bm{I}_2$ are the inertial moments and $\mass$ is the mass matrix of the plate.

The external virtual work accounts for body forces, boundary traction forces, and thermal loads, in the form:
\begin{equation}
\delta W_e = \int_\Omega \delta \bm{u}^T \bar{\bm{b}} \de \Omega + \sum_j \int_{\boundc_j} \delta \bm{u}^T\bar{\bm{t}}_j \de \boundc_j+ \int_\Omega \strain \plbr{\delta \bm{u}}^T \hat{\bm{R}} \de \Omega,
\label{eq:ext_virt_work_plate}
\end{equation}
where $\bar{\bm{b}}$ and $\bar{\bm{t}}_j$ are the prescribed body and the traction forces on the generic edge $\boundc_j$, respectively. Prescribed displacements are enforced as in standard finite element strategies.

\section{Virtual Element Method} \label{section:VEM} 

The numerical approximation is carried out via the Virtual Element Method (VEM). The method allows the use of arbitrary polygonal elements, which is a desirable feature to simplify the mesh generation on complex domains. Moreover, the VEM naturally handles hanging nodes, allowing the treatment of nonconforming discretizations as conforming ones. 

By following the approach in~\cite{daveiga2019virtual}, the standard VEM framework is extended here to account for curved edges, effectively enabling the treatment of arbitrary curved boundaries. Moreover, the $\kord$-version of VEM is considered, where $\kord \geq 1$ is the order of accuracy. The combined effect of curved edges representation capability and higher-order approximations is particularly effective to achieve simple yet accurate models, eliminating the need to operate geometry-induced mesh refinement strategies.

In the following, the VEM discretization, space, and associated degrees of freedom (DOFs), are recalled. Then, the discrete forms and corresponding projector operators of the energy contributions are presented, for both standard and $\mathrm{VC}$-VEM~\cite{foligno2026benchmarking}. The reader is referred to the \hyperref[chapter:appendix]{Supporting material} Section for more details on the steps to build the relevant matrices, as well as for the definitions of the different operators and quantities not otherwise specified.

\subsection{Discretization}

The domain $\surom$ is decomposed into a finite set $\globdom$ of non-overlapping star-shaped polygons with possibly curved edges $\el \in \globdom$~\cite{daveiga2013basic}, with $\diammax=\max_{\el \in \globdom} \diam$ and $\diam$ the diameter of $\el$. The boundary of $\el$ is denoted by $\boundel = \partial \el$, which can be partitioned into straight segments $\bound$ with $e=1,\dots,\nedges$ and curved edges $\boundcurv$ with $e=1,\dots,\nedgescurv$, the latter satisfying the regularity assumptions introduced earlier.

The generic continuous bilinear form $\bila \plbr{\cdot,\cdot}$ is now introduced. It is obtained as the sum of the contributions of the elements $\el$ in $\globdom$ as:
\begin{equation}
    \bila \plbr{\delta \disp,\disp} = \sum_{\el \in \globdom} \bilael\plbr{\delta \disp,\disp} \quad \forall \delta \disp,\disp \in \Vspacevec,
    \label{eq:continuos_bilinearform}
\end{equation}
where $\Vspacevec \equiv \sqbr{ \Vspace^u \times \Vspace^v \times \Vspace^w \times \Vspace^{\theta_x} \times \Vspace^{\theta_y} }$ is the continuous space in which the generalized displacement components lie.
The discretization uses the space 
$\Vspacedvec \equiv \sqbr{ \Vspaced^u \times \Vspaced^v \times \Vspaced^w \times \Vspaced^{\theta_x} \times \Vspaced^{\theta_y} }$, with
$\Vspacedvec \subset \Vspacevec$. Consequently, the discrete version of \eq{continuos_bilinearform} translates into:
\begin{equation}
    \bilad \plbr{\testdispd,\dispd} = \sum_{\el \in \globdom} \biladel \plbr{\testdispd,\dispd} \quad \forall \testdispd,\dispd \in \Vspacedvec,
    \label{eq:discrete_bilinearform}
\end{equation}
where $\dispd$ is the discrete counterpart of $\disp$.

\subsection{Virtual Element space and degrees of freedom}\label{subsection:VEMspacedofs}

The local virtual element spaces, on a generic polygon $\el$, with order of accuracy $\kord \geq 1$ and $m \geq \kord$, are hereafter defined.
For a generic $\phi$ representing a generalized displacement, the space is defined as:
\begin{equation}
\begin{aligned}
    \Vspaced^{r} \plbr{\el} = \left\{ \right. &{\phi}_h \in H^1 \plbr{\el} \cap C^0 \plbr{\el} : \bm{L} \sqbr{\cost \strain \plbr{{\phi}_h}} \big|_\el \in \sqbr{\Pspace_{\kord}\plbr{\el}}^5 , \\
    &\left. {\phi}_h \big|_{\bound} \in \Pspace_\kord \plbr{\bound} \quad\forall e=1,\dots,\nedges, \; {\phi}_h \big|_{\boundcurv} \in \Pspacecurv_\kord \plbr{\boundcurv} \quad \forall e=1,\dots,\nedgescurv, \right.\\
    &\left. \int_\el {\phi}_h \; \qpoly \; \de \el =  \int_\el {\PiKenhascalar}^\phi {\phi}_h \; \qpoly \; \de \el \quad\forall \qpoly \in \Pspace_{\kord} \plbr{\el} \setminus \Pspace_{\kord-r} \plbr{\el}
    \right\},
\end{aligned}
    \label{eq:localvirtual_spaces_general}
\end{equation}
where $\strain \plbr{{\phi}_h}$ is the strain operator obtained by considering only the corresponding displacement component. $\Pspace_{\kord} \plbr{\el}$ is the polynomial space of degree less than or equal to $\kord$ and $\Pspacecurv_\kord \plbr{\boundcurv}$ is the polynomial space on a curved edge, which, following~\cite{daveiga2019virtual}, can be expressed as:
\begin{equation}
    \Pspacecurv_\kord \plbr{\boundcurv} = \cubr{ \Tilde{\qpoly} = \qpoly \circ \gamma_e^{-1}: \qpoly \in \Pspace_\kord \plbr{I_e}}.
    \label{eq:curvedpoly_space}
\end{equation}
Hence, the functions on the curved edges are polynomials with respect to the parametrization $\gamma_e$.
The operator ${\PiKenhascalar}^\phi$ is the scalar counterpart, referred to the generic generalized displacement $\phi$, of the operator $\PiKenha$ that projects the trial functions from the VEM to the polynomial space, which will be defined later. The spaces for the different displacement components are then obtained as:
\begin{equation}
\Vspaced^{u} \plbr{\el} = \Vspaced^{v} \plbr{\el} = \Vspaced^{2} \plbr{\el}, \qquad \Vspaced^{w} \plbr{\el} = \Vspaced^{1} \plbr{\el}, \qquad \Vspaced^{\theta_x} \plbr{\el} = \Vspaced^{\theta_y} \plbr{\el} = \Vspaced^{0} \plbr{\el}.
\label{eq:localvirtual_spaces}
\end{equation}

The total local virtual element space for a generic element $\el$ of dimension $\Vspacedim$ reads:
\begin{equation}
    \Vspacedvec \plbr{\el} = \sqbr{ \Vspaced^u \times \Vspaced^v \times \Vspaced^w \times \Vspaced^{\theta_x} \times \Vspaced^{\theta_y} }.
    \label{eq:totlocalvirtual_space}
\end{equation}

The unknown $\dispd \in \Vspacedvec \plbr{\el}$ is uniquely identified by the following set of degrees of freedom ($\textbf{DOF}_i$):
\begin{itemize}
    \item $\textbf{DOF}_1$: the values of $\dispd$ at the $\nv$ vertices of $\el$,
    \item $\textbf{DOF}_2$: for $\kord > 1$, the values of $\dispd$ at the $\kord-1$ internal Gauss-Lobatto quadrature points on each straight edge $\bound$, and the values of $\dispd$ at the $\kord-1$ points on each curved edge $\boundcurv$, which are images through $\gamma_e$ of the $\kord-1$ internal Gauss-Lobatto quadrature points on $I_e$~\cite{daveiga2019virtual},
    \item $\textbf{DOF}_{3,4}$: for $\kord > 1$, the internal moments of $\Vspaced^{u} \plbr{\el}$ and $\Vspaced^{v} \plbr{\el}$ up to order $\kord-2$:
        \begin{equation}
            \frac{1}{\shbr{\el}} \int_\el \qpoly \plbr{x,y} u_h\plbr{x,y} \de \el \quad \forall \qpoly \in \Qpolyspace_{\kord-2}\plbr{\el}, \quad \frac{1}{\shbr{\el}} \int_\el \qpoly \plbr{x,y} v_h\plbr{x,y} \de \el \quad \forall \qpoly \in \Qpolyspace_{\kord-2}\plbr{\el}.
            \label{eq:internalmomentes_uv}
        \end{equation}   
        \item $\textbf{DOF}_5$: the internal moments of $\Vspaced^{w} \plbr{\el}$ up to order $\kord-1$:
        \begin{equation}
            \frac{1}{\shbr{\el}} \int_\el \qpoly \plbr{x,y} w_h\plbr{x,y} \de \el \quad \forall \qpoly \in \Qpolyspace_{\kord-1}\plbr{\el}.
            \label{eq:internalmomentes_w}
        \end{equation}
        \item $\textbf{DOF}_{6,7}$: the internal moments of $\Vspaced^{\theta_x} \plbr{\el}$ and $\Vspaced^{\theta_y} \plbr{\el}$ up to order $\kord$:
        \begin{equation}
            \frac{1}{\shbr{\el}} \int_\el \qpoly \plbr{x,y} {\theta_x}_h\plbr{x,y} \de \el \quad \forall \qpoly \in \Qpolyspace_{\kord}\plbr{\el}, \quad \frac{1}{\shbr{\el}} \int_\el \qpoly \plbr{x,y} {\theta_y}_h\plbr{x,y} \de \el \quad \forall \qpoly \in \Qpolyspace_{\kord}\plbr{\el},
            \label{eq:internalmomentes_xy}
        \end{equation}
\end{itemize}
where $\shbr{\el}$ is the area of the element $\el$ and $\Qpolyspace_{\kord}\plbr{\el}$ is a basis for the polynomial space $\Pspace_{\kord}\plbr{\el}$.
The displacements are polynomials of order $\kord$ on the edges. Conversely, in the element interior they differ for the FSDT formulation~\cite{daltri2022first}, and are associated to polynomials of degree $\kord-2$, $\kord-1$ and $\kord$ for the in-plane displacements, out-of-plane displacement and in-plane rotations, respectively.

\begin{figure}[!htbp]
    \centering
    \subfigure[$u$ and $v$.\label{fig:dof_k2_uv}]{
        \includegraphics[width=0.22\textwidth]{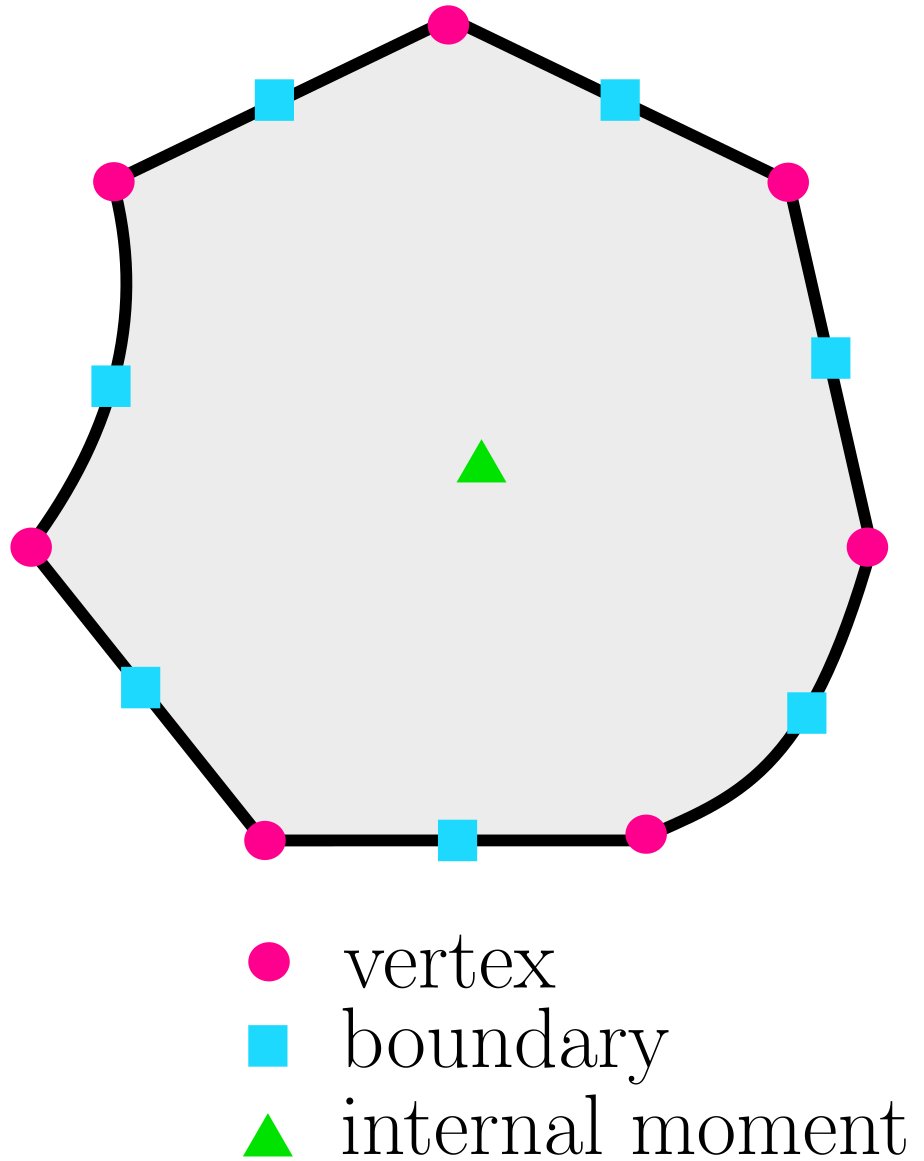}
    }
    \hfill
    \subfigure[$w$.\label{fig:dof_k2_w}]{
        \includegraphics[width=0.22\textwidth]{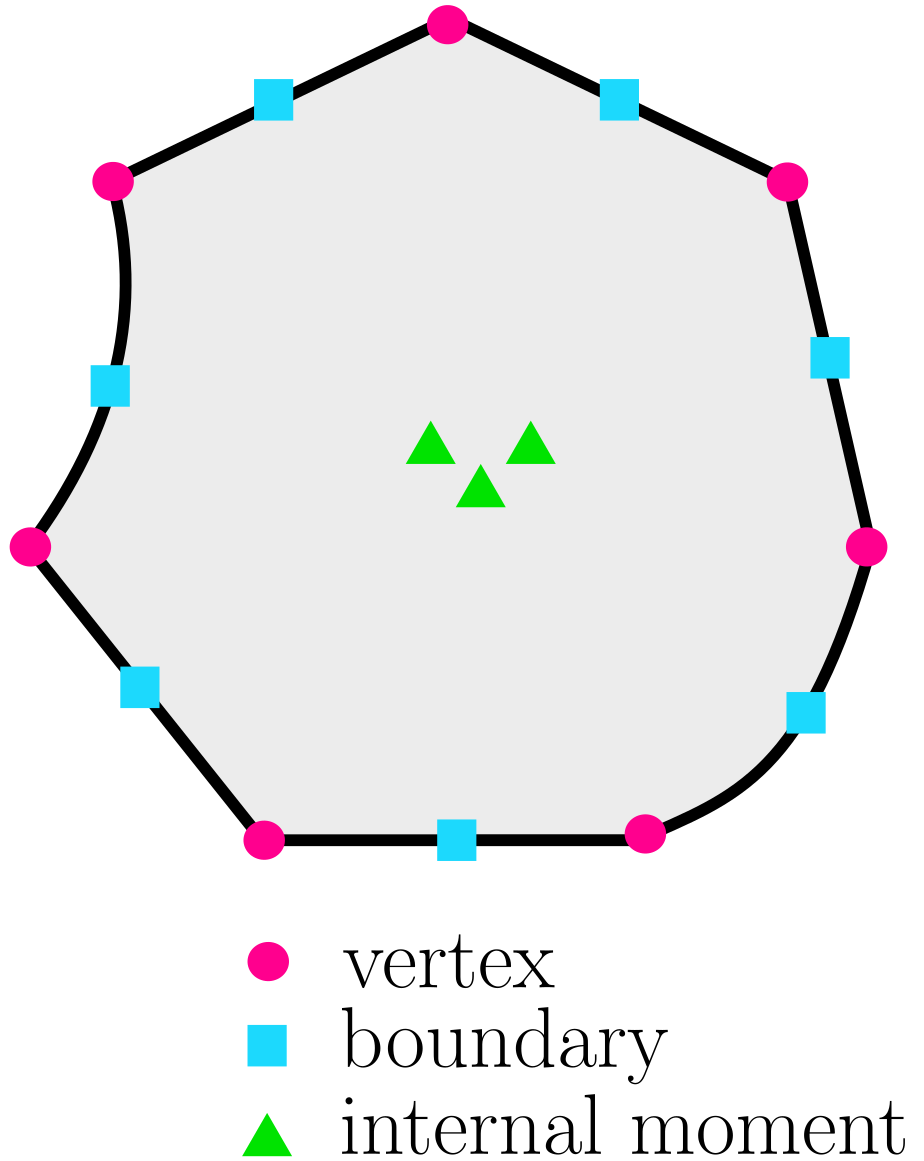}
    }
    \hfill
    \subfigure[$\theta_x$ and $\theta_y$.\label{fig:dof_k2_xy}]{
        \includegraphics[width=0.22\textwidth]{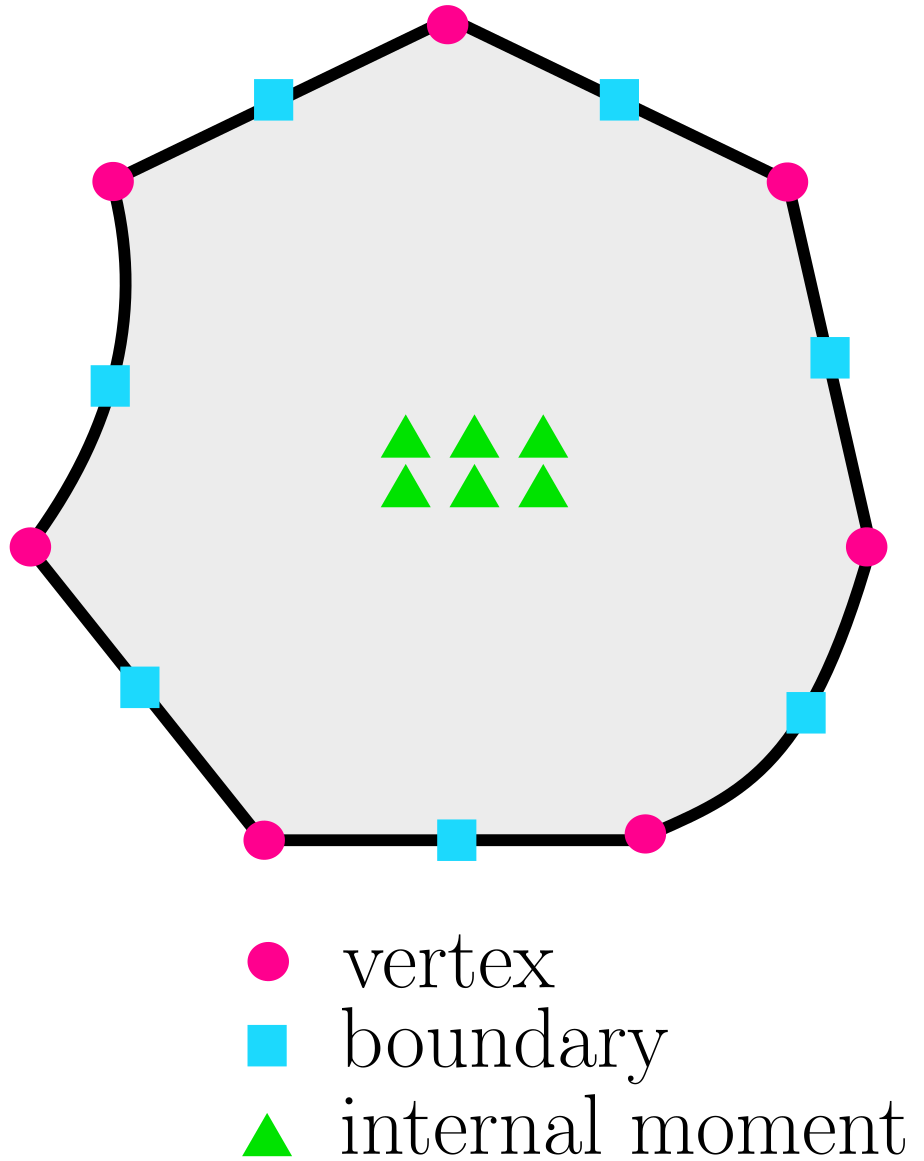}
    }
    \caption{Degrees of freedom for $\kord=2$.\label{fig:dof_k2}}
\end{figure}

A schematic representation of the degrees of freedom (DOFs) for $\kord=2$ is shown in \fig{dof_k2}, where circles and squares represent the vertex and edge DOFs, respectively, while triangles correspond to the internal moments.

The total global virtual element space is then obtained as:
\begin{equation}
    \Vspacedvec = \cubr{ \disp \in \sqbr{H_0^1\plbr{\surom}}^5 : \disp \in \Vspacedvec \plbr{\el} \quad \forall \el \in \globdom} = \sqbr{\Vspaced^u \times \Vspaced^v \times \Vspaced^w \times \Vspaced^{\theta_x} \times \Vspaced^{\theta_y} }.
    \label{eq:totglobalvirtual_space}
\end{equation}

The dependence of the constitutive law $\cost$ on the spatial coordinates $x,y$ is here neglected. Otherwise, the definition of the discrete spaces is inevitably more involved. The extension to this case is addressed at the end of this section.

\subsection{Discrete bilinear forms and projection operators} \label{subsection:VEMstandard}

As usual with VEM, the trial functions are solutions of a PDE inside each element and they are not explicitly computed.
Consequently, only a projection from the VEM onto the polynomial space of order $\kord$ is computable through the degrees of freedom:
\begin{equation}
    \bm{\Pi}_\el : \Vspacedvec \plbr{\el} \rightarrow \Pspacevec_\kord \plbr{\el}.
    \label{eq:projection_map}
\end{equation}

For the generic bilinear form, the projection is the solution of:
\begin{equation}
    \bilael\plbr{\dispd,\qpolyvec} = \bilael\plbr{\bm{\Pi}_\el\dispd,\qpolyvec} \quad \forall \qpolyvec \in \Pspacevec_\kord \plbr{\el}.
    \label{eq:projection_generic}
\end{equation}

The displacement field can be written using the Lagrangian-type interpolation identity as:
\begin{equation}
    \dispd = \sum_{i=1}^{\Vspacedim} \text{dof}_i \plbr{\dispd} \trialfcnvec_i \quad \forall \dispd \in \Vspacedvec \plbr{\el},
    \label{eq:lagrangian_interpolation}
\end{equation}
with $\text{dof}_j \plbr{\trialfcnvec_i} = \bm{\delta}_{ij}$ for all $i, j = 1,\dots, \Vspacedim$, where $\bm{\delta}_{ij}$ is the $5 \times 1 $ Kronecker-$\delta$ vector having unit value on the corresponding displacement component when $i=j$ and zeros elsewhere.
Substitution of \eq{lagrangian_interpolation} into \eq{projection_generic} yields:
\begin{equation}
    \bilael\plbr{\trialfcnvec_i,\qpolyvec} = \bilael\plbr{\bm{\Pi}_\el\trialfcnvec_i,\qpolyvec} \quad \forall \qpolyvec \in \Pspacevec_\kord \plbr{\el}.
    \label{eq:projection_generic_trialfunctions}
\end{equation}

In the following, the projections are defined for the stiffness matrix, the geometric stiffness matrix, the mass matrix, and the body and thermal force vectors. 

\subsubsection{Stiffness matrix}

The evaluation of the stiffness matrix is first presented by assuming that the constitutive tensor is constant within each element, with a value equal to its average at the integration points. This approximation facilitates the computation of the projection using the degrees of freedom. 
The generalization to spatially varying constitutive tensors $\cost \plbr{x,y}$ is outlined at the end of the section.

Accordingly, the bilinear form associated with the stiffness matrix is of elliptic type and is defined as:
\begin{equation}
    \begin{aligned}
        &\bilKenha \plbr{\trialfcnvec_i, \trialfcnvec_j}= \int_\el \strain\plbr{\trialfcnvec_i}^T \cost \strain\plbr{\trialfcnvec_j} \de \el,\\
    \end{aligned}
    \label{eq:stiffness_discretebilinearforms}
\end{equation}
and the projection operator $\PiKenha \trialfcnvec_i$ is the solution of:
\begin{equation}
\begin{cases}
    \begin{aligned}
        &\bilKenha \plbr{\trialfcnvec_i, \qpolyvec} = \bilKenha \plbr{\PiKenha\trialfcnvec_i, \qpolyvec} && \forall \qpolyvec \in \Pspacevec_\kord \plbr{\el}\\
        &P_0 \plbr{\trialfcnvec_i, \qpolyvec} = P_0 \plbr{\PiKenha\trialfcnvec_i, \qpolyvec} && \forall \qpolyvec \in \Pspacevec_0 \plbr{\el} \oplus \hat{\bm{P}} \plbr{\el},
    \end{aligned}
\end{cases}    
    \label{eq:stiffness_projections}
\end{equation}
for all $ i=1,\dots,\Vspacedim$; $\hat{\bm{P}} \plbr{\el}$ corresponds to $\begin{Bmatrix} x \; 0 \; 0 \; 0 \; 0 \end{Bmatrix}^T$ and $P_0 \plbr{\trialfcnvec_i, \qpolyvec}$, introduced to take care of the rigid body motions, is defined as:
\begin{equation}
    \begin{aligned}
        &P_0 \plbr{\trialfcnvec_i, \qpolyvec} = \frac{1}{\nv} \sum_{j=1}^{5\nv} \text{dof}_j \plbr{\trialfcnvec_i} \text{dof}_j \plbr{\qpolyvec}.
    \end{aligned}
    \label{eq:rigidbodymotions}
\end{equation}

The evaluation of the right-hand side of \eq{stiffness_projections} is straightforward as it involves the product of known polynomials. Regarding the left-hand side, integration by parts is applied. Subsequently, line integrals are easily computed as the trial functions are known on the boundary. The surface integrals are evaluated from the internal degrees of freedom.

Once the projection operator is available, following~\cite{daveiga2014hitchhiker}, a generic virtual function can be written as:
\begin{equation}
    \trialfcnvec_i = \PiKenha\trialfcnvec_i + \plbr{\bm{I}-\PiKenha}\trialfcnvec_i.
    \label{eq:VEM_function_polyandnot}
\end{equation}

Substituting \eq{VEM_function_polyandnot} into the stiffness bilinear form, and omitting the mixed terms, yields:
\begin{equation}
        \bilKVC \plbr{\trialfcnvec_i, \trialfcnvec_j} = \bilKVC \plbr{\PiKenha\trialfcnvec_i, \PiKenha\trialfcnvec_j} + \bilKVC \plbr{\plbr{\bm{I}-\PiKenha}\trialfcnvec_i, \plbr{\bm{I}-\PiKenha}\trialfcnvec_j},
    \label{eq:stiffness_bilinearform_cs}
\end{equation}
where:
\begin{equation}
    \begin{aligned}
        &\bilKVC \plbr{\trialfcnvec_i, \trialfcnvec_j}= \int_\el \strain\plbr{\trialfcnvec_i}^T \cost \plbr{x,y} \strain\plbr{\trialfcnvec_j} \de \el.\\
    \end{aligned}
    \label{eq:stiffness_discretebilinearforms_VC}
\end{equation}
In \eq{stiffness_bilinearform_cs}, if $\cost$ is constant within the elements, the two mixed terms vanish by the definition of $\PiKenha$; otherwise, they are omitted. 

The first term on the right-hand side of \eq{stiffness_bilinearform_cs} is the consistency term that ensures the accuracy of the solution. The second term is a contribution required to guarantee the stability of the solution; however, since the trial functions are not explicitly known inside the element, it cannot be computed explicitly.

\subsubsection*{Stabilized VEM}

Following~\cite{daveiga2013basic}, the stabilization term is a symmetric and semi-positive definite bilinear form defined in such a way that there exist two positive constants $\alpha_* \plbr{\kord}$ and $\alpha^* \plbr{\kord}$ independent of $\diammax$:
\begin{equation}
    \alpha_* \bilael\plbr{\testdispd,\testdispd} \leq \KEScalar \plbr{\testdispd,\testdispd} \leq \alpha^* \bilael\plbr{\testdispd,\testdispd} \quad \forall \testdispd \in \Vspacedvec \plbr{\el} \quad \text{with} \quad \PiKenha \testdispd = 0.
    \label{eq:stability_requirement}
\end{equation}

Consequently, the discrete version of \eq{stiffness_bilinearform_cs} is defined as:
\begin{equation}
        \bilKVC_h \plbr{\trialfcnvec_i, \trialfcnvec_j} = \bilKVC \plbr{\PiKenha\trialfcnvec_i, \PiKenha\trialfcnvec_j} + \KEScalar \plbr{\plbr{\bm{I}-\PiKenha}\trialfcnvec_i, \plbr{\bm{I}-\PiKenha}\trialfcnvec_j}.
    \label{eq:stiffness_bilinearform_cs_approx}
\end{equation}

In~\cite{foligno2026benchmarking}, several stabilization recipes from the literature have been compared for different PDEs and approximation orders to identify the most robust one. Relying upon the outcomes of the mentioned study, the stabilization term is selected as in~\cite{mascotto2018ill}.
In particular, the stabilization is defined block-wise according to the different components:
\begin{equation}
\left\{
\begin{aligned}
    & \KEScalar  =  \tau \sum_{r \in \mathcal{I}^{u,v}} \text{dof}_r \plbr{\plbr{\bm{I}-\PiKenha}\trialfcnvec_i} \cdot \plbr{s^E}_{r} \cdot \text{dof}_r \plbr{\plbr{\bm{I}-\PiKenha}\trialfcnvec_j}, && \text{if} \quad i,j \in \mathcal{I}^{u,v} \\
    & \KEScalar =  \tau \sum_{r \in \mathcal{I}^{w}} \text{dof}_r \plbr{\plbr{\bm{I}-\PiKenha}\trialfcnvec_i} \cdot \plbr{s^E}_{r} \cdot \text{dof}_r \plbr{\plbr{\bm{I}-\PiKenha}\trialfcnvec_j}, && \text{if} \quad i,j \in \mathcal{I}^{w} \\
    & \KEScalar =  \tau \sum_{r \in \mathcal{I}^{\theta_x,\theta_y}} \text{dof}_r \plbr{\plbr{\bm{I}-\PiKenha}\trialfcnvec_i} \cdot \plbr{s^E}_{r} \cdot \text{dof}_r \plbr{\plbr{\bm{I}-\PiKenha}\trialfcnvec_j}, && \text{if} \quad i,j \in \mathcal{I}^{\theta_x,\theta_y}.
\end{aligned}
\right.
    \label{eq:K_s_mascotto_diag}
\end{equation}
where $\mathcal{I}^{u,v}$, $\mathcal{I}^{w}$ and $\mathcal{I}^{\theta_x,\theta_y}$ refer to the subset of degrees of freedom associated to the in-plane displacements, out-of-plane displacements and in-plane rotations. Therefore, the stabilization term is block diagonal.

Moreover, $\tau$ is a user-defined parameter, taken as $0.5$ in linear elasticity~\cite{mengolini2019engineering}, and $\plbr{s^E}_{r} = \max \plbr{1, \bilKVC \plbr{\trialfcnvec_r, \trialfcnvec_r}}$, which scales the stabilization term as the diagonal of the consistency term. 

Consequently, in the standard stabilized VEM, the second term in \eq{stiffness_bilinearform_cs} is substituted with \eq{K_s_mascotto_diag}. Since this heuristic choice may spoil the accuracy of the method, self-stabilized formulations have been introduced in the literature. These formulations, by adopting higher-order polynomial projections, do not rely on the arbitrary stabilization term. Owing to the results in~\cite{foligno2026benchmarking,foligno2026novel}, in which several stabilized and self-stabilized formulations from the literature have been compared in terms of accuracy and conditioning for different classes of PDEs, one stabilized and one self-stabilized strategy are adopted in the present work, corresponding to the most robust and suitable for elasticity problems.

\subsubsection*{Self-stabilized VEM}

Self-stabilized VEM formulations do not rely on ad hoc stabilization terms. In contrast, they employ higher-order polynomial projections to ensure the stability of the linear system.

In~\cite{daveiga2016virtual,foligno2026benchmarking}, it was demonstrated that strategies based on $L^2$ projections are generally more robust in the presence of variable coefficients. Owing to these outcomes, just one self-stabilized formulation based on the $L^2$ projection is used in the present work. Specifically, the one introduced in~\cite{berrone2025lowest} for the lowest-order VEM in the context of the Laplace problem is here adapted to linear elasticity and generic order $\kord$.

In order to compute higher-order polynomial projections, the spaces in \eq{localvirtual_spaces} must be enlarged.
For a generic $\phi$, it holds that:
\begin{equation}
\begin{aligned}
    \VspacedenhV^{r} \plbr{\el} = \left\{ \right. &{\phi}_h \in H^1 \plbr{\el} \cap C^0 \plbr{\el} : \bm{L} \sqbr{\cost \strain \plbr{{\phi}_h}} \big|_\el \in \sqbr{\Pspace_{\kord+\laug}\plbr{\el}}^5 , \\
    &\left. {\phi}_h \big|_{\bound} \in \Pspace_\kord \plbr{\bound} \quad\forall e=1,\dots,\nedges, \; {\phi}_h \big|_{\boundcurv} \in \Pspacecurv_\kord \plbr{\boundcurv} \quad\forall e=1,\dots,\nedgescurv \right.,\\
    &\int_\el {\phi}_h \; \qpoly \; \de \el =  \int_\el {\PiKenhascalar}^\phi {\phi}_h \; \qpoly \; \de \el \quad\forall \qpoly \in \Pspace_{\kord} \plbr{\el} \setminus \Pspace_{\kord-r} \plbr{\el},\\
    &\left. \int_\el {\phi}_h \; \qpoly \; \de \el =  \int_\el {\Piok}^\phi {\phi}_h \; \qpoly \; \de \el \quad\forall \qpoly \in \Pspace_{\kord+\laug} \plbr{\el} \setminus \Pspace_{\kord} \plbr{\el} \right\},
\end{aligned}
    \label{eq:localvirtual_spaces_selfstabenhanced_V4_general}
\end{equation}
and the spaces for the different displacement components are obtained as:
\begin{equation}
\VspacedenhV^{u} \plbr{\el} = \VspacedenhV^{v} \plbr{\el} = \VspacedenhV^{2} \plbr{\el}, \qquad \VspacedenhV^{w} \plbr{\el} = \VspacedenhV^{1} \plbr{\el}, \qquad \VspacedenhV^{\theta_x} \plbr{\el} = \VspacedenhV^{\theta_y} \plbr{\el} = \VspacedenhV^0 \plbr{\el}.
\label{eq:localvirtual_spaces_selfstabenhanced_V4}
\end{equation}
The subscript $\laug$ denotes the required increase in polynomial order to obtain a non-singular system. The operators ${\PiKenhascalar}^\phi$ and ${\Piok}^\phi$ are the scalar versions of $\PiKenha$ and $\PiokK$ associated with the generic displacement component. 

The $L^2$ projection for the self-stabilized VEM reads:
\begin{equation}
    \begin{aligned}
        \plbr{\strain\plbr{\trialfcnvec_i}, \cost \qpolyvecunc}_\el = \plbr{ \PiKenhaSone\strain\plbr{\trialfcnvec_i}^T, \cost \qpolyvecunc}_\el    \quad \forall \qpolyvecunc \in \Pspacevecunc_{\kordaug-1} \plbr{\el},
    \end{aligned}  
    \label{eq:stiffness_projections_ss}
\end{equation}
for all $ i=1,\dots,\Vspacedim$, where $\qpolyvecunc$ and $\Pspacevecunc_{\kordaug-1} \plbr{\el}$ are the polynomial vector and polynomial space referred to the self-stabilized formulation. The projector operator $\PiKenhaSone$ maps the strains and not the functions themselves, as in the stabilized VEM, in the polynomial space.

The procedure to obtain the projection is the same as in standard stabilized VEM; hence, integration by parts is applied to the left-hand side. For what concerns line integrals, they are evaluated as in stabilized VEM. Conversely, surface integrals cannot be evaluated solely from the internal degrees of freedom, as in stabilized formulations. Therefore, the enlarged spaces in \eq{localvirtual_spaces_selfstabenhanced_V4} are used to evaluate higher-order polynomials. Once the projector is found, the bilinear form is evaluated as in \eq{stiffness_bilinearform_cs}, by considering only the consistency term and the self-stabilized projection operator.

To determine the augmented polynomial order $\laug$, two conditions are verified. First, the number of generalized strain modes must be greater than or equal to the number of degrees of freedom minus the number of rigid body motions. 
Second, a full-rank condition is enforced. 

With this purpose, an algorithm is implemented that iteratively increases $\laug$ until this condition is satisfied. More details can be found in~\cite{foligno2026benchmarking,foligno2026novel}.

\subsubsection{Remaining forms}

The projection procedure and the resulting discrete forms for the geometric stiffness and mass matrices, as well as for the body and thermal load vectors, follow the same construction adopted for the stiffness matrix. The stabilization term is required only for the stiffness matrix to ensure invertibility of the linear system. For all the other contributions, only the consistency term is retained. A summary of the discrete forms is provided in \tab{summary_discretebilinearforms}.

\begin{sidewaystable}[p]
\centering
{\footnotesize
\renewcommand{\arraystretch}{1.20}
\setlength{\tabcolsep}{4pt}
\begin{tabularx}{\textheight}{
                               >{\raggedright\arraybackslash}p{6.5cm}
                               >{\raggedright\arraybackslash}p{7.5cm}
                               >{\raggedright\arraybackslash}p{6.5cm}}
\toprule
\textbf{Continuous form} & \textbf{Projection} & \textbf{Discrete form} \\
\midrule

\multicolumn{3}{l}{\textit{Stiffness matrix (stabilized VEM)}}\\
\cline{1-1}
\parbox[t]{\linewidth}{
$\begin{aligned}
\bilKenha \plbr{\trialfcnvec_i, \trialfcnvec_j}= \int_\el \strain\plbr{\trialfcnvec_i}^T \cost \strain\plbr{\trialfcnvec_j} \de \el
\end{aligned}$

$\begin{aligned}
\bilKVC \plbr{\trialfcnvec_i, \trialfcnvec_j}= \int_\el \strain\plbr{\trialfcnvec_i}^T \cost \plbr{x,y} \strain\plbr{\trialfcnvec_j} \de \el
\end{aligned}$
}
&
for all $ i=1,\dots,\Vspacedim$:
\parbox[t]{\linewidth}{
$\left\{
\begin{aligned}
&\bilKenha \plbr{\trialfcnvec_i, \qpolyvec} = \bilKenha \plbr{\PiKenha\trialfcnvec_i, \qpolyvec} && \forall \qpolyvec \in \Pspacevec_\kord \plbr{\el}\\
&P_0 \plbr{\trialfcnvec_i, \qpolyvec} = P_0 \plbr{\PiKenha\trialfcnvec_i, \qpolyvec} && \forall \qpolyvec \in \Pspacevec_0 \plbr{\el} \oplus \hat{\bm{P}} \plbr{\el}
\end{aligned}
\right.$
}
&
$\aligned \bilKVC_h &\plbr{\trialfcnvec_i, \trialfcnvec_j} =\bilKVC \plbr{\PiKenha\trialfcnvec_i, \PiKenha\trialfcnvec_j}\\ &+ \KEScalar \plbr{\plbr{\bm{I}-\PiKenha}\trialfcnvec_i, \plbr{\bm{I}-\PiKenha}\trialfcnvec_j} \endaligned$
\\

\midrule

\multicolumn{3}{l}{\textit{Stiffness matrix (self-stabilized VEM)}}\\
\cline{1-1}
\parbox[t]{\linewidth}{

$\begin{aligned}
\bilKVC \plbr{\trialfcnvec_i, \trialfcnvec_j}= \int_\el \strain\plbr{\trialfcnvec_i}^T \cost \plbr{x,y} \strain\plbr{\trialfcnvec_j} \de \el
\end{aligned}$
}
&
for all $ i=1,\dots,\Vspacedim$:
\parbox[t]{\linewidth}{
$\begin{aligned}
\plbr{\strain\plbr{\trialfcnvec_i}, \cost \qpolyvecunc}_\el = \plbr{ \PiKenhaSone\strain\plbr{\trialfcnvec_i}^T, \cost \qpolyvecunc}_\el
\end{aligned}$

$\begin{aligned}
\forall \qpolyvecunc \in \Pspacevecunc_{\kordaug-1} \plbr{\el}
\end{aligned}$
}
&
$\bilKVC_h \plbr{\trialfcnvec_i, \trialfcnvec_j} = \plbr{ \PiKenhaSone\strain\plbr{\trialfcnvec_i}^T, \cost \plbr{x,y}\PiKenhaSone\strain\plbr{\trialfcnvec_j}}_\el$
\\

\midrule

\multicolumn{3}{l}{\textit{Geometric stiffness matrix}}\\
\cline{1-1}
\parbox[t]{\linewidth}{
$\begin{aligned}
\bilGenha \plbr{\trialfcn_i, \trialfcn_j}= \int_\el \strain_\buck\plbr{\trialfcn_i}^T \mathbb{N} \strain_\buck\plbr{\trialfcn_j} \de \el
\end{aligned}$

$\begin{aligned}
\bilGVC \plbr{\trialfcn_i, \trialfcn_j}= \int_\el \strain_\buck\plbr{\trialfcn_i}^T \mathbb{N} \plbr{x,y} \strain_\buck\plbr{\trialfcn_j} \de \el
\end{aligned}$
}
&
for all $ i=1,\dots,\Vspacedim^w$:
\parbox[t]{\linewidth}{
$\left\{
\begin{aligned}
&\bilGenha \plbr{\trialfcn_i, \qpoly} = \bilGenha \plbr{\PiGenha \trialfcn_i, \qpoly} && \forall \qpoly \in \Pspace_\kord \plbr{\el}\\
&P_0 \plbr{\trialfcn_i, \qpoly} = P_0 \plbr{\PiGenha\trialfcn_i, \qpoly} && \forall \qpoly \in \Pspace_0 \plbr{\el}
\end{aligned}
\right.$
}
&
$\bilGVC_h \plbr{\trialfcn_i, \trialfcn_j} = \bilGVC \plbr{\PiGenha\trialfcn_i, \PiGenha\trialfcn_j}$
\\

\midrule

\multicolumn{3}{l}{\textit{Mass matrix}}\\
\cline{1-1}
$\begin{aligned}
\bilMenha \plbr{\trialfcnvec_i, \trialfcnvec_j}= \int_\el \trialfcnvec_i^T \mathbb{M} \trialfcnvec_j \de \el
\end{aligned}$
&
for all $ i=1,\dots,\Vspacedim$:
$\begin{aligned}
&\bilMenha \plbr{\trialfcnvec_i, \qpolyvec} = \bilMenha \plbr{\PiMenha\trialfcnvec_i, \qpolyvec} && \forall \qpolyvec \in \Pspacevec_\kord \plbr{\el}
\end{aligned}$
&
$\begin{aligned}
\bilMenha_h \plbr{\trialfcnvec_i, \trialfcnvec_j} =\bilMenha \plbr{\PiMenha\trialfcnvec_i, \PiMenha\trialfcnvec_j}
\end{aligned}$
\\

\midrule

\multicolumn{3}{l}{\textit{Body forces vector}}\\
\cline{1-1}
$\bilF \plbr{\trialfcnvec_i} = \int_\el \trialfcnvec_i^T \bodyforcevec \de \el$
&
for all $ i=1,\dots,\Vspacedim$:
\parbox[t]{\linewidth}{
$\plbr{\trialfcnvec_i^T, \qpolyvec}_\el = \plbr{\PiokK\trialfcnvec_i^T, \qpolyvec}_\el \quad \forall \qpolyvec \in \Pspacevec_\kord \plbr{\el}$
}
&
$\bilF_h \plbr{\trialfcnvec_i} = \bilF \plbr{\PiokK\trialfcnvec_i}$
\\

\midrule

\multicolumn{3}{l}{\textit{Thermal forces vector}}\\
\cline{1-1}
\parbox[t]{\linewidth}{
$\begin{aligned}
\bilTenha \plbr{\trialfcnvecther_i, \trialfcnvecther_j}= \int_\el \strain_\ther\plbr{\trialfcnvecther_i}^T \check{\bm{R}} \strain_\ther\plbr{\trialfcnvecther_j} \de \el
\end{aligned}$

$\begin{aligned}
\bilTVC \plbr{\trialfcnvecther_i}= \int_\el \check{\bm{R}} \plbr{x,y}\strain_\ther\plbr{\trialfcnvecther_i} \de \el
\end{aligned}$
}
&
for all $ i=1,\dots,\Vspacedim^*$:
\parbox[t]{\linewidth}{
$\left\{
\begin{aligned}
&\bilTenha \plbr{\trialfcnvecther_i, \qpolyvecther} = \bilTenha \plbr{\PiTenha\trialfcnvecther_i, \qpolyvecther} && \forall \qpolyvecther \in \Pspacevecther_\kord \plbr{\el}\\
&P_0 \plbr{\trialfcnvecther_i, \qpolyvecther} = P_0 \plbr{\PiTenha\trialfcnvecther_i, \qpolyvecther} && \forall \qpolyvecther \in \Pspacevecther_0 \plbr{\el} \oplus \hat{\bm{P}}^* \plbr{\el}
\end{aligned}
\right.$
}
&
$\bilTVC_h \plbr{\trialfcnvecther_i} = \bilTVC \plbr{\PiTenha\trialfcnvecther_i}$
\\

\bottomrule
\end{tabularx}
}
\caption{Summary of remaining forms.}
\label{tab:summary_discretebilinearforms}
\end{sidewaystable}

The projection is of elliptic type for the stabilized stiffness matrix, the geometric stiffness matrix and the thermal load vector, and of $L^2$ type for the self-stabilized stiffness matrix, the mass matrix and the body force vector. For the geometric stiffness matrix, only the contributions associated with the out-of-plane displacement $w$ are retained, as the remaining terms are neglected in the buckling formulation and would otherwise lead to a singular system.

For the mass matrix, the density is assumed constant. So, the inertial coefficient matrix $\mass$ is independent of $x$ and $y$. For the thermal load vector, the out-of-plane displacement $w$ is excluded, as it is assumed not to contribute to thermal effects. 

The line load vector is constructed following the standard finite element procedure, since the trial functions are polynomials, or polynomial images, along the element edges. Prescribed displacements are enforced at the assembled system level, as in standard finite elements.

\subsection{Discrete forms and projection operators with variable coefficients} \label{subsection:VEMVC}

In variable stiffness laminates, the elastic properties are not constant over the domain, so proper handling is required for the projection operators. In the previous section, the elastic coefficients were approximated as a constant within each element, which is appropriate for low-order approaches based on $h$-refinement, but is, in general, not suitable within a $\kord$-refinement framework.

To overcome this issue, the Variable Coefficients-VEM approach ($\mathrm{VC}$-VEM) proposed by the authors in the recent work~\cite{foligno2026benchmarking} is here employed.
As opposed to standard VEM practice, the coefficient variability is directly included in the projection, and non-computable terms are approximated using the $L^2$ VEM projector $\PiokK$. 
The $\mathrm{VC}$-VEM applies to the forms with a non-constant coefficient, i.e. stiffness matrix, geometric stiffness matrix, and thermal forces vector. A summary of the corresponding projections is provided in \tab{vcvem_projections}.

\begin{table}[!htbp]
\centering
{\small
\renewcommand{\arraystretch}{1.2}
\setlength{\tabcolsep}{4pt}
\begin{tabularx}{\columnwidth}{>{\raggedright\arraybackslash}p{2.2cm}
                               >{\raggedright\arraybackslash}X}
\toprule
\textbf{Term} & \textbf{Projection} \\
\midrule

Stiffness matrix (stabilized)
&
for all $ i=1,\dots,\Vspacedim$:
\parbox[t]{\linewidth}{
$\left\{
\begin{aligned}
&\bilKVCstar \plbr{\trialfcnvec_i, \qpolyvec} = \bilKVC \plbr{\PiKVC\trialfcnvec_i, \qpolyvec} && \forall \qpolyvec \in \Pspacevec_\kord \plbr{\el}\\
&P_0 \plbr{\trialfcnvec_i, \qpolyvec} = P_0 \plbr{\PiKVC\trialfcnvec_i, \qpolyvec} && \forall \qpolyvec \in \Pspacevec_0 \plbr{\el} \oplus \hat{\bm{P}} \plbr{\el}
\end{aligned}
\right.$
}
\\

\midrule

Stiffness matrix (self-stabilized)
&
for all $ i=1,\dots,\Vspacedim$:
$\begin{aligned}
\plbr{\strain^*\plbr{\trialfcnvec_i}, \cost \plbr{x,y} \qpolyvecunc}_{*,\el}
=
\plbr{ \PiKVCSone\strain\plbr{\trialfcnvec_i}^T, \cost \plbr{x,y}\qpolyvecunc}_\el
\quad
\forall \qpolyvecunc \in \Pspacevecunc_{\kordaug-1} \plbr{\el}
\end{aligned}$
\\

\midrule

Geometric stiffness matrix
&
for all $ i=1,\dots,\Vspacedim^w$:
\parbox[t]{\linewidth}{
$\left\{
\begin{aligned}
&\bilGVCstar \plbr{\trialfcn_i, \qpoly} = \bilGVC \plbr{\PiGVC \trialfcn_i, \qpoly} && \forall \qpoly \in \Pspace_\kord \plbr{\el}\\
&P_0 \plbr{\trialfcn_i, \qpoly} = P_0 \plbr{\PiGVC\trialfcn_i, \qpoly} && \forall \qpoly \in \Pspace_0 \plbr{\el}
\end{aligned}
\right.$
}
\\

\midrule

Thermal forces vector
&
for all $ i=1,\dots,\Vspacedim^*$:
\parbox[t]{\linewidth}{
$\left\{
\begin{aligned}
&\bilTVCstar \plbr{\trialfcnvecther_i, \qpolyvecther} = \bilTVC \plbr{\PiTVC\trialfcnvecther_i, \qpolyvecther} && \forall \qpolyvecther \in \Pspacevecther_\kord \plbr{\el}\\
&P_0 \plbr{\trialfcnvecther_i, \qpolyvecther} = P_0 \plbr{\PiTVC\trialfcnvecther_i, \qpolyvecther} && \forall \qpolyvecther \in \Pspacevecther_0 \plbr{\el} \oplus \hat{\bm{P}}^* \plbr{\el}
\end{aligned}
\right.$
}
\\

\bottomrule
\end{tabularx}
}
\caption{Projection for VC-VEM.}
\label{tab:vcvem_projections}
\end{table}

For the stiffness matrix, the procedure described above is extended to the variable stiffness case, i.e., spatially varying constitutive laws $\cost\plbr{x,y}$.
The idea behind $\mathrm{VC}$-VEM is to compute the polynomial projection $\PiKVC$ as the solution of the following problem, which directly involves the bilinear form:
\begin{equation}\label{eq:vc_ideal}
    \left\{
\begin{aligned}
&\bilKVC \plbr{\trialfcnvec_i, \qpolyvec} = \bilKVC \plbr{\PiKVC\trialfcnvec_i, \qpolyvec} && \forall \qpolyvec \in \Pspacevec_\kord \plbr{\el}\\
&P_0 \plbr{\trialfcnvec_i, \qpolyvec} = P_0 \plbr{\PiKVC\trialfcnvec_i, \qpolyvec} && \forall \qpolyvec \in \Pspacevec_0 \plbr{\el} \oplus \hat{\bm{P}} \plbr{\el}.
\end{aligned}
\right.
\end{equation}
However, while the right-hand side is fully computable as it involves only polynomials (and the coefficients), the left-hand side is not. Indeed, integration by parts yields the following identity:
\begin{equation}
    \begin{aligned}
    \bilKVC \plbr{\trialfcnvec_i, \qpolyvec} = &- \int_\el \trialfcnvec_i \cdot \bm{L}_d \sqbr{\cost \plbr{x,y} \strain\plbr{\qpolyvec}}  \de \el + \int_{\boundel} \trialfcnvec_i \cdot  \hat{\stress}\plbr{\qpolyvec} \normedgeel  \de \boundel +\\
    &+ \int_\el \strain^\const\plbr{\trialfcnvec_i}^T \cost \plbr{x,y} \strain\plbr{\qpolyvec} \de \el,
    \end{aligned}
    \label{eq:stiffness_VC_projection}
\end{equation}
where it is clear that the surface integrals are not computable through the degrees of freedom. \eq{vc_ideal} is then modified as follows:
\begin{equation}
    \left\{
\begin{aligned}
&\bilKVCstar \plbr{\trialfcnvec_i, \qpolyvec} = \bilKVC \plbr{\PiKVC\trialfcnvec_i, \qpolyvec} && \forall \qpolyvec \in \Pspacevec_\kord \plbr{\el}\\
&P_0 \plbr{\trialfcnvec_i, \qpolyvec} = P_0 \plbr{\PiKVC\trialfcnvec_i, \qpolyvec} && \forall \qpolyvec \in \Pspacevec_0 \plbr{\el} \oplus \hat{\bm{P}} \plbr{\el},
\end{aligned}
\right.
\end{equation}
where the left-hand side is the \textit{modified} form defined as:
\begin{equation}
    \begin{aligned}
    \bilKVCstar \plbr{\trialfcnvec_i, \qpolyvec} := &- \int_\el \PiokK \trialfcnvec_i \cdot \bm{L}_d \sqbr{\cost \plbr{x,y} \strain\plbr{\qpolyvec}}  \de \el + \int_{\boundel} \trialfcnvec_i \cdot  \hat{\stress}\plbr{\qpolyvec} \normedgeel  \de \boundel +\\
    &+ \int_\el \strain^\const\plbr{\PiokK\trialfcnvec_i}^T \cost \plbr{x,y} \strain\plbr{\qpolyvec} \de \el.
    \end{aligned}
    \label{eq:stiffness_VC_projection_2}
\end{equation}
$\bilKVCstar$ is a computable approximation of the original bilinear form through the $L^2$ projection operator $\PiokK$, and $\hat{\stress}$ contains the spatial variation of $\cost \plbr{x,y}$. The evaluation of the line integral does not require special care, although the non-polynomial nature of  $\cost \plbr{x,y}$ leads to a non-exact integration. 

As implied by \eq{stiffness_VC_projection_2}, the coefficient $\cost \plbr{x,y}$ should be at least of class $\mathcal{C}^1$ within each element. This is a suitable assumption for the problems under consideration, as the constitutive law $\cost \plbr{x,y}$ is obtained through standard laminate stiffness matrices assembly and the fiber orientations are interpolated via Lagrange polynomials.

The same strategy applies to the self-stabilized VEM formulation, as well as to the geometric stiffness matrix and the thermal load vector. The related $\mathrm{VC}$-VEM projections are collected in \tab{vcvem_projections}. The subscript ``$*$'' denotes the \textit{modified} forms defined by replacing the trial functions with their polynomial projection in the non-computable integrals, after having applied integration by parts. In addition, $\mathbb{N}\plbr{x,y}$ and $\check{\bm{R}}\plbr{x,y}$ are required to be of class $\mathcal{C}^1$ within each element, a condition easily satisfied for the class of problems considered in this investigation.

The spatial variation of $\mathbb{N}\plbr{x,y}$ depends not only on the constitutive law $\cost \plbr{x,y}$, but also on the strain field. Therefore, even for constant fiber orientation, $\mathbb{N}\plbr{x,y}$ remains spatially dependent.

The final discrete bilinear forms are constructed following the standard VEM formulation, with the stabilization term of the stiffness matrix defined using the standard projector $\PiKenha$, which does not include the spatial variation of the constitutive tensor in the projection.

\section{Results} \label{section:results}

In this section, the accuracy and robustness of the proposed VEM formulation are assessed by comparison with analytic solutions and benchmark problems from the literature, as well as with commercial finite element simulations conducted using Abaqus. First, the convergence of the method is assessed through three representative test cases in a $\kord$-refinement setting. Specifically, a first test case serves to investigate the performance of standard and $\mathrm{VC}$-VEM, in their stabilized and self-stabilized versions, in the presence of different layups, ranging from constant to curvilinear fiber orientations. Moreover, the static and free-vibration responses of problems featuring high-gradient solutions are investigated to demonstrate the effectiveness of performing local mesh refinements. Lastly, the method is validated against results from the literature and Abaqus for both free-vibration and buckling analyses. Geometries of varying complexity are considered to illustrate the ability of the method to handle general geometries.
Different meshes and approximation orders are considered throughout the section. In all the cases, the number of integration points is selected to ensure sufficient accuracy of the results. For each element, all curve types are interpolated as B\'ezier curves, which are selected in this work for their robustness and simplicity of implementation.

\subsection{Test case 1}

The first test case investigates the convergence properties of the standard VEM and its $\mathrm{VC}$-version in a $\kord$-refinement framework. In particular, this analysis extends the results of~\cite{foligno2026benchmarking} to a more complex domain with a cutout, curvilinear fiber orientations, and both membrane and bending behavior.

\begin{figure}[!htbp]
    \centering
        \includegraphics[width=0.35\textwidth]{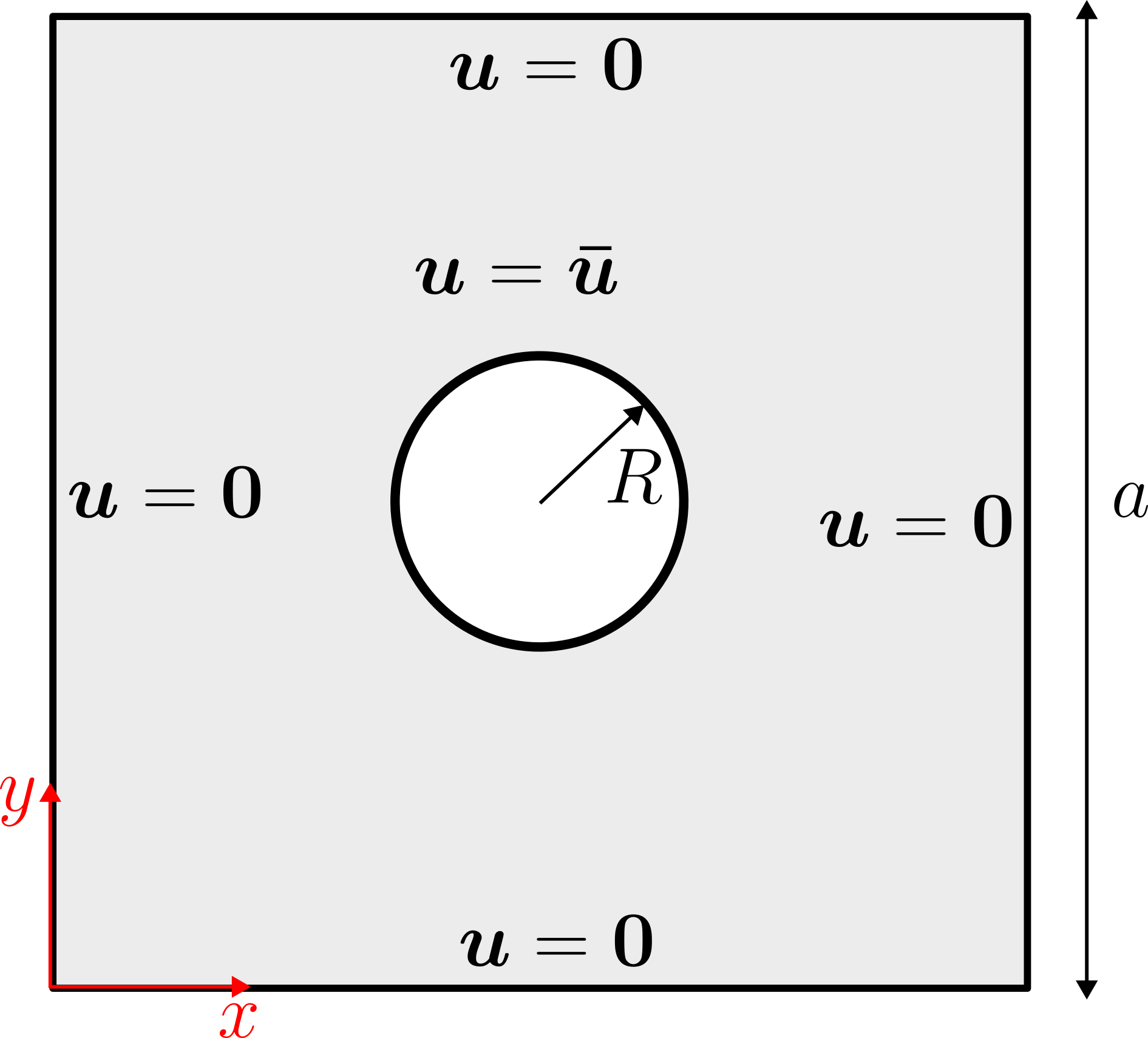}
    \caption{Test case 1: configuration. \label{fig:testcase0_configuration}}
\end{figure}

The plate is square with side $a=1000$~mm, and a circular cutout of radius $R=150$~mm is centered in the middle of the plate. A schematic representation is available in \fig{testcase0_configuration}, where the essential boundary conditions, prescribed to all the displacement components $\bm{u} = \{ u \, v \, w \, \theta_x \, \theta_y \}^{\mathrm{T}}$, are reported. 

An orthotropic material with properties $E_{11}=181000$~MPa, $E_{22}=10273$~MPa, $G_{12}=G_{13}=G_{23}=7170.5$~MPa and $\nu_{12}=0.28$ is considered. The laminate has sixteen plies with thickness $h=0.1272$~mm each. Two different layups are investigated:
\begin{equation}
    \begin{aligned}
     \text{layup 1}: [\pm90\langle 45|45\rangle]_{4s}, \qquad
     \text{layup 2}: [\langle 90|0\rangle_8, \langle -90|0\rangle_8].
    \end{aligned}
    \label{eq:testcase1_layups}
\end{equation}

These layups feature an increasing level of complexity. The first one displays uniform properties over the domain, while the second accounts for variable properties through a linear variation of the fiber angle and adds membrane-bending coupling due to the asymmetry of the stacking sequence. 

The errors of the VEM solution are evaluated with respect to the exact solution of the problem, formulated via an in inverse approach: the exact solution is postulated in advance and the corresponding loading conditions are retrieved. In particular, the solution is imposed to be: 
\begin{equation}
    \disp \plbr{x,y} = 
    \sin\plbr{\pi x} \sin\plbr{\pi y}
    \begin{Bmatrix}
    1 &
    1 &
    1 &
    1 &
    1
    \end{Bmatrix}^T,
    \label{eq:testcase1_exactsolution}
\end{equation}
where fixed boundary conditions are assumed along the outer edges.

Error estimates are conducted using the energy norm:
\begin{equation}
\begin{aligned}
    e_{\bm{\varepsilon}}^{\disp} =  \frac{ \plbr{\sum_{\el \in \globdom}  \norm{ \sqrt{\cost\plbr{x,y}} \strain \plbr{ \disp - \PiokK \dispd } }_{0,\el}^2  }^{\frac{1}{2}} }{ \norm{ \sqrt{\cost\plbr{x,y}} \strain \plbr{\disp}}_{0,\el} }.
\end{aligned}    
    \label{eq:linearelasticity_l2error_dispenergy}
\end{equation}
 
A $\kord$-refinement strategy is adopted, with $\kord=1,\dots,10$. Both standard and $\mathrm{VC}$-VEM, as well as stabilized and self-stabilized VEM, are employed, thereby allowing a comprehensive comparison of the different strategies in the presence of curved edges, high approximation orders, and complex layup configurations.

\begin{figure}[!htbp]
        \centering
	\subfigure[layup 1. \label{fig:testcase1_layup1_quad}]{
		\includegraphics[width=0.47\textwidth]{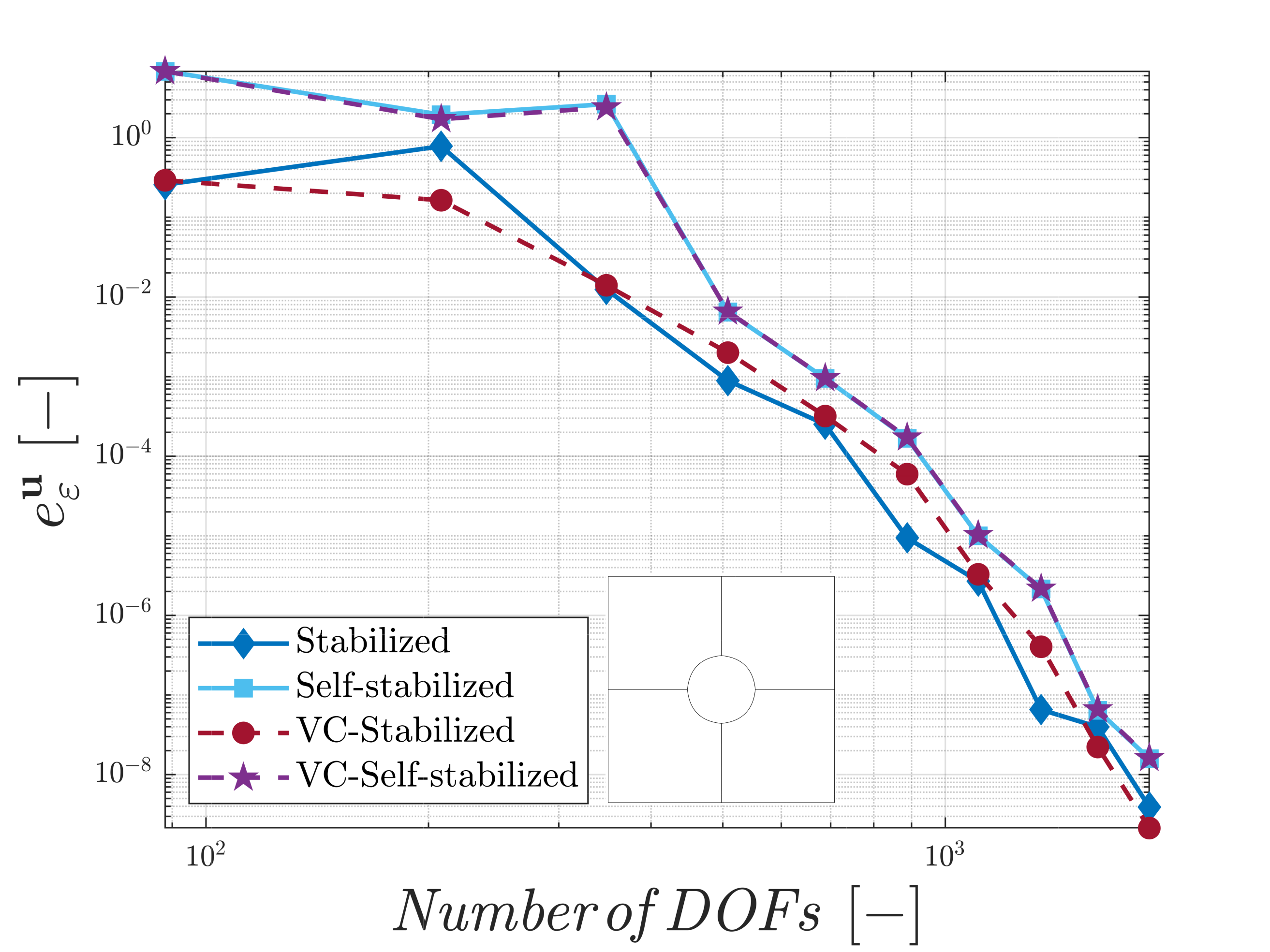}
	}
        \centering
	\subfigure[layup 2. \label{fig:testcase1_layup3_quad}]{
		\includegraphics[width=0.47\textwidth]{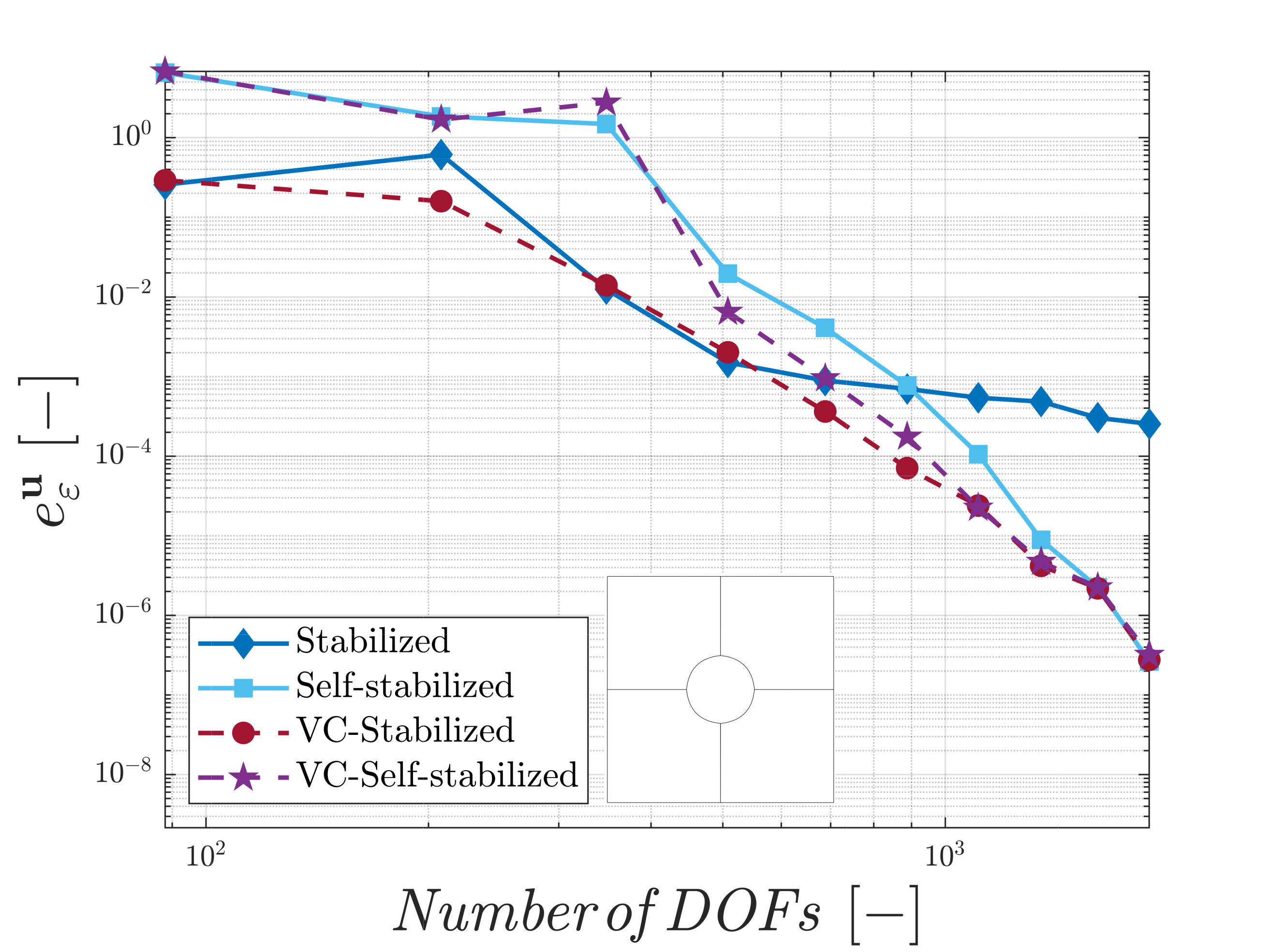}
	}
    \caption{Test case 1: energy norm error.\label{fig:testcase1_quad}}
 \end{figure}

The error estimates are reported in \fig{testcase1_quad} for the two layups at hand. The mesh, which features curved edges, is plotted in the same figures next to the legend.

For the first layup, \fig{testcase1_layup1_quad}, the results demonstrate equal accuracy across the different strategies. In contrast, when the fiber orientation is no longer constant (layup 2 in \fig{testcase1_layup3_quad}), the standard stabilized VEM loses its accuracy as the order $\kord$ increases. Instead, both the stabilized and self-stabilized $\mathrm{VC}$-VEM, and the standard self-stabilized VEM, which adopts a $L^2$ projection, maintain a good level of accuracy, despite the increased complexity given by the non-uniform stiffness distribution and membrane-bending coupling. Across all layups, some oscillations are present at lower orders; while, as the order increases, a smoother trend is observed. For layup 2, the error is slightly higher than for layup 1. This behavior is ascribed to the complexity introduced by the presence of the curvilinear fiber orientations.

These findings are in agreement with those presented in~\cite{daveiga2016virtual}, in which it is shown that $L^2$ projections typically perform better in the presence of variable coefficients, and with those in~\cite{foligno2026benchmarking}, in which it is demonstrated that $\mathrm{VC}$-VEM maintains optimal accuracy even with elliptic projections. The present test case broadens the ones analyzed in~\cite{foligno2026benchmarking} to a more complex scenario. Indeed, while~\cite{foligno2026benchmarking} investigates the performance of $\mathrm{VC}$-VEM by analyzing the membrane behavior of a square domain with a polynomial variation of the elastic properties, the present work assesses the effectiveness of the method by considering a plate with a cutout, curvilinear fiber orientations, and membrane-bending coupling.

\subsection{Test case 2}

The second test case regards a cracked panel under tension load, with a singular stress state at the crack tip. This is a classical benchmark to assess the convergence properties of numerical methods, as the presence of the singularity makes the convergence particularly challenging~\cite{zander2017multi}. In this work, this test case is of interest to demonstrate the potential of VEM to perform mesh refinements where desired, thereby enhancing the convergence properties in the presence of singularities. Hanging nodes, i.e. nodes located along the edges of an element, are used to preserve mesh conformity in the refined areas.   

\begin{figure}[!htbp]
    \centering
        \includegraphics[width=0.65\textwidth]{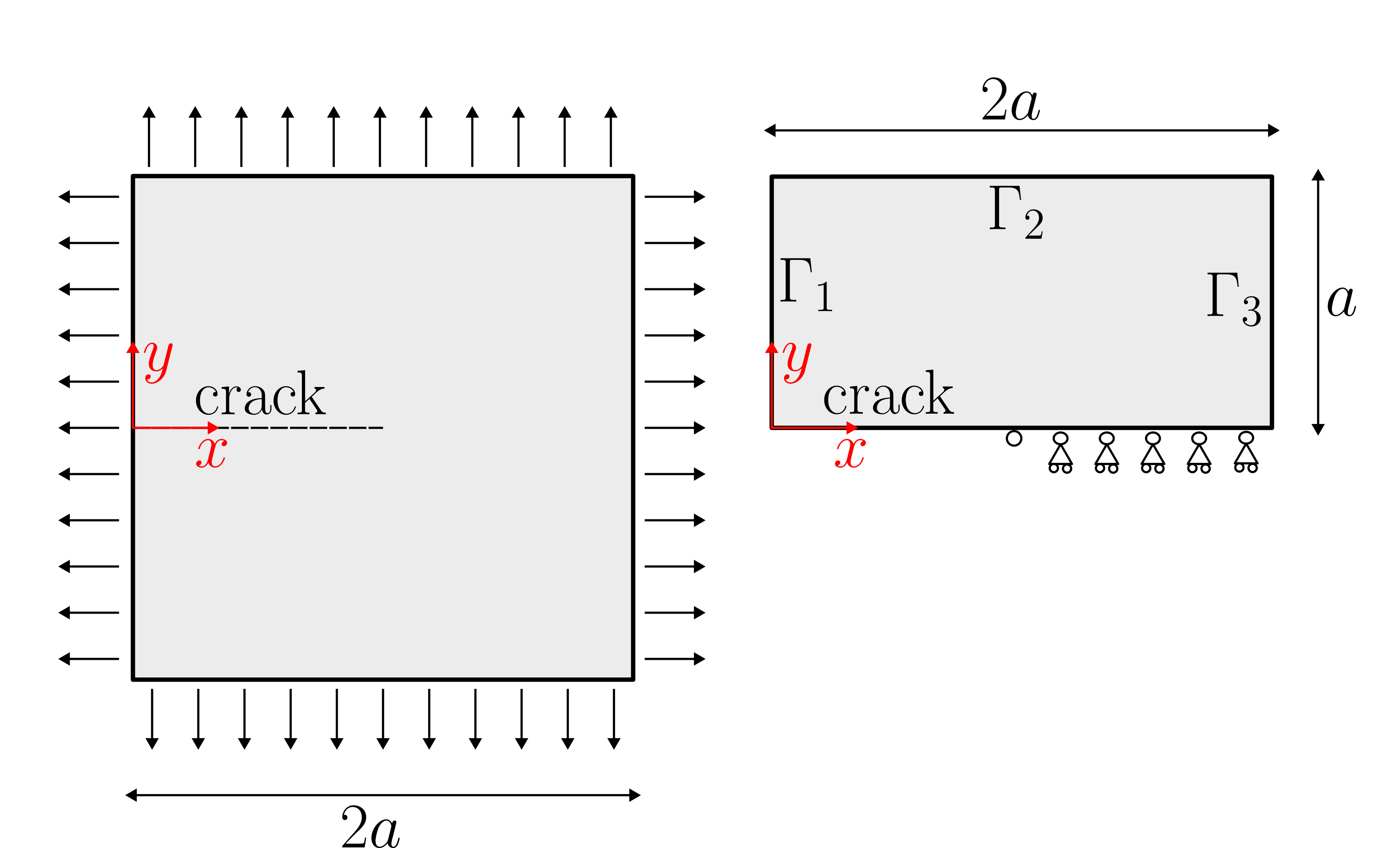}
    \caption{Test case 2: configuration.\label{fig:testcase1_configuration}}
\end{figure}

The configuration is reported in \fig{testcase1_configuration}. By exploiting the symmetry of the problem, only half of the panel is considered and the crack is simulated via Dirichlet boundary conditions. The panel has a half-length $a=2$~mm and is made of isotropic material with $E=10$~MPa and $\nu=0.3$. Plane strain conditions are assumed.

The analytical stress field in polar coordinates $\plbr{\theta,r}$ is given by~\cite{zander2017multi}:
\begin{equation}
    \begin{aligned}
        &\sigma_{xx} = \frac{K_1}{\sqrt{2\pi r}} \cos{\frac{\theta}{2}} \plbr{1-\sin{\frac{\theta}{2}}\sin{\frac{3\theta}{2}}},&&
        \sigma_{yy} = \frac{K_1}{\sqrt{2\pi r}} \cos{\frac{\theta}{2}} \plbr{1+\sin{\frac{\theta}{2}}\sin{\frac{3\theta}{2}}},\\
        &\sigma_{xy} = \frac{K_1}{\sqrt{2\pi r}} \sin{\frac{\theta}{2}} \cos{\frac{\theta}{2}}\cos{\frac{3\theta}{2}},
    \end{aligned}
    \label{eq:testcase1_stress_analytical}
\end{equation}
where $K_1=1$ is the stress intensity factor. Neumann boundary conditions are applied along the free edges~\cite{zander2017multi}:
\begin{equation}
    \begin{aligned}
        &\bm{t}_1 = \bm{\sigma} \cdot \bm{n} = - 
        \begin{bmatrix}
            \sigma_{xx}\\
            \sigma_{xy}
        \end{bmatrix}\quad
        \forall \bm{x} \in \Gamma_1,&&
        \bm{t}_2 = \bm{\sigma} \cdot \bm{n} = + 
        \begin{bmatrix}
            \sigma_{xy}\\
            \sigma_{yy}
        \end{bmatrix}\quad
        \forall \bm{x} \in \Gamma_2,\\
        &\bm{t}_3 = \bm{\sigma} \cdot \bm{n} = + 
        \begin{bmatrix}
            \sigma_{xx}\\
            \sigma_{xy}
        \end{bmatrix}\quad
        \forall \bm{x} \in \Gamma_3.
    \end{aligned}
    \label{eq:testcase1_neumann_bcs}
\end{equation}

The convergence properties are evaluated using the energy norm error:
\begin{equation}
    e_{\bm{\varepsilon}}^{\disp} =  \frac{ \sum_{\el \in \globdom}  \norm{ \sqrt{\bm{A}} \strain^0 \plbr{ \disp - \PiokK \dispd } }_{0,\el}^2  } { \norm{ \sqrt{\bm{A}} \strain^0 \plbr{\disp}}_{0} }.
    \label{eq:testcase1_error}
\end{equation}

In the following, different refinement strategies are employed to assess the convergence properties of the method. Firstly, uniform $h$-refinement and $\kord$-refinement are adopted. Then, local $h$-refinement is performed at the crack tip. The effect of mesh distortion is investigated, too. With this purpose, both structured quadrilateral and Voronoi meshes are considered.

All the simulations reported below are based on standard stabilized VEM, as preliminary analyses showed no significant differences compared to the self-stabilized variant. 

\subsubsection*{Uniform $h$-refinement and $\kord$-refinement}

\begin{figure}[!htbp]
    \centering
    \subfigure[$h=4 \times 2$.\label{fig:testcase1_h_ref_4_2}]{
        \includegraphics[width=0.3\textwidth]{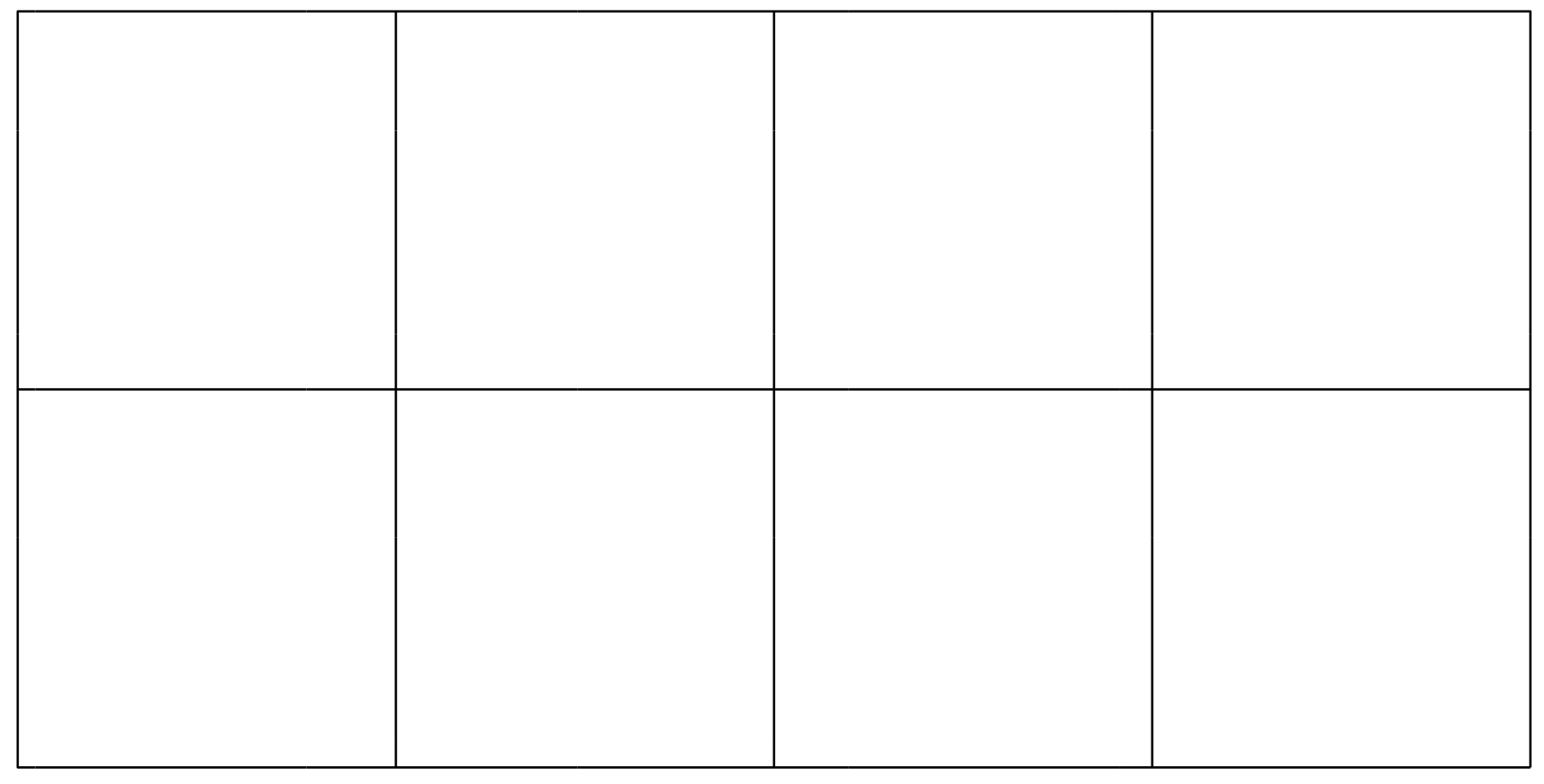}
    }
    \subfigure[$h=8 \times 4$.\label{fig:testcase1_h_ref_8_4}]{
        \includegraphics[width=0.3\textwidth]{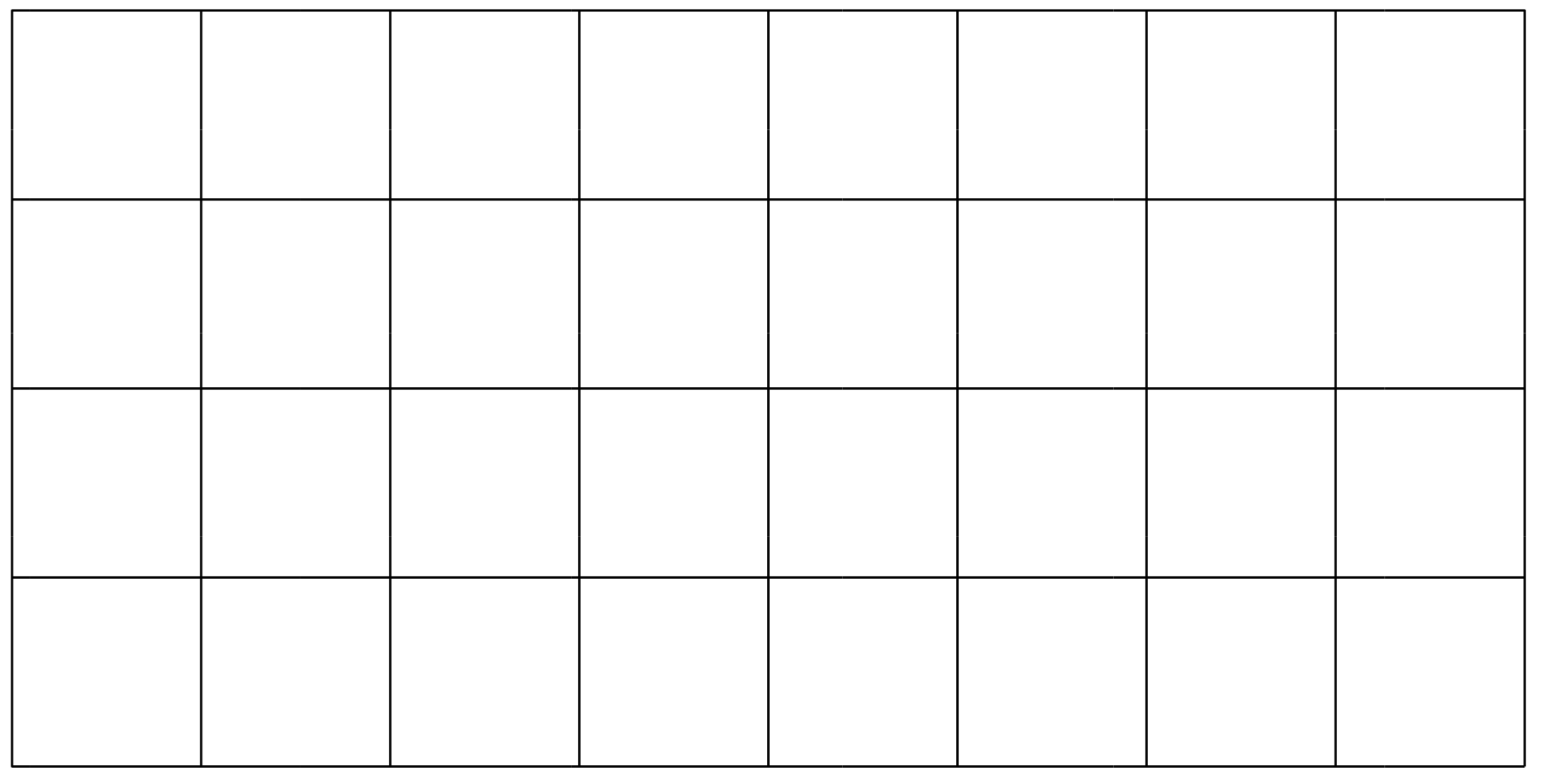}
    }
    \subfigure[$h=12 \times 6$.\label{fig:testcase1_h_ref_12_6}]{
        \includegraphics[width=0.3\textwidth]{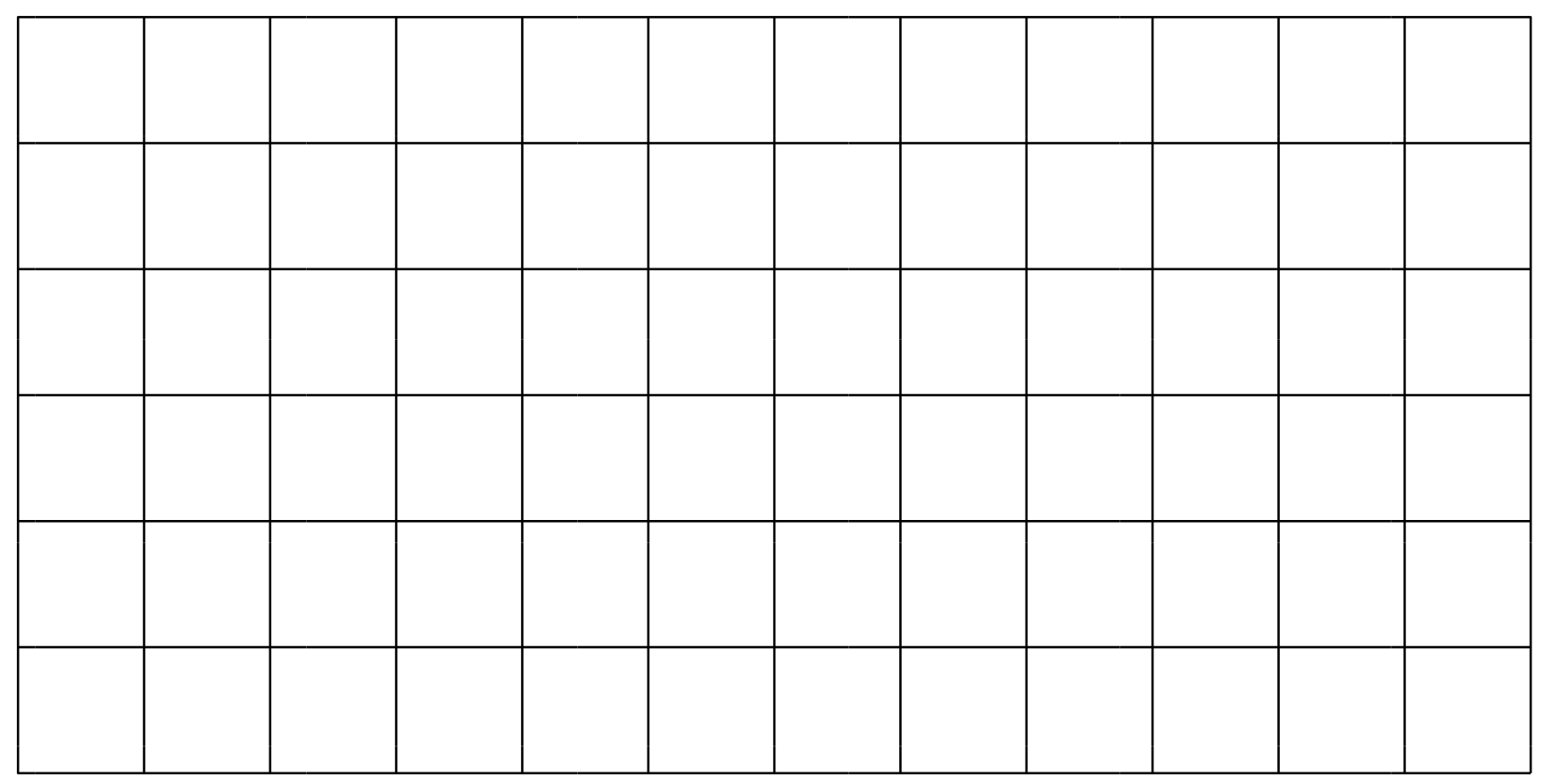}
    }
    \caption{Test case 2: uniform $h$-refinement.\label{fig:testcase1_h_ref}}
\end{figure}

The first investigation deals with the $h$- and $\kord$-refinement strategies. For the latter, the order $\kord$ is progressively increased from $1$ to $4$. The meshes used in the simulations are presented in \fig{testcase1_h_ref}. They consist of rectangular elements of increasing density. 

\begin{figure}[!htbp]
    \centering
    \subfigure[uniform $h$-refinement.\label{fig:testcase1_error_h}]{
        \includegraphics[width=0.47\textwidth]{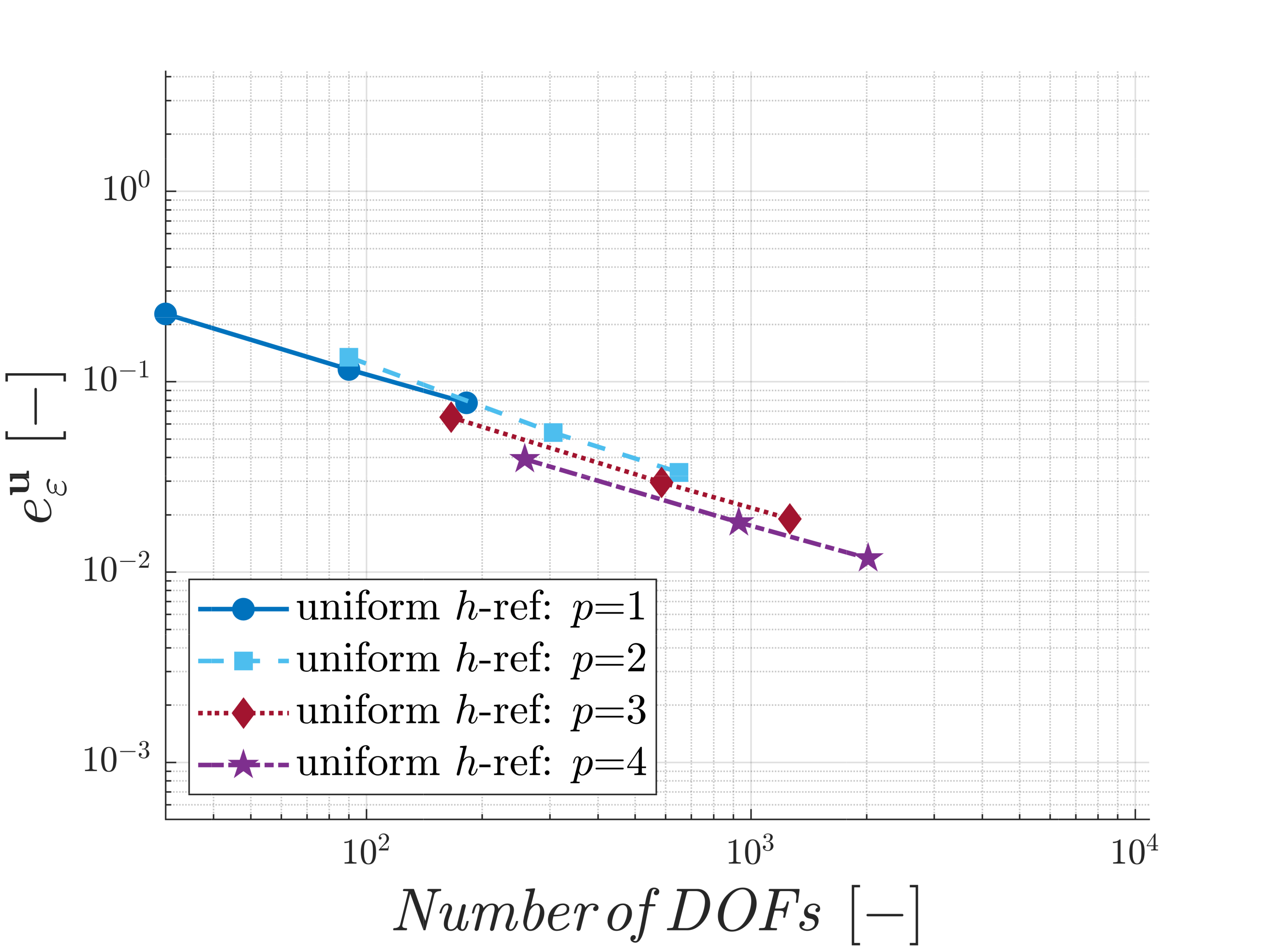}
    }
    \subfigure[$\kord$-refinement.\label{fig:testcase1_error_p}]{
        \includegraphics[width=0.47\textwidth]{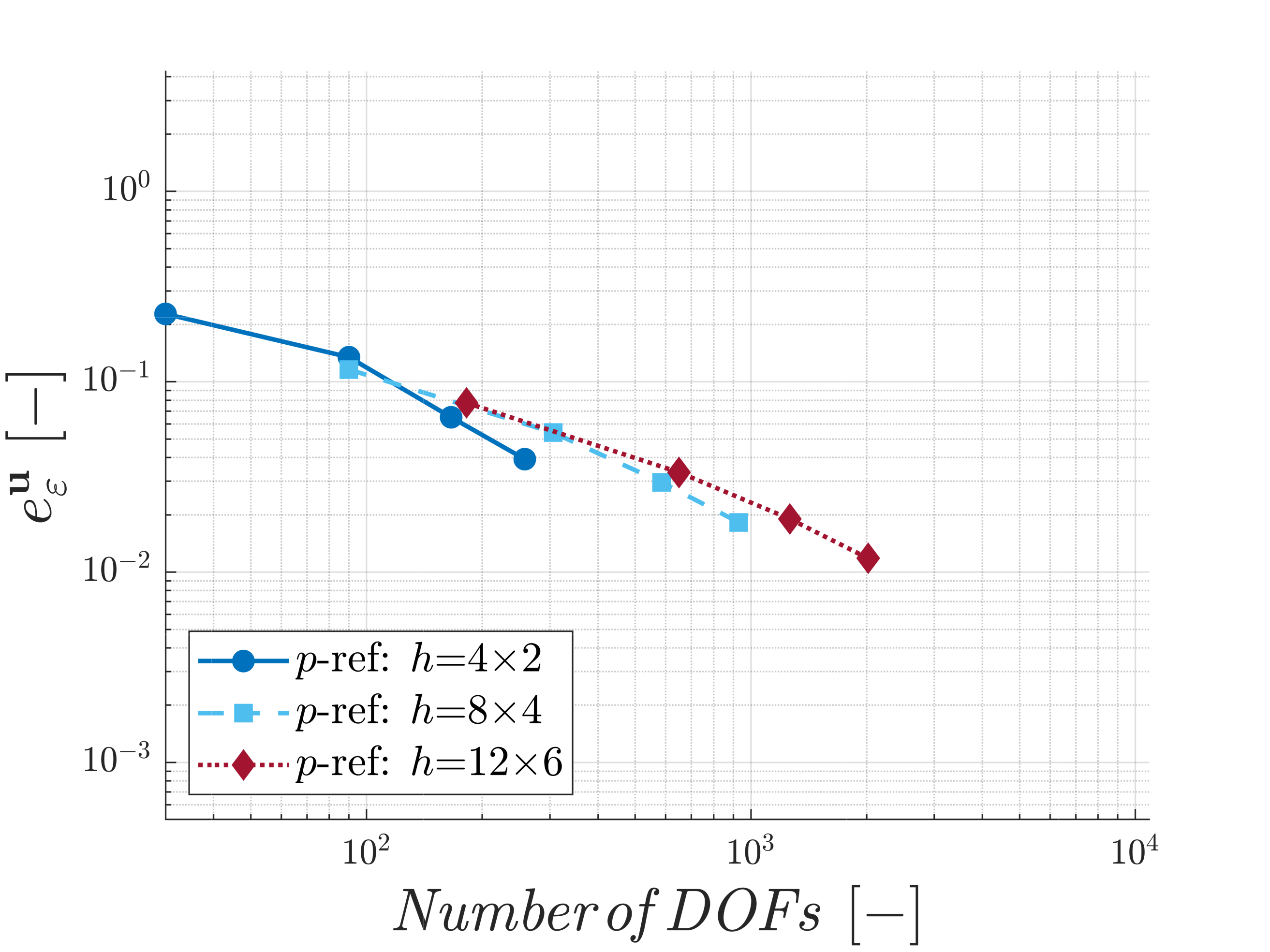}
    }
    \caption{Test case 2: error curves for uniform $h$ and $\kord$ refinements.\label{fig:testcase1_error_hp}}
\end{figure}

A summary of the results is available in \fig{testcase1_error_hp}, where the errors in the energy norm are plotted against the total number of degrees of freedom.

As shown in \fig{testcase1_error_h}, the convergence rate observed by $h$-refinement is essentially independent of the order $\kord$.
This behavior is due to the presence of the singularity in the stress field, hence any increase in the order $\kord$ does not provide any benefit, and the convergence remains limited by this fixed rate.

Similar conclusions are drawn by inspection of \fig{testcase1_error_p}, where the results of the $\kord$-refinement strategy are reported. In particular, an algebraic convergence rate with respect to $\kord$ can be observed. The convergence is not exponential due to the presence of the singularity at the crack tip. Furthermore, a similar convergence rate is achieved for the three different meshes considered here. 

\begin{figure}[!htbp]
    \subfigure[uni. $h$-ref: $h=12\times6$, $\kord=1$.\label{fig:testcase1_stress_profile_1}]{
        \includegraphics[width=0.31\textwidth]{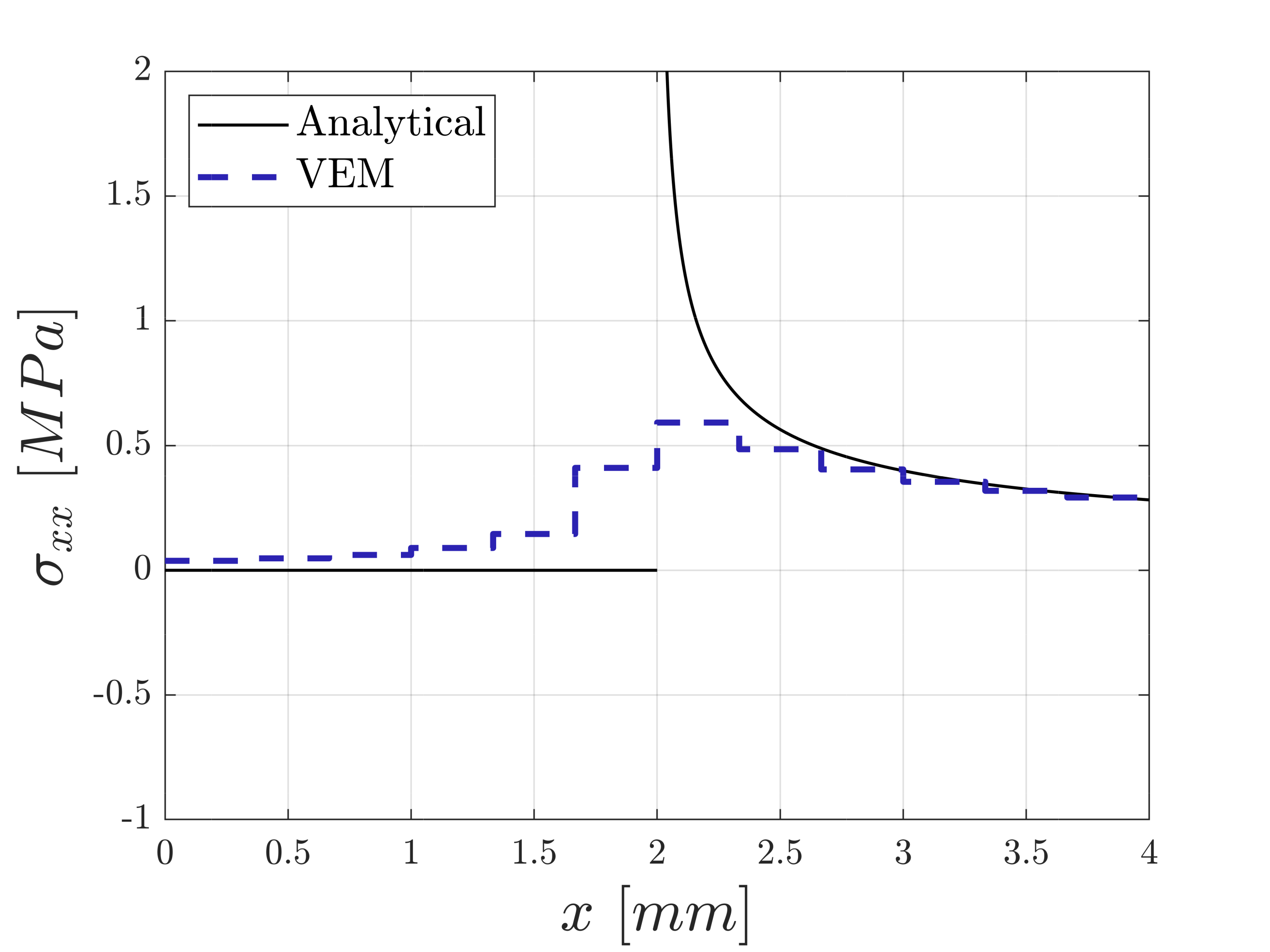}
    }
    \subfigure[$\kord$-ref: $h=4\times2$, $\kord=4$.\label{fig:testcase1_stress_profile_2}]{
        \includegraphics[width=0.31\textwidth]{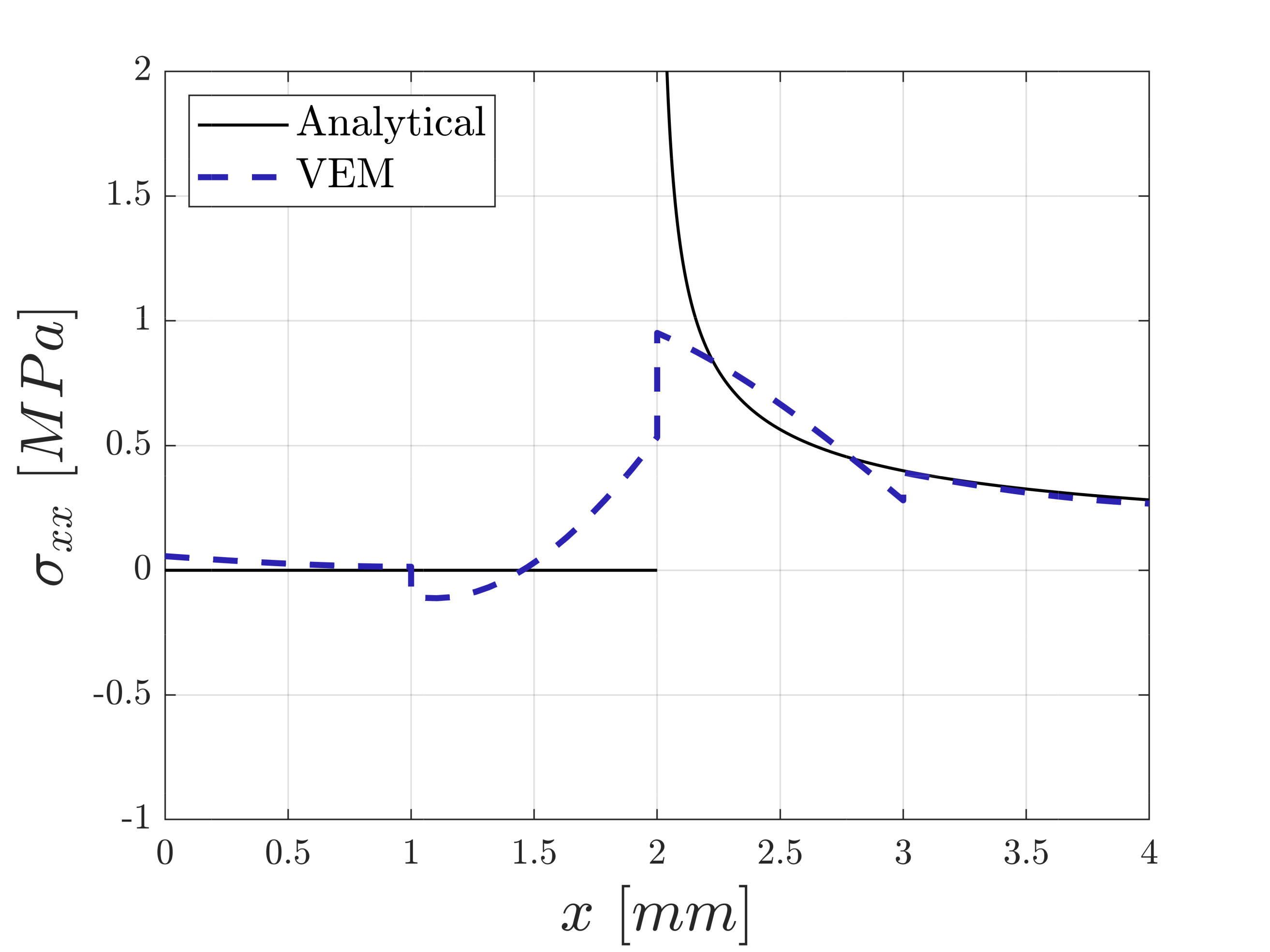}
    }
    \subfigure[uni. $h\kord$-ref: $h=12\times6$, $\kord=4$.\label{fig:testcase1_stress_profile_3}]{
        \includegraphics[width=0.31\textwidth]{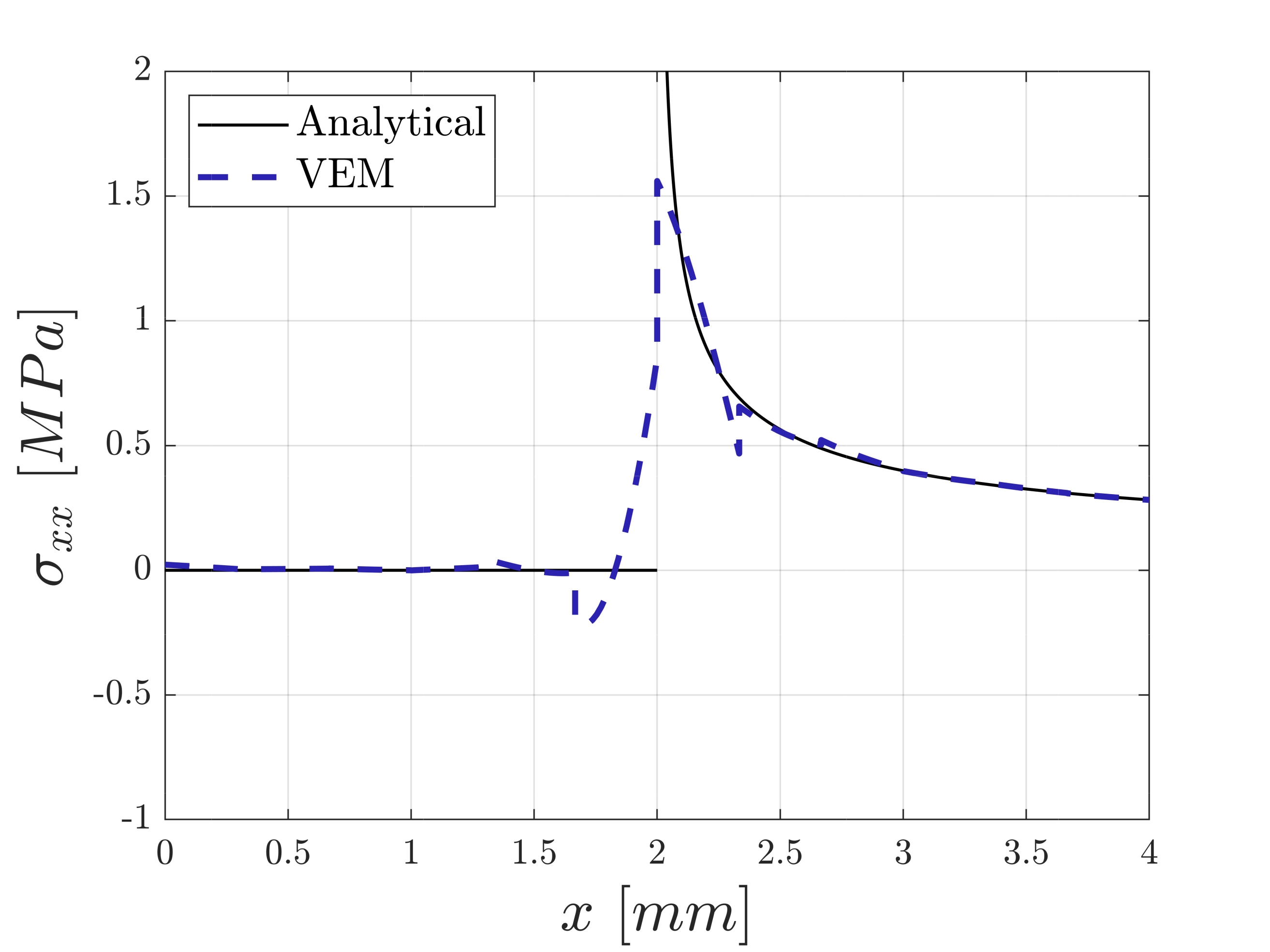}
    }
    \caption{Test case 2: stress profile at $y=0$ for uniform $h$ and $\kord$ refinements.\label{fig:stress_profile_hp}}
\end{figure}

As a further assessment, the stress profile at the coordinate $y=0$ is reported for the different refinement strategies in \fig{stress_profile_hp}. This comparison aims at illustrating the quality of the predictions in the proximity of the crack tip. Three different models are considered for this purpose. The first one corresponds to the $h$-refined mesh with $12\times6$ elements, the second one is obtained via $\kord$-refinement of the $4\times2$ mesh, and the third one by combining these two strategies. The VEM stress is recovered from the polynomial projection via $\PiokK$ of the numerical solution.

As shown in \fig{testcase1_stress_profile_1}, the first strategy provides a relatively inaccurate piecewise constant description in the proximity of the stress peak. An improved stress profile prediction is observed with the second one, as revealed by \fig{testcase1_stress_profile_2}. The $h\kord$-refined strategy, whose results are presented in \fig{testcase1_stress_profile_3}, further improves the quality of the predictions. However, noticeable discrepancies are still present.

\subsubsection*{Local $h$-refinement}

Moving from the results obtained in the previous section, the potential of VEM is now exploited to locally refine the grid in the proximity of the stress singularity. In particular, the size of the elements is progressively halved and mesh conformity is preserved by employing hanging nodes at the center of the edge of the interface elements with non-matching dimensions. This strategy is defined as \textit{local $h$-refinement}.

\begin{figure}[!htbp]
    \centering
    \subfigure[Refinement level $2$.\label{fig:testcase1_d_ref_2}]{
        \includegraphics[width=0.3\textwidth]{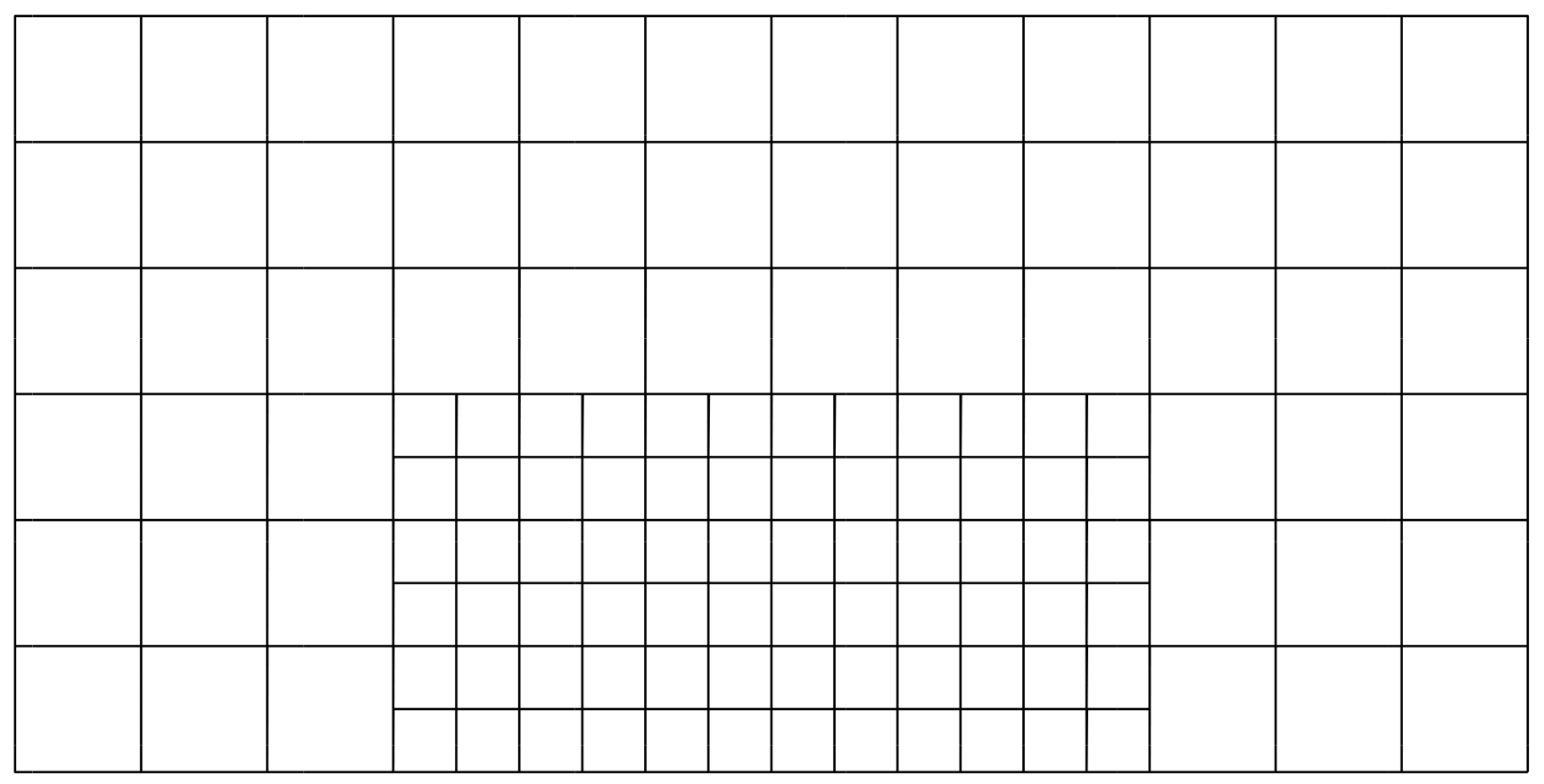}
    }
    \subfigure[Refinement level $4$.\label{fig:testcase1_d_ref_4}]{
        \includegraphics[width=0.3\textwidth]{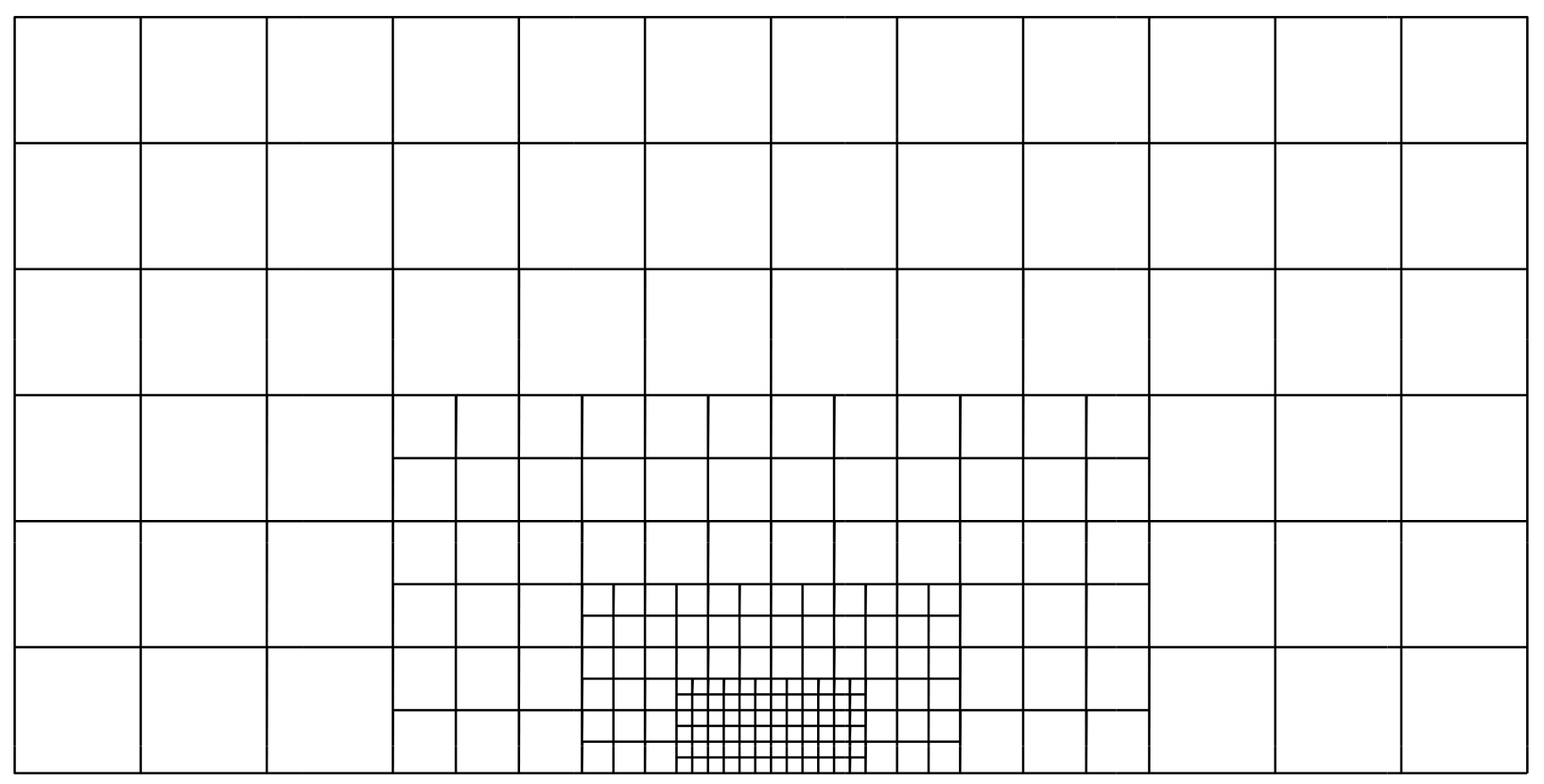}
    }
    \subfigure[Refinement level $7$.\label{fig:testcase1_d_ref_7}]{
        \includegraphics[width=0.3\textwidth]{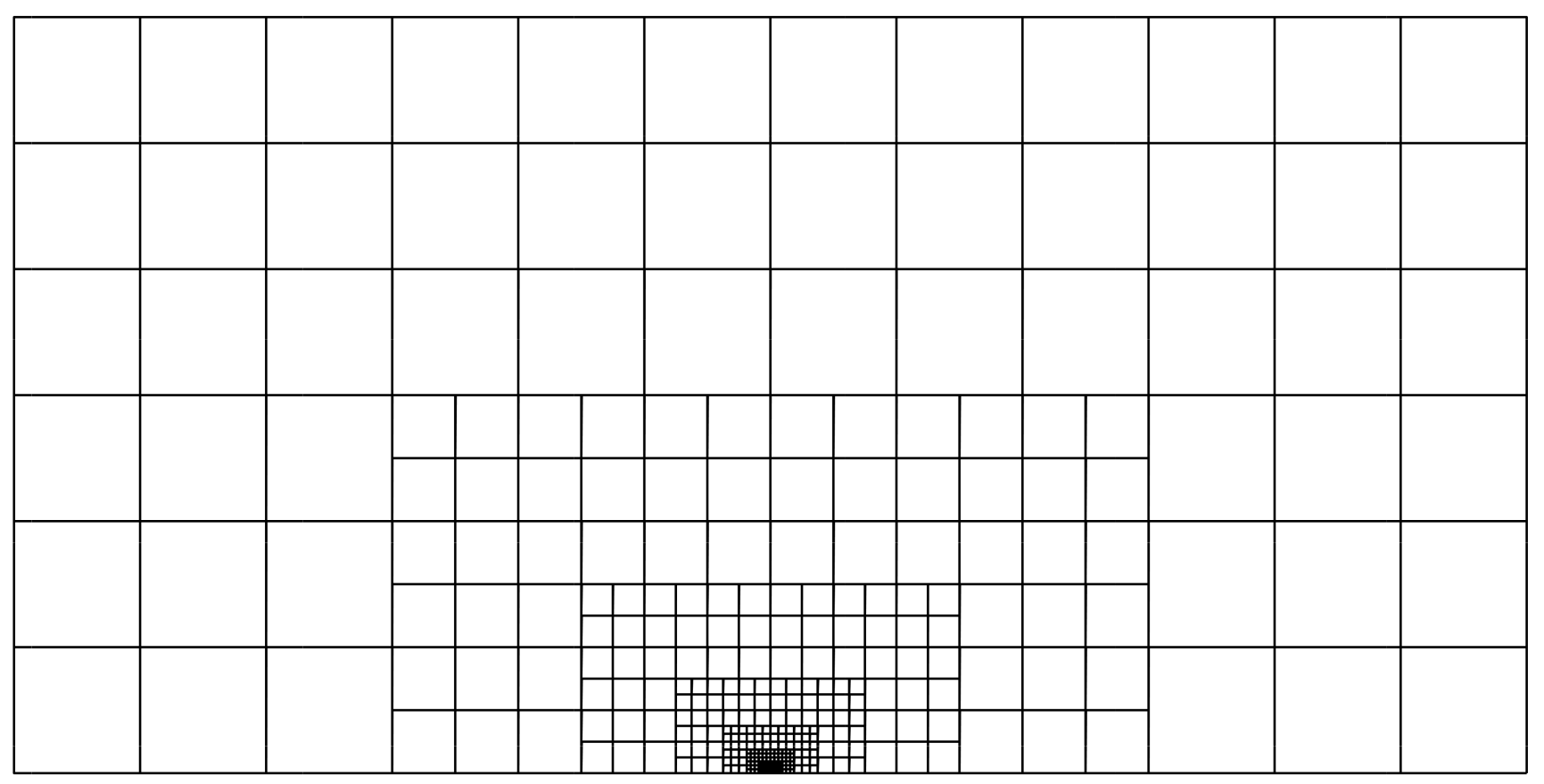}
    }
    \caption{Test case 2: local $h$-refinement for $h=12 \times 6$.\label{fig:testcase1_d_ref}}
\end{figure}

By defining as level 1 the uniform mesh, six more refinements are considered, ranging from level 2 up to level 7. In particular, level 2 corresponds to a mesh where two different sizes of $h$ coexist: one for the broad field and another in correspondence of the crack tip. Similarly, the other refinement levels correspond to an increasing number of progressively refined meshes. For clarity, different levels in the case of $h=12\times6$ are presented in \fig{testcase1_d_ref}.

\begin{figure}[!htbp]
    \centering
    \subfigure[local $h$-refinement (fixed $\kord=1$).\label{fig:testcase1_error_d_h}]{
        \includegraphics[width=0.47\textwidth]{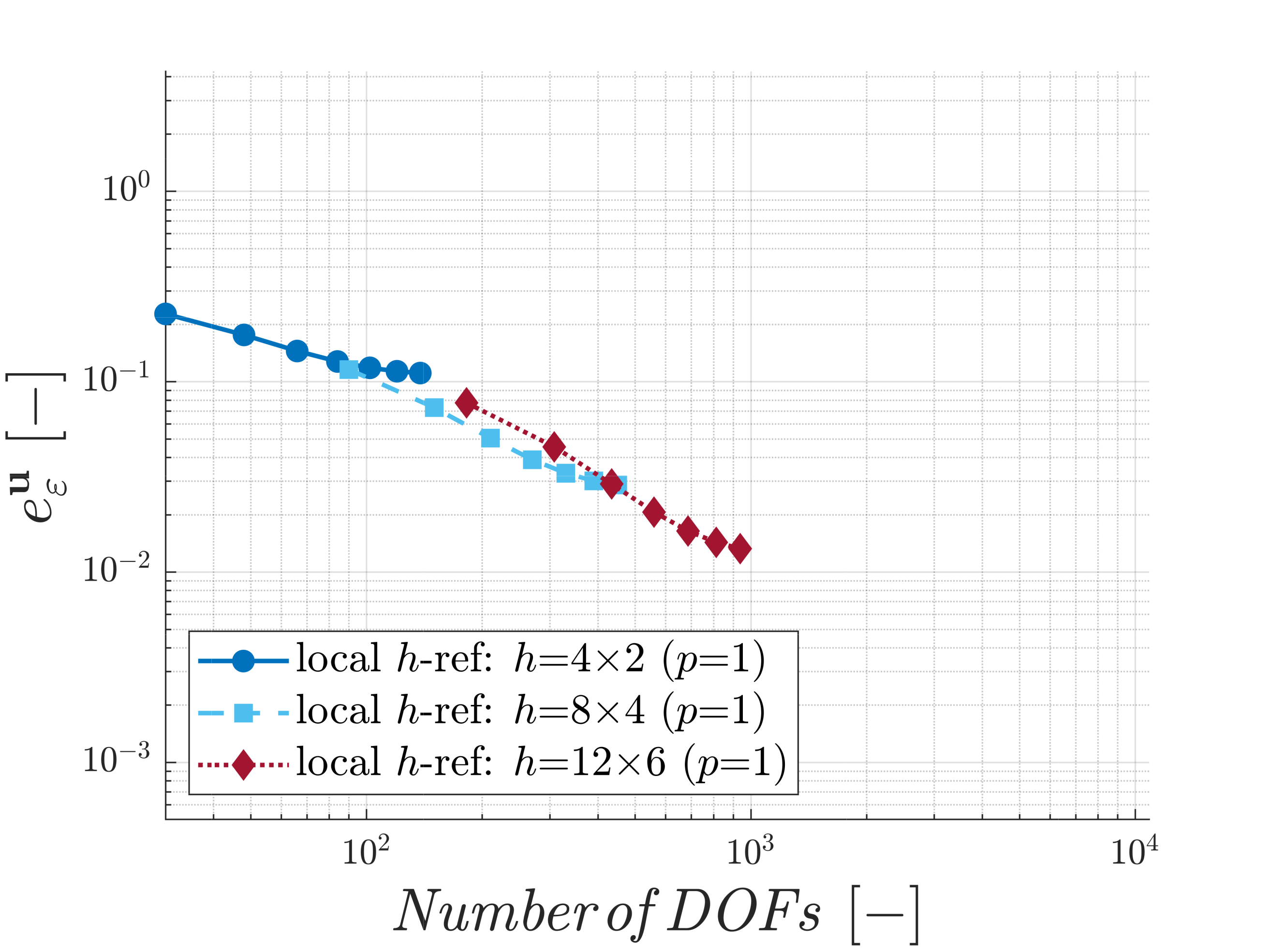}
    }
    \subfigure[local $h$-refinement (fixed $h=12\times6$).\label{fig:testcase1_error_d_p}]{
        \includegraphics[width=0.47\textwidth]{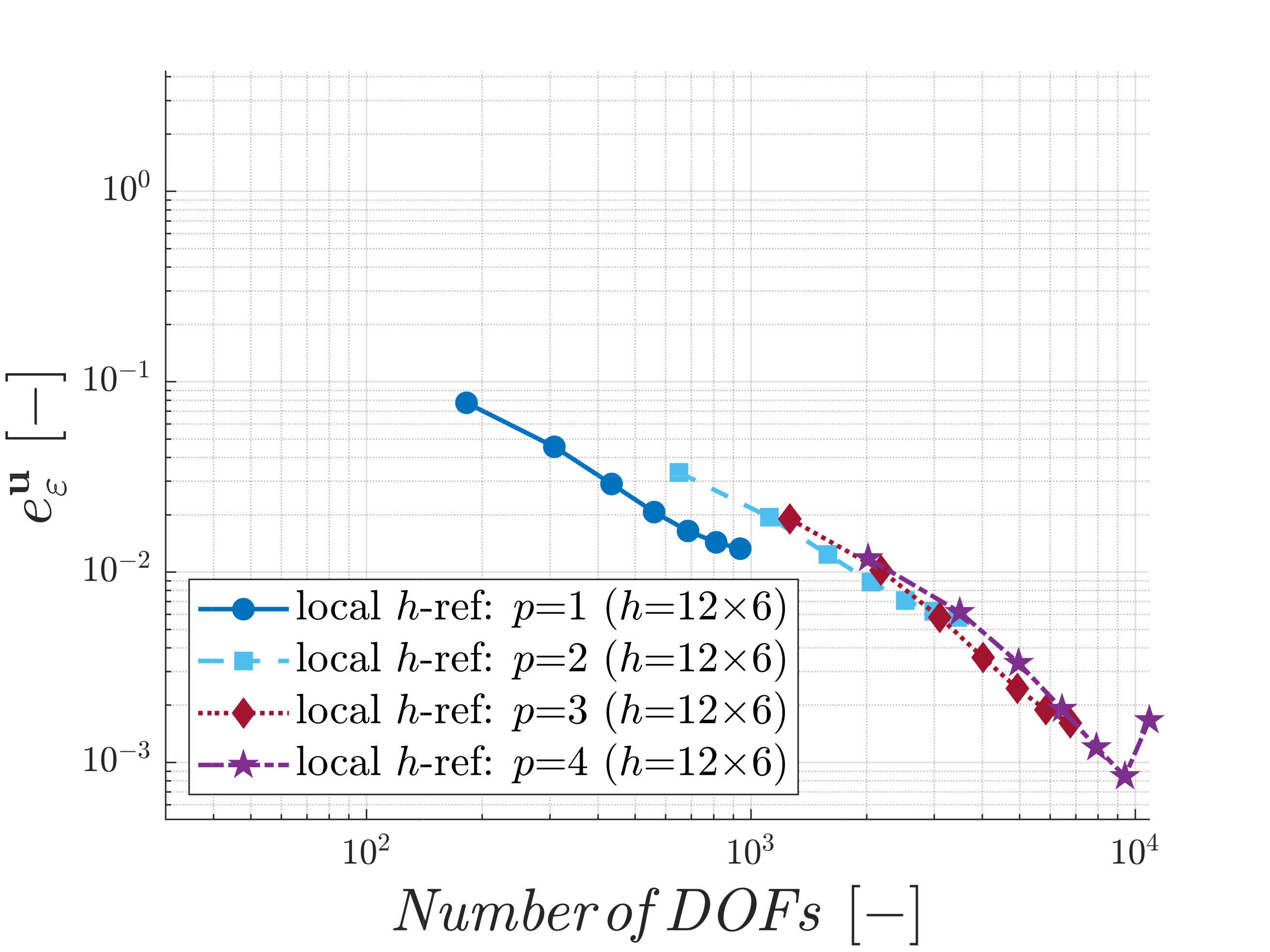}
    }
    \caption{Test case 2: error curves for local $h$-refinement.\label{fig:testcase1_error_d}}
\end{figure}

The error curves are shown in \fig{testcase1_error_d}. Within the same local $h$-refinement paradigm, two different studies are conducted. In the first case, three different mesh densities are considered, and refinement is operated while keeping fixed the order at $\kord=1$. In the second case, the initial mesh density is fixed at $12\times6$ elements. Different orders $\kord=1,\ldots,4$ are considered and each model is progressively refined locally.

 Compared to the uniform refinement investigated earlier, the local $h$-refinement of \fig{testcase1_error_d_h} illustrates an improved convergence rate.

\begin{figure}[!htbp]
    \centering
    \subfigure[Stress contour ($\sigma_{xx}$~$\lbrack$MPa$\rbrack$).\label{fig:testcase1_stress_contour}]{
        \includegraphics[width=0.4\textwidth]{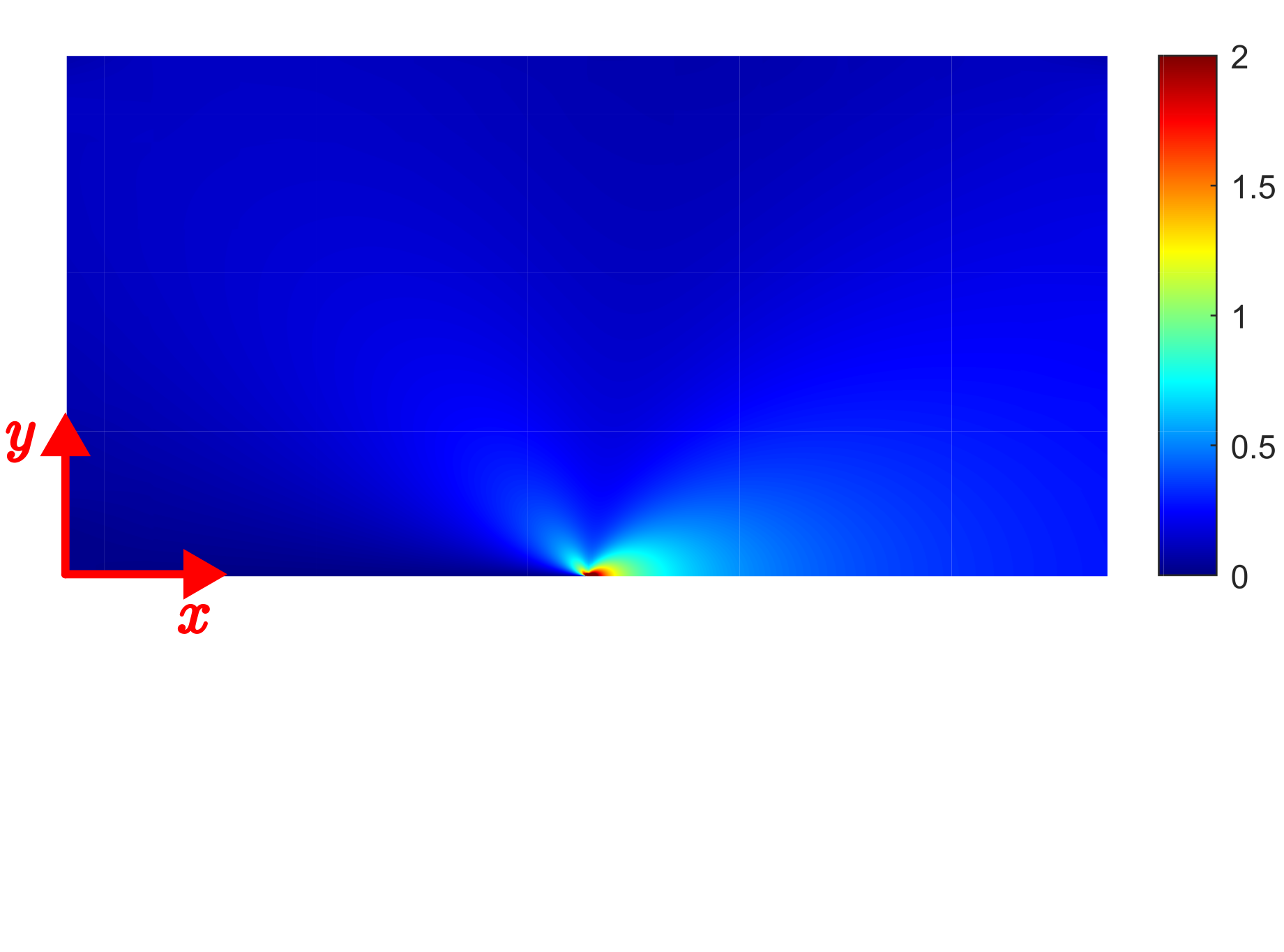}
    }
    \subfigure[Stress profile at $y=0$.\label{fig:testcase1_stress_profile_4}]{
        \includegraphics[width=0.4\textwidth]{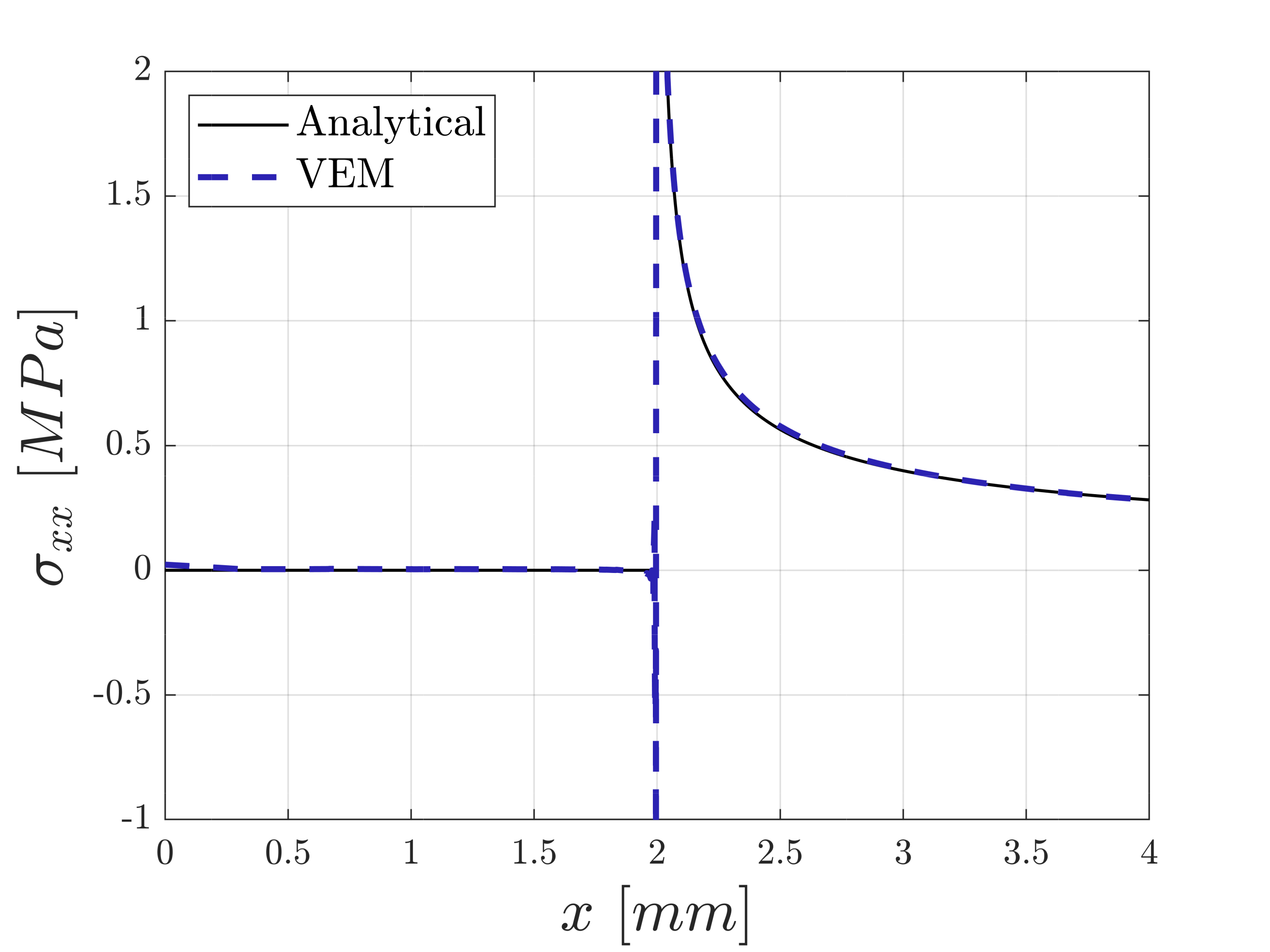}
    }
    \caption{Test case 2: local $h$-refinement: $h=12 \times 6$, $\kord=4$, refinement level $7$.\label{fig:stress_profile_d}}
\end{figure}

 While the results above refer to the global energy response, it is interesting to address the local behavior in terms of stress distribution in the most critical region, i.e. the crack tip. For this reason, the stress profile and the contour of the stress component $\sigma_{xx}$ are reported in \fig{stress_profile_d} for the mesh with $h=12\times6$, $\kord=4$, and $7$ refinement levels. 

The quality of the stress prediction is excellent, as seen from  \fig{testcase1_stress_profile_4}. A slight overshoot of the solution can be noted at the crack tip. This is a common effect of high-order approximations, but it does not affect the overall accuracy of the solution. 

To conclude, the local $h$-refinement strategy guarantees improved convergence and local stress prediction capability when compared to uniform $h$- and $\kord$-refinements approaches presented earlier. Owing to the VEM inherent ability to handle hanging nodes, ease of modeling and excellent accuracy-to-degrees of freedom ratio are guaranteed.

\subsubsection*{Non-structured mesh}

\begin{figure}[!htbp]
    \centering
    \subfigure[Refinement level $1$.\label{fig:testcase1_d_ref_1_voronoi}]{
        \includegraphics[width=0.3\textwidth]{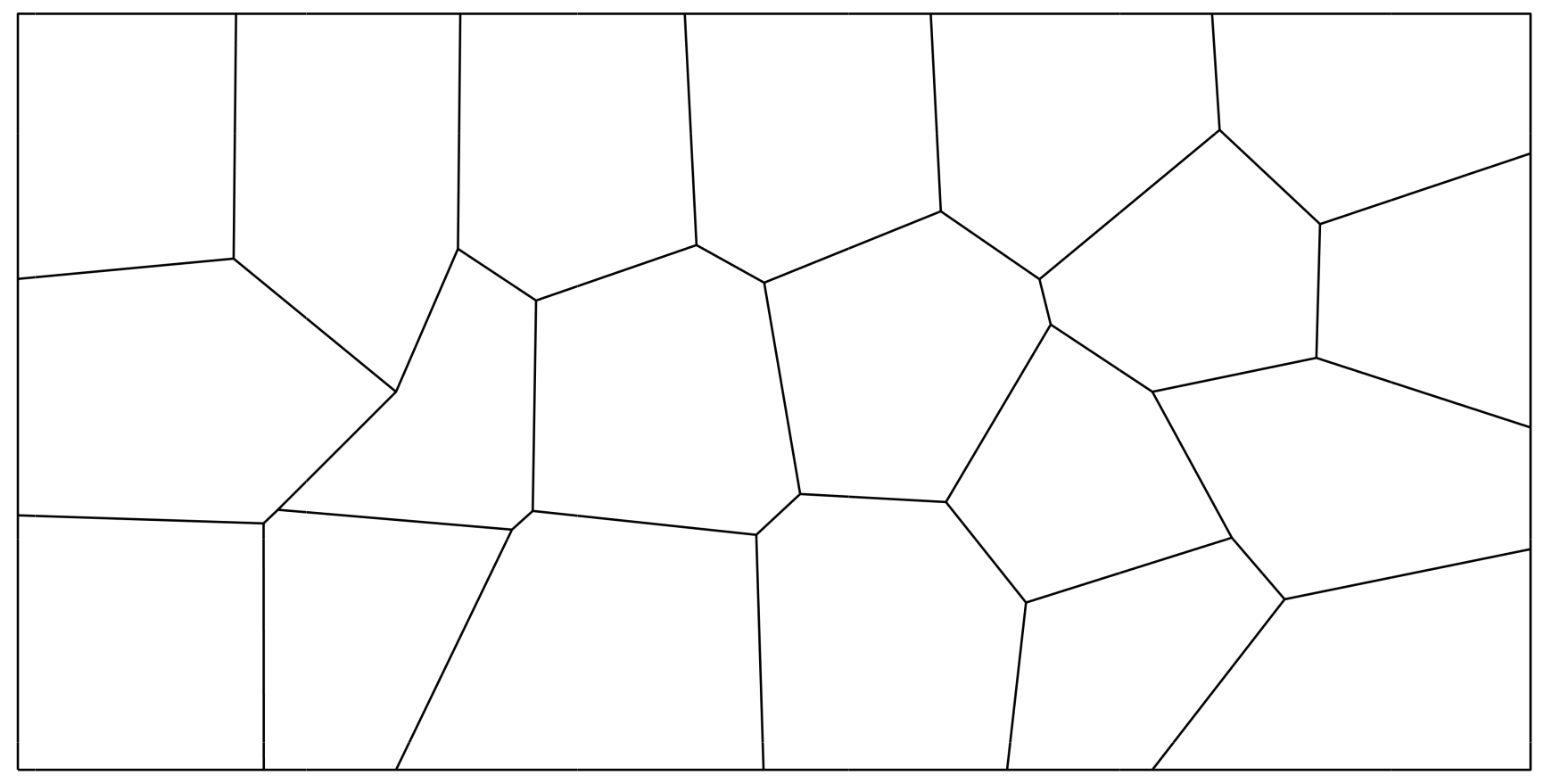}
    }
    \subfigure[Refinement level $4$.\label{fig:testcase1_d_ref_4_voronoi}]{
        \includegraphics[width=0.3\textwidth]{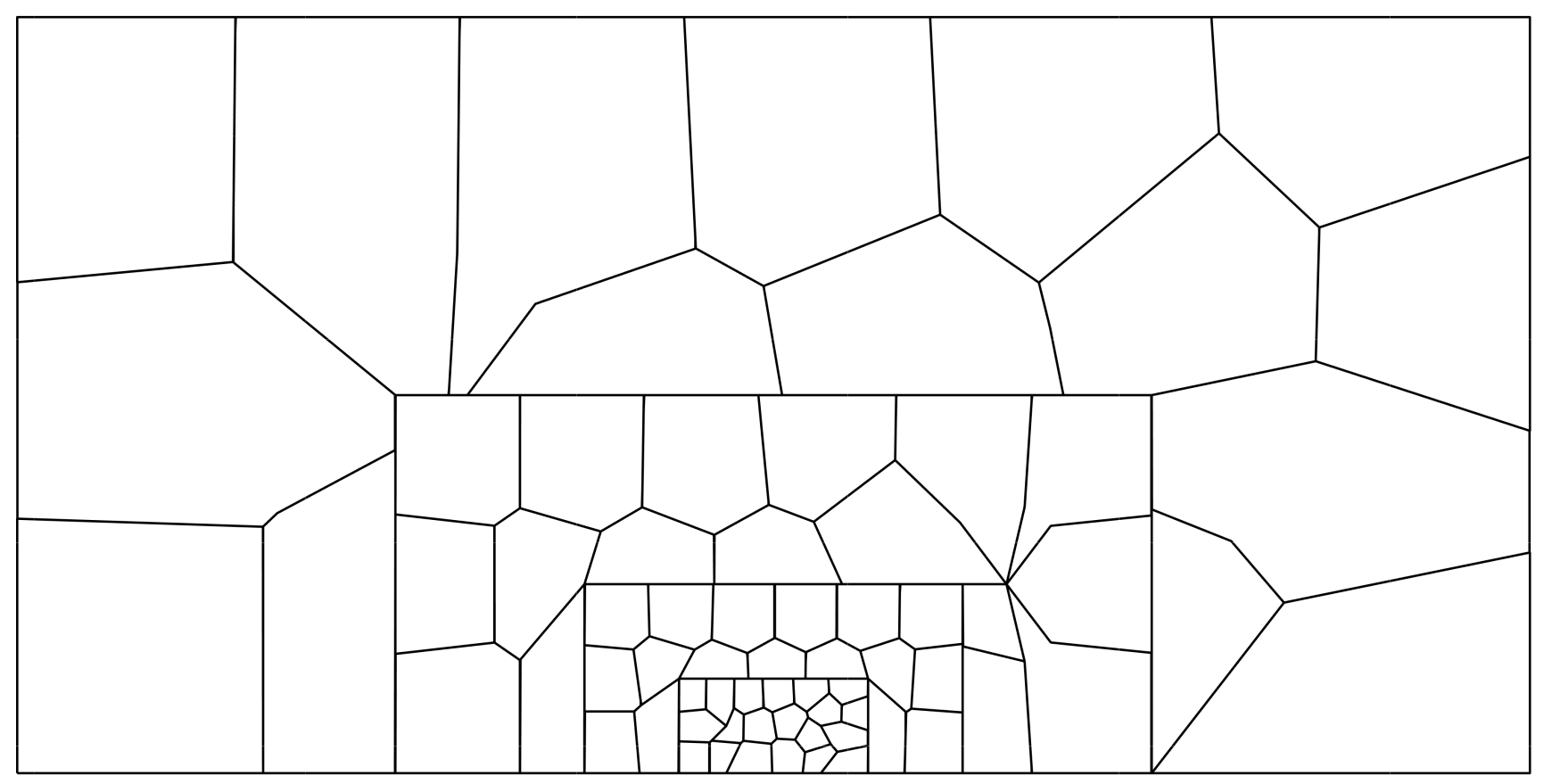}
    }
    \subfigure[Refinement level $7$.\label{fig:testcase1_d_ref_7_voronoi}]{
        \includegraphics[width=0.3\textwidth]{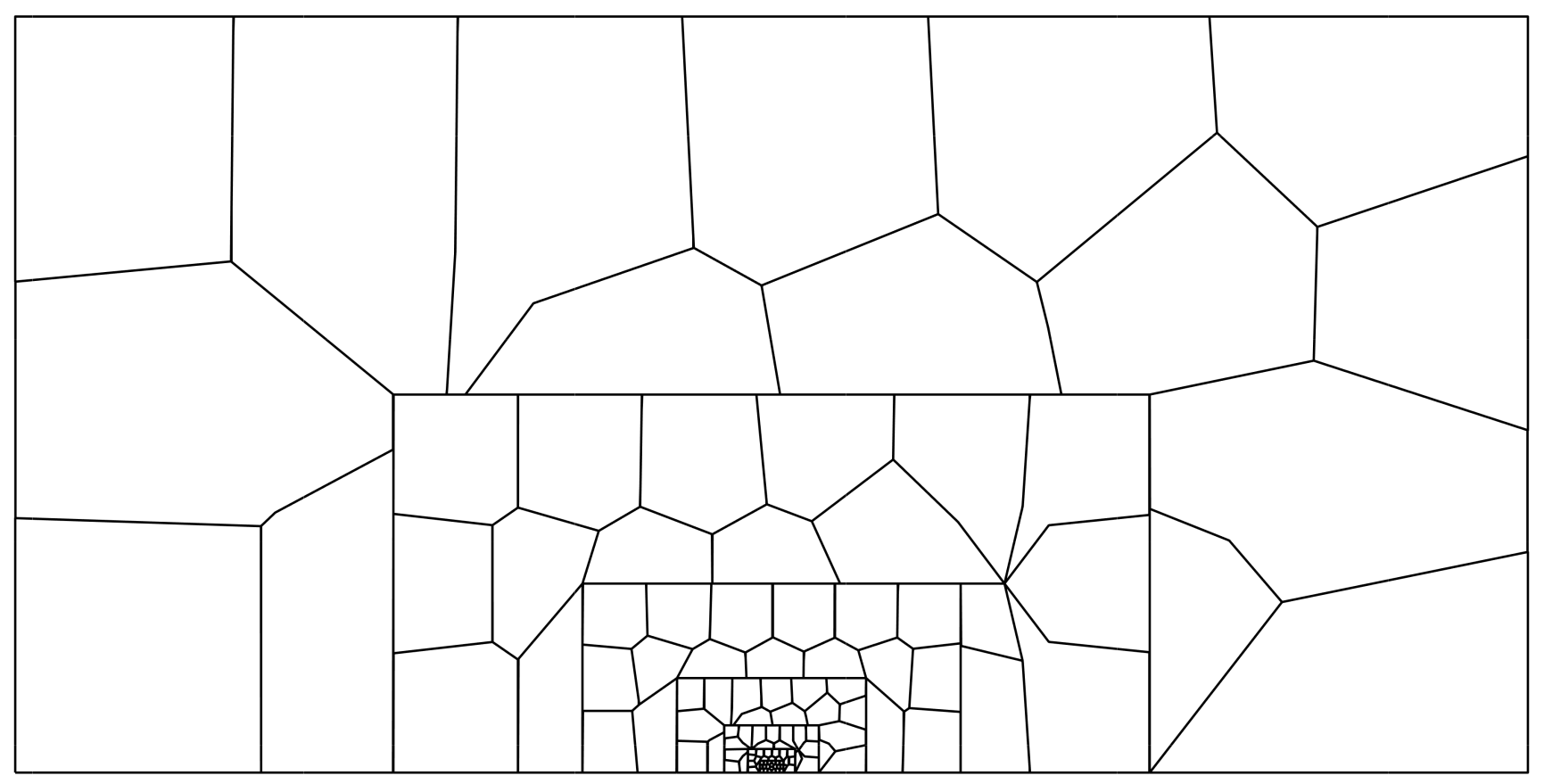}
    }
    \caption{Test case 2: local $h$-refinement for Voronoi mesh.\label{fig:testcase1_d_ref_voronoi}}
\end{figure}

\begin{figure}[!htbp]
    \centering
    \subfigure[$\kord$-refinement.\label{fig:testcase1_error_voronoi_pref}]{
        \includegraphics[width=0.47\textwidth]{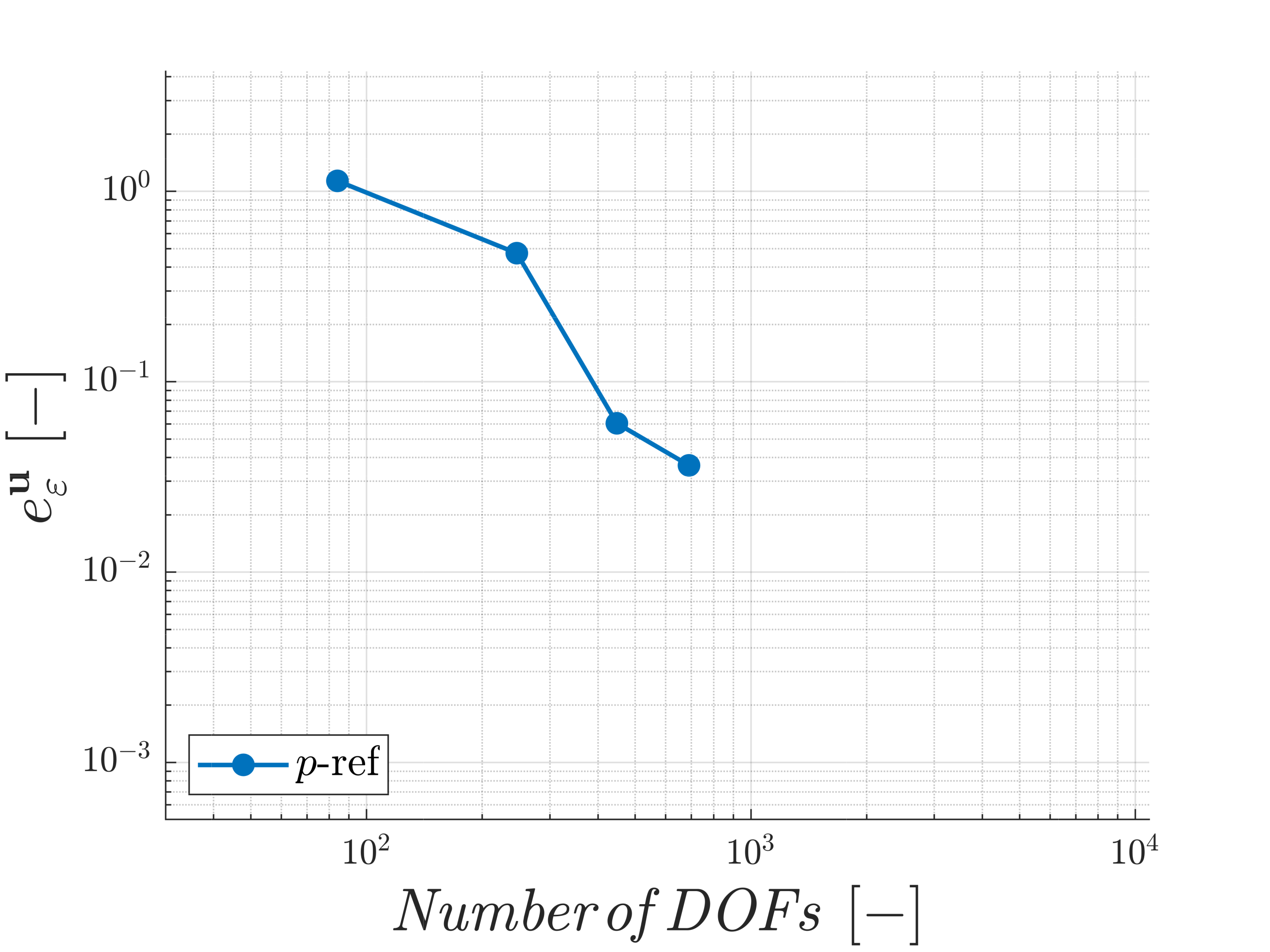}
    }
    \subfigure[local $h$-refinement.\label{fig:testcase1_error_voronoi_localhref}]{
        \includegraphics[width=0.47\textwidth]{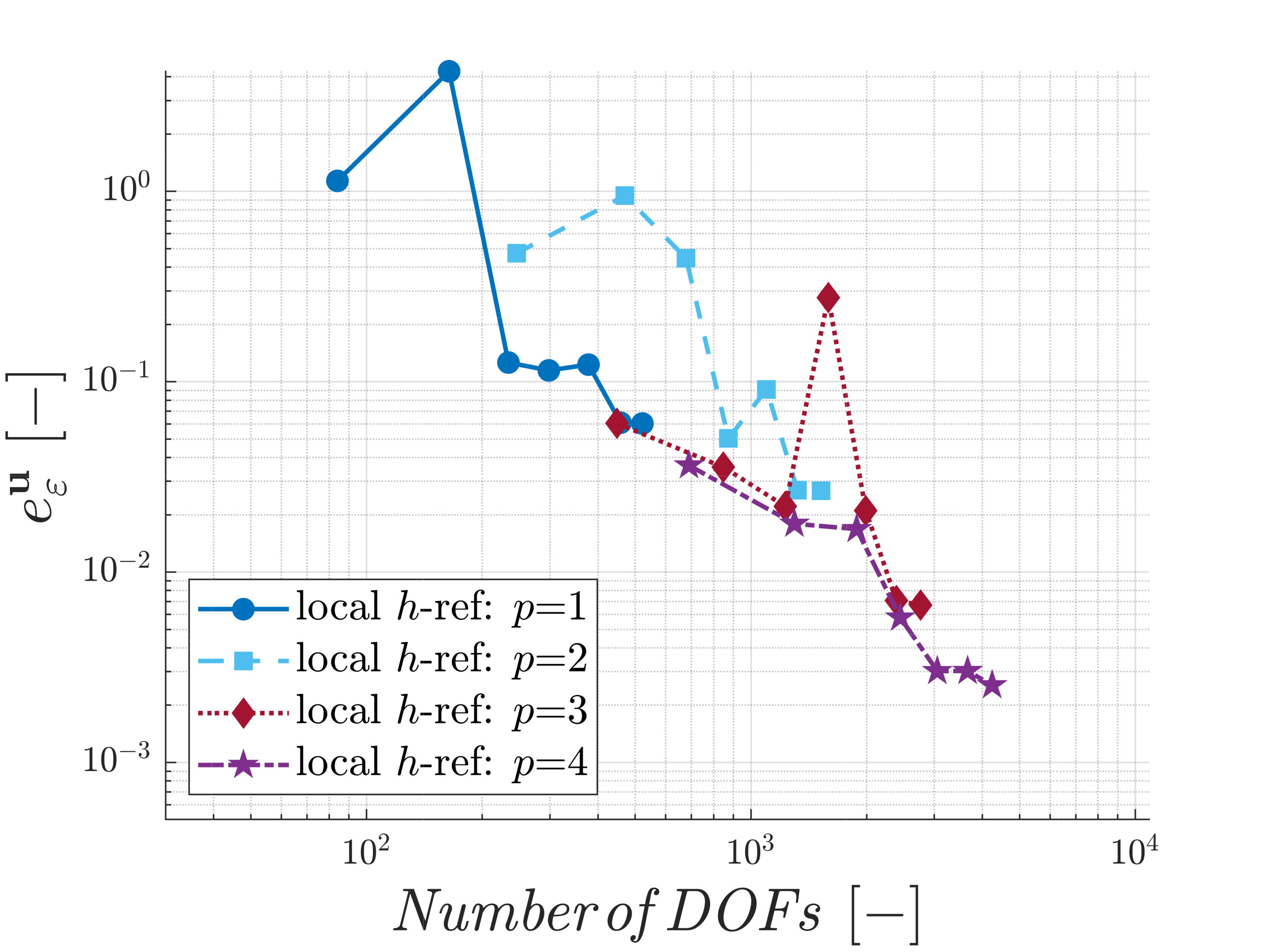}
    }
    \caption{Test case 2: error curves for $\kord$ and local $h$ refinements.\label{fig:testcase1_error_voronoi}}
\end{figure}

The investigation is further extended to verify how the element regularity may affect the conclusions drawn in the previous section. So, a distorted Voronoi mesh is used. 
Two refinements are considered: $\kord$-refinement performed on the level $1$ mesh, and local $h$-refinement up to seven levels. Refinement levels $1$, $4$, and $7$ are shown in \fig{testcase1_d_ref_voronoi}.
On the contrary, uniform $h$-refinement is excluded in this analysis owing to inherent restrictions associated with the type of mesh at hand, since non-structured meshes may lose accuracy under uniform refinements due to increasing mesh distortion. The errors curves for $\kord$ and local $h$-refinement are reported in \figs{testcase1_error_voronoi_pref}{testcase1_error_voronoi_localhref}, respectively.

The curves in \fig{testcase1_error_voronoi} illustrate that diffused distortion of the elements has a detrimental effect on the convergence of the method. On the one hand, the results demonstrate that convergence is achieved despite the high degree of distortion. This robustness is a key feature of VEM. On the other hand, no benefits are now experienced by the local $h$-refinement strategy against the pure $\kord$-refinement. Indeed, the rate of convergence is similar in both cases. Moreover, some oscillations are present, especially for lower orders.

\subsubsection*{Remarks}
The results obtained in this test case demonstrate the difficulties of the standard $h$- and $\kord$-refinement techniques to handle scenarios characterized by drastic stress gradients. In contrast, local $h$-refinements are more suitable for these cases, as shown by noticeable benefits on the convergence plots. Furthermore, this example demonstrates the effectiveness of employing general polygonal elements: local mesh refinements are easily introduced where needed, still preserving the mesh conformity. Moreover, the effects of mesh distortion are investigated by adopting a non-structured Voronoi mesh, illustrating that the presence of distorted elements affects the convergence of the method to some extent, leading to oscillations while preserving the overall trend.

\subsection{Test case 3}

In this third test case, the convergence of the method is assessed for the free-vibration response of a highly anisotropic plate. The benchmark has been studied in~\cite{wu2012comparison,vescovini2018application} and proves to be useful to investigate the VEM ability to capture localization induced by drastic material anisotropy.  

\begin{figure}[!htbp]
    \centering
        \includegraphics[width=0.35\textwidth]{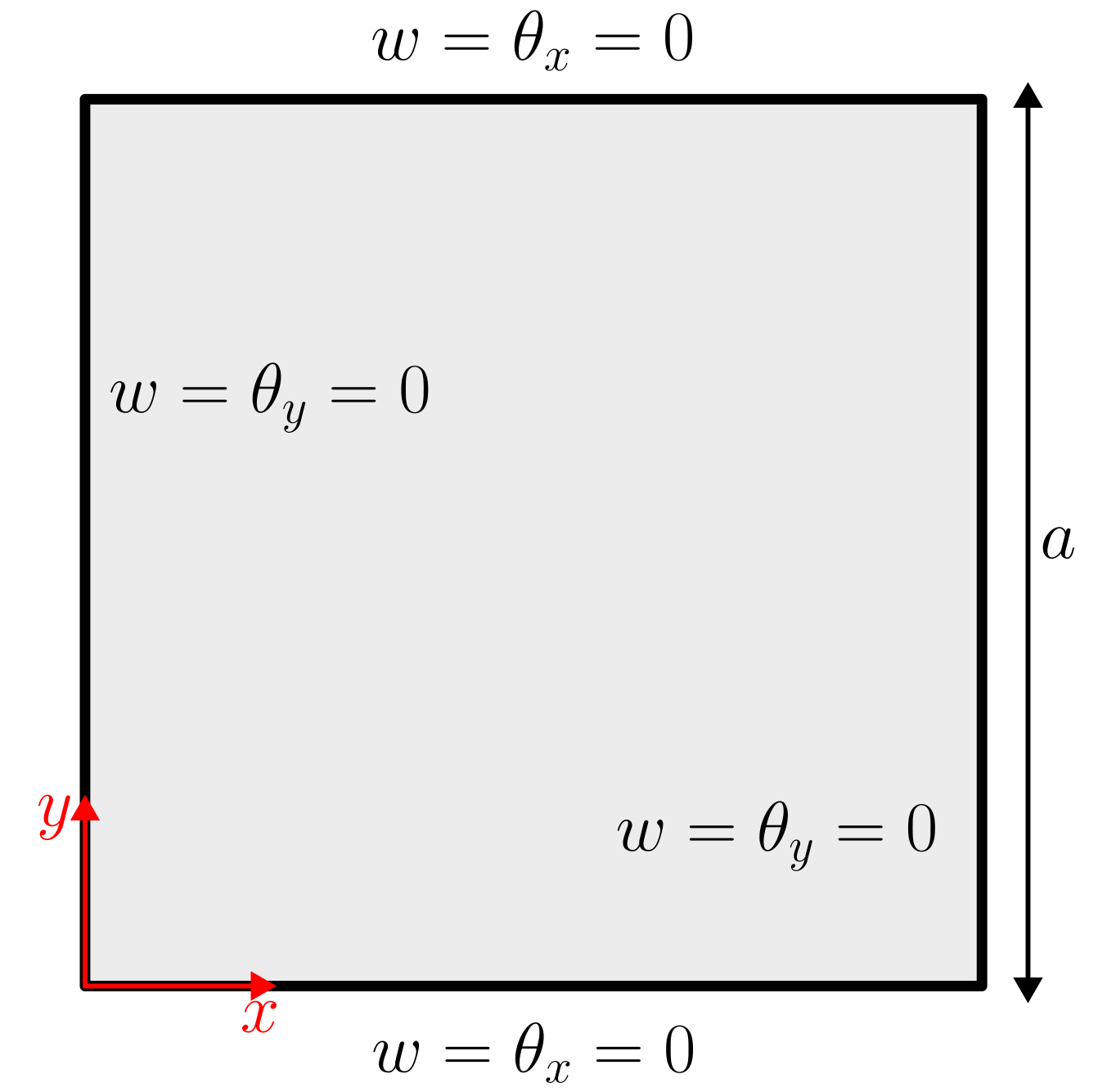}
    \caption{Test case 3: configuration.\label{fig:testcase2_configuration}}
\end{figure}

A simply supported square plate of dimension $a=100$~mm and thickness $h=0.01$~mm is considered. The domain and the corresponding boundary conditions at the edges are specified in \fig{testcase2_configuration}.

The material is a highly anisotropic pre-preg P100/AS3501 carbon/epoxy, with elastic properties: $E_{11}=369000$~MPa, $E_{22}=5030$~MPa, $G_{12}=G_{13}=G_{23}=5240$~MPa, $\nu_{12}=\nu_{13}=\nu_{23}=0.31$ and $\rho=1.5\times 10^{-6}$~kg/mm$^3$, exhibiting an artificially high orthotropic ratio. The laminate is composed of a single ply oriented at $45^\circ$. So, membrane and flexural anisotropy are maximized. The combination of these material and layup features exacerbates the convergence challenges, making this benchmark of special interest.

The free vibration response of the plate is investigated in terms of the first nondimensional circular frequency, defined as:
\begin{equation}
    \bar{\omega} = \omega \frac{a^2}{h} \sqrt{\frac{\rho}{E_{22}}}.
    \label{eq:testcase2_nondim_w}
\end{equation}
The error is measured as:
\begin{equation}
    e_{\bar{\omega}} = \frac{\shbr{\bar{\omega}_{\text{VEM}}-\bar{\omega}_{\text{REF}}}}{\bar{\omega}_{\text{REF}}}.
    \label{eq:testcase2_error}
\end{equation}
The subscripts VEM and REF refer to the present VEM and the reference solution reported in~\cite{yan2023application}, $\bar{\omega}=21.9263$, obtained by application of a refined $ps$-FEM based model.

As done previously, convergence is studied in terms of uniform $h$-, $\kord$- and local $h$-refinement. Moreover, both the standard stabilized and self-stabilized VEM are employed to assess the influence of the stabilization term on the accuracy of the solution. This aspect is particularly relevant, as eigenvalue problems are sensitive to the choice of the stabilization term.

\subsubsection*{Uniform $h$-refinement and $\kord$-refinement}

\begin{figure}[!htbp]
    \centering
    \subfigure[$h=4 \times 4$.\label{fig:testcase2_h_ref_4_4}]{
        \includegraphics[width=0.28\textwidth]{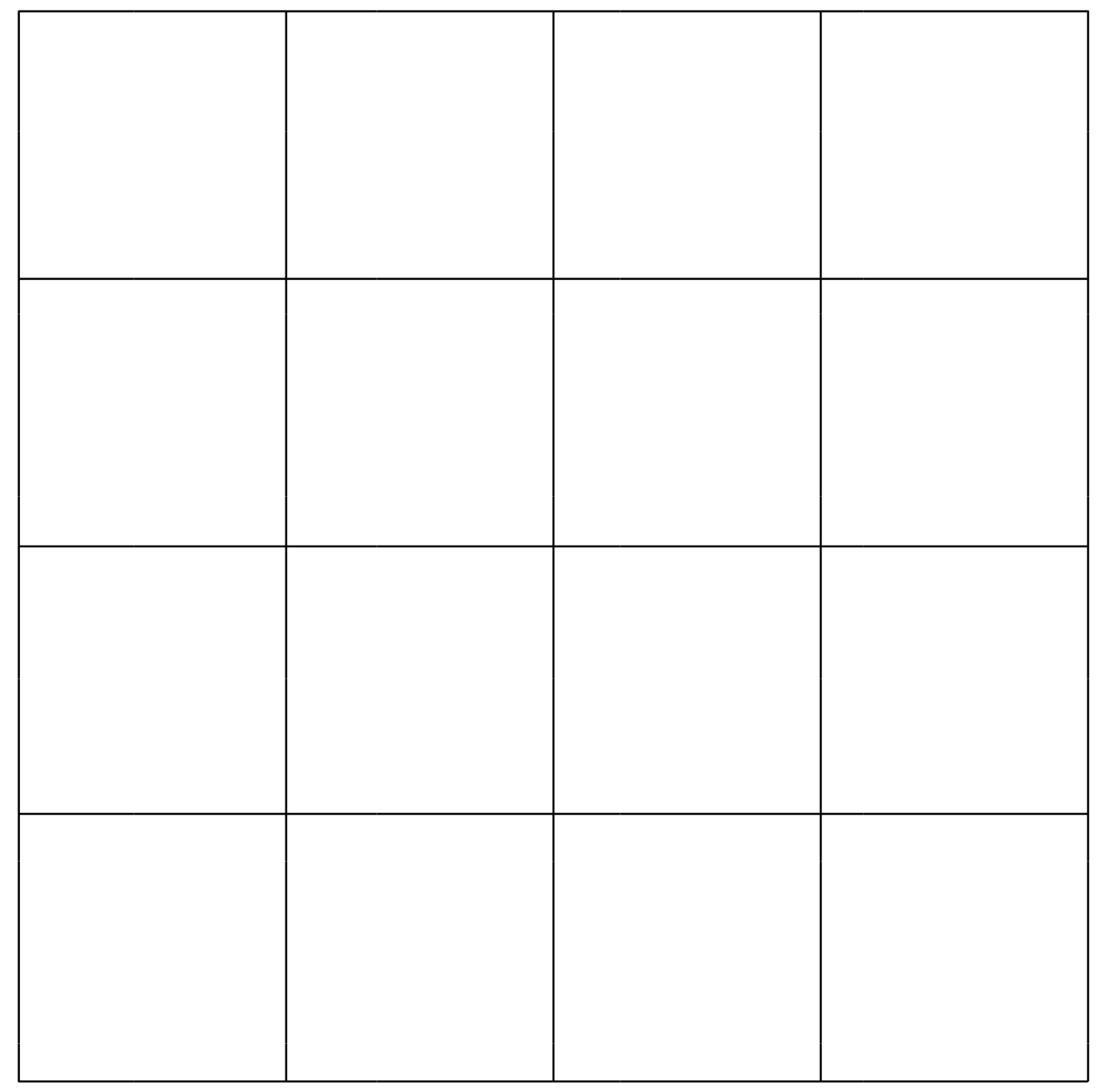}
    }
    \subfigure[$h=8 \times 8$.\label{fig:testcase2_h_ref_8_8}]{
        \includegraphics[width=0.28\textwidth]{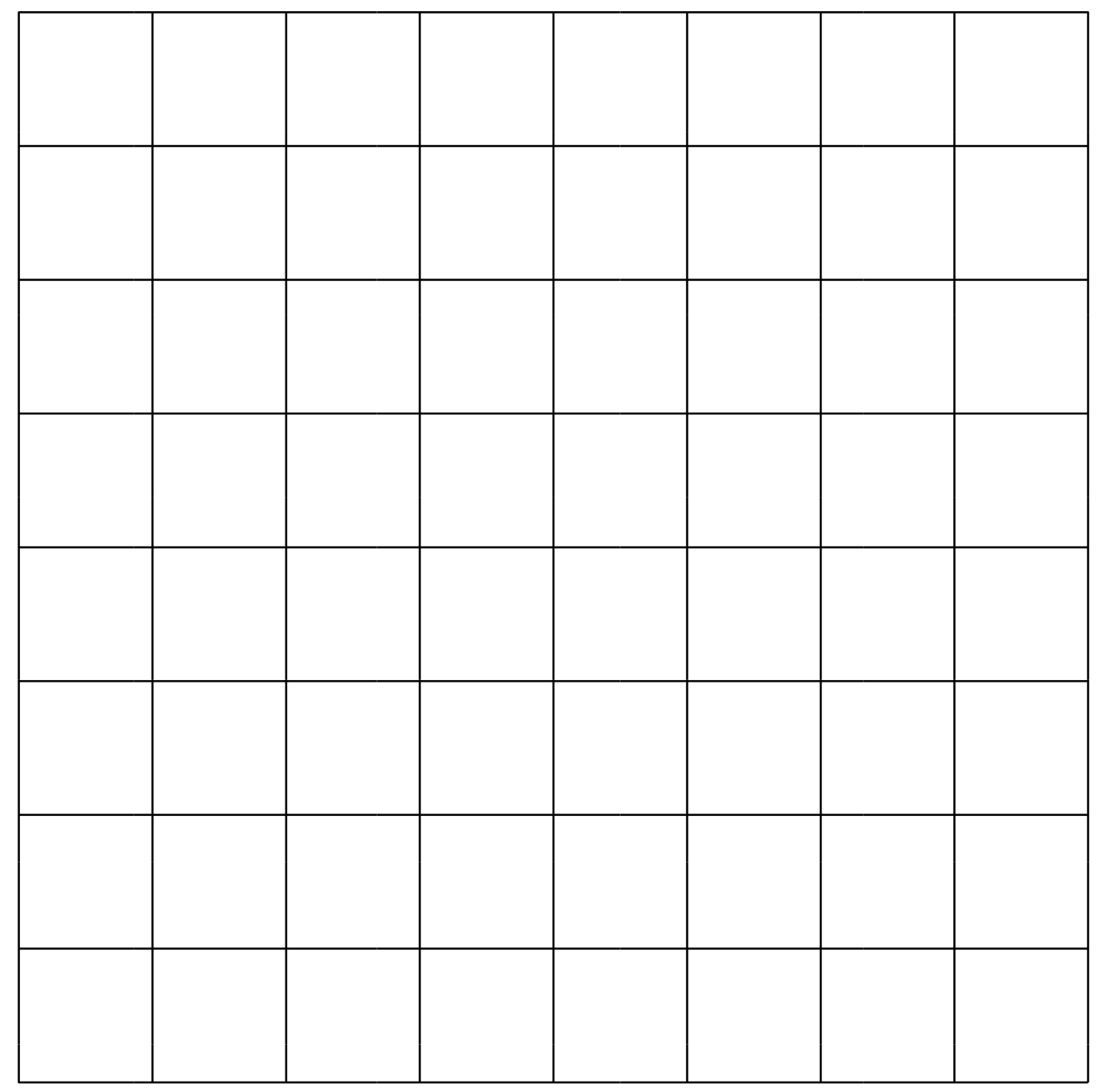}
    }
    \subfigure[$h=12 \times 12$.\label{fig:testcase2_h_ref_12_12}]{
        \includegraphics[width=0.28\textwidth]{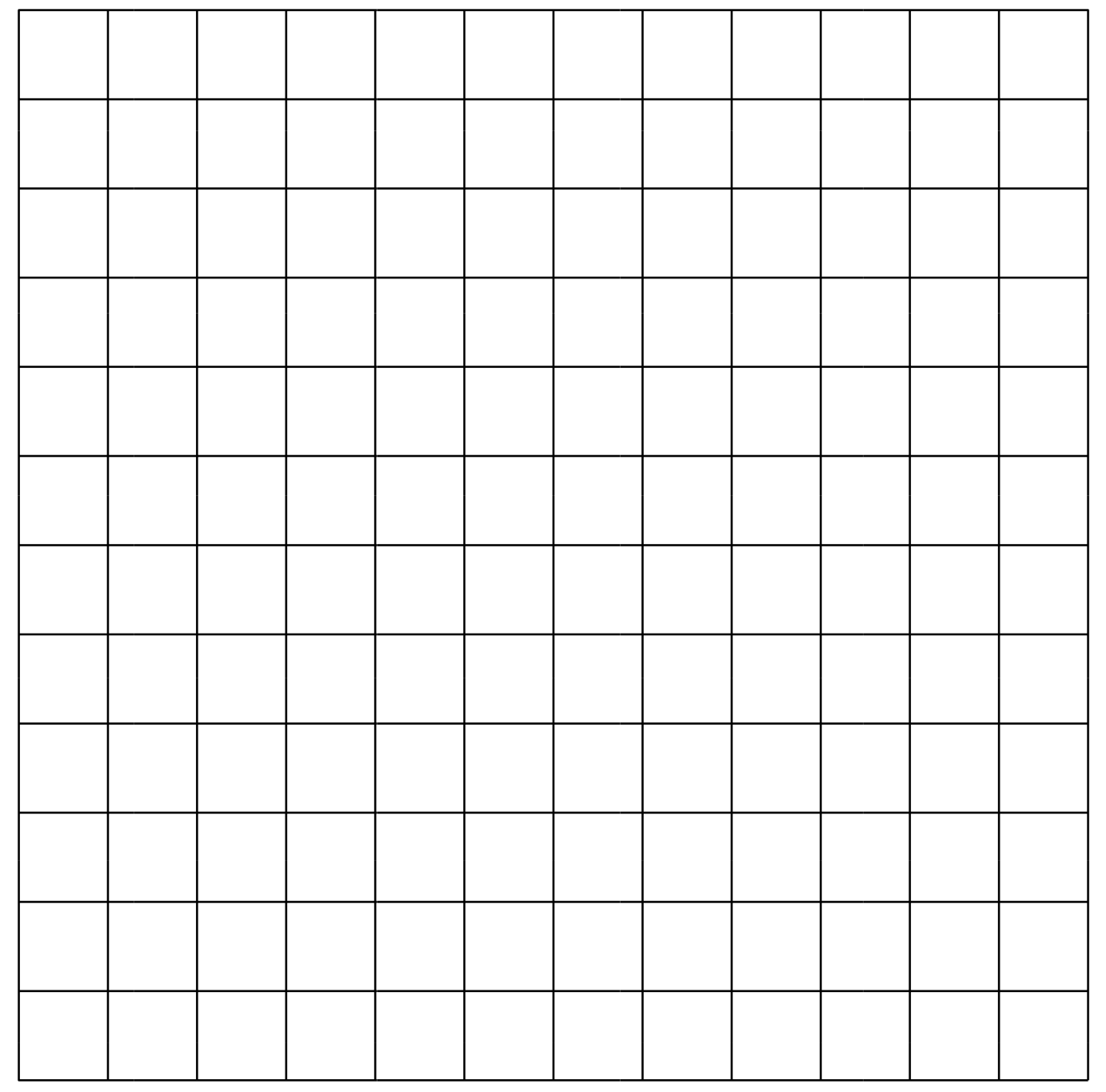}
    }
    \caption{Test case 3: uniform $h$-refinement.\label{fig:testcase2_h_ref}}
\end{figure}

 The meshes used for the $h$-refinement are shown in \fig{testcase2_h_ref}. The $\kord$-refinement is conducted for orders $\kord=2,\dots,5$.

\begin{figure}[!ht]
    \centering
    \subfigure[Stabilized VEM: uniform $h$-refinement.\label{fig:testcase2_error_h_stab}]{
        \includegraphics[width=0.47\textwidth]{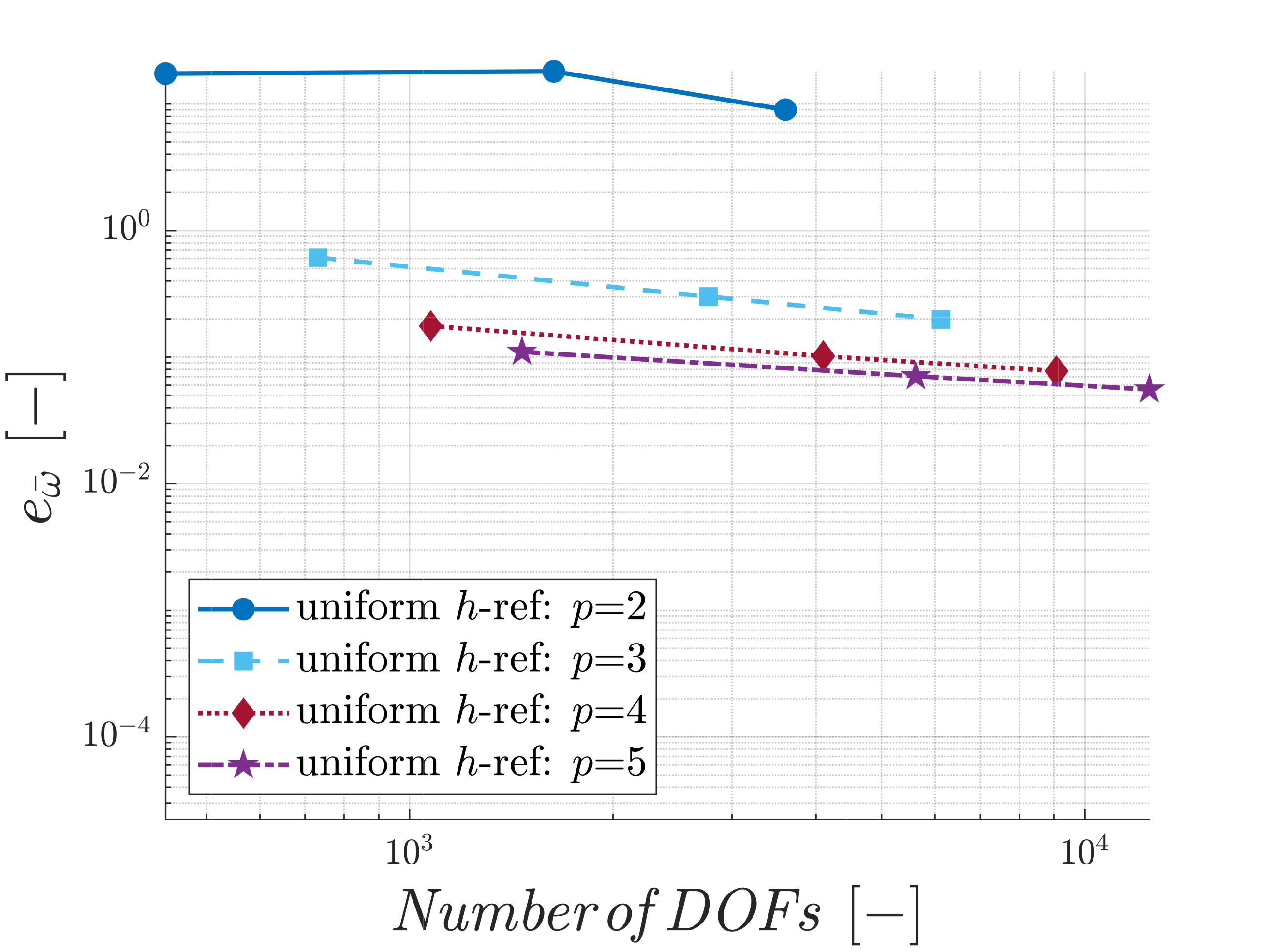}
    }
    \subfigure[Self-stabilized VEM: uniform $h$-refinement.\label{fig:testcase2_error_h_selfstab}]{
        \includegraphics[width=0.47\textwidth]{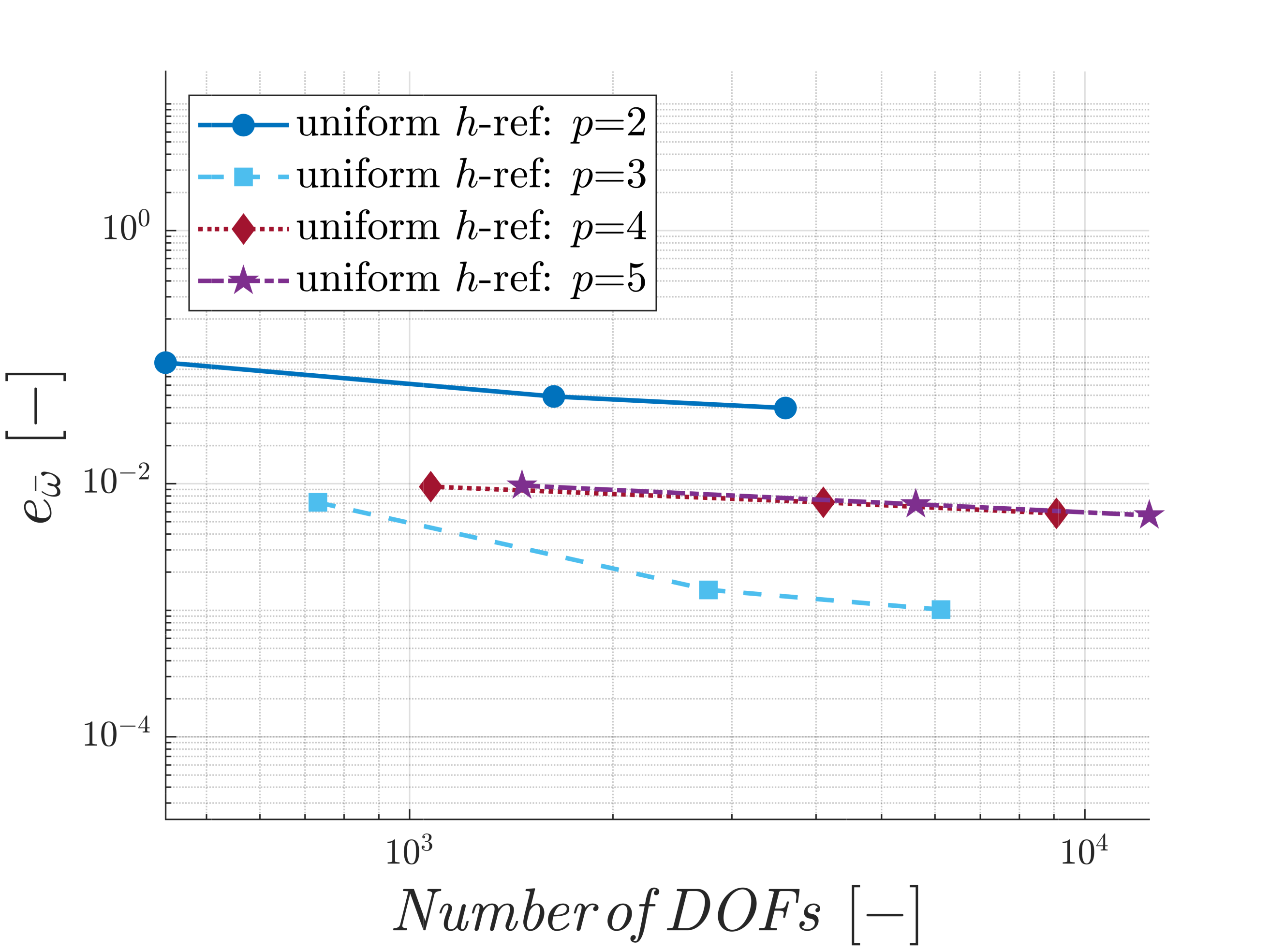}
    }\\
    \subfigure[Stabilized VEM: $\kord$-refinement.\label{fig:testcase2_error_p_stab}]{
        \includegraphics[width=0.47\textwidth]{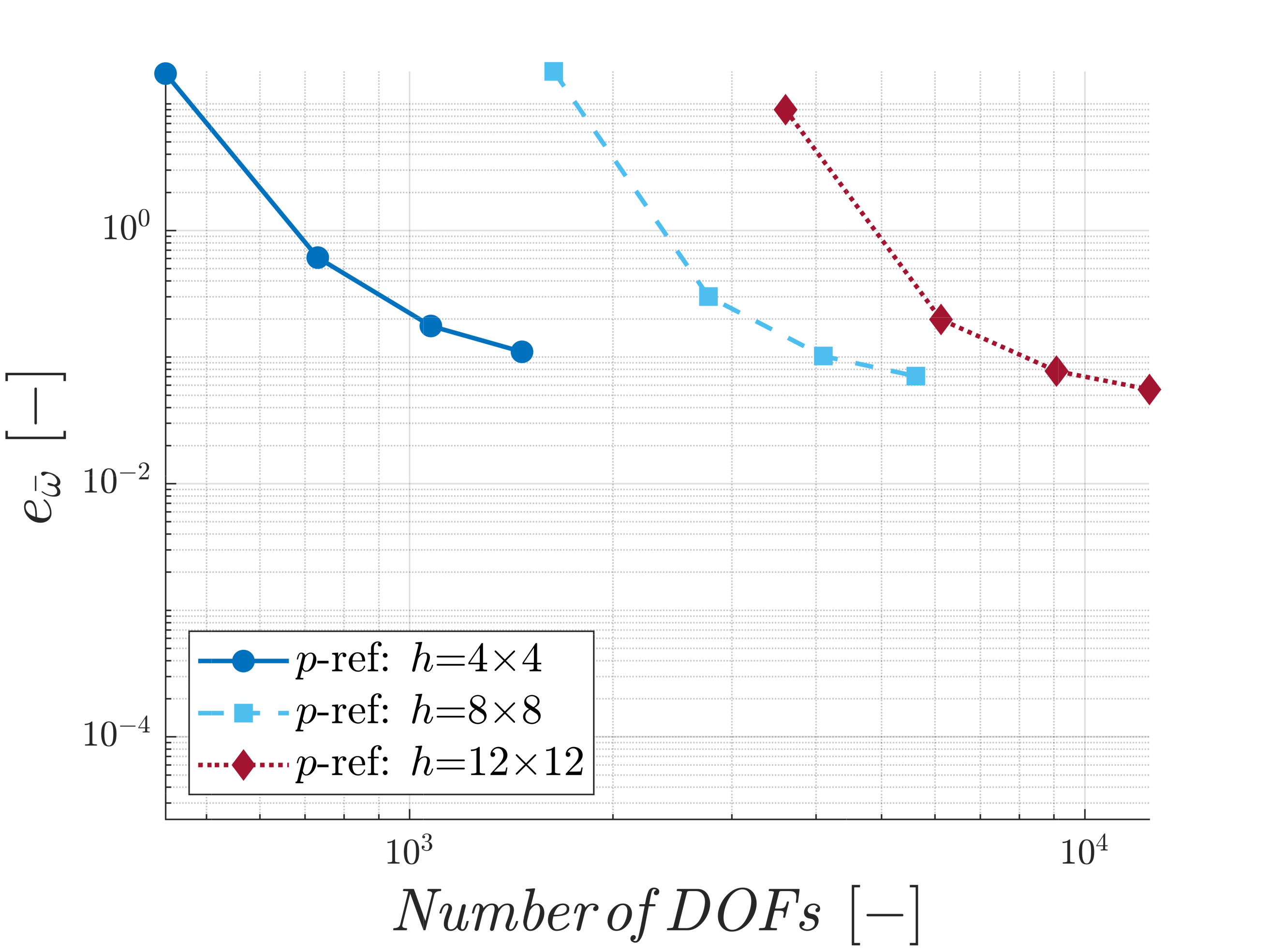}
    }
    \subfigure[Self-stabilized VEM: $\kord$-refinement.\label{fig:testcase2_error_p_selfstab}]{
        \includegraphics[width=0.47\textwidth]{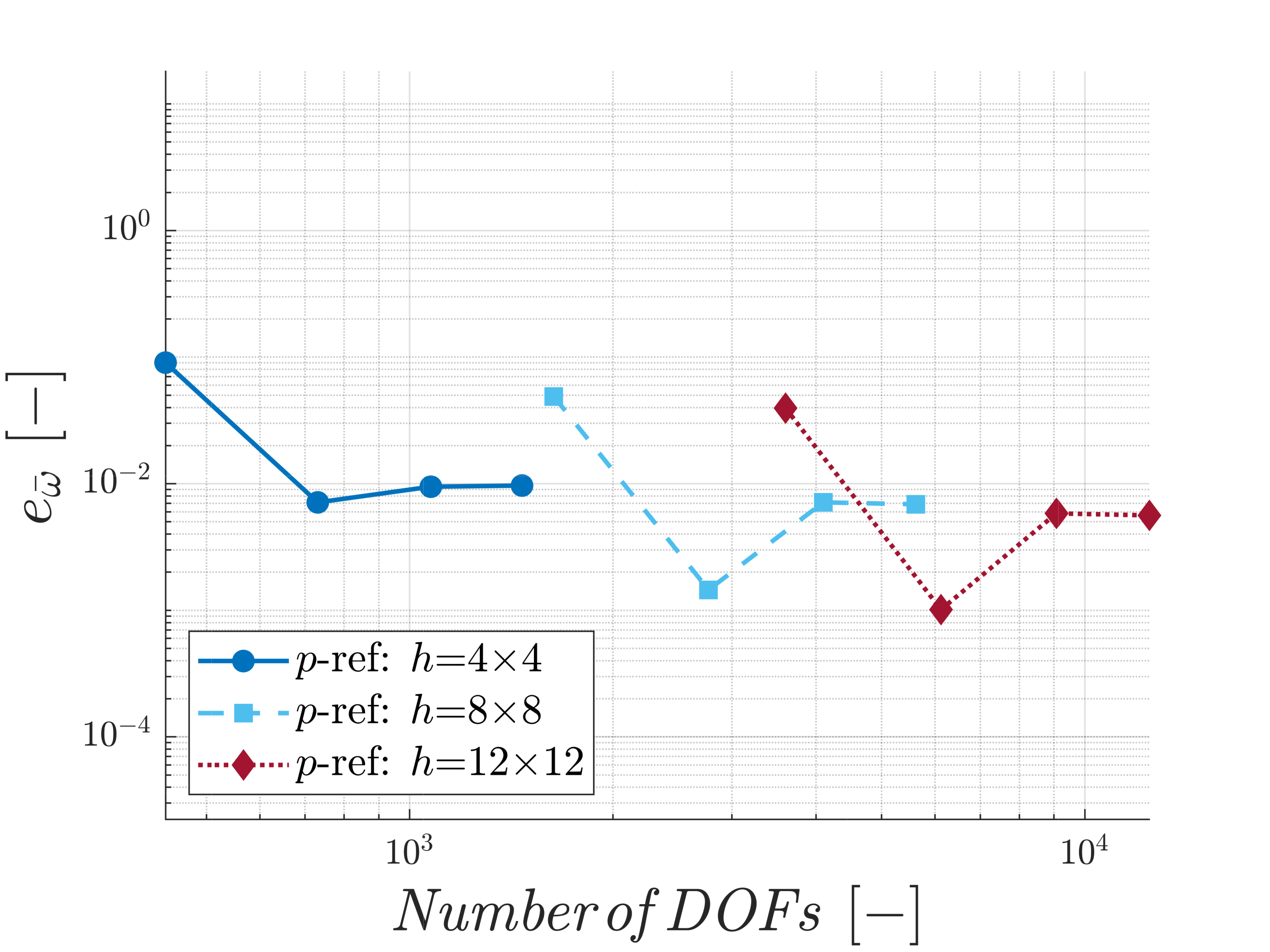}
    }
    \caption{Test case 3: error curves for uniform $h$- and $\kord$- refinements. \label{fig:testcase2_error_hp}}
\end{figure}

For both the uniform $h$ and $\kord$-refinement, the error plots are reported in \fig{testcase2_error_hp}. The number of degrees of freedom refers solely to the flexural ones.

By adopting a uniform $h$-refinement strategy, both the stabilized and the self-stabilized formulation exhibit a relatively flat trend, as revealed by \figs{testcase2_error_h_stab}{testcase2_error_h_selfstab}. This trend is ascribed to the presence of strong anisotropy-induced shear gradients in correspondence of the corners. Compared to the self-stabilized variant, the stabilized VEM yields higher errors. This behavior is motivated by the sensitivity of eigenvalue problems to the choice of the stabilization term. Therefore, self-stabilized formulations are preferable in these cases. 

Regarding the $\kord$-refinement in \figs{testcase2_error_p_stab}{testcase2_error_p_selfstab}, the stabilized VEM exhibits a faster error decay, but larger errors. Conversely, the self-stabilized formulation features some oscillations, but lower overall errors.  

\subsubsection*{Local $h$-refinement}

\begin{figure}[!htbp]
    \centering
    \subfigure[Refinement level $2$.\label{fig:testcase2_d_ref_2}]{
        \includegraphics[width=0.28\textwidth]{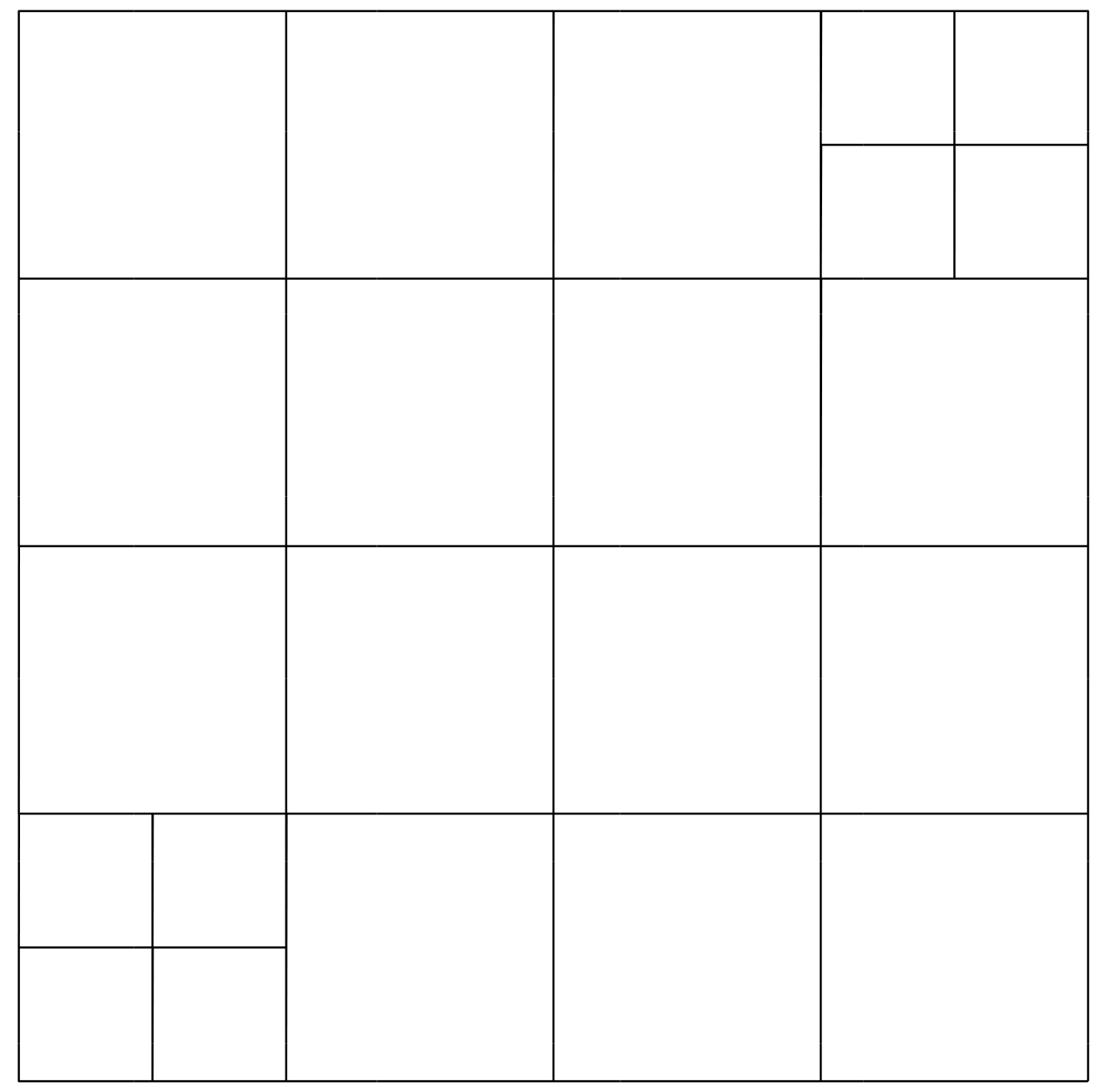}
    }
    \subfigure[Refinement level $5$.\label{fig:testcase2_d_ref_5}]{
        \includegraphics[width=0.28\textwidth]{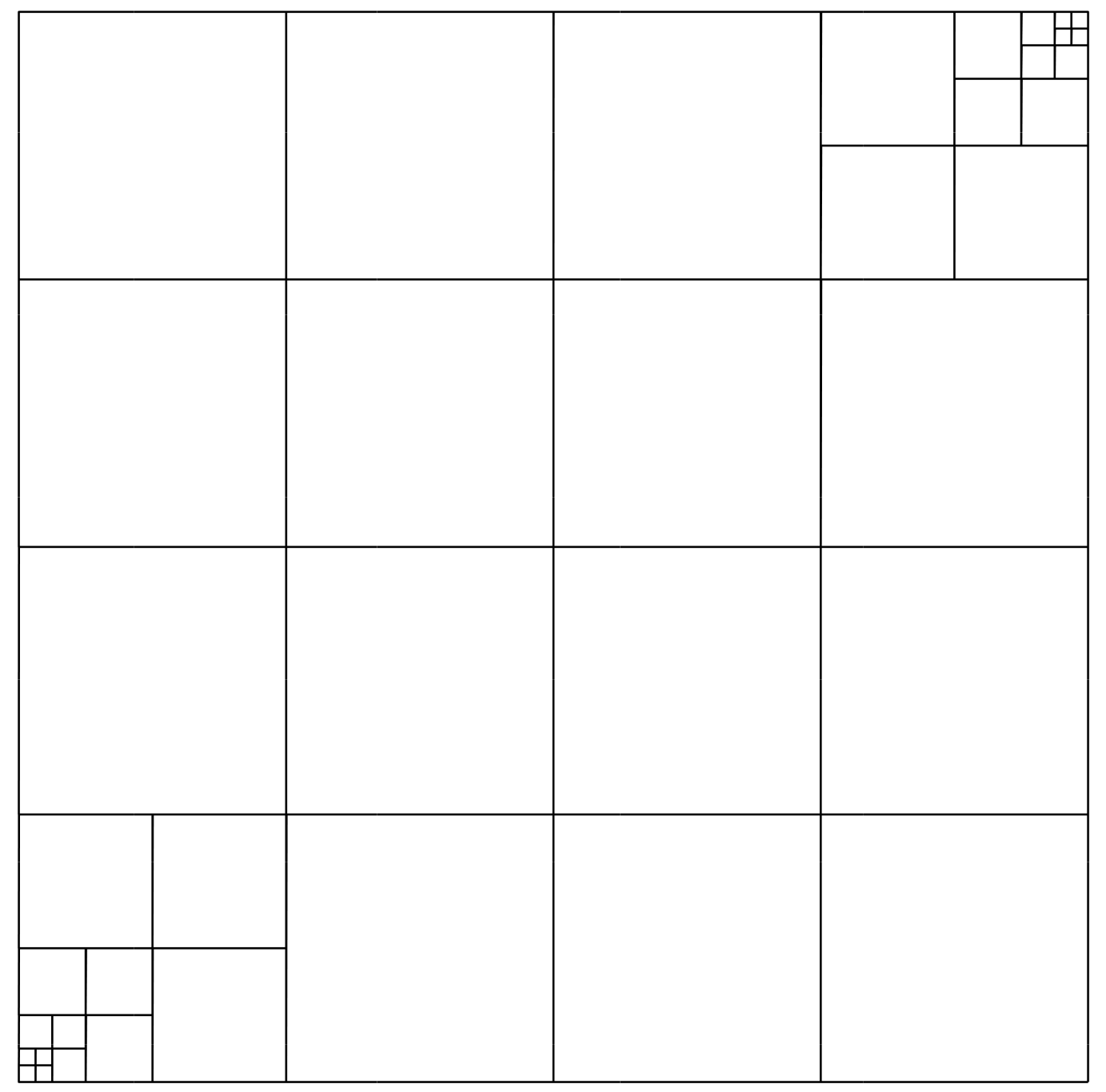}
    }
    \subfigure[Refinement level $10$.\label{fig:testcase2_d_ref_10}]{
        \includegraphics[width=0.28\textwidth]{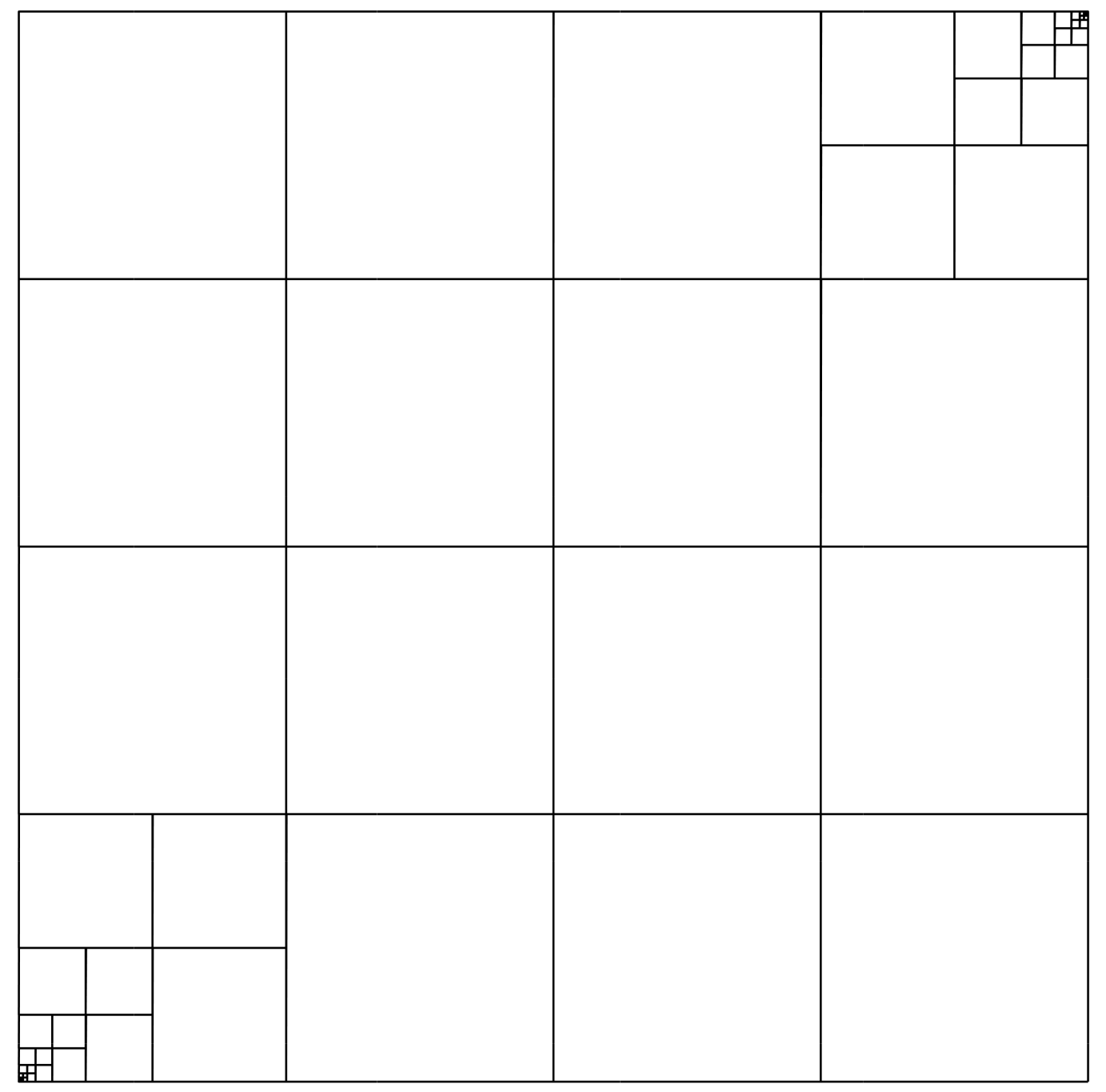}
    }
    \caption{Test case 3: local $h$-refinement for $h=4 \times 4$.}\label{fig:testcase2_d_ref}
\end{figure}

Improved results can be obtained by application of the local $h$-refinement strategy already presented in the previous test cases. For the problem at hand, local refinement is beneficial due to internal shear gradients in the proximity of the corners. For this reason, ten refinement levels are considered with increasing mesh density toward the corners. The different levels of refinement are shown in \fig{testcase2_d_ref}, where the mesh $4\times4$ is used as the initial reference discretization. 

\begin{figure}[!htbp]
    \centering
    \subfigure[Stabilized VEM.\label{fig:testcase2_error_pd_stab}]{
        \includegraphics[width=0.47\textwidth]{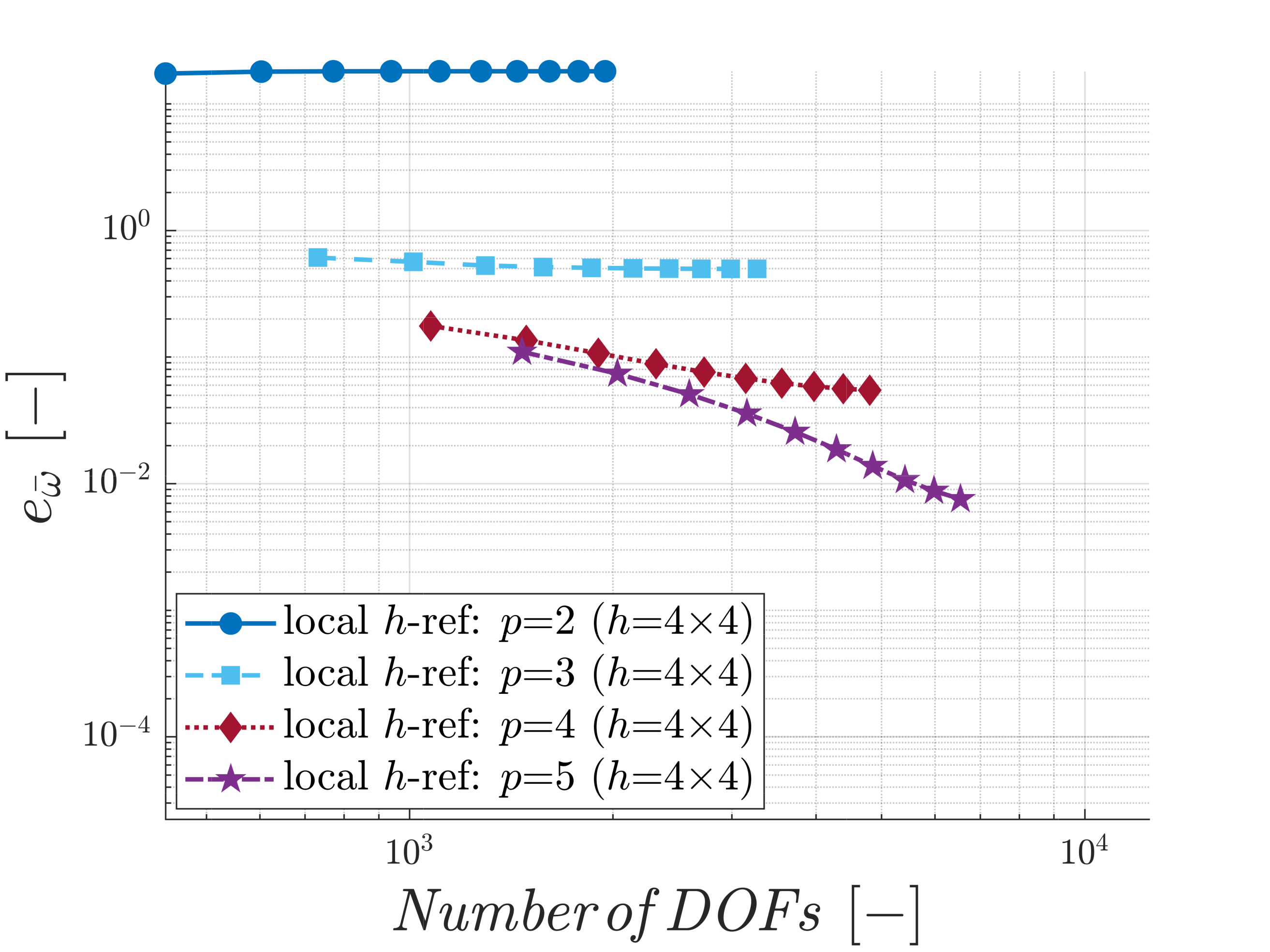}
    }
    \subfigure[Self-stabilized VEM.\label{fig:testcase2_error_pd_selfstab}]{
        \includegraphics[width=0.47\textwidth]{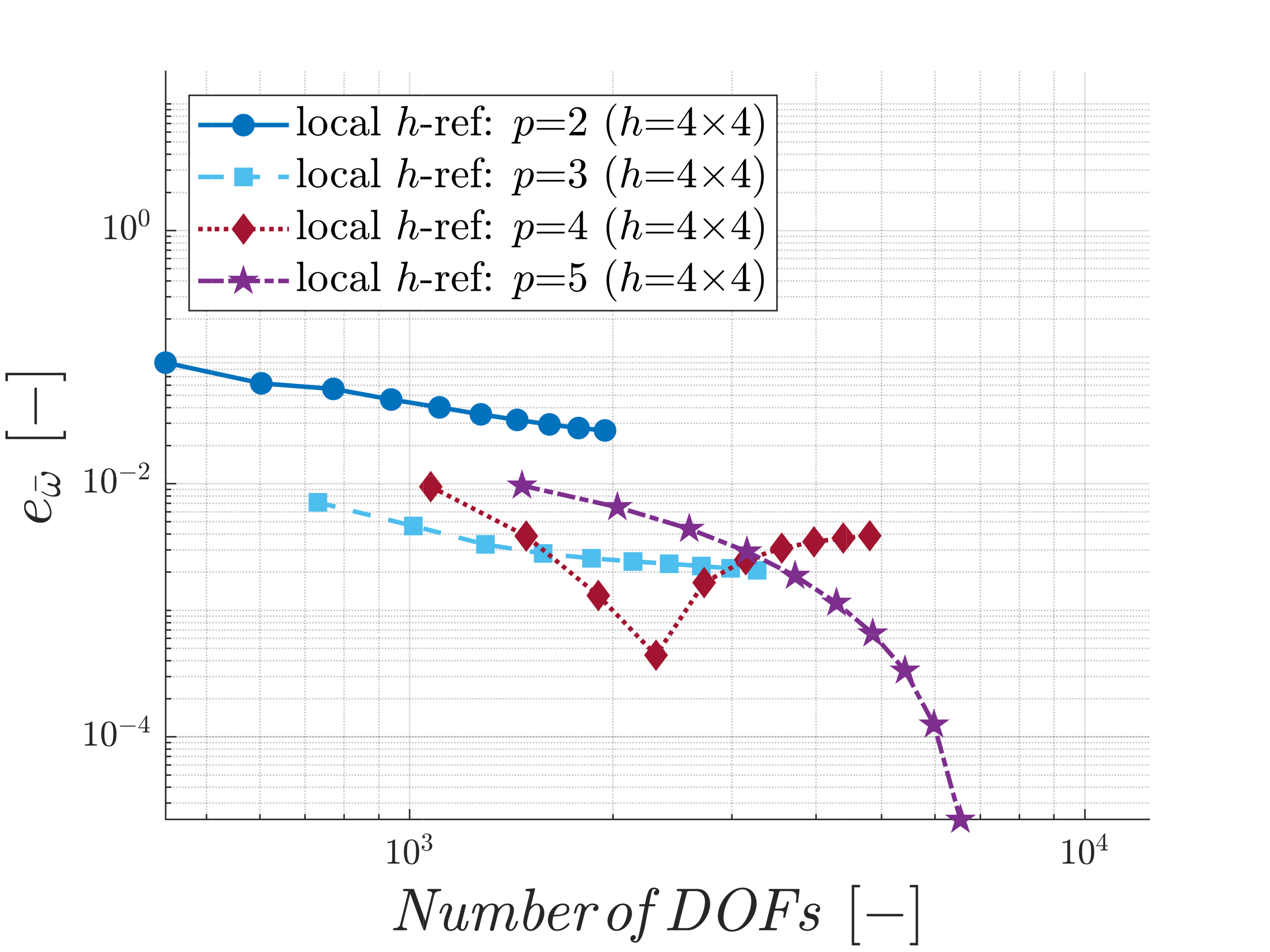}
    }
    \caption{Test case 3: error curves for local $h$-refinement: (a)-(b) fixed $h=4\times4$. \label{fig:testcase2_error_d}}
\end{figure}

The error curves are reported in \fig{testcase2_error_d}, where a fixed mesh $h=4\times4$ and different polynomial orders up to $\kord=5$ are considered.

The stabilized VEM features a flattening trend, indicating a slow convergence rate, as shown in \fig{testcase2_error_pd_stab}. Conversely, self-stabilized VEM achieves lower errors overall, but with oscillations and occasional error increases, as seen by inspection of \fig{testcase2_error_pd_selfstab}. The convergence rate improves significantly as the order is increased up to $\kord=5$. This is consistent with the trend observed for progressively refined models, which tend to yield lower frequency estimates, suggesting convergence toward a more accurate solution. 

\subsubsection*{Remarks}
Consistently with the previous test case, this example further demonstrates the effectiveness of the virtual element method in enabling local $h$-refinement to capture localized behavior. Since eigenvalue problems are particularly sensitive to the choice of the stabilization term, some differences can be observed between the stabilized and self-stabilized formulation. In particular, the stabilized VEM shows smoother convergence but with higher errors, whereas the self-stabilized VEM achieves lower errors although it exhibits some oscillations.

\subsection{Test case 4}

Having established the convergence properties of the method, its numerical performance is further assessed by investigating the response of plates with complex geometries. This benchmark, taken from~\cite{devarajan2021free,devarajan2020thermal}, considers the free-vibration and thermal buckling analysis of a square plate with a heart-shaped cutout. 

The plate is square with side length $a=10$~m. The heart-shaped cutout is defined as a NURBS curve in the reference, whereas in the present work it is constructed by combining three simple geometries: a square of side $a_1=4$~m, centered at $\plbr{x_s,y_s}=\plbr{4,6}$~m, and two circles of radius $R_1=R_2=a_1/2$, centered at $\plbr{x_{c1},y_{c1}}=\plbr{4,4}$~m and $\plbr{x_{c2},y_{c2}}=\plbr{6,6}$~m, respectively. 

\begin{figure}[!htbp]
    \centering
    \subfigure[Free-vibration.\label{fig:testcase3_configuration_freevib}]{
        \includegraphics[width=0.4\textwidth]{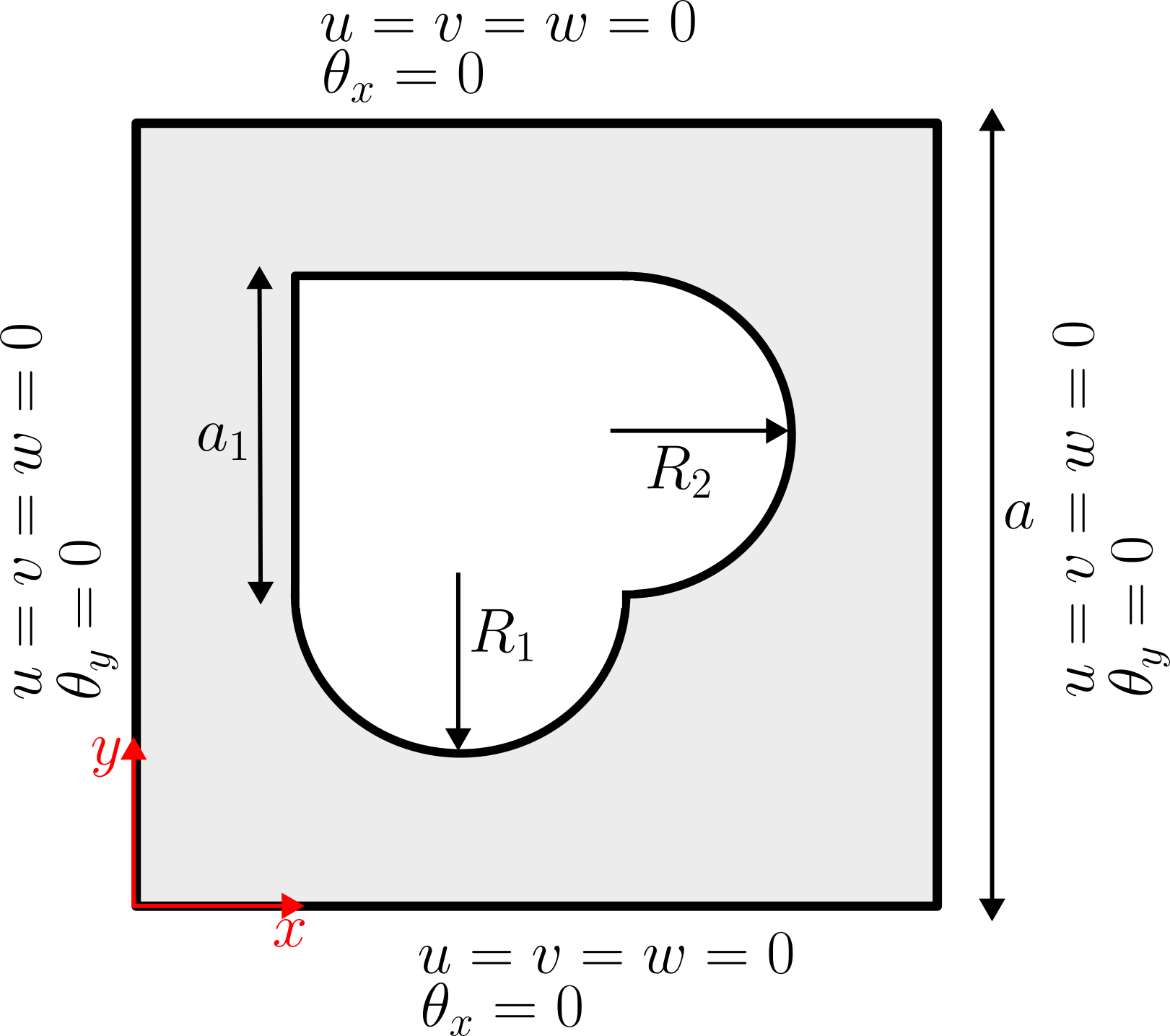}
    }
    \subfigure[Buckling.\label{fig:testcase3_configuration_buck}]{
        \includegraphics[width=0.4\textwidth]{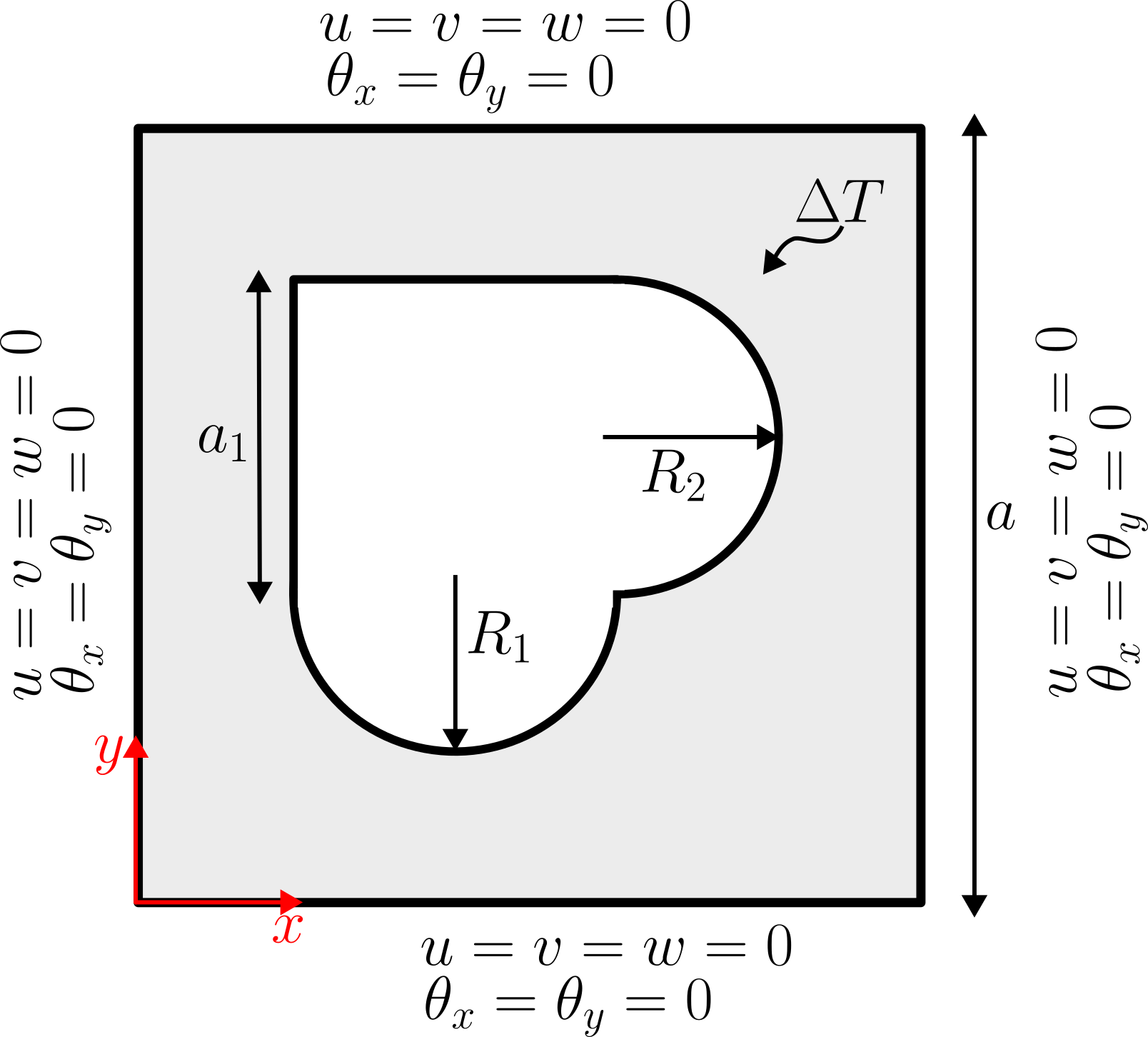}
    }
    \caption{Test case 4: configuration.\label{fig:testcase3_configuration}}
\end{figure}

Two configurations are considered with equal geometry, but different boundary conditions depending on the analysis type. 
In particular, the plate is simply supported or fully clamped, dependently on whether free-vibration or buckling analysis are considered. The dimensions and boundary conditions are summarized in \fig{testcase3_configuration}.

The material properties and layups are specified in the respective subsections.

\begin{figure}[!htbp]
    \centering
    \subfigure[Mesh 1.\label{fig:testcase3_heart_mesh_less}]{
        \includegraphics[width=0.28\textwidth]{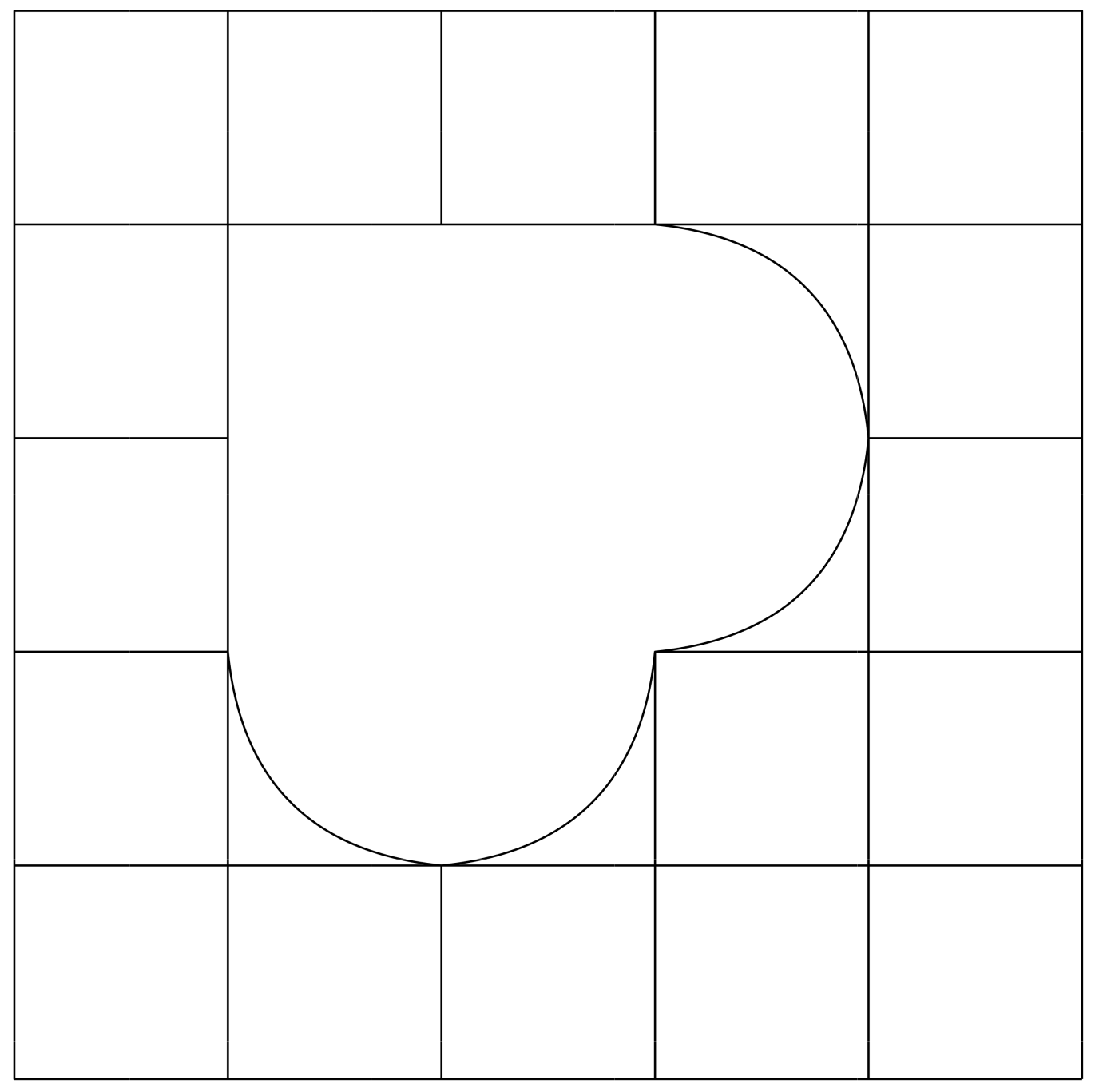}
    }
    \subfigure[Mesh 2.\label{fig:testcase3_heart_mesh_more}]{
        \includegraphics[width=0.28\textwidth]{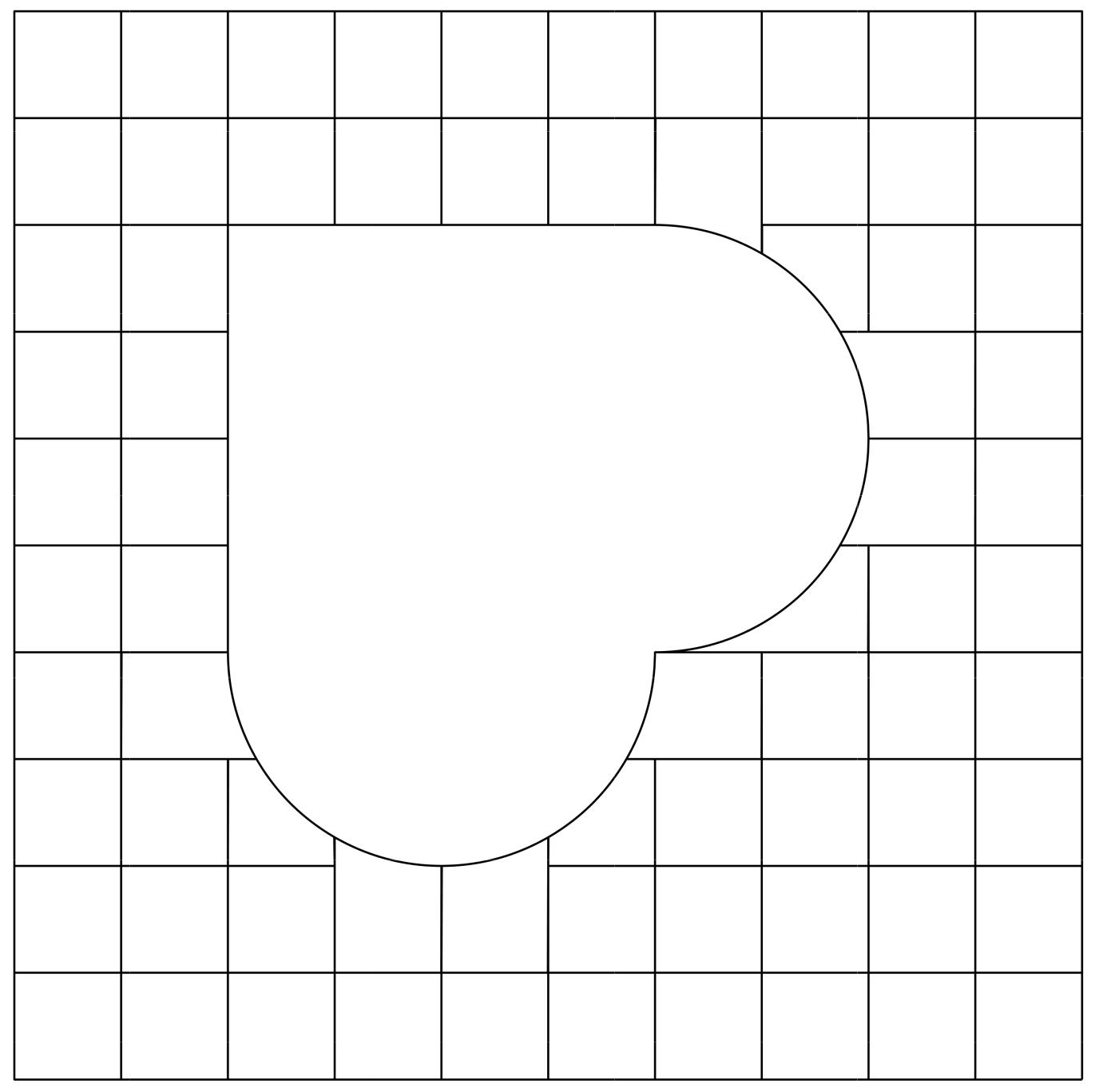}
    }
    \caption{Test case 4: meshes. The elements' circular arcs are interpolated as B\'ezier curves..\label{fig:testcase3_heart_mesh}}
\end{figure}

Two structured meshes are employed, as shown in \fig{testcase3_heart_mesh}. The first mesh is relatively coarse and is employed in combination with high approximation orders, whereas the second is finer and coupled with lower approximation orders. This setup enables an assessment of the influence of both mesh density and approximation order $\kord$ on the accuracy and computational efficiency of the method.
Particular attention should be paid to the presence of elements with curved edges. Their treatment is particularly straightforward within the VEM framework: in the coarse mesh, such elements are relatively irregular, yet accurately captured thanks to the proposed curved-edge formulation; in the finer mesh, a different scenario is observed, where hanging nodes prove useful during the meshing process.

\subsubsection*{Free-vibrations}

For the free-vibration analysis, a composite material with properties $E_{11}/E_{22}=2.45$, $G_{12}/E_{22}=0.48$, $G_{13}/E_{22}=G_{23}/E_{22}=0.2$, $\nu_{12}=0.23$, $\rho=8000$~kg/m$^3$, and thickness $h=0.06$~m, is considered. Three different layups are analyzed:
\begin{equation}
    \begin{aligned}
     &\text{layup 1}: [15^\circ/-15^\circ/15^\circ], 
     &&\text{layup 2}: [30^\circ/-30^\circ/30^\circ], 
     &\text{layup 3}: [45^\circ/-45^\circ/45^\circ] .
    \end{aligned}
    \label{eq:testcase3_layups}
\end{equation}
The frequencies are normalized as~\cite{devarajan2021free}:
\begin{equation}
    \bar{\omega}_i = \plbr{\frac{ \rho \, h \, \omega_i \, a^{4} } { \overline{D} }}^{\frac{1}{2}}, \quad \text{with} \quad \overline{D} = \frac{ E_{11} \, h^{3} }{ 12 \plbr{ 1 - \nu_{12} \, \frac{E_{22}}{E_{11}} \, \nu_{12} } }.
    \label{eq:testcase3_nornalizefreq_comp}
\end{equation}

In order to simplify the comparison across the three layups, only Mesh 2 with approximation order $\kord=5$ is employed. The standard stabilized VEM is employed.

The normalized natural frequencies for the three layups are reported in \tab{freevib_modes_comp_heart} and compared with those obtained via isogeometric analysis (IGA) in~\cite{devarajan2021free}. 
\begin{table}[!htbp]
\centering
\begin{tabular}{cccccccccccc}
\hline
\multirow{2}{*}{\textbf{Mode}} 
& \multicolumn{3}{c}{{layup 1}} 
&& \multicolumn{3}{c}{{layup 2}} 
&& \multicolumn{3}{c}{{layup 3}} \\
\cline{2-4} \cline{6-8} \cline{10-12}
& {VEM} & {IGA}~\cite{devarajan2021free} & {Error [\%]}
&& {VEM} & {IGA}~\cite{devarajan2021free} & {Error [\%]}
&& {VEM} & {IGA}~\cite{devarajan2021free} & {Error [\%]} \\
\hline
1 & 19.04 & 18.91 & 0.69 && 20.53 & 20.40 & 0.64 && 21.24 & 21.10 & 0.66 \\
2 & 31.70 & 31.83 & 0.41 && 33.51 & 33.66 & 0.45 && 34.30 & 34.45 & 0.44 \\
3 & 35.69 & 36.09 & 1.11 && 36.82 & 37.23 & 1.10 && 37.49 & 37.91 & 1.11 \\
4 & 56.98 & 57.00 & 0.04 && 59.19 & 59.20 & 0.02 && 60.29 & 60.29 & 0.00 \\
5 & 61.93 & 62.73 & 1.28 && 64.09 & 65.03 & 1.45 && 65.20 & 66.22 & 1.54 \\
6 & 83.78 & 83.93 & 0.18 && 87.84 & 87.92 & 0.09 && 91.02 & 90.68 & 0.37 \\
\hline
\end{tabular}
\caption{Test case 4: nondimensional natural frequencies $\bar{\omega}$ and normalized percentage errors with respect to the reference values.}
\label{tab:freevib_modes_comp_heart}
\end{table}
In all cases, the results obtained with the present VEM formulation show good agreement with the reference results from~\cite{devarajan2021free}, with a relative error below $1\%$ in most cases. In particular, for all layups, the highest error is observed for Mode 5, with a value slightly higher than $1\%$.

\subsubsection*{Thermal buckling}

For the thermal buckling benchmark~\cite{devarajan2020thermal}, the material properties are $E_{11}/E_{22}=15$, $G_{12}/E_{22}=0.5$, $G_{13}/E_{22}=G_{23}/E_{22}=0.3356$, $\nu_{12}=0.3$, with $E_{22}=1$~GPa, and thermal expansion coefficients $\alpha_{11}/\alpha_{0}=0.015$ and $\alpha_{22}/\alpha_{0}=1$, with $\alpha_{0}=0.025$ $1/^\circ$C. The thickness is $h=0.1$~m, and the laminate layup is $[0^\circ /90^\circ /90^\circ /0^\circ]$. The plate is loaded with a temperature gradient $\Delta T$.

Both meshes shown in \fig{testcase3_heart_mesh} are employed in combination with the stabilized VEM. The approximation orders are selected to balance accuracy and computational efficiency, and the two following models are considered: Mesh 1 with a higher-order approximation ($\kord=8$), and Mesh 2 with a lower-order approximation ($\kord=5$). 

The comparison is performed with an Abaqus model consisting of $2032$ S4R shell elements. 

The nondimensional critical temperatures corresponding to the first five buckling modes are summarized in \tab{buckling_heart}. 
The results are obtained with the standard stabilized VEM approach. So, the spatial dependency of the pre-stress distribution is neglected in the computation of the geometric stiffness matrix projection operator. 
 
\begin{table}[!htbp]
\centering
\begin{tabular}{cccccc}
\hline
\multirow{2}{*}{\textbf{Mode}} 
& \multicolumn{2}{c}{{Mesh 1 and {$\kord=8$}}} 
& \multicolumn{2}{c}{{Mesh 2 and {$\kord=5$}}} 
& \multirow{2}{*}{{Abaqus}} \\
\cline{2-5}
& {VEM} & {Error [\%]} & {VEM} & {Error [\%]} & \\
\hline
1 & 0.01197 & 0.50 & 0.01206 & 0.25 & 0.01203 \\
2 & 0.01242 & 0.64 & 0.01254 & 0.32 & 0.01250 \\
3 & 0.01308 & 0.15 & 0.01315 & 0.38 & 0.01310 \\
4 & 0.01320 & 0.00 & 0.01324 & 0.30 & 0.01320 \\
5 & 0.01438 & 0.69 & 0.01451 & 0.21 & 0.01448 \\
\hline
\end{tabular}
\caption{Test case 4: nondimensional critical temperatures $\alpha_0 \Delta T^{cr}$ and normalized percentage errors with respect to the reference values.}
\label{tab:buckling_heart}
\end{table}
The results obtained with the present VEM formulation closely match the Abaqus ones. In all cases, the relative error in the nondimensional critical temperature is below $1 \%$ for every mode. Comparing the two meshes, the combination of the coarser Mesh 1 with the higher approximation order $\kord=8$ yields slightly lower nondimensional critical temperatures than those obtained with the finer Mesh 2 with $\kord=5$. This highlights that coarse discretization with high-order approximations can provide accurate results with a reduced number of degrees of freedom. Indeed, for Mesh 2 with $\kord=5$, the number of degrees of freedom is $9684$, whereas for Mesh 1 with $\kord=8$ it is $5877$. 

\begin{figure}[!htbp]
    \centering
    \subfigure[Mode 1 VEM.\label{fig:buck_heart_1_VEM}]{
        \includegraphics[width=0.175\textwidth]{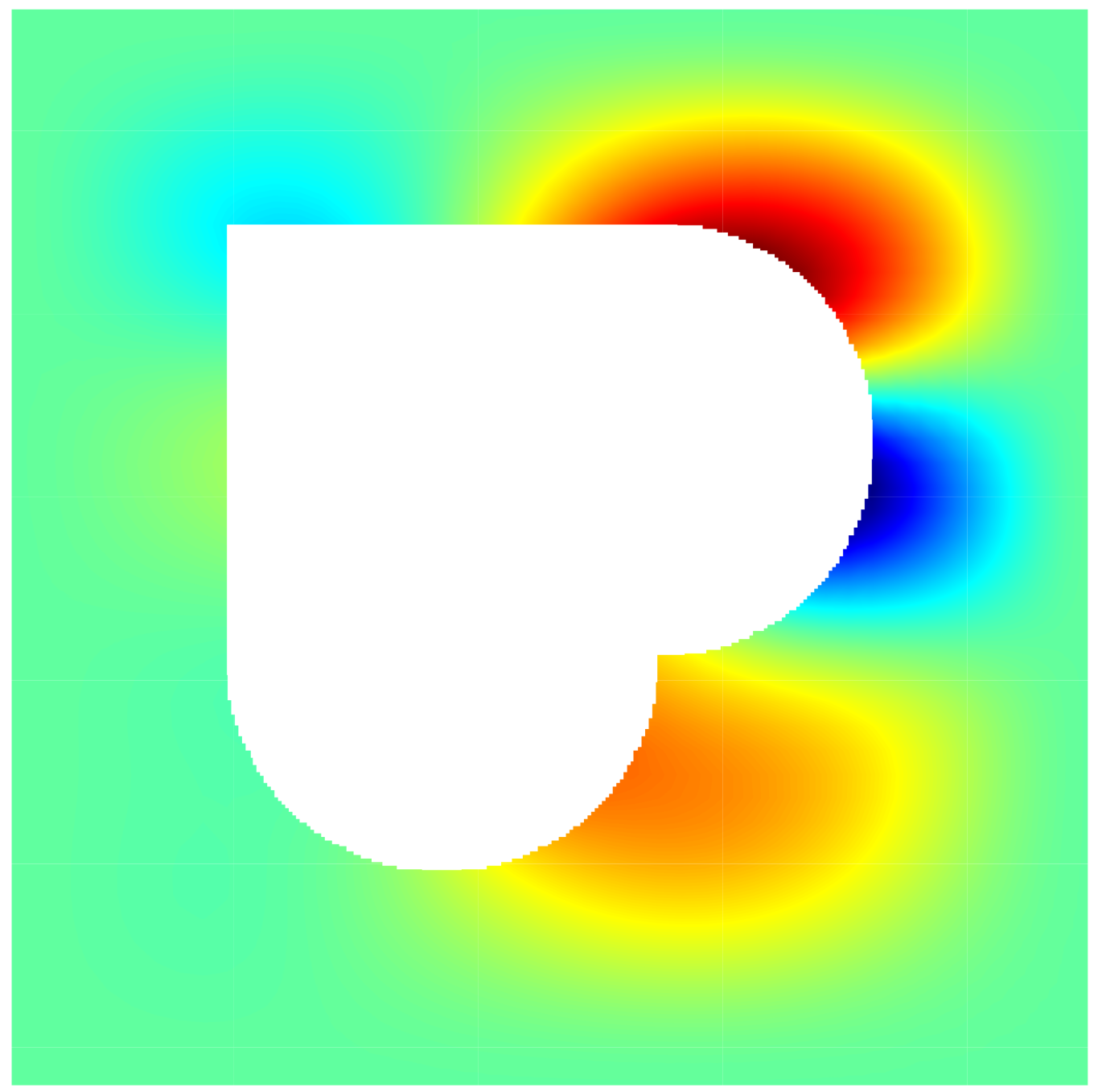}
    }    
    \subfigure[Mode 2 VEM.\label{fig:buck_heart_2_VEM}]{
        \includegraphics[width=0.175\textwidth]{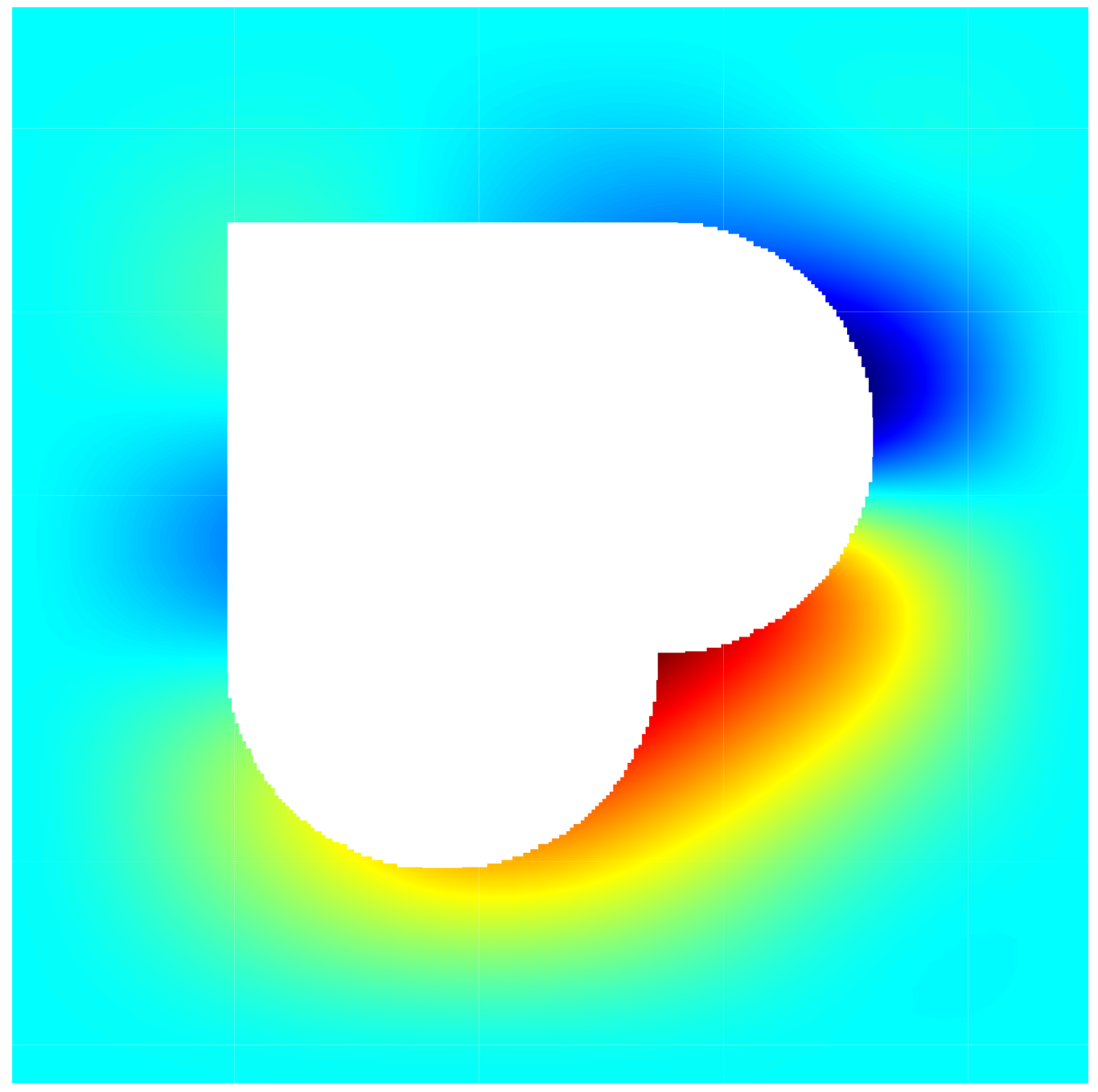}
    }    
    \subfigure[Mode 3 VEM.\label{fig:buck_heart_3_VEM}]{
        \includegraphics[width=0.175\textwidth]{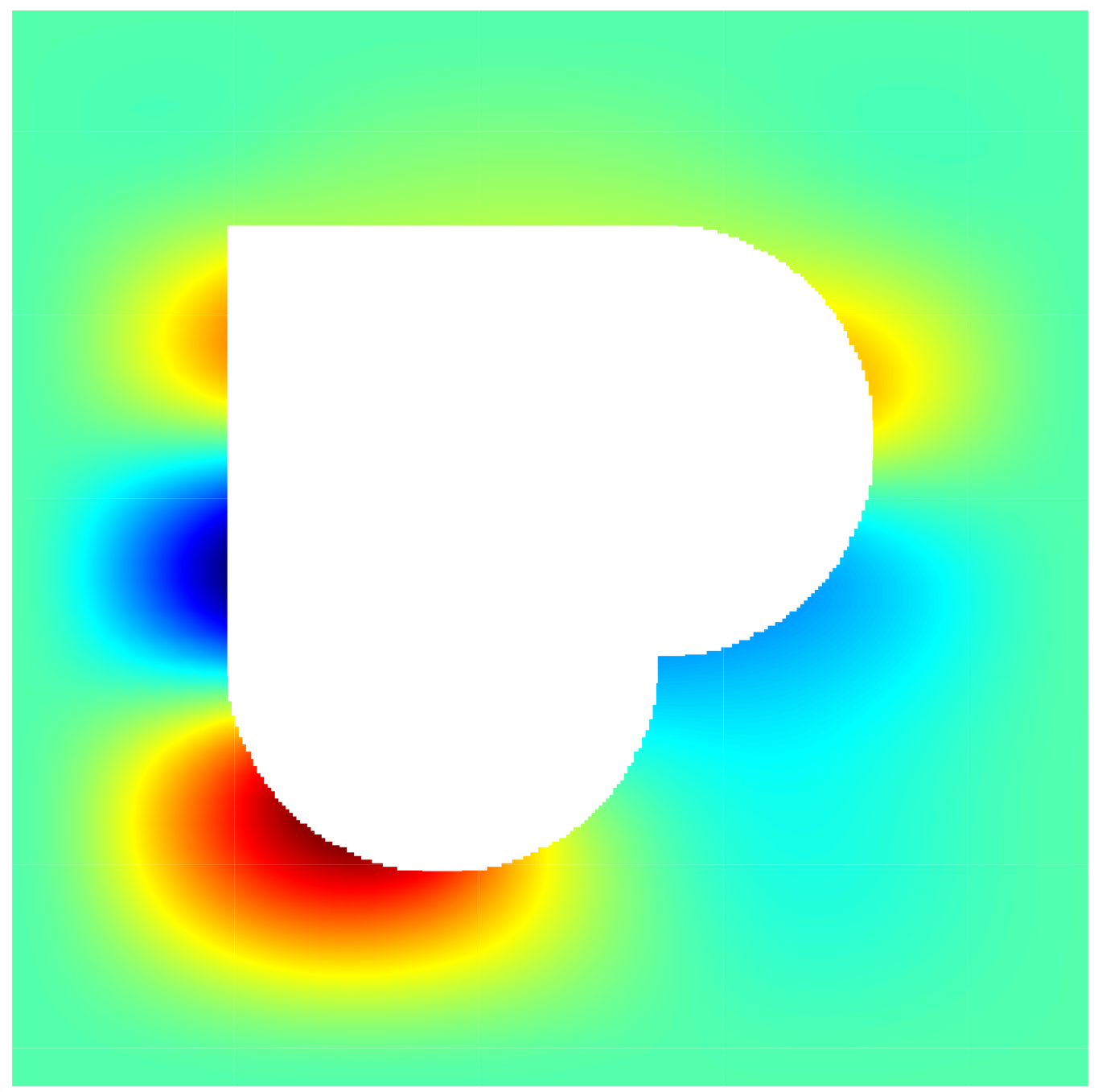}
    }    
    \subfigure[Mode 4 VEM.\label{fig:buck_heart_4_VEM}]{
        \includegraphics[width=0.175\textwidth]{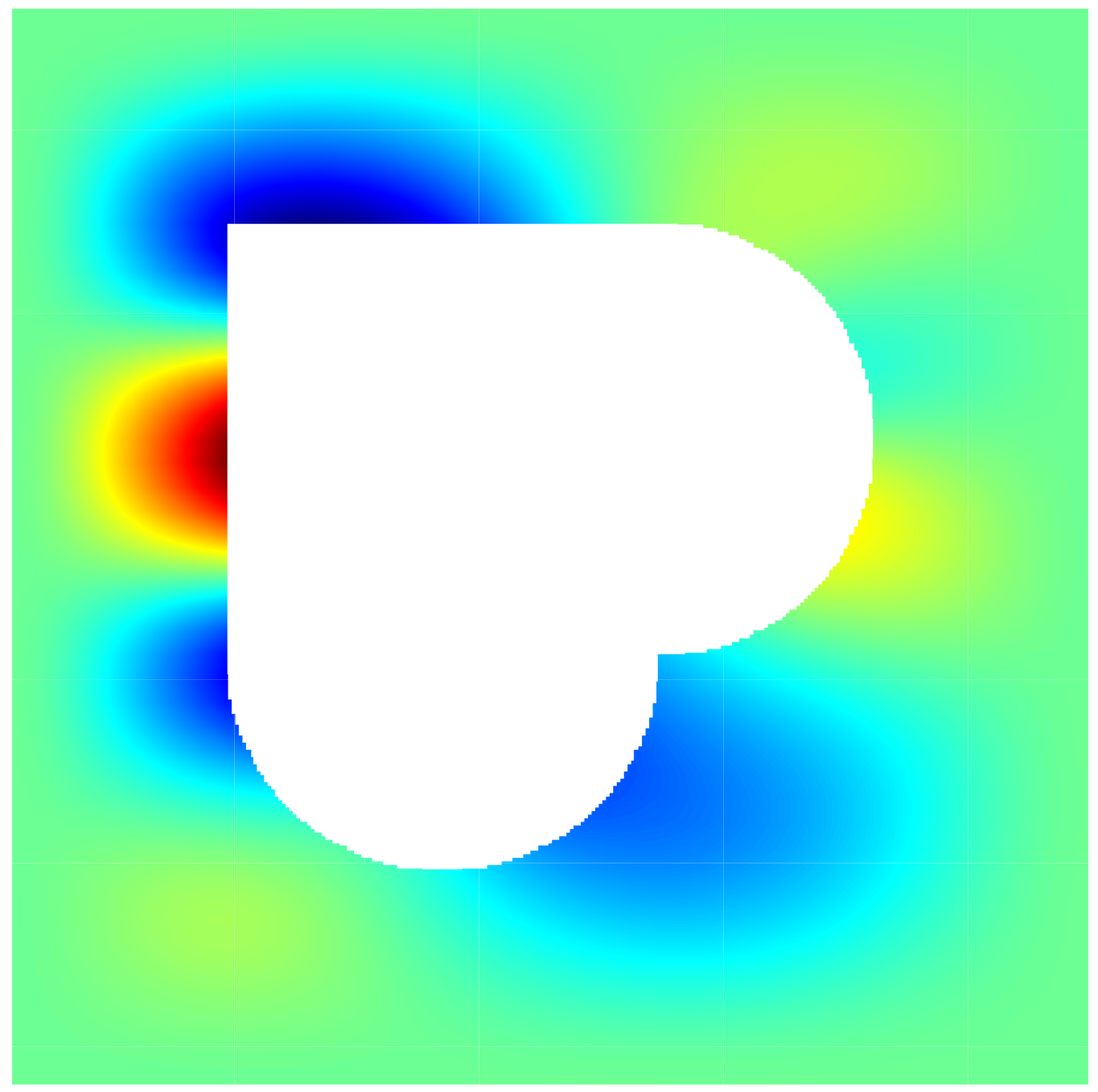}
    }    
    \subfigure[Mode 5 VEM.\label{fig:buck_heart_5_VEM}]{
        \includegraphics[width=0.175\textwidth]{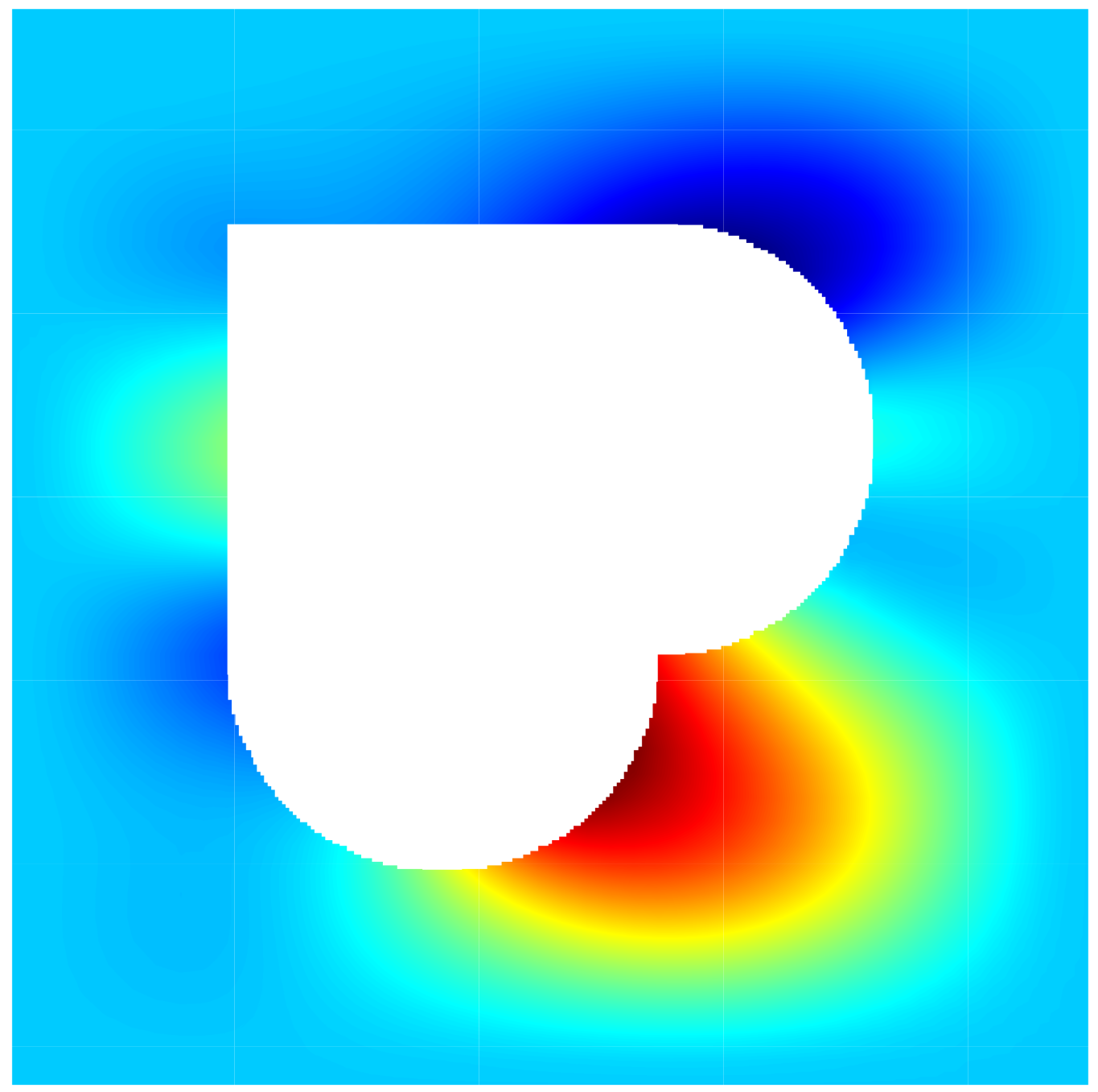}
    }    
    \\
    \subfigure[Mode 1 Abaqus.\label{fig:buck_heart_1_Abaqus}]{
        \includegraphics[width=0.175\textwidth]{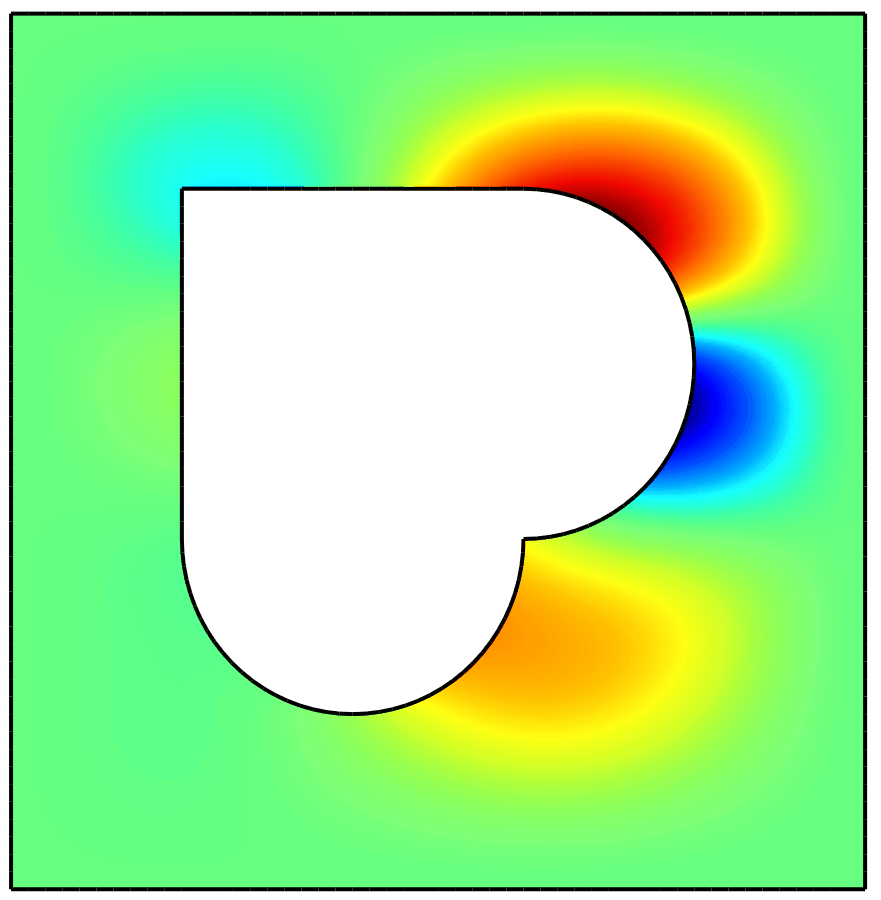}
    }    
    \subfigure[Mode 2 Abaqus.\label{fig:buck_heart_2_Abaqus}]{
        \includegraphics[width=0.175\textwidth]{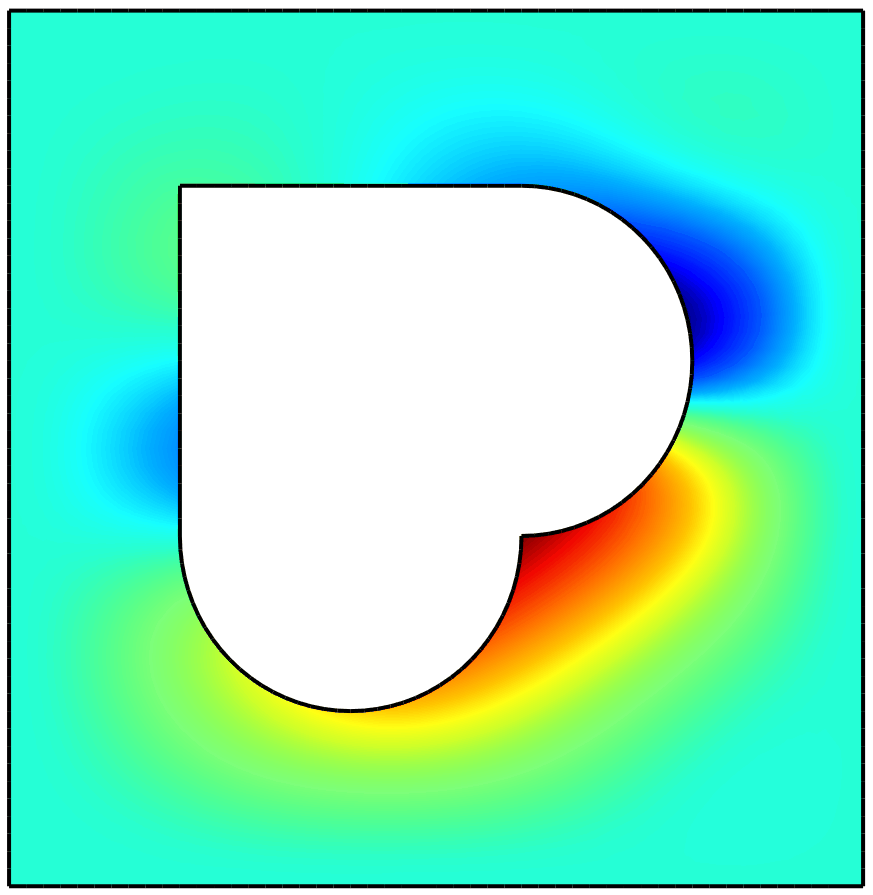}
    }    
    \subfigure[Mode 3 Abaqus.\label{fig:buck_heart_3_Abaqus}]{
        \includegraphics[width=0.175\textwidth]{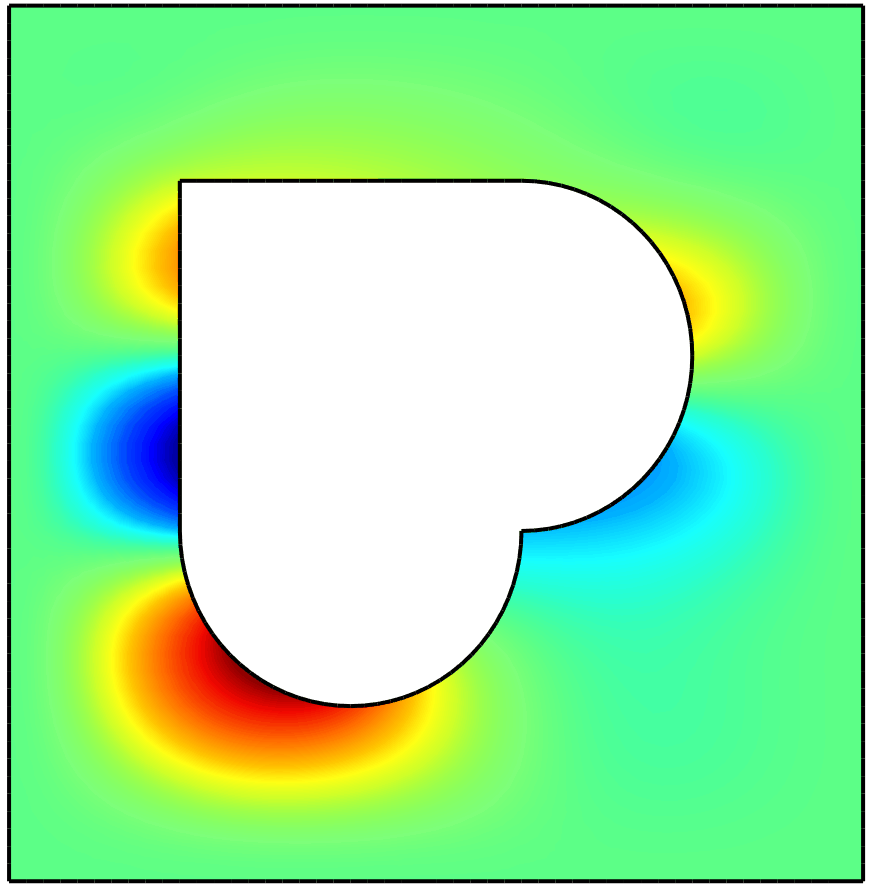}
    }    
    \subfigure[Mode 4 Abaqus.\label{fig:buck_heart_4_Abaqus}]{
        \includegraphics[width=0.175\textwidth]{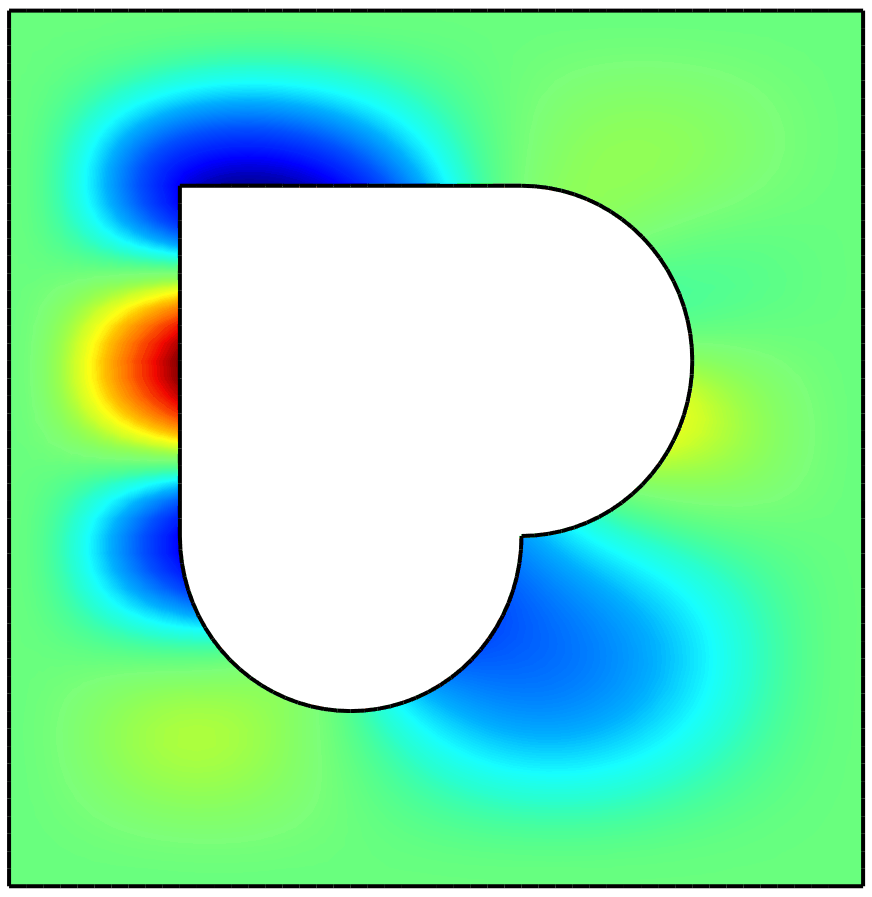}
    }    
    \subfigure[Mode 5 Abaqus.\label{fig:buck_heart_5_Abaqus}]{
        \includegraphics[width=0.175\textwidth]{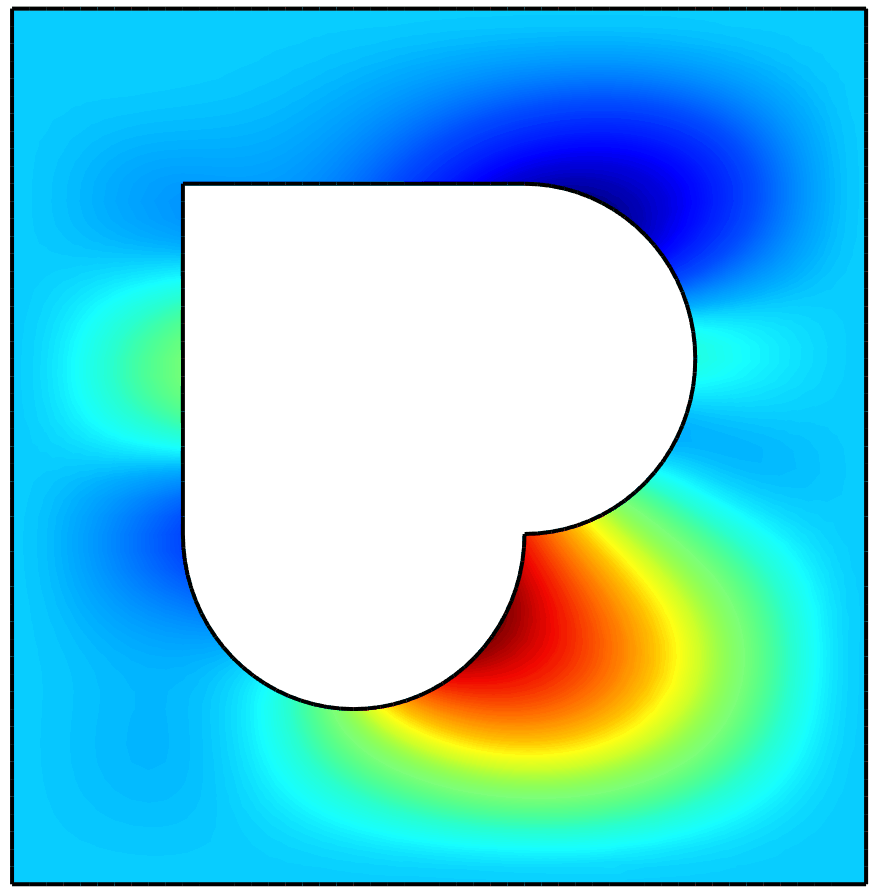}
    }    
    \caption{Test case 4: buckling modes ($w$).\label{fig:buck_heart}}
\end{figure}

The VEM-Abaqus comparison in terms of buckling modes is provided in \fig{buck_heart}. Mesh 1 with $\kord=8$ is considered, although similar results are obtained with Mesh 1, but they are omitted here for the sake of brevity.

The contours illustrate the close agreement between the two computational strategies, further demonstrating the correctness of the proposed implementation. 

\subsection{Test case 5}

\begin{figure}[!htbp]
    \centering
        \includegraphics[width=0.45\textwidth]{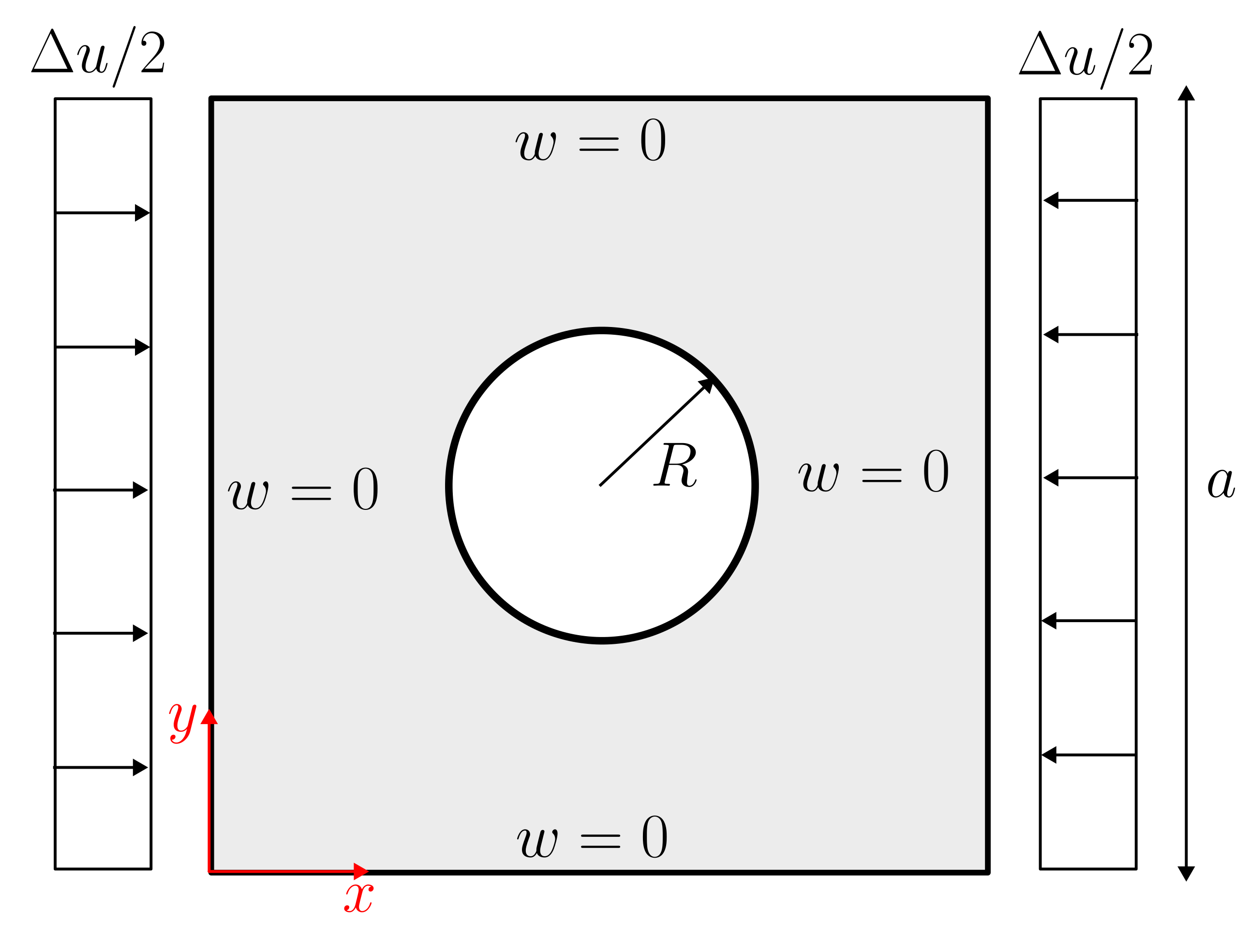}
    \caption{Test case 5: configuration.\label{fig:testcase5_configuration}}
\end{figure}

To further assess the effectiveness of the developed tool, the validation is now extended to variable stiffness plates with complex geometries. The test case is taken from~\cite{milazzo2023buckling}, where the Ritz method is employed to investigate the buckling of a square plate with a circular cutout. The plate is square with dimensions $a=254$~mm, with a circular cutout of radius $R=50.8$~mm centered in the middle of the plate. 
The plate is simply supported: the out-of-plane deflections are prevented, while the in-plane conditions are free apart from the normal component on the loaded edges, as shown in \fig{testcase5_configuration}.

The material properties are $E_{11}=181000$~MPa, $E_{22}=10273$~MPa, $G_{12}=G_{13}=G_{23}=7170.5$~MPa and $\nu_{12}=0.28$. The laminate consists of sixteen plies, each with thickness $h=0.1272$~mm. Four different layups are considered. In particular:
\begin{equation}
    \begin{aligned}
     \text{layup 1}: [\pm90\langle 0|75\rangle]_{4s}, \;\;
     \text{layup 2}: [\pm0\langle 60|30\rangle]_{4s}, \;\;
     \text{layup 3}: [\pm45\langle 15|0\rangle]_{4s}, \;\;
     \text{layup 4}: [\pm90\langle 45|45\rangle]_{4s}.
    \end{aligned}
    \label{eq:testcase5_layups}
\end{equation}

The first three layups are characterized by stiffness variability, while the last one is a constant stiffness configuration. The presence of both variable and constant stiffness plates is useful to compare standard and $\mathrm{VC}$-VEM.

\begin{figure}[!htbp]
    \centering
    \subfigure[Mesh 1.\label{fig:testcase5_mesh_quad_less}]{
        \includegraphics[width=0.225\textwidth]{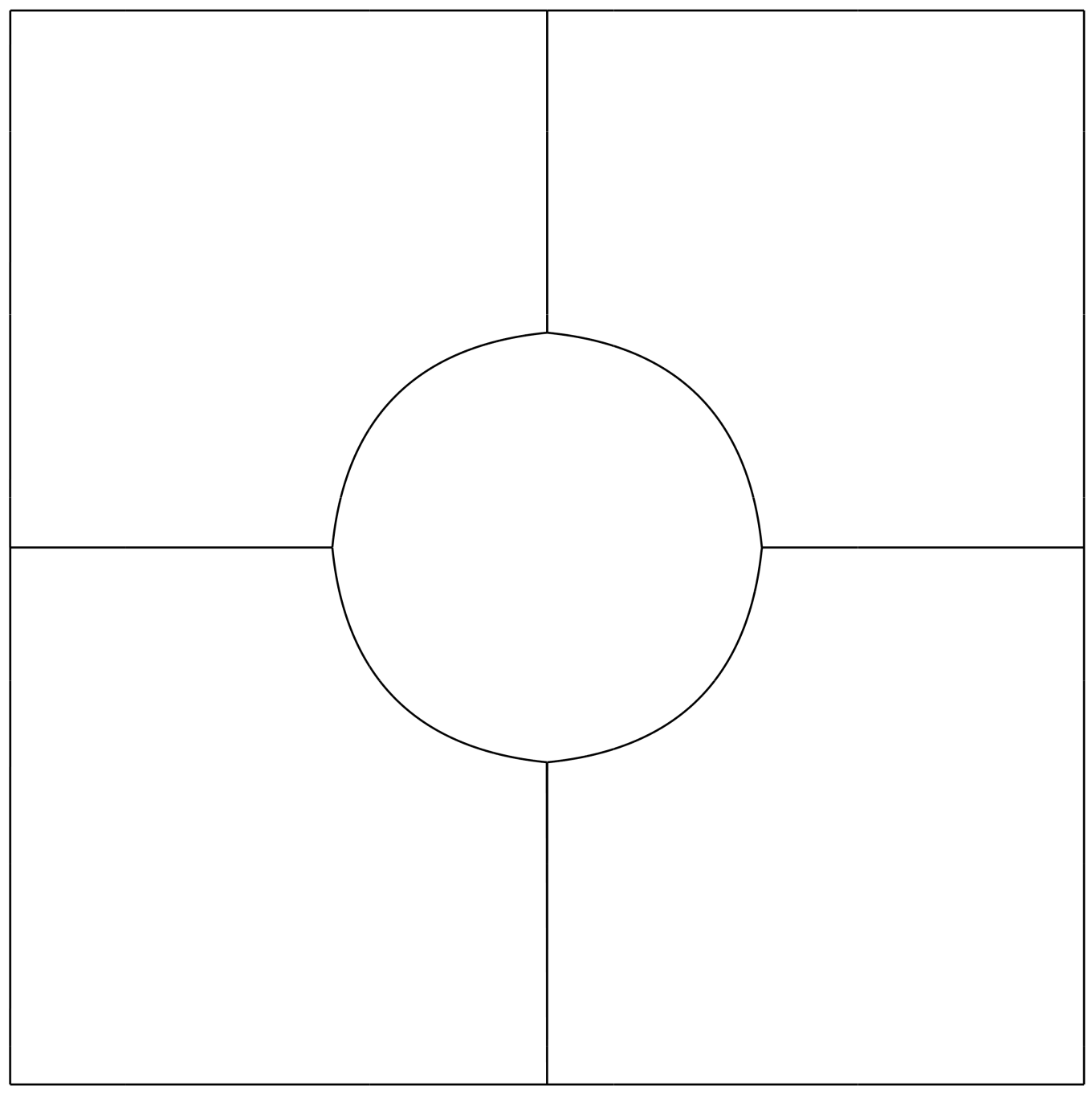}
    }
    \subfigure[Mesh 2.\label{fig:testcase5_mesh_quad_more}]{
        \includegraphics[width=0.225\textwidth]{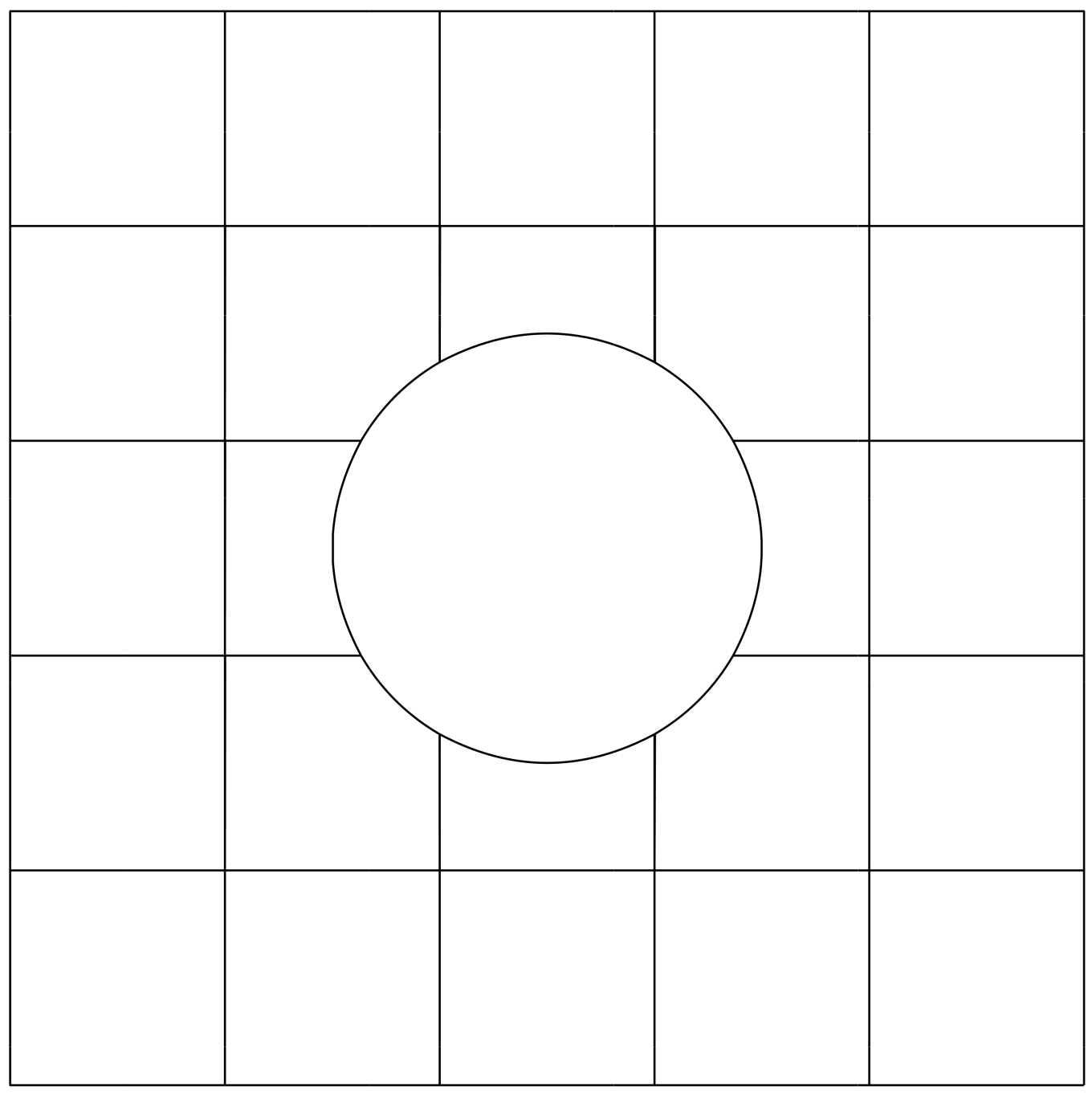}
    }
    \subfigure[Mesh 3.\label{fig:testcase5_mesh_voronoi_less}]{
        \includegraphics[width=0.225\textwidth]{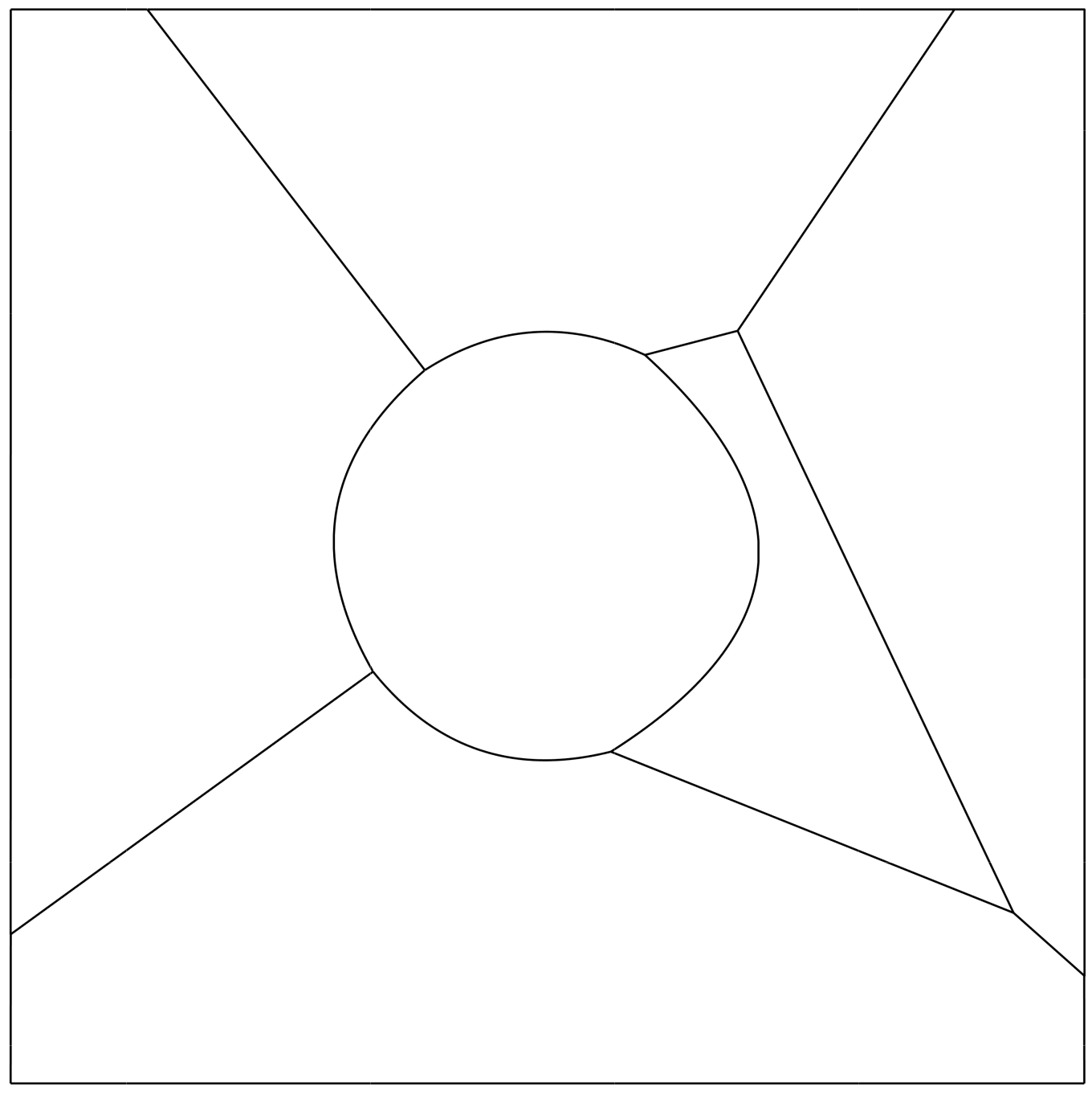}
    }
    \subfigure[Mesh 4.\label{fig:testcase5_mesh_voronoi_more}]{
        \includegraphics[width=0.225\textwidth]{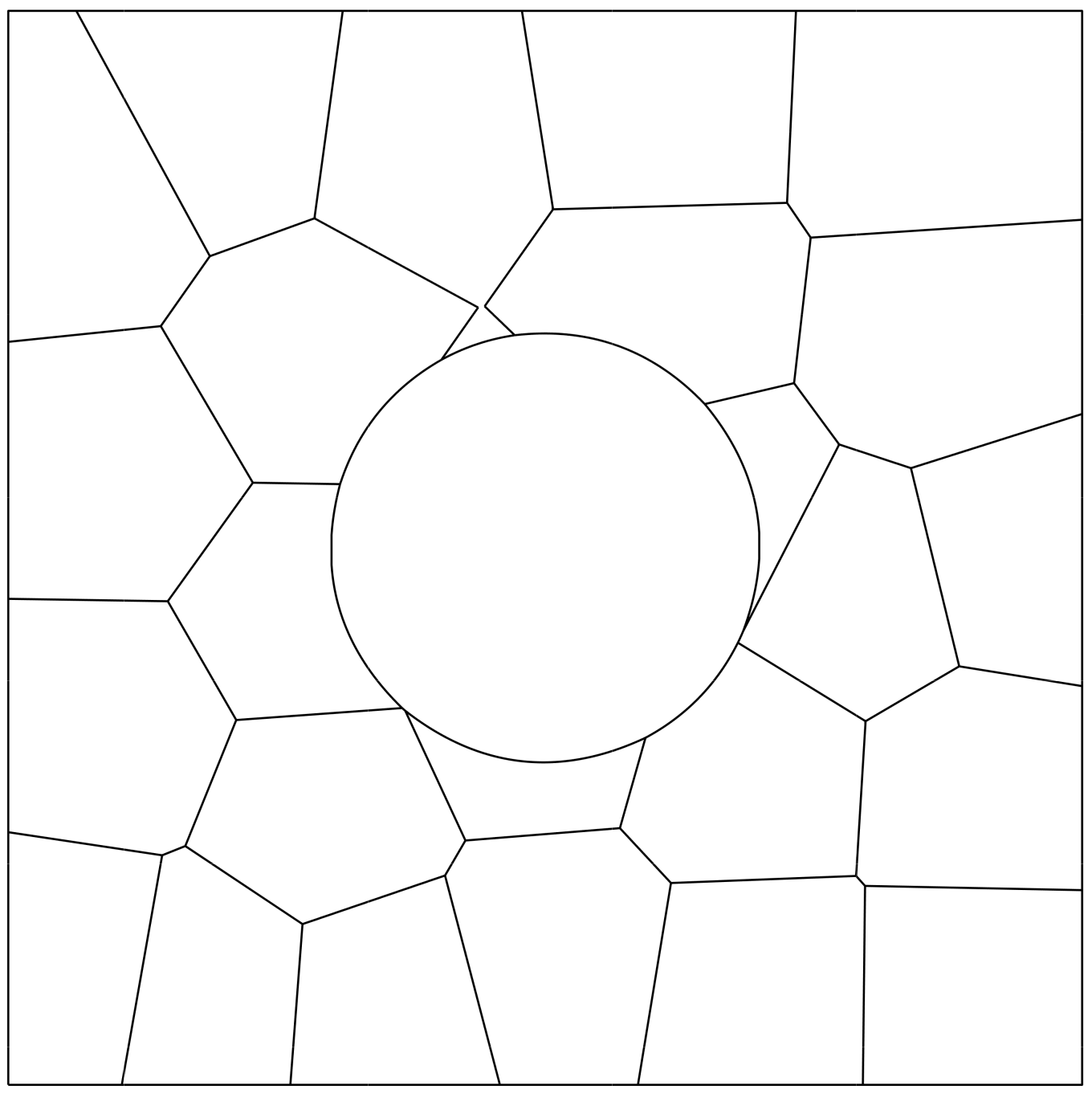}
    }
    \caption{Test case 5: meshes. The elements' circular arcs are interpolated as B\'ezier curves.\label{fig:testcase5_mesh}}
\end{figure}

Four different meshes are chosen, as shown in \fig{testcase5_mesh}. 

In particular, both structured and non-structured meshes are employed. The first two meshes are regular, while the third and fourth are characterized by the presence of elements with a relatively high degree of distortion. 
In all the cases, the curved edge feature is exploited to accurately represent the cutout.

The four meshes are associated with different orders $\kord$ to account for the different refinement levels. Mesh 1 is used with $\kord=7$, Mesh 2 with $\kord=4$, Mesh 3 with $\kord=6$, and Mesh 4 with $\kord=4$. 
Both standard and $\mathrm{VC}$-VEM are employed, in their stabilized and self-stabilized versions. 

A summary of the results is provided in \tab{testcase5_buckling} for the different meshes and stabilization techniques, and the comparison is presented against Ref.~\cite{milazzo2023buckling}. In particular, the results are presented in terms of the nondimensional ratio $\bar{N}_{xx}^{cr}/\bar{N}_{iso,xx}^{cr}$, the former being the buckling resultant of the plate under investigation, the latter the buckling resultant of a corresponding isotropic plate without cutout and elastic properties $E=69668$~MPa and $\nu=0.296$.

\begin{table}[!htbp]
\centering
\begin{threeparttable}
\begin{tabular}{ccccccc}
\hline
 & & & {layup 1} & {layup 2} & {layup 3} & {layup 4} \\
 \hline
 & & & \multicolumn{4}{c}{${\bar{N}_{xx}^{cr}/\bar{N}_{xx,iso}^{cr}}$ [-]}\\
\hline
\multirow{4}{*}{{Mesh 1 \& ${p=7}$}}
 & \multirow{2}{*}{{Stabilized}} 
   & {Standard} & 2.03 & 1.06 & 1.07 & 1.07 \\
   \cline{3-3}
  & & {VC} & 2.06 & 1.03 & 1.08 & 1.06 \\
  \cline{2-3}
 & \multirow{2}{*}{{Self-Stabilized}} 
   & {Standard} & 2.02 & 1.02 & 1.02 & 1.04 \\
   \cline{3-3}
 &  & {VC} & 2.01 & 1.01 & 1.04 & 1.04 \\
 \cline{1-3}
\multirow{4}{*}{{Mesh 2 \& ${p=4}$}}
 & \multirow{2}{*}{{Stabilized}} 
   & {Standard} & 2.01 & 1.04 & 1.04 & 1.06 \\
   \cline{3-3}
  & & {VC} & 2.05 & 1.04 & 1.05 & 1.07 \\
  \cline{2-3}
 & \multirow{2}{*}{{Self-Stabilized}} 
   & {Standard} & 2.03 & 1.03 & 1.04 & 1.06 \\
   \cline{3-3}
 &  & {VC} & 2.02 & 1.02 & 1.04 & 1.06 \\
 \cline{1-3}
\multirow{4}{*}{{Mesh 3 \& ${p=6}$}}
 & \multirow{2}{*}{{Stabilized}} 
   & {Standard} & 1.60 & 1.33 & 0.46 & 1.51 \\
   \cline{3-3}
 &  & {VC} & 2.05 & 1.05 & 1.12 & 1.05 \\
 \cline{2-3}
 & \multirow{2}{*}{{Self-Stabilized}} 
   & {Standard} & 2.04 & 1.10 & 1.07 & 1.08 \\
   \cline{3-3}
 &  & {VC} & 2.04 & 1.09\tnote{a} & 1.13\tnote{a} & 1.08 \\
 \cline{1-3}
\multirow{4}{*}{{Mesh 4 \& ${p=4}$}}
& \multirow{2}{*}{{Stabilized}} 
   & {Standard} & 1.81\tnote{a} & 1.18\tnote{a} & 1.04 & 0.96\tnote{a} \\
   \cline{3-3}
 &  & {VC} & 2.03 & 0.99 & 1.04 & 1.06 \\
 \cline{2-3}
 & \multirow{2}{*}{{Self-Stabilized}} 
   & {Standard} & 2.05\tnote{a} & 1.04 & 1.05 & 1.07\tnote{a}\\
   \cline{3-3}
 & & {VC} & 2.04\tnote{a} & 1.03 & 1.05 & 1.07 \\
\hline 
{Ref.}~\cite{milazzo2023buckling} & & & 2.06 & 1.02 & 1.05 & 1.08 \\
\hline
\end{tabular}
\begin{tablenotes}
\footnotesize
\item[a] spurious modes were removed
\end{tablenotes}
\caption{Test case 5: normalized buckling resultant $\bar{N}_{xx}^{cr}/\bar{N}_{xx,iso}^{cr}$ for different meshes and stabilization techniques.}
\label{tab:testcase5_buckling}
\end{threeparttable}
\end{table}

A first consideration regards the substantial agreement between the results obtained with the two structured meshes, i.e. Mesh 1 and Mesh 2, and the reference results.
No significant discrepancies are observed in these cases, irrespective of the layup and adopted VEM strategy. 

On the other hand, some noticeable deviations can be seen for the distorted Meshes 3 and 4. In particular, the standard stabilized VEM exhibits incorrect results for some configurations, both in terms of buckling load and mode shape. 
A first observation regards the modes predicted by Mesh 3 for the standard stabilized VEM: for all the layups, they deviate from the expected shape. 
Instead, regarding Mesh 4, several buckling loads are slightly underestimated and the corresponding mode shapes feature a milder amplitude, hence an increase in the approximation order $\kord$ may be necessary. Nonetheless, the order is not further augmented since $\kord=4$ provides sufficiently accurate results for the $\mathrm{VC}$-VEM, as well as for both self-stabilized variants. Therefore, the same order $\kord=4$ is retained to ensure a proper comparison of the different formulations. 
In contrast, both stabilized and self-stabilized $\mathrm{VC}$-VEM approaches, as well as the standard self-stabilized VEM, are able to correctly predict the buckling load even in the presence of mesh distortion. These results are in agreement with those presented in the first test case, further highlighting the superior robustness of the $L^2$ projection and $\mathrm{VC}$-VEM in the presence of more complex scenarios.

\begin{figure}[!htbp]
    \centering
    \subfigure[layup 1 VEM.\label{fig:buck_milazzo_1_VEM}]{
        \includegraphics[width=0.22\textwidth]{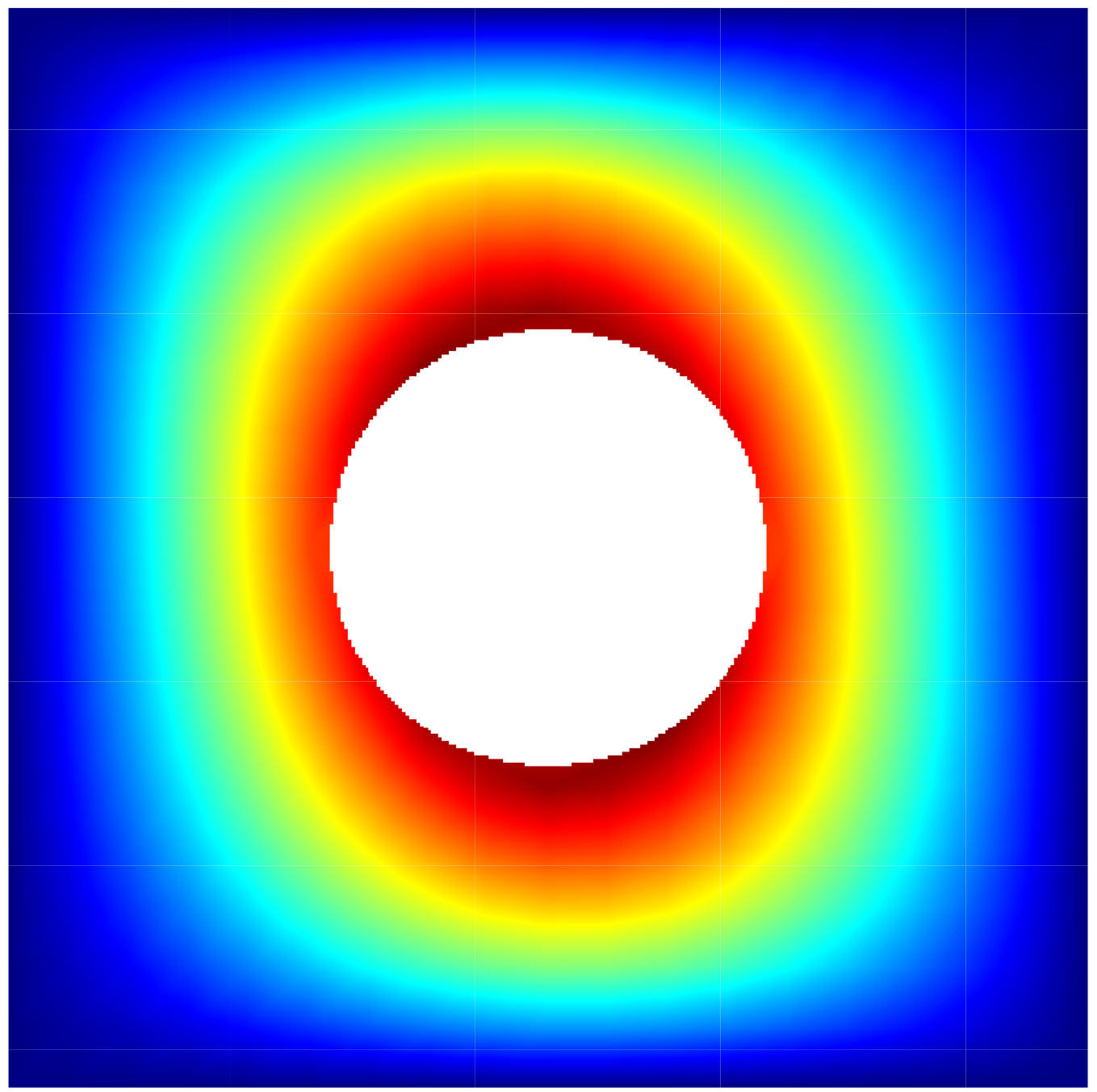}
    }    
    \subfigure[layup 2 VEM.\label{fig:buck_milazzo_2_VEM}]{
        \includegraphics[width=0.22\textwidth]{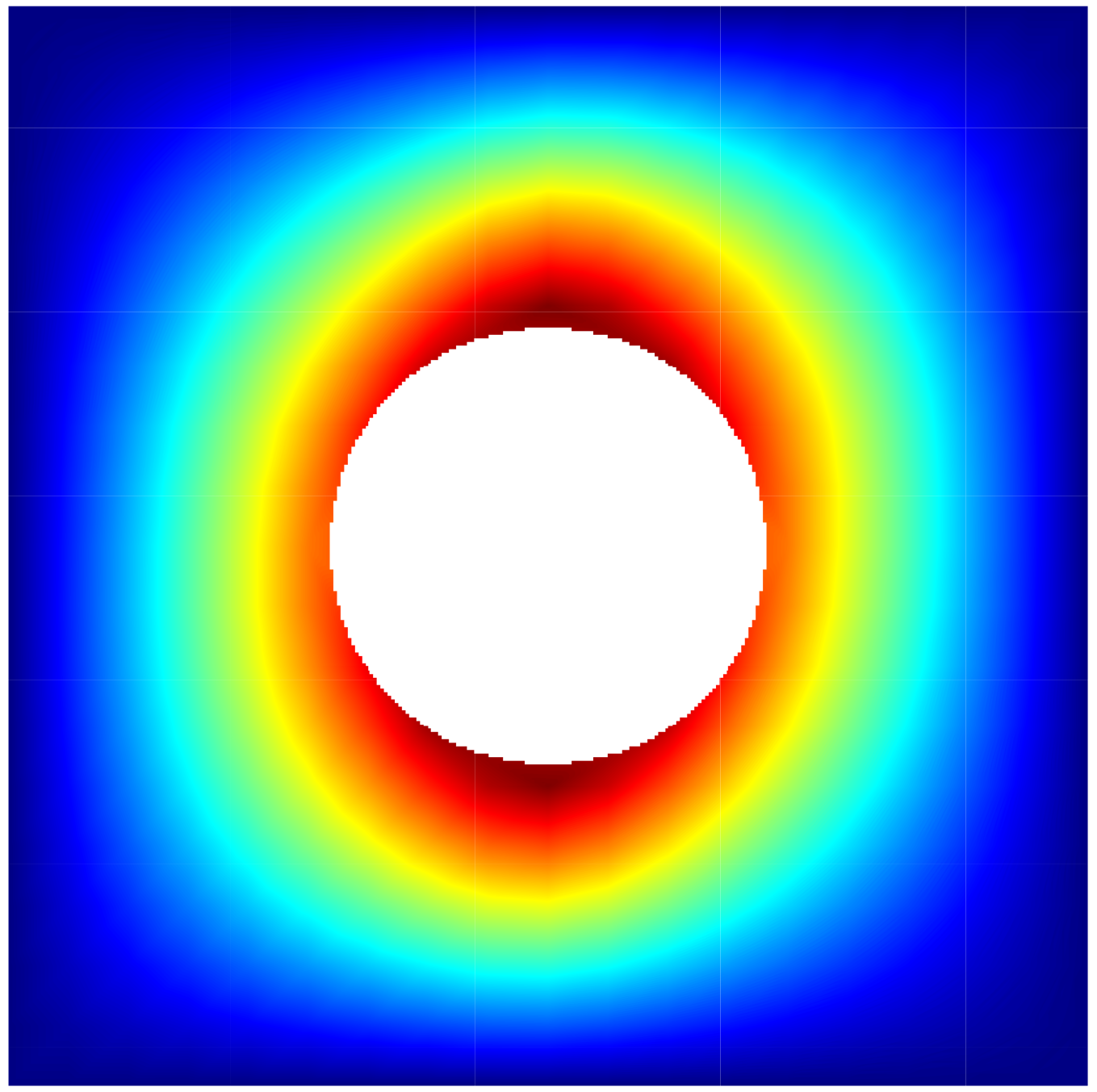}
    }    
    \subfigure[layup 3 VEM.\label{fig:buck_milazzo_3_VEM}]{
        \includegraphics[width=0.22\textwidth]{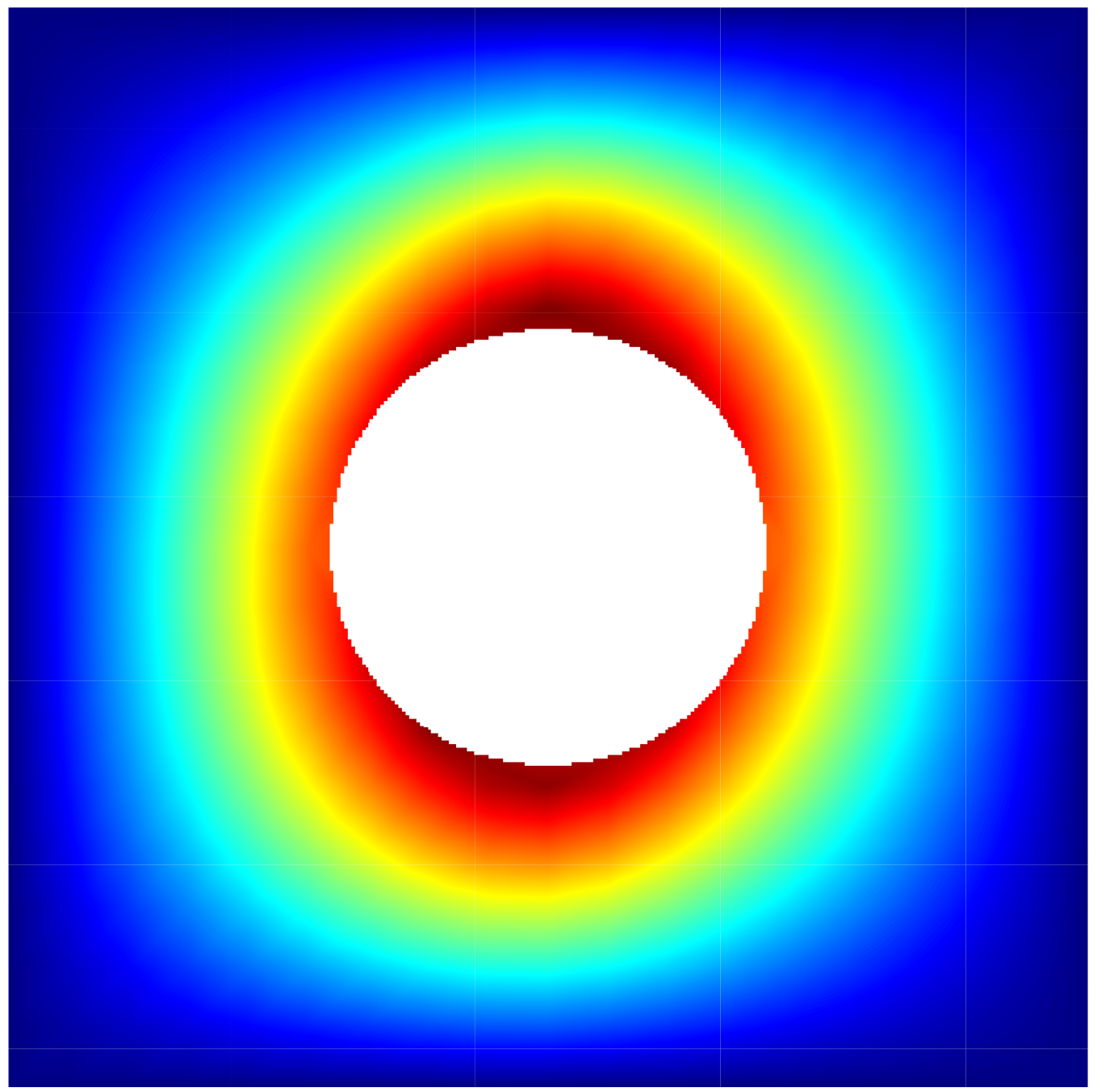}
    }    
    \subfigure[layup 4 VEM.\label{fig:buck_milazzo_4_VEM}]{
        \includegraphics[width=0.22\textwidth]{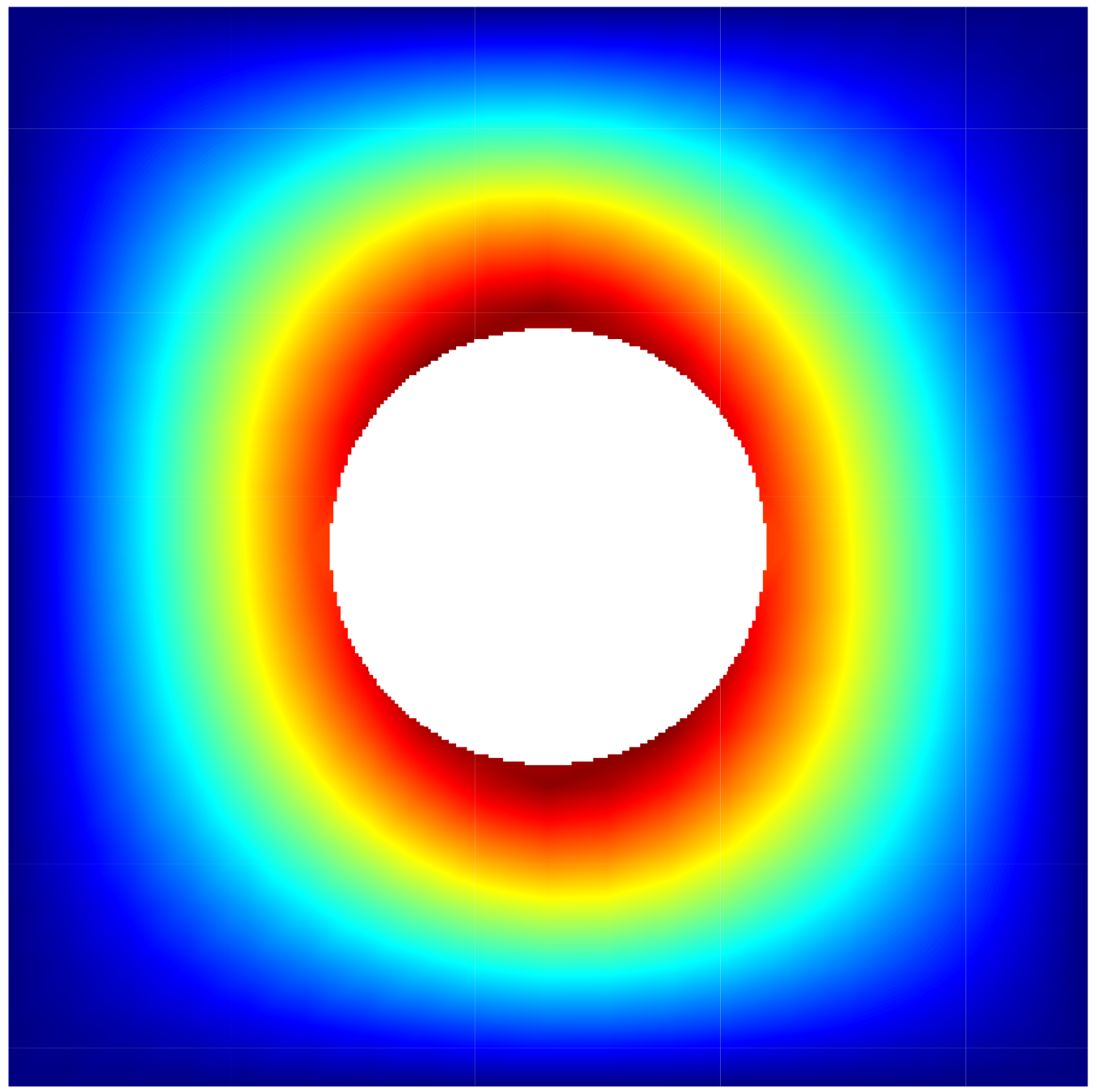}
    }     
    \\
    \subfigure[layup 1 Ritz~\cite{milazzo2023buckling}.\label{fig:buck_milazzo_1_Ritz}]{
        \includegraphics[width=0.22\textwidth]{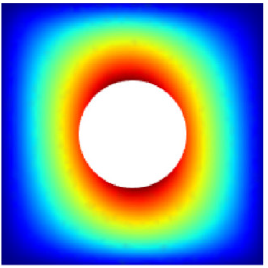}
    }    
    \subfigure[layup 2 Ritz~\cite{milazzo2023buckling}.\label{fig:buck_milazzo_2_Ritz}]{
        \includegraphics[width=0.22\textwidth]{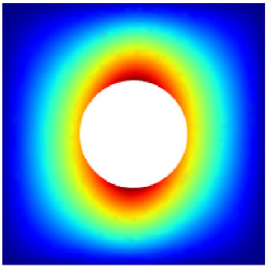}
    }    
    \subfigure[layup 3 Ritz~\cite{milazzo2023buckling}.\label{fig:buck_milazzo_3_Ritz}]{
        \includegraphics[width=0.22\textwidth]{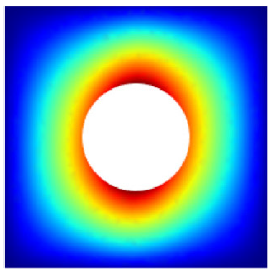}
    }    
    \subfigure[layup 4 Ritz~\cite{milazzo2023buckling}.\label{fig:buck_milazzo_4_Ritz}]{
        \includegraphics[width=0.22\textwidth]{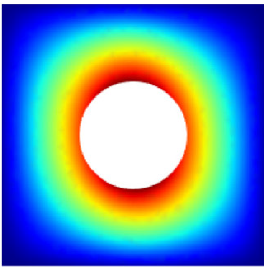}
    }      
    \caption{Test case 5: buckling modes ($w$).\label{fig:buck_milazzo} }
\end{figure}

Lastly, the buckling modes of the different layups are shown in \fig{buck_milazzo}. The VEM results are obtained using Mesh 1 with stabilized $\mathrm{VC}$-VEM. 

As shown in \fig{buck_milazzo}, the VEM predictions are in good agreement with the reference ones for all the layups considered. Overall, this test case further proves the effectiveness of the developed tool, and in particular of the $\mathrm{VC}$-VEM, in predicting the buckling behavior of complex plate domains with variable stiffness properties.

\section{Conclusions} \label{section:conclusions}

This work presented a comprehensive high-order Virtual Element Method ($\kord$-VEM) framework for the static, free-vibration, and buckling analysis of variable stiffness plates with arbitrary shapes. The proposed formulation aimed at bridging the gap between existing VEM mathematical formulations and engineering applications to plate problems in structural mechanics. For this purpose, an implementation-oriented formulation based on standard finite element notation was developed. Furthermore, a number of advanced VEM capabilities have been integrated within a unified computational framework: arbitrary polygonal elements with curved edges, hanging nodes, stabilized and self-stabilized formulations, and a newly proposed Variable Coefficients VEM (VC-VEM) approach.

Five test cases were presented for problems involving both standard and innovative variable stiffness configurations. The numerical investigations demonstrate the accuracy and robustness of the method for a wide range of structural applications. In particular, the combination of high-order approximation with local mesh refinement is effective for problems characterized by localized phenomena, such as stress concentration; the use of arbitrary polygonal elements, curved edges, and hanging nodes offers an excellent potential to simplify the discretization of complex geometries. Furthermore, owing to the $\kord$-VEM capabilities, accurate solutions are obtained on relatively coarse and distorted meshes, highlighting the robustness of the formulation and the effectiveness of $\kord$-refinement.
For variable stiffness laminates, the proposed $\mathrm{VC}$-VEM formulation and the standard self-stabilized VEM, which employs a $L^2$ projection, provide improved treatment of the spatially varying constitutive properties, leading to higher accuracy than the standard stabilized VEM formulation, particularly for distorted meshes and higher approximation orders.

Overall, the test cases highlight both the general capabilities of the framework and its specific advantages for the analysis of variable stiffness structures. The proposed framework successfully combines the geometric flexibility of VEM with an accurate treatment of spatially varying constitutive properties, while maintaining an implementation-oriented formulation suitable for practical engineering applications. The resulting methodology therefore represents a versatile computational tool for the analysis of advanced composite structures characterized by complex geometries and non-uniform stiffness distributions.

\bibliography{Bibliography}
\bibliographystyle{unsrt}
\newpage

\newpage
\section*{Supporting material} \label{chapter:appendix}

Hereafter, the construction of the relevant quantities of the VEM formulation is presented in detail. A more in-depth derivation is available in~\cite{foligno2026novel}.

\subsection*{Supporting material for \subsect{VEMspacedofs}: Spaces definition}

\subsubsection*{Virtual Element space}

The differential operator $\bm{L} \sqbr{\cdot}$ defines the problem in its strong form, and its expression reads:
\begin{equation}
    \bm{L} \sqbr{\cdot} = 
    \begin{bmatrix}
    ,x &  0 & ,y & 0  & 0  & 0  & 0  & 0\\ 
    0  & ,y & ,x & 0  & 0  & 0  & 0  & 0\\ 
    0  &  0 & 0  & 0  & 0  & 0  & ,x & ,y\\ 
    0  &  0 & 0  & ,x & 0  & ,y & -1  & 0\\ 
    0  &  0 & 0  & 0  & ,y & ,x & 0  & -1    
    \end{bmatrix}.
    \label{eq:differntialmatrix_L}
\end{equation}

The dimensions of the local spaces defined in \eq{localvirtual_spaces} are:
\begin{equation}
    \begin{aligned}
        &\Vspacedim^s = \text{dim} \plbr{\Vspaced^s \plbr{\el}} = \kord \nv + \frac{\plbr{\kord-1}\kord }{2}\ \text{for }  s=u,v,\\
        &\Vspacedim^w = \text{dim} \plbr{\Vspaced^w \plbr{\el}} = \kord \nv + \frac{\kord\plbr{\kord+1} }{2},\\
        &\Vspacedim^{s} = \text{dim} \plbr{\Vspaced^{s} \plbr{\el}} = \kord \nv + \frac{\plbr{\kord+1}\plbr{\kord+2} }{2}  \ \text{for }  s=\theta_x,\theta_y.
    \end{aligned}
    \label{eq:localvirtual_dofs}
\end{equation}
The dimension of the total local space of \eq{totlocalvirtual_space} is:
\begin{equation}
    \Vspacedim = \text{dim} \plbr{\Vspacedvec \plbr{\el}} =  5 \kord \nv + 2 \frac{\plbr{\kord-1}\kord }{2} + \frac{\kord\plbr{\kord+1} }{2} + 2\frac{\plbr{\kord+1}\plbr{\kord+2} }{2}.
    \label{eq:totlocalvirtual_dofs}
\end{equation}

\subsubsection*{Polynomial space}

The most common choice for $\Qpolyspace_{\kord}\plbr{\el}$, which serves as a basis for $\Pspace_{\kord}\plbr{\el}$, is the set of scaled monomials, here denoted as $\bar{\qpoly}_{\bm{n}} \plbr{x,y}$. By denoting with $\plbr{x_C,y_C}$ the coordinates of the centroid of $\el$, the basis is defined as~\cite{mascotto2018ill}:
\begin{equation}
    \bar{\qpoly}_{\bm{n}} \plbr{x,y} \big|_\el = \plbr{\frac{x-x_C}{\diam}}^{n_1}\plbr{\frac{y-y_C}{\diam}}^{n_2} \quad \forall \bm{n}=\sqbr{n_1,n_2} = \sqbr{\plbr{0,\dots,\kord},\plbr{0,\dots,\kord}}.
    \label{eq:scaled_monomials}
\end{equation}
However, scaled monomials suffer from numerical instability for higher values of $\kord$. To improve stability, $L^2\plbr{\el}$ orthonormal polynomials $\qpoly_{\bm{n}} \plbr{x,y}$ are employed in this work, obtained via Modified Gram-Schmidt (MGS) orthonormalization~\cite{mascotto2018ill}:
\begin{equation}
    \qpoly_{\bm{n}} \plbr{x,y} \big|_\el = \sum_{\bm{\beta}=1}^{\bm{n}} \bm{GS}_{\bm{n},\bm{\beta}} \big|_\el \bar{\qpoly}_{\bm{n}} \plbr{x,y} \big|_\el \quad \forall \bm{n}=\sqbr{n_1,n_2} = \sqbr{\plbr{0,\dots,\kord},\plbr{0,\dots,\kord}},
    \label{eq:orthonorm_monomials}
\end{equation}
where $\bm{GS} \big|_\el$ is the matrix of the orthonormalization coefficients, which are obtained for each element $\el$~\cite{bassi2012flexibility}.

To build a complete polynomial basis of order $\kord$, $\frac{\plbr{\kord+1}\plbr{\kord+2} }{2}$ independent polynomials are required. Hence, for the five displacement components, the dimension of the polynomial space is:
\begin{equation}
    \pspacedim = \text{dim} \plbr{\Pspacevec_\kord \plbr{\el}} =  5\frac{\plbr{\kord+1}\plbr{\kord+2} }{2}.
    \label{eq:localpoly_spacedim}
\end{equation}

\subsection*{Supporting material for \subsect{VEMstandard}: Matrices and vectors construction}

\subsubsection*{Stiffness matrix}

By using \eq{lagrangian_interpolation}, the unknowns $\dispd$ can be compactly written as a function of the VEM trial functions as:
\begin{equation}
    \dispd = 
    \begin{Bmatrix}
        u_h\quad
        v_h\quad
        w_h\quad
        {\theta_x}_h\quad
        {\theta_y}_h
    \end{Bmatrix}^T
    = 
    \begin{Bmatrix}
        {\trialfcnm^u}\quad
        {\trialfcnm^v}\quad
        {\trialfcnm^w}\quad
        {\trialfcnm^{\theta_x}}\quad
        {\trialfcnm^{\theta_y}}
    \end{Bmatrix}^T
    \dofs
    = \sum_{i=1}^{\Vspacedim} \trialfcnvec_i {\dofs}_i
    = \trialfcnm \dofs.
    \label{eq:VEM_trialfunctions}
\end{equation}
where $\trialfcnm^k$ are $\Vspacedim \times 1$ column vectors associated to the $k$-th displacement component, and $\dofs$ is the vector of dimension $\Vspacedim \times 1$, collecting the degrees of freedom.
The projection of the generic VEM trial functions $\trialfcnvec_i \big|_{i=1,\dots,\Vspacedim}$ onto the polynomial space is defined as~\cite{mengolini2019engineering}:
\begin{equation}
    \PiKenha \trialfcnvec_i = \sum_{\beta=1}^{\pspacedim} s_{i,\beta} \qpolyvec_{\beta} \quad \forall i=1,\dots,\Vspacedim.
    \label{eq:trialfcn_projection}
\end{equation}
By assembling the VEM trial functions and the polynomial basis functions column-wise, \eq{trialfcn_projection} becomes:
\begin{equation}
    \PiKenha \trialfcnm = \qpolymtrx \PiKenhatilde ,
    \label{eq:trialfcn_projection_matrix}
\end{equation}
where $\PiKenha \trialfcnm$ is the matrix whose columns contain the projections of the VEM trial functions, and $\qpolymtrx$ contains the polynomial basis functions. Matrix $\PiKenha \trialfcnm$ has dimension $5 \times \Vspacedim$, while matrix $\qpolymtrx$ has dimension $5 \times \pspacedim$. Accordingly, matrix $\PiKenhatilde$ has dimension $\pspacedim \times \Vspacedim$ and its $i$-th column contains the polynomial coefficients of the projection of the $i$-th VEM trial function.
The first six polynomials, for which the strain $\strain \plbr{\qpolyvec}=0$, correspond to the rigid body motions. Therefore, the following invertible augmented system is considered~\cite{mengolini2019engineering}:
\begin{equation}
\left\{
    \begin{aligned}
        &\int_\el \strain\plbr{\trialfcnvec_i}^T \cost \strain\plbr{\qpolyvec} \de \el = \int_\el \strain\plbr{\PiKenha\trialfcnvec_i}^T \cost \strain\plbr{\qpolyvec} \de \el    && \forall \qpolyvec \in \Pspacevec_\kord \plbr{\el}\\
        &\frac{1}{\nv} \sum_{j=1}^{5\nv} \text{dof}_j \plbr{\trialfcnvec_i} \text{dof}_j \plbr{\qpolyvec_{\alpha}} = \frac{1}{\nv} \sum_{j=1}^{5\nv} \text{dof}_j \plbr{\PiKenha\trialfcnvec_i} \text{dof}_j \plbr{\qpolyvec_{\alpha}} && \forall \alpha = 1,\dots,6,
    \end{aligned}
\right.
    \label{eq:system_projection}
\end{equation}
for all $ i=1,\dots,\Vspacedim$. Here, $\text{dof}_j \plbr{\qpolyvec_{\alpha}}$ represents the value of the polynomial at the vertex associated with the degree of freedom $j$.
The system in \eq{system_projection} can be expressed in matrix form as~\cite{mengolini2019engineering}:
\begin{equation}
    \BmKenha = \GmKenhatilde \PiKenhatilde \quad \longrightarrow \quad \PiKenhatilde = \plbr{ \GmKenhatilde }^{-1} \BmKenha.
    \label{eq:matrix_projection}
\end{equation}
Matrices $\BmKenha$ and $\GmKenhatilde$ have dimensions $\pspacedim \times \Vspacedim$ and $\pspacedim \times \pspacedim$, respectively, and are defined as:
\begin{equation}
    \begin{aligned}
        \BmKenha_{\alpha,i} = \int_\el \strain\plbr{\trialfcnvec_i}^T \cost \strain\plbr{\qpolyvec_\alpha} \de \el,\qquad
        \GmKenhatilde_{\alpha,\beta} = \int_\el \strain\plbr{\qpolyvec_{\beta}}^T \cost \strain\plbr{\qpolyvec_\alpha} \de \el.
    \end{aligned}
    \label{eq:BG_stiffness}
\end{equation}
Notice that, according to the second line of \eq{system_projection}, additional rows are included in $\GmKenhatilde$ and $\BmKenha$ to account for the rigid body motions. The same applies to the subsequent matrices and vectors that require the augmented rows.
To construct these matrices, it is convenient first to define polynomial vectors. Let $\qpoly_n \in \Pspace_\kord \plbr{\el}$ with $n=1,\dots,\frac{\pspacedim}{5}$ denote the scalar polynomial basis. The corresponding vector is defined as:
\begin{equation}
    \qpolyvec_{\alpha=5\plbr{n-1}+l} = \bm{e}_l \qpoly_n, \quad \forall l = 1,\dots,5, \quad \forall n=1,\dots,\frac{\pspacedim}{5}, \quad \forall \alpha = 1,\dots,\pspacedim,
    \label{eq:poly_space_versors}
\end{equation}
where $\bm{e}_l$ is the $l$-th versor in $\mathbb{R}^5$.
The construction of matrix $\GmKenhatilde_{\alpha,\beta}$ in \eq{BG_stiffness} is straightforward, as it involves only polynomial terms. For $\alpha=1,\dots,6$, the rigid body motions defined in \eq{system_projection} must also be included.
The construction of $\BmKenha$ is more involved because the VEM trial functions are not explicitly known inside the element. 
First, the strain operator $\strain \plbr{\cdot}$ can be decomposed as:
\begin{equation}
    \begin{aligned}
    \strain \plbr{\cdot} &=
    \begin{bmatrix}
        ,x & 0  & 0  & 0  & 0\\
        0  & ,y & 0  & 0  & 0\\
        ,y & ,x & 0  & 0  & 0\\
        0  & 0  & 0  & ,x & 0\\
        0  & 0  & 0  & 0  & ,y\\
        0  & 0  & 0  & ,y & ,x\\
        0 & 0   & ,x & 1  & 0\\
        0 & 0   & ,y & 0  & 1
    \end{bmatrix}
    = 
    \begin{bmatrix}
        ,x & 0  & 0  & 0  & 0\\
        0  & ,y & 0  & 0  & 0\\
        ,y & ,x & 0  & 0  & 0\\
        0  & 0  & 0  & ,x & 0\\
        0  & 0  & 0  & 0  & ,y\\
        0  & 0  & 0  & ,y & ,x\\
        0 & 0   & ,x & 0  & 0\\
        0 & 0   & ,y & 0  & 0
    \end{bmatrix}
    +
    \begin{bmatrix}
        0 & 0 & 0 & 0 & 0\\
        0 & 0 & 0 & 0 & 0\\
        0 & 0 & 0 & 0 & 0\\
        0 & 0 & 0 & 0 & 0\\
        0 & 0 & 0 & 0 & 0\\
        0 & 0 & 0 & 0 & 0\\
        0 & 0 & 0 & 1 & 0\\
        0 & 0 & 0 & 0 & 1
    \end{bmatrix}    
    = \strain^\diff \plbr{\cdot}+ \strain^\const \plbr{\cdot}.
    \end{aligned}
    \label{eq:espilon_operator}
\end{equation}
Using this decomposition, the matrix $\BmKenha$ can be written as:
\begin{equation}
        \begin{aligned}
        \BmKenha_{\alpha,i} = \int_\el \strain\plbr{\trialfcnvec_i}^T \cost \strain\plbr{\qpolyvec_\alpha} \de \el &= \int_\el \sqbr{\strain^\diff\plbr{\trialfcnvec_i}+\strain^\const\plbr{\trialfcnvec_i}}^T \cost \strain\plbr{\qpolyvec_\alpha} \de \el\\
        & = \int_\el \strain^\diff\plbr{\trialfcnvec_i}^T \cost \strain\plbr{\qpolyvec_\alpha} \de \el + \int_\el \strain^\const\plbr{\trialfcnvec_i}^T \cost \strain\plbr{\qpolyvec_\alpha} \de \el.
        \end{aligned}
    \label{eq:B_stiffness}
\end{equation}
Integration by parts can now be applied to the first term, obtaining:
\begin{equation}
    \begin{aligned}
    \BmKenha_{\alpha,i} &= - \int_\el \trialfcnvec_i \cdot \bm{L}_d \sqbr{ \cost \strain\plbr{\qpolyvec_\alpha}}  \de \el + \int_{\boundel} \trialfcnvec_i \cdot  \hat{\stress}\plbr{\qpolyvec_\alpha} \normedgeel  \de \boundel 
    + \int_\el \strain^\const\plbr{\trialfcnvec_i}^T \cost \strain\plbr{\qpolyvec_\alpha} \de \el,
    \end{aligned}
    \label{eq:B_stiffness_integrationbyparts}
\end{equation}
where $\normedgeel$ is the outward unit normal on the element boundary, and $\bm{L}_d \sqbr{\cdot}$ and $\hat{\stress}$ are defined as:
\begin{equation}
\begin{aligned}
    \bm{L}_d \sqbr{\cdot}= 
    \begin{bmatrix}
    ,x &  0 & ,y & 0  & 0  & 0  & 0  & 0\\ 
    0  & ,y & ,x & 0  & 0  & 0  & 0  & 0\\ 
    0  &  0 & 0  & 0  & 0  & 0  & ,x & ,y\\ 
    0  &  0 & 0  & ,x & 0  & ,y & 0  & 0\\ 
    0  &  0 & 0  & 0  & ,y & ,x & 0  & 0     
    \end{bmatrix}, \qquad
    \hat{\stress} =
   \begin{bmatrix}
       N_{xx} & N_{xy}\\
       N_{xy} & N_{yy}\\
       Q_{xx} & Q_{yy}\\
       M_{xx} & M_{xy}\\
       M_{xy} & M_{yy}
   \end{bmatrix}.
\end{aligned}
    \label{eq:differntialmatrix_Ld}
\end{equation}
The three terms defined in \eq{B_stiffness_integrationbyparts} can be computed entirely from the degrees of freedom. Starting from the first term and approximating $\cost$ as constant within the element, with a value equal to its average at the integration points, it follows that:
\begin{equation}
\left\{
    \begin{aligned}
    &\bm{L}_d \sqbr{\cost \strain\plbr{\qpolyvec_\alpha}} \in \Pspacevec_{\kord-1} \plbr{\el} && \forall \alpha = 3 + 5r, \;\; \text{with} \;\; r = 0, \dots, \frac{\pspacedim-5}{5}\\
    &\bm{L}_d \sqbr{\cost \strain\plbr{\qpolyvec_\alpha}} \in \Pspacevec_{\kord-2} \plbr{\el} && \text{otherwise}.
    \end{aligned}
\right.    
    \label{eq:B1_poly_stiffness}
\end{equation}
This reduces to an integral of a VEM trial function multiplied by a polynomial of degree $\kord-1$ for $w$, and of degree $\kord-2$ for the other displacement components. Thus, the first term depends only on the internal degrees of freedom of $\Vspacedvec$. 
Following~\cite{mengolini2019engineering}, the coefficients can be expressed in terms of the polynomial basis as:
\begin{equation}
\left\{
    \begin{aligned}
        &\bm{L}_d \sqbr{\cost \strain\plbr{\qpolyvec_\alpha}} = \sum_{\beta=1}^{\pspacedim-1} d_{\alpha,\beta} \qpolyvec_{\beta} && \forall \alpha = 3 + 5r, \;\; \text{with} \;\; r = 0, \dots, \frac{\pspacedim-5}{5}\\
        &\bm{L}_d \sqbr{\cost \strain\plbr{\qpolyvec_\alpha}} = \sum_{\beta=1}^{\pspacedim-2} d_{\alpha,\beta} \qpolyvec_{\beta} && \text{otherwise}.
    \end{aligned}
\right.    
    \label{eq:B1_coeff_stiffness}
\end{equation}
For a scaled monomial basis, these coefficients are trivial to obtain. In the case of the orthonormal basis adopted in this work, the coefficients are computed using the orthonormality property~\cite{mascotto2018ill}:
\begin{equation}
    d_{\alpha,\beta} = \int_\el \bm{L}_d \sqbr{\cost \strain\plbr{\qpolyvec_\alpha}} \cdot \qpolyvec_{\beta} \de \el.
    \label{eq:d_alphabeta}
\end{equation}
Regarding the third term, it holds that:
\begin{equation}
    \cost \strain\plbr{\qpolyvec_\alpha} \in \Pspacevec_{\kord} \plbr{\el} \quad \forall \alpha = 4 + s + 5r, \;\; \text{with} \;\; r=0,\dots,\frac{\pspacedim-5}{5}, \; \; s=0,1.
    \label{eq:B3_poly_stiffness}
\end{equation}
Therefore, this matrix also depends solely on the internal degrees of freedom, as it corresponds to the integral of a VEM trial function times a polynomial of degree $\kord$ for $\theta_x$ and $\theta_y$. The coefficients can be written as:
\begin{equation}
        \cost \strain\plbr{\qpolyvec_\alpha} = \sum_{\beta=1}^{\pspacedim} d_{\alpha,\beta} \qpolyvec_{\beta} \quad \forall \alpha = 4 + s + 5r, \;\; \text{with} \;\; r=0,\dots,\frac{\pspacedim-5}{5}, \; \; s=0,1.
    \label{eq:B3_coeff_stiffness}
\end{equation}
This step is straightforward, and the coefficients are obtained in the same way as in \eq{B1_coeff_stiffness}.
Lastly, the boundary integrals can be easily evaluated because the trial functions are known along the element boundary:
\begin{equation}
    \begin{aligned}
    \int_{\boundel} \trialfcnvec_i \cdot \hat{\stress}\plbr{\qpolyvec_\alpha} \normedgeel  \de \boundel &= \sum_{j=1}^{\nedges} \int_{\bound} \trialfcnvec_i \cdot \hat{\stress}\plbr{\qpolyvec_\alpha} \normedge  \de \bound 
    +\sum_{j=1}^{\nedgescurv} \int_{\boundcurv} \trialfcnvec_i \cdot \hat{\stress}\plbr{\qpolyvec_\alpha} \normedgecurv  \de \boundcurv.
    \end{aligned}
    \label{eq:B2_stiffness}
\end{equation}
For straight edges, the integrals are computed using the Gauss-Lobatto quadrature:
\begin{equation}
    \sum_{j=1}^{\nedges} \int_{\bound} \trialfcnvec_i \cdot \hat{\stress}\plbr{\qpolyvec_\alpha} \normedge  \de \bound = \sum_{j=1}^{\nedges} \big|e_j\big| \sum_{r=1}^{n_{int}} w_r \trialfcnvec_i \big|_{x_r,y_r} \cdot \hat{\stress}\plbr{\qpolyvec_\alpha}\big|_{x_r,y_r} \normedge,
    \label{eq:B2_stiffness_straight}
\end{equation}
where $\big|e_j\big|$ is the edge length, $w_r$ and $x_r,y_r$ are the quadrature weights and points, and $n_{int}$ is the number of integration points required for exact integration. 
For curved edges, following~\cite{daveiga2019virtual}, the mapping is exploited:
\begin{equation}
    \begin{aligned}
        \sum_{j=1}^{\nedgescurv} \int_{\boundcurv} \trialfcnvec_i \cdot \hat{\stress}\plbr{\qpolyvec_\alpha} \normedgecurv  \de \boundcurv  &= \sum_{j=1}^{\nedgescurv} \int_{\boundcurv} \sqbr{\tilde{\trialfcnvec}_i \cdot \bm{\tilde{\hat{S}}\plbr{\qpolyvec_\alpha}}\normedgecurvtilde} \circ \gamma_e^{-1}   \de \boundcurv \\
        &= \sum_{j=1}^{\nedgescurv} \int_{I_e} \sqbr{\tilde{\trialfcnvec}_i \cdot \bm{\tilde{\hat{S}}\plbr{\qpolyvec_\alpha}}\normedgecurvtilde} \norm{\gamma_e^{'}}   \de I_e\\
        &= \sum_{j=1}^{\nedgescurv} \sum_{r=1}^{n_{int}} w_r \sqbr{\tilde{\trialfcnvec}_i \cdot \bm{\tilde{\hat{S}}\plbr{\qpolyvec_\alpha}}\normedgecurvtilde}\big|_{x_r,y_r} \norm{\gamma_e^{'}}\big|_{x_r,y_r}\\
        &= \sum_{j=1}^{\nedgescurv} \sum_{r=1}^{n_{int}} w_r \tilde{\trialfcnvec}_i\big|_{x_r,y_r} \cdot \sqbr{\hat{\stress}\plbr{\qpolyvec_\alpha}\normedgecurv}\big|_{\gamma_e\plbr{x_r,y_r}} \norm{\gamma_e^{'}}\big|_{x_r,y_r}.
    \end{aligned}
    \label{eq:B2_stiffness_curve}
\end{equation}
Here, $\tilde{\trialfcnvec}_i$ is the image of $\trialfcnvec_i$ through $\gamma_e$ and is a polynomial, see \eq{curvedpoly_space}; specifically, Lagrange polynomials are employed. 
Therefore, matrix $\BmKenha_{\alpha,i}$ is fully computable. 
For $\alpha=1,\dots,6$, the rigid body motions defined in \eq{system_projection} must also be included.
It is convenient to also define the matrix:
\begin{equation}
    \Dm_{i,\alpha} = \text{dof}_i \plbr{\qpolyvec_{\alpha}}.
    \label{eq:D_matrix}
\end{equation}
For the boundary degrees of freedom, this corresponds to evaluating the polynomial at the given node, while for internal moments it reduces to integrating the product of known polynomials. Numerical integration over a generic polygon with curved B\'ezier edges is carried out using the open-access code from~\cite{sommariva2023low,sommariva2007product}.

By defining the same operator as in \eq{trialfcn_projection}, this time expressed in terms of the trial functions themselves~\cite{daveiga2014hitchhiker}, the following expression is obtained:
\begin{equation}
\begin{aligned}
    \PiKenha \trialfcnvec_i = \sum_{j=1}^{\Vspacedim} \pi_{i,j} \trialfcnvec_{j} \quad \forall i=1,\dots,\Vspacedim, && \text{with} \quad \pi_{i,j} = \sum_{\beta=1}^{\pspacedim} s_{i,\beta} \text{dof}_j \plbr{\qpolyvec_{\beta}}.
\end{aligned}    
    \label{eq:trialfcn_projection_trialfcns}
\end{equation}
In matrix form:
\begin{equation}
    \bm{\Pi} = \Dm \PiKenhatilde.
    \label{eq:stiffness_pimatrix}
\end{equation}

The consistency part of the stiffness matrix is constructed as:
\begin{equation}
    \KEc = \plbr{\PiKenhatilde}^T \GmKVC\PiKenhatilde . 
    \label{eq:K_matrix}
\end{equation}
Notice that, to construct the stiffness matrix, $\GmKVC$ accounts for the variable stiffness:
\begin{equation}
    \GmKVC_{\alpha,\beta} = \int_\el \strain\plbr{\qpolyvec_{\beta}}^T \cost \plbr{x,y} \strain\plbr{\qpolyvec_\alpha} \de \el.
    \label{eq:stiffness_KG}
\end{equation}
To evaluate \eq{stiffness_KG}, for a spatially varying $\cost \plbr{x,y}$, a higher-order quadrature rule is required compared to the constant case.

\subsubsection*{Stabilized VEM}
In matrix form, \eq{K_s_mascotto_diag} can be written as:
\begin{equation}
\KES{}
=
\begin{bmatrix}
{\KES{}}^{u,v} & \bm{0} & \bm{0} \\
\bm{0} & {\KES{}}^{w} & \bm{0} \\
\bm{0} & \bm{0} & {\KES{}}^{\theta_x,\theta_y}
\end{bmatrix}, \qquad \text{with} \quad {\KES{}}^{\alpha}
=
\tau
\plbr{\bm{I}^{\alpha}-\Pim^{\alpha}}^T
{\bm{s}^{E}}^\alpha
\plbr{\bm{I}^{\alpha}-\Pim^{\alpha}},
\label{eq:K_matrix_s_mascotto_diag}
\end{equation}
where $\Pim^{\alpha}$ and ${\bm{s}^{E}}^\alpha$ are the sub-matrices associated to the degrees of freedom of the sub-block $\alpha$, see Eq.~\ref{eq:K_s_mascotto_diag}.

So that, in the case of stabilized VEM, the stiffness matrix can be written as:
\begin{equation}
    \KE = \KEc + \KES{}.
    \label{eq:K_stiffness_stabilizedVEM}
\end{equation}

\subsubsection*{Self-stabilized VEM}

The VEM space has the same dimension as the classical stabilized formulation, i.e. $\spacedimaugV = \Vspacedim$.
Conversely, the dimension of the enlarged polynomial space is given by:
\begin{equation}
    \pspacedimaugV = \text{dim} \plbr{\Pspacevecunc_{\kordaug-1} \plbr{\el}} = 8\frac{\plbr{\kordaug}\plbr{\kordaug+1} }{2}.
    \label{eq:localpoly_spacedim_V4}
\end{equation}

The projection in the $L^2$ norm is obtained by solving the system:
\begin{equation}
        \int_\el \strain\plbr{\trialfcnvec_i}^T \cost \qpolyvecunc \de \el = \int_\el \PiKenhaSone\strain\plbr{\trialfcnvec_i}^T \cost \qpolyvecunc \de \el    \quad \forall \qpolyvecunc \in \Pspacevecunc_{\kordaug-1} \plbr{\el},
    \label{eq:system_projection_selfstab_nointdofs_pproj}
\end{equation}
for all $ i=1,\dots,\Vspacedim$, where $\qpolyvecunc$ is the polynomial vector defined as:
\begin{equation}
    \qpolyvecunc_{\alpha=8\plbr{n-1}+l} = \bm{e}_l \qpoly_n, \quad \forall l = 1,\dots,8, \quad \forall n=1,\dots,\frac{\pspacedim}{5}, \quad \alpha = 1,\dots,\pspacedimaugV,
    \label{eq:uncoupledpoly_space}
\end{equation}
with $\bm{e}_l$ being the $l-th$ versor in $\mathbb{R}^8$ and $\Pspacevecunc_{\kordaug-1} \plbr{\el} = \sqbr{ \Pspace_{\kordaug-1} \plbr{\el} }^8$. The projector $\PiKenhaSone$ maps the strain $\strain\plbr{\trialfcnvec_i}$, rather than the trial function $\trialfcnvec_i$ itself, into this new polynomial space. 

\eq{system_projection_selfstab_nointdofs_pproj} in matrix form becomes:
\begin{equation}
    \BmKenhaS = \GmKenhaS \PiKenhatildeSone \quad \longrightarrow \quad \PiKenhatildeSone = \plbr{ \GmKenhaS }^{-1} \BmKenhaS.
    \label{eq:matrix_projection_selfstab_nointdofs_pprojection}
\end{equation}
The matrix $\BmKenhaS$ has dimension $\pspacedimaugV \times \Vspacedim$ and the matrix $\GmKenhaS$ has dimension $\pspacedimaugV \times \pspacedimaugV$. The projection operator $\PiKenhatildeSone$ has dimension $\pspacedimaugV \times \Vspacedim$. The matrices $\BmKenhaS$ and $\GmKenhaS$ are defined as:
\begin{equation}
    \begin{aligned}
        \BmKenhaS_{\alpha,i} = \int_\el \strain\plbr{\trialfcnvec_i}^T \cost \qpolyvecunc_\alpha \de \el, \qquad
        \GmKenhaS_{\alpha,\beta} = \int_\el {\qpolyvecunc_{\beta}}^T \cost \qpolyvecunc_\alpha \de \el.
    \end{aligned}
    \label{eq:BG_stiffness_selfstab_nointdofs_pproj}
\end{equation}
These matrices are constructed following the same procedure as in the stabilized VEM.
After applying integration by parts, to construct matrix $\BmKenhaS_{\alpha,i}$, the enlarged enhanced space in \eq{localvirtual_spaces_selfstabenhanced_V4} is used, giving:
\begin{equation}
    \begin{aligned}
    \BmKenhaS_{\alpha,i} &= - \int_\el \PiokK \trialfcnvec_i \cdot \bm{L}_d \sqbr{\cost \qpolyvecunc_\alpha}  \de \el + \int_{\boundel} \trialfcnvec_i \cdot  \hat{\stress}\plbr{\qpolyvecunc_\alpha} \normedgeel  \de \boundel 
    + \int_\el \strain^\const\plbr{\PiokK\trialfcnvec_i}^T \cost \qpolyvecunc_\alpha \de \el,
    \end{aligned}
    \label{eq:B_stiffness_integrationbyparts_selfstab_nointdofs_projection_pproj}
\end{equation}
where $\hat{\stress}\plbr{\qpolyvecunc_\alpha}$ is constructed starting from $\cost \qpolyvecunc_\alpha$. The operator $\PiokK$ is the $L^2$ projector and is defined later in \eq{body_L2_projection_extended}.
The stiffness matrix is then obtained in the form:
\begin{equation}
        \KE = \plbr{\PiKenhatildeSone}^T \GmKVCS \PiKenhatildeSone,
    \label{eq:K_stiffness_selfstabilizedVEM}
\end{equation}
where $\GmKVCS$ is the variable coefficient counterpart of $\GmKenhaS$.

\subsubsection*{Geometric stiffness matrix}

Following the same logical flow adopted for the stiffness matrix, the following system is obtained:
\begin{equation}
\left\{
    \begin{aligned}
        &\int_\el \strain_\buck\plbr{\trialfcn_i}^T \mathbb{N} \strain_\buck\plbr{\qpoly} \de \el = \int_\el \strain_\buck\plbr{\PiGenha\trialfcn_i}^T \mathbb{N} \strain_\buck\plbr{\qpoly} \de \el    && \forall \qpoly \in \Pspace_\kord \plbr{\el}\\
        &\frac{1}{\nv} \sum_{j=1}^{\nv} \text{dof}_j \plbr{\trialfcn_i} \text{dof}_j \plbr{\qpoly_{\alpha}} = \frac{1}{\nv} \sum_{j=1}^{\nv} \text{dof}_j \plbr{\PiGenha\trialfcn_i} \text{dof}_j \plbr{\qpoly_{\alpha}} && \alpha = 1,
    \end{aligned}
    \label{eq:bucklingsystem_projection}
\right.    
\end{equation}
for all $ i=1,\dots,\Vspacedim^w$.
Since only the terms related to $w$ are retained, the rigid body motion corresponds to $\alpha=1$ and $\trialfcn_i$ and $\qpoly$ are both scalars. 
The strain operator is defined as $\strain_\buck = 
\begin{bmatrix}
    ,x, \quad
    ,y 
\end{bmatrix}^T$, corresponding to the scalar version of \eq{plate_epsilonbuck}.
The system in \eq{bucklingsystem_projection} can be written in matrix form as:
\begin{equation}
    \BmGenha = \GmGenhatilde \PiGenhatilde \quad \longrightarrow \quad \PiGenhatilde = \plbr{ \GmGenhatilde }^{-1} \BmGenha.
    \label{eq:bucklingmatrix_projection}
\end{equation}
Matrix $\GmGenhatilde$ has dimension $\frac{\pspacedim}{5} \times \frac{\pspacedim}{5}$, while matrix $\BmGenha$ has dimension $\frac{\pspacedim}{5} \times \Vspacedim^w$, and they read:
\begin{equation}
    \begin{aligned}
        &\BmGenha_{\alpha,i} = \int_\el \strain_\buck\plbr{\trialfcn_i}^T \mathbb{N}\strain_\buck\plbr{\qpoly_\alpha} \de \el,
        &\GmGenhatilde_{\alpha,\beta} = \int_\el \strain_\buck\plbr{\qpoly_{\beta}}^T  \mathbb{N}\strain_\buck\plbr{\qpoly_\alpha} \de \el.
    \end{aligned}
    \label{eq:BG_buckling}
\end{equation}
In this context, $\qpoly_\alpha \in \Pspace_\kord \plbr{\el}$ for $\alpha=1,\dots,\frac{\pspacedim}{5}$.
Regarding matrix $\BmGenha$, integration by parts yields:
\begin{equation}
    \BmGenha_{\alpha,i} = - \int_\el \trialfcn_i \cdot \bm{L}_b \sqbr{\mathbb{N}\strain_\buck\plbr{\qpoly_\alpha}}  \de \el + \int_{\boundel} \trialfcn_i \cdot \mathbb{N}\strain_\buck \plbr{\qpoly_\alpha} \normedgeel  \de \boundel,
    \label{eq:B_buckling_integrationbyparts}
\end{equation}
where $\bm{L}_b \sqbr{\cdot}= 
\begin{bmatrix}
    ,x & ,y
\end{bmatrix}$.
The steps to compute \eq{B_buckling_integrationbyparts} follow the same procedure detailed for the stiffness matrix, requiring only the degrees of freedom. The construction of matrix $\GmGenhatilde$ is straightforward, allowing the projection in \eq{bucklingmatrix_projection} to be computed.  
Therefore, the matrix is constructed as:
\begin{equation}
    \KEg  = \plbr{\PiGenhatilde}^T \GmGVC\PiGenhatilde ,
    \label{eq:Kg_matrix}
\end{equation}
with:
\begin{equation}
    \begin{aligned}
      &\GmGVC_{\alpha,\beta} = \int_\el \strain_\buck\plbr{\qpoly_{\beta}}^T \mathbb{N} \plbr{x,y} \strain_\buck\plbr{\qpoly_\alpha} \de \el.
      \end{aligned}
    \label{eq:buckling_GG}
\end{equation}

\subsubsection*{Mass matrix}

The steps for the construction of the mass matrix are detailed hereafter. The projection is of $L^2$ type and its extended form reads:
\begin{equation}
    \int_\el \trialfcnvec_i^T \mass \qpolyvec \de \el = \int_\el \PiMenha\trialfcnvec_i^T \mass\qpolyvec \de \el    \quad \forall \qpolyvec \in \Pspacevec_\kord \plbr{\el},
    \label{eq:mass_L2_projection}
\end{equation}
for all $ i=1,\dots,\Vspacedim$. In matrix form:
\begin{equation}
    \QmMenha = \HmMenha \PiMenhatilde \quad \longrightarrow \quad \PiMenhatilde = \plbr{ \HmMenha }^{-1} \QmMenha,
    \label{eq:massmatrix_projection}
\end{equation}
where matrix $\HmMenha$ has dimension $\pspacedim \times \pspacedim$ and matrix $\QmMenha$ has dimension $\pspacedim \times \Vspacedim$, with entries:
\begin{equation}
    \begin{aligned}
        \QmMenha_{\alpha,i} = \int_\el \trialfcnvec_i^T \mass \qpolyvec_\alpha \de \el, \qquad
        \HmMenha_{\alpha,\beta} = \int_\el \qpolyvec_{\beta}^T \mass \qpolyvec_\alpha \de \el.
    \end{aligned}
    \label{eq:QH_mass}
\end{equation}

By using the enhanced space, matrix $\QmMenha$ is constructed as follows:
\begin{equation}
    \QmMenha_{\alpha,i} = 
    \begin{cases}
        \int_\el \trialfcnvec_i^T \mass \qpolyvec_\alpha \de \el \qquad\;\, \forall \alpha = 1 + s + 5r, \; \text{with} \; r=0,\dots,\frac{\pspacedimtwo-5}{5}, \; s=0,1\\
        \int_\el \trialfcnvec_i^T \mass \qpolyvec_\alpha \de \el \qquad\;\, \forall \alpha = 3 + 5r, \; \text{with} \; r=0,\dots,\frac{\pspacedimone-5}{5} \\
        \int_\el \trialfcnvec_i^T \mass \qpolyvec_\alpha \de \el \qquad\;\, \forall \alpha = 4 + s + 5r, \; \text{with} \; r=0,\dots,\frac{\pspacedim-5}{5}, \; s=0,1 \\
        \int_\el \PiKenhatilde \trialfcnvec_i^T \mass \qpolyvec_\alpha \de \el \quad \text{otherwise}.
    \end{cases}
    \label{eq:Qmassmatrix}
\end{equation}
The first three rows correspond to the internal degrees of freedom and are computed exactly. The remaining terms are evaluated using the enhanced condition defined in \eqs{localvirtual_spaces_general}{localvirtual_spaces}.
The construction of the matrix $\HmMenha$ is straightforward. 
The mass matrix is then obtained as:
\begin{equation}
    \ME = \plbr{\PiMenhatilde}^T \HmMenha\PiMenhatilde .
    \label{eq:M_matrix}
\end{equation}

\subsubsection*{Body forces vector}

The $L^2$ projection for the body forces vector reads:
\begin{equation}
    \int_\el \trialfcnvec_i^T \qpolyvec \de \el = \int_\el \PiokK\trialfcnvec_i^T \qpolyvec \de \el    \quad \forall \qpolyvec \in \Pspacevec_\kord \plbr{\el},
    \label{eq:body_L2_projection_extended}
\end{equation}
for all $ i=1,\dots,\Vspacedim$. In matrix form:
\begin{equation}
    \Qmok= \Hmok \PiokKtilde \quad \longrightarrow \quad \PiokKtilde = \plbr{ \Hmok }^{-1} \Qmok,
    \label{eq:bodymatrix_projection}
\end{equation}
where:
\begin{equation}
    \begin{aligned}
        {\Qmok}_{\alpha,i} = \int_\el \trialfcnvec_i^T \qpolyvec_\alpha \de \el, \qquad
        {\Hmok}_{\alpha,\beta} = \int_\el \qpolyvec_{\beta}^T \qpolyvec_\alpha \de \el.
    \end{aligned}
    \label{eq:QH_body}
\end{equation}
Matrix $\Qmok$ is straightforward as it is the product of known polynomials, whereas matrix $\Hmok$ is built by using the enhanced condition in \eq{localvirtual_spaces}.
Using this projection, the body forces vector can be computed as:
\begin{equation}
    \plbr{\bodyfE}_i = \int_\el \PiokK\trialfcnvec_i^T \bodyforcevec \de \el ,
    \label{eq:bodyforces_vector}
\end{equation}
where $\bodyforcevec$ are the body forces.

\subsubsection*{Thermal forces vector}

The projection for the thermal forces vector can be written as:
\begin{equation}
\left\{
    \begin{aligned}
        &\int_\el \strain_\ther\plbr{\trialfcnvecther_i}^T \check{\bm{R}} \strain_\ther\plbr{\qpolyvecther} \de \el = \int_\el \strain_\ther\plbr{\PiTenha\trialfcnvecther_i}^T \check{\bm{R}} \strain_\ther\plbr{\qpolyvecther} \de \el    && \; \qpolyvecther \in \Pspacevecther_\kord \plbr{\el}\\
        &\frac{1}{\nv} \sum_{j=1}^{\nv} \text{dof}_j \plbr{\trialfcnvecther_i} \text{dof}_j \plbr{\qpolyvecther_{\alpha}} = \frac{1}{\nv} \sum_{j=1}^{\nv} \text{dof}_j \plbr{\PiTenha\trialfcnvecther_i} \text{dof}_j \plbr{\qpolyvecther_{\alpha}} && \; \alpha = 1,\dots,5,
    \end{aligned}
\right.    
    \label{eq:thermalsystem_projection}
\end{equation}
for all $i=1,\dots,\Vspacedim^*$, with $\Vspacedim^*=\Vspacedim^u + \Vspacedim^v + \Vspacedim^{\theta_x} + \Vspacedim^{\theta_y}$, and $\Pspacevecther_\kord = \sqbr{\Pspace_\kord}^4$. In $\trialfcnvecther$, the out-of-plane displacement $w$ is excluded, as the shear contribution is not accounted for in the thermal forces. Notice that in \tab{summary_discretebilinearforms}, $\hat{\bm{P}}^* \plbr{\el}$ corresponds to $\begin{Bmatrix} x \; 0 \; 0 \; 0 \end{Bmatrix}^T$. The strain operator $\strain_\ther \plbr{\cdot}$ and matrix $\check{\bm{R}}$ are defined as:
\begin{equation}
\begin{aligned}
    \strain_\ther \plbr{\cdot} =
    \begin{bmatrix}
        ,x & 0  & 0  & 0\\
        0  & ,y & 0  & 0\\
        ,y & ,x & 0  & 0\\
        0  & 0  & ,x & 0\\
        0  & 0  & 0  & ,y\\
        0  & 0  & ,y & ,x\\
        0 & 0   & 0  & 0\\
        0 & 0   & 0  & 0
    \end{bmatrix}, \qquad
    \check{\bm{R}}\plbr{x,y} = 
    \begin{bmatrix}
        \check{\bm{N}}\plbr{x,y} & \bm{0} & \bm{0}\\
        \bm{0} & \check{\bm{M}}\plbr{x,y} & \bm{0}\\
         \bm{0} & \bm{0} & \check{\bm{Q}}\plbr{x,y}
    \end{bmatrix}.
\end{aligned}
    \label{eq:epsilon_th}
\end{equation}
The operator $\strain_\ther \plbr{\cdot}$ corresponds to a reduced version of $\strain \plbr{\cdot}$, in which shear contributions and the out-of-plane displacement $w$ are omitted. In $\check{\bm{R}}\plbr{x,y}$:
\begin{equation}
    \check{\bm{N}}\plbr{x,y} =
    \begin{bmatrix}
        \hat{N}_{xx} & \hat{N}_{xy}\\
        \hat{N}_{xy} & \hat{N}_{yy}      
    \end{bmatrix},
    \qquad
    \check{\bm{M}}\plbr{x,y} =
    \begin{bmatrix}
        \hat{M}_{xx} & \hat{M}_{xy}\\
        \hat{M}_{xy} &  \hat{M}_{yy}
    \end{bmatrix},
    \qquad
    \check{\bm{Q}} =
    \begin{bmatrix}
        0 & 0\\
        0 & 0
    \end{bmatrix}.
    \label{eq:NMQ_thermal}
\end{equation}

In matrix form, the projection is expressed as:
\begin{equation}
    \BmTenha = \GmTenhatilde \PiTenhatilde \quad \longrightarrow \quad \PiTenhatilde = \plbr{ \GmTenhatilde }^{-1} \BmTenha.
    \label{eq:thermalmatrix_projection}
\end{equation}
Matrices $\BmTenha$ and $\GmTenhatilde$ have dimensions $\frac{4\pspacedim}{5} \times \Vspacedim^*$ and $\frac{4\pspacedim}{5} \times \frac{4\pspacedim}{5}$, respectively, and are defined as:
\begin{equation}
    \begin{aligned}
        \BmTenha_{\alpha,i} = \int_\el \strain_\ther\plbr{\trialfcnvecther_i}^T \check{\bm{R}} \strain_\ther\plbr{\qpolyvecther_\alpha} \de \el, \qquad
        \GmTenhatilde_{\alpha,\beta} = \int_\el \strain_\ther\plbr{\qpolyvecther_{\beta}}^T \check{\bm{R}} \strain_\ther\plbr{\qpolyvecther_\alpha} \de \el.
    \end{aligned}
    \label{eq:BG_thermal}
\end{equation}

In this case, the polynomial space is defined as:
\begin{equation}
    \qpolyvecther_{\alpha=4\plbr{n-1}+l} = \bm{e}_l \qpoly_n, \quad \forall l = 1,\dots,4, \quad \forall n=1,\dots,\frac{\pspacedim}{5}, \quad \forall \alpha = 1,\dots,\frac{4\pspacedim}{5},
    \label{eq:poly_space_versors_ther}
\end{equation}
where $\bm{e}_l$ is the $l-th$ versor in $\mathbb{R}^4$. 
The matrix $\BmTenha_{\alpha,i}$ is constructed following the same steps detailed for the stiffness matrix. Specifically, by applying integration by parts, it can be expressed as:
\begin{equation}
    \BmTenha_{\alpha,i} = - \int_\el \trialfcnvecther_i \cdot \bm{L}_\ther \sqbr{\check{\bm{R}} \strain_\ther\plbr{\qpolyvecther_\alpha}}  \de \el + \int_{\boundel} \trialfcnvecther_i \cdot  \tilde{\bm{R}}\plbr{\qpolyvecther_\alpha} \normedgeel  \de \boundel ,
    \label{eq:B_thermal_integrationbyparts}
\end{equation}
where the operator $\bm{L}_\ther \sqbr{\cdot}$ and the matrix $\tilde{\bm{R}}$ are defined as:
\begin{equation}
\begin{aligned}
    \bm{L}_\ther \sqbr{\cdot}= 
    \begin{bmatrix}
    ,x &  0 & ,y & 0  & 0  & 0  & 0  & 0\\ 
    0  & ,y & ,x & 0  & 0  & 0  & 0  & 0\\  
    0  &  0 & 0  & ,x & 0  & ,y & 0  & 0\\ 
    0  &  0 & 0  & 0  & ,y & ,x & 0  & 0     
    \end{bmatrix}, \qquad
    \tilde{\bm{R}} =
   \begin{bmatrix}
       \hat{N}_{xx} & \hat{N}_{xy}\\
       \hat{N}_{xy} & \hat{N}_{yy}\\
       \hat{M}_{xx} & \hat{M}_{xy}\\
       \hat{M}_{xy} & \hat{M}_{yy}
   \end{bmatrix}.
\end{aligned}   
    \label{eq:differntialmatrix_Lth_Shat_operator_th}
\end{equation}
Therefore, the matrix $\BmTenha_{\alpha,i}$ is fully computable from the degrees of freedom. The matrix $\GmTenhatilde_{\alpha,\beta}$ is easily evaluated since it involves only polynomial components. 
The thermal load vector is then obtained as:
\begin{equation}
    \plbr{\therfE}_i = \int_\el \check{\bm{R}} \plbr{x,y} \strain_\ther \plbr{\PiTenha{\trialfcnvecther_i}^T}  \de \el .
    \label{eq:thermalforces_vector}
\end{equation}

\subsubsection*{Line loads and prescribed displacements}
The construction of the line load vector follows the standard FEM procedure, as trial functions are polynomials along the element edges. The line load vector is expressed as:
\begin{equation}
    \plbr{\linefE}_i = \sum_j \int_{\boundel_j} \trialfcnvec_i \lineforcevec_j \de \boundel_j,
    \label{eq:lineforces_vector}
\end{equation}
where $\lineforcevec_j$ is the prescribed traction on the generic edge $\boundel_j$. The prescribed displacements are also imposed as in standard FEM. 

\subsection*{Supporting material for \subsect{VEMVC}: Matrix and vectors construction}

\subsubsection*{Stiffness matrix}

To obtain a projection that accounts for the spatial variation of the fiber orientation,
the following strategy is proposed. The objective is to construct a projection operator $\PiKVC$ such that:

\begin{equation}
\left\{
    \begin{aligned}
    &\left[\int_\el \strain\plbr{\trialfcnvec_i}^T \cost \plbr{x,y} \strain\plbr{\qpolyvec} \de \el \right]^{\mathrm{VC}}=
    \int_\el \strain\plbr{\PiKVC\trialfcnvec_i}^T \cost \plbr{x,y} \strain\plbr{\qpolyvec} \de \el    && \forall \qpolyvec \in \Pspacevec_\kord \plbr{\el}\\
    &\frac{1}{\nv} \sum_{j=1}^{5\nv} \text{dof}_j \plbr{\trialfcnvec_i} \text{dof}_j \plbr{\qpolyvec_{\alpha}} =
    \frac{1}{\nv} \sum_{j=1}^{5\nv} \text{dof}_j \plbr{\PiKVC\trialfcnvec_i} \text{dof}_j \plbr{\qpolyvec_{\alpha}} && \forall \alpha = 1,\dots,6,
    \end{aligned}
\right. 
    \label{eq:energynorm_projection_VC}
\end{equation}
for all $ i=1,\dots,\Vspacedim$. This leads to the standard matrix form:
\begin{equation}
    \BmKVC = \GmKVCtilde \PiKVCtilde \quad \longrightarrow \quad \PiKVCtilde = \plbr{ \GmKVCtilde }^{-1} \BmKVC,
    \label{eq:projectionVC_stiffness}
\end{equation}
where matrix $\GmKVCtilde$ reads:

\begin{equation}
    \begin{aligned}
        \GmKVCtilde_{\alpha,\beta} = \int_\el \strain\plbr{\qpolyvec_{\beta}}^T \cost \plbr{x,y} \strain\plbr{\qpolyvec_\alpha} \de \el,
    \end{aligned}
    \label{eq:BG_stiffness_VC}
\end{equation}
and its computation only requires a sufficient number of integration points to accurately evaluate the non-constant constitutive law. The notation $[\cdot]^\mathrm{VC}$ used at the left-hand side of \eq{energynorm_projection_VC} refers to the $\mathrm{VC}$-VEM technique and hides the manipulations described in Section~\ref{subsection:VEMVC}. The left-hand side of \eq{energynorm_projection_VC} is thus defined as:
\begin{equation}
    \begin{aligned}
    \BmKVC_{\alpha,i} = &- \int_\el \PiokK \trialfcnvec_i \cdot \bm{L}_d \sqbr{\cost \plbr{x,y} \strain\plbr{\qpolyvec_\alpha}}  \de \el + \int_{\boundel} \trialfcnvec_i \cdot  \hat{\stress}\plbr{\qpolyvec_\alpha} \normedgeel  \de \boundel\\
    &+ \int_\el \strain^\const\plbr{\PiokK\trialfcnvec_i}^T \cost \plbr{x,y} \strain\plbr{\qpolyvec_\alpha} \de \el,
    \end{aligned}
    \label{eq:B_stiffness_integrationbyparts_VC_projection}
\end{equation}
where $\hat{\stress}$ also contains the spatial variation of $\cost \plbr{x,y}$. 

The subsequent steps for constructing the final stiffness matrix are identical to those of the standard VEM, with the stabilization term built using the standard projector.

\subsubsection*{Self-stabilized VEM}

In the case of self-stabilized VEM, the corresponding projection is obtained:

\begin{equation}
        \left[\int_\el \strain\plbr{\trialfcnvec_i}^T \cost \plbr{x,y} \qpolyvecunc \de \el \right]^{\mathrm{VC}} = \int_\el \PiKVCSone\strain\plbr{\trialfcnvec_i}^T \cost \plbr{x,y} \qpolyvecunc \de \el \quad \forall \qpolyvecunc \in \Pspacevecunc_{\kordaug-1} \plbr{\el},
    \label{eq:system_projection_VC_selfstab_nointdofs_pproj}
\end{equation}
for all $ i=1,\dots,\Vspacedim$. In matrix form:
\begin{equation}
    \BmKVCS = \GmKVCS \PiKVCtildeSone \quad \longrightarrow \quad \PiKVCtildeSone = \plbr{ \GmKVCS }^{-1} \BmKVCS.
    \label{eq:matrix_projection_VC_selfstab_nointdofs_pprojection}
\end{equation}
The matrix $\BmKVCS$ is computed as:
 \begin{equation}
    \begin{aligned}
    \BmKVCS_{\alpha,i} &= - \int_\el \PiokK \trialfcnvec_i \cdot \bm{L}_d \sqbr{\cost \plbr{x,y} \qpolyvecunc_\alpha}  \de \el + \int_{\boundel} \trialfcnvec_i \cdot  \hat{\stress}\plbr{\qpolyvecunc_\alpha} \normedgeel  \de \boundel +\\
    &+\int_\el \strain^\const\plbr{\PiokK\trialfcnvec_i}^T \cost \plbr{x,y}\qpolyvecunc_\alpha \de \el,
    \end{aligned}
    \label{eq:B_stiffness_integrationbyparts_selfstab_nointdofs_projection_VC_pproj}
\end{equation}
where $\hat{\stress}$ also contains the spatial variation of $\cost \plbr{x,y}$. 
The matrix $\GmKVCS$ is defined as:
\begin{equation}
    \GmKVCS_{\alpha,\beta} = \int_\el {\qpolyvecunc_{\beta}}^T \cost \plbr{x,y} \qpolyvecunc_\alpha \de \el.
    \label{eq:G_VC_selfstab_nointdofs_pproj}
\end{equation}

\subsubsection*{Geometric stiffness matrix}

The projection that accounts for the spatial variability of $\mathbb{N} \plbr{x,y}$ is defined such that:

\begin{equation}
\left\{
    \begin{aligned}
    &\left[\int_\el \strain_\buck\plbr{\trialfcn_i}^T \mathbb{N} \plbr{x,y} \strain_\buck\plbr{\qpoly} \de \el \right]^{\mathrm{VC}} =
    \int_\el \strain_\buck\plbr{\PiGVC\trialfcn_i}^T \mathbb{N} \plbr{x,y} \strain_\buck\plbr{\qpoly} \de \el && \forall \qpoly \in \Pspace_\kord \plbr{\el}\\
    &\frac{1}{\nv} \sum_{j=1}^{\nv} \text{dof}_j \plbr{\trialfcn_i} \text{dof}_j \plbr{\qpoly_{\alpha}} =
    \frac{1}{\nv} \sum_{j=1}^{\nv} \text{dof}_j \plbr{\PiGVC\trialfcn_i} \text{dof}_j \plbr{\qpoly_{\alpha}} &&  \alpha = 1  ,
    \end{aligned}
\right.    
    \label{eq:bucklinhenergynorm_projection_VC}
\end{equation}
for all $ i=1,\dots,\Vspacedim^w$. In matrix form, this can be written as:
\begin{equation}
    \BmGVC = \GmGVCtilde \PiGVCtilde \quad \longrightarrow \quad \PiGVCtilde = \plbr{ \GmGVCtilde }^{-1} \BmGVC,
    \label{eq:projectionVC_buckling}
\end{equation}
with:

\begin{equation}
    \begin{aligned}
        \GmGVCtilde_{\alpha,\beta} = \int_\el \strain_\buck\plbr{\qpoly_{\beta}}^T \mathbb{N} \plbr{x,y} \strain_\buck \plbr{\qpoly_\alpha} \de \el.
    \end{aligned}
    \label{eq:BG_buckling_VC}
\end{equation}

To handle the surface integrals involving unknown terms, the matrix $\BmGVC_{\alpha,i}$ is written as:
\begin{equation}
    \BmGVC_{\alpha,i} = - \int_\el \Piok \trialfcn_i \cdot \bm{L}_b \sqbr{\mathbb{N} \plbr{x,y} \strain_\buck \plbr{\qpoly_\alpha}} \de \el + \int_{\boundel} \trialfcn_i \cdot \mathbb{N} \plbr{x,y} \strain_\buck \plbr{\qpoly_\alpha} \normedgeel  \de \boundel.
    \label{eq:B_buckling_integrationbyparts_VC_projection}
\end{equation}

\subsubsection*{Thermal forces vector}

For the thermal forces vector, the projection operator is defined such that:

\begin{equation}
\left\{
    \begin{aligned}
    &\left[\int_\el \strain_\ther\plbr{\trialfcnvecther_i}^T \check{\bm{R}} \plbr{x,y} \strain_\ther\plbr{\qpolyvecther} \de \el \right]^{\mathrm{VC}} =
    \int_\el \strain_\ther\plbr{\PiTVC\trialfcnvecther_i}^T \check{\bm{R}} \plbr{x,y} \strain_\ther\plbr{\qpolyvecther} \de \el && \forall \qpolyvecther \in \Pspacevecther_\kord \plbr{\el}\\
    &\frac{1}{\nv} \sum_{j=1}^{\nv} \text{dof}_j \plbr{\trialfcnvecther_i} \text{dof}_j \plbr{\qpolyvecther_{\alpha}} =
    \frac{1}{\nv} \sum_{j=1}^{\nv} \text{dof}_j \plbr{\PiTVC\trialfcnvecther_i} \text{dof}_j \plbr{\qpolyvecther_{\alpha}} && \forall \alpha = 1,\dots,5,
    \end{aligned}
\right.    
    \label{eq:thermalenergynorm_projection_VC}
\end{equation}
for all $ i=1,\dots,\Vspacedim^*$. This can be compactly rewritten as:
\begin{equation}
    \BmTVC = \GmTVCtilde \PiTVCtilde \quad \longrightarrow \quad \PiTVCtilde = \plbr{ \GmTVCtilde }^{-1} \BmTVC.
    \label{eq:projectionVC_thermal}
\end{equation}
In order to compute $\PiTVCtilde$, matrix $\GmGVCtilde_{\alpha,\beta}$ is constructed as:

\begin{equation}
    \begin{aligned}
        \GmGVCtilde_{\alpha,\beta} = \int_\el \strain_\ther\plbr{\qpolyvecther_{\beta}}^T \check{\bm{R}} \plbr{x,y} \strain_\ther\plbr{\qpolyvecther_\alpha} \de \el.
    \end{aligned}
    \label{eq:BG_thermal_VC}
\end{equation}

To make the expression of $\BmTVC_{\alpha,i}$ computable, the matrix is rewritten as:
\begin{equation}
    \BmTVC_{\alpha,i} = - \int_\el \PiokKther \trialfcnvecther_i \cdot \bm{L}_\ther \sqbr{\check{\bm{R}} \plbr{x,y} \strain_\ther\plbr{\qpolyvecther_\alpha}}  \de \el + \int_{\boundel} \trialfcnvecther_i \cdot  \tilde{\bm{R}}\plbr{\qpolyvecther_\alpha} \normedgeel  \de \boundel ,
    \label{eq:B_thermal_integrationbyparts_VC_projection}
\end{equation}
where $\tilde{\bm{R}}$ accounts for the dependency on $x,y$ and $\PiokKther$ is $\PiokK$ referred to $\trialfcnvecther$. 

\end{document}